\documentclass[aos,preprint]{imsart}

\usepackage{hyperref}
\usepackage{amsfonts}
\usepackage{mathtools}
\usepackage{amsgen,amsmath,amstext,amsbsy,amsopn,amssymb,subfigure,stmaryrd}
\usepackage[dvips]{graphicx}
\usepackage{comment}
\usepackage{enumitem}
\usepackage[longnamesfirst]{natbib}
\usepackage{booktabs}
\usepackage{blkarray}
\usepackage{algorithm}
\usepackage{algorithmic}
\usepackage{url}
\usepackage{longtable}
\usepackage{multirow}
\usepackage{multicol}
\usepackage[textsize=tiny]{todonotes}

\makeatletter
\def\maketag@@@#1{\hbox{\m@th\normalfont\normalsize#1}}
\makeatother

\usepackage[normalem]{ulem}

\arxiv{arXiv:2002.11255}

\startlocaldefs

\numberwithin{equation}{section}

\newtheorem{Theorem}{Theorem}[section]
\newtheorem{Definition}{Definition}[section]
\newtheorem{Assumption}{Assumption}[section]

\newtheorem{Lemma}{Lemma}[section]
\newtheorem{Remark}{Remark}[section]
\newtheorem{Corollary}{Corollary}[section]

\newtheorem{Proposition}{Proposition}[section]

\newcommand{\subscript}[2]{$#1 _ #2$}

\newcommand{\cC}{\mathcal{C}}

\newcommand{\cS}{\mathcal{S}}
\newcommand{\cT}{\mathcal{T}}
\newcommand{\cP}{\mathcal{P}}

\newcommand{\tF}{{\rm F}}

\newcommand{\cI}{{\mathcal{I}}}

\newcommand{\cA}{{\mathcal{A}}}

\newcommand{\cF}{{\mathcal{F}}}

\newcommand{\cX}{{\mathcal{X}}}
\newcommand{\cY}{{\mathcal{Y}}}
\newcommand{\cZ}{{\mathcal{Z}}}
\newcommand{\cH}{{\mathcal{H}}}

\renewcommand{\r}{{\boldsymbol{r}}}
\newcommand{\p}{{\boldsymbol{p}}}
\newcommand{\x}{{\mathbf{x}}}

\newcommand{\X}{{\mathbf{X}}}
\newcommand{\B}{{\mathbf{B}}}

\newcommand{\I}{{\mathbf{I}}}

\newcommand{\N}{{\mathbf{N}}}

\newcommand{\R}{{\mathbf{R}}}
\renewcommand{\S}{{\mathbf{S}}}
\newcommand{\U}{{\mathbf{U}}}
\newcommand{\V}{{\mathbf{V}}}
\renewcommand{\H}{{\mathbf{H}}}

\newcommand{\W}{{\mathbf{W}}}
\newcommand{\A}{{\mathbf{A}}}
\newcommand{\Y}{{\mathbf{Y}}}
\newcommand{\Z}{{\mathbf{Z}}}
\newcommand{\T}{{\mathbf{T}}}

\newcommand{\Proj}{{\mathbb{P}}}
\newcommand{\Poisson}{{\rm Poisson}}

\newcommand{\cN}{{\cal N}}

\newcommand{\rank}{{\rm rank}}

\newcommand{\tr}{{\rm tr}}

\newcommand{\diag}{{\rm diag}}

\newcommand{\SVD}{{\rm SVD}}
\newcommand{\VEC}{{\rm Vec}}

\newcommand{\HeteroPCA}{\rm HeteroPCA}

\newcommand{\argmin}{\mathop{\rm arg\min}}
\newcommand{\argmax}{\mathop{\rm arg\max}}
\newcommand{\bbR}{\mathbb{R}}
\newcommand{\bbN}{\mathbb{N}}
\newcommand{\bbO}{\mathbb{O}}
\newcommand{\bbP}{\mathbb{P}}
\newcommand{\bbE}{\mathbb{E}}
\newcommand{\cM}{\mathcal{M}}

\endlocaldefs
\begin{document}

\begin{frontmatter}

\title{An Optimal Statistical and Computational Framework for Generalized Tensor Estimation}
\runtitle{Generalized Tensor Estimation}


\begin{aug}
\author[A]{\fnms{Rungang} \snm{Han}\ead[label=e1,mark]{rhan32@stat.wisc.edu}},
\author[B]{\fnms{Rebecca} \snm{Willett}\ead[label=e2]{willett@uchicago.edu}},
\and
\author[A]{\fnms{Anru R.} \snm{Zhang}\ead[label=e3,mark]{anruzhang@stat.wisc.edu}}

\runauthor{R. Han, R. Willett, AND A. R. Zhang}

\address[A]{Department of Statistics, University of Wisconsin-Madison, \printead{e1,e3}}
\address[B]{Departments of Statistics and Computer Science, University of Chicago,\printead{e2}}          
\end{aug}

\begin{abstract}
This paper describes a flexible framework for generalized low-rank tensor estimation problems that includes many important instances arising from applications in computational imaging, genomics, and network analysis. The proposed estimator consists of finding a low-rank tensor fit to the data under generalized parametric models. To overcome the difficulty of non-convexity in these problems, we introduce a unified approach of projected gradient descent that adapts to the underlying low-rank structure. Under mild conditions on the loss function, we establish both an upper bound on statistical error and the linear rate of computational convergence through a general deterministic analysis. Then we further consider a suite of generalized tensor estimation problems, including sub-Gaussian tensor PCA, tensor regression, and Poisson and binomial tensor PCA. We prove that the proposed algorithm achieves the minimax optimal rate of convergence in estimation error. Finally, we demonstrate the superiority of the proposed framework via extensive experiments on both simulated and real data.
\end{abstract}



\end{frontmatter}

\section{Introduction}\label{sec:intro}
In recent years, the analysis of tensors or high-order arrays has emerged as an active topic in statistics, applied mathematics, machine learning, and data science. 
Datasets in the form of tensors arise from various scientific applications \citep{kroonenberg2008applied}, such as collaborative filtering \citep{bi2018multilayer,shah2019iterative}, neuroimaging analysis \citep{zhou2013tensor,li2018tucker}, hyperspectral imaging \citep{li2010tensor}, longitudinal data analysis \citep{hoff2015multilinear}, and more. In many of these problems, although the tensor of interest is high-dimensional in the sense that the ambient dimension of the dataset is substantially greater than the sample size, there is often hidden low-dimensional structures in the tensor that can be exploited to facilitate the data analysis. In particular, the low-rank condition renders convenient decomposable structure and has been proposed and widely used in the analysis of tensor data \citep{kroonenberg2008applied,kolda2009tensor}. However, leveraging these hidden low-rank structures in estimation and inference can pose great statistical and computational challenges in real practice. 
\subsection{Generalized Tensor Estimation}\label{sec:generalized-tensor-estimation}
In this paper, we consider a statistical and optimization framework for generalized tensor estimation. Suppose we observe a random sample $D$ drawn from some distribution parametrized by an unknown low-rank tensor parameter $\cX^*\in \mathbb{R}^{p_1\times p_2\times p_3}$. A straightforward idea to estimate $\cX^*$ is via optimization:
\begin{equation}\label{eq:minimization}
    \hat\cX = \argmin_{\cX \text{ is low-rank}} L(\cX; D).
\end{equation}
Here, $L(\cX; D)$ can be taken as the negative log-likelihood function (then $\hat\cX$ becomes the maximum likelihood estimator (MLE)) or any more general loss function. We can even broaden the scope of this framework to a deterministic setting: suppose we observe $D$ that is ``associated" with an unknown tensor parameter $\cX^*\in\bbR^{p_1\times p_2 \times p_3}$; to estimate $\cX^*$, we try to minimize the loss function $L(\cX; D)$ that is specified by the problem scenario. This general framework includes many important instances arising in real applications. For example:
\begin{itemize}[leftmargin=*]
	\item \emph{Computational imaging}. Photon-limited imaging appears in signal processing \citep{salmon2014poisson}, material science \citep{yankovich2016non}, astronomy \citep{timmerman1999multiscale,willett2007multiscale}, and often involves arrays with non-negative photon counts contaminated by substantial noise. Data from photon-limited imaging are often in the form of tensors (e.g., stacks of spectral images in which each image corresponds to a different wavelength of light). How to denoise these images is often crucial for the subsequent analysis. To this end, Poisson tensor PCA serves as a prototypical model for tensor photon-limited imaging analysis; see Sections \ref{sec:poisson} and \ref{sec:real-data} for more details.
	\item \emph{Analysis of multilayer network data}. In network analysis, one often observes multiple snapshots of static or dynamic networks \citep{sewell2015latent,lei2019consistent,arroyo2019inference,pensky2019spectral}. How to perform an integrative analysis for the network structure using multilayer network data has become an important problem in practice. By stacking adjacency matrices from multiple snapshots to an adjacency tensor, the hidden community structure of network can be transformed to the low-rankness of adjacency tensor, and the generalized tensor learning framework can provide a new perspective on
	the analysis of multilayer network data.	
	\item \emph{Biological sequencing data analysis}. Tensor data also commonly appear in biological sequencing data analysis \citep{faust2012microbial,flores2014temporal,wang2017three}. The identification of significant triclusters or modules, i.e., coexpressions of different genes or coexistence of different microbes, often has significant biological meanings \citep{henriques2019triclustering}. From a statistical perspective, these modules often correspond to low-rank tensor structure, so the generalized tensor learning framework could be naturally applied.

	\item \emph{Online-click through Prediction.} Online click-through data analysis in e-commerce has become an increasingly important tool in building the online recommendation system \citep{mcmahan2013ad,sun2016sparse,shan2016predicting}. There are three major entities: users, items, and time, and the data can be organized as a tensor, where each entry represents the click times of one user on a specific category of items in a time period (e.g., noon or evening). Then generalized tensor estimation could be applied to study the implicit features of users and items for better prediction of user behaviors. 
\end{itemize}
Additional applications include neuroimaging analysis \citep{zhou2013tensor}, collaborative filtering \citep{yu2018recovery}, mortality rate analysis \citep{wilmoth2016human}, and more. We also elucidate specific model setups and real data examples in detail later in Sections \ref{sec:instances} and \ref{sec:real-data}, respectively. 

The central tasks of generalized tensor estimation problems include two elements. 
From a statistical perspective, it is important to investigate how well one can estimate the target tensor parameter $\cX^*$ and  the optimal rates of estimation error. From an optimization perspective, it is crucial to develop a computationally efficient procedure for estimating $\cX^*$ with provable theoretical guarantees. To estimate the low-rank tensor parameter $\cX^*$, a straightforward idea is to perform the rank constrained minimization on the loss function $L(\cX; D)$ in \eqref{eq:minimization}. Since the low-rank constraint is highly non-convex, the direct implementation of \eqref{eq:minimization} is computationally infeasible in practice. If $\cX^*$ is a sparse vector or low-rank matrix, common substitutions often involve convex regularization methods, such as M-estimators with an $\ell_1$ penalty or matrix nuclear norm penalty for estimating sparse or low-rank structure \citep{tibshirani1996regression,fazel2002matrix}. These methods enjoy great empirical and theoretical success for vector/matrix estimators, but it is unclear whether they can achieve good performance on generalized tensor estimation problems. First, different from the matrix nuclear norm, tensor nuclear norm is generally NP-hard to even approximate \citep{friedland2018nuclear}, so that the tensor nuclear norm regularization approach can be computationally intractable. Second, other computationally feasible convex regularization methods, such as the overlapped nuclear norm minimization \citep{tomioka2011statistical,tomioka2013convex}, may be statistical sub-optimal based on the theory of simultaneously structured model estimation \citep{oymak2015simultaneously}.

\sloppypar In contrast, we focus on a unified non-convex approach for generalized tensor estimation problems in this paper. Our central idea is to decompose the low-rank tensor into $\cX = \llbracket\cS; \U_1, \U_2, \U_3\rrbracket$ (see Section \ref{sec:notation} for explanations of tensor algebra) and reformulate the original problem to
\small
\begin{equation}\label{eq:model2}
\begin{split}
& (\hat \cS, \hat \U_1, \hat \U_2, \hat \U_3) =  \argmin_{\cS, \U_1, \U_2, \U_3} \left\{L(\llbracket\cS; \U_1, \U_2, \U_3\rrbracket; D) + \frac{a}{2} \sum_{k=1}^3\left\|\U_k^\top \U_k - b^2\I_{r_k}\right\|_\tF^2\right\},
\end{split}
\end{equation}
\normalsize
which can be efficiently solved by (projected) gradient descent on all components. The resulting $\hat\cX = \llbracket \hat \cS; \hat \U_1, \hat \U_2 , \hat \U_3\rrbracket$ naturally admits a low-rank structure. The auxiliary regularizers $\left\|\U_k^\top \U_k - b^2\I_{r_k}\right\|_\tF^2$ in \eqref{eq:model2} can keep $\hat\U_k$ from being singular. It is actually easy to check that \eqref{eq:minimization} and  \eqref{eq:model2} are exactly equivalent.

We provide strong theoretical guarantees for the proposed procedure on generalized tensor estimation problems.  In particular, we establish the linear rate of local convergence for gradient descent methods under a general deterministic setting with the \emph{Restricted Correlated Gradient} condition (see Section \ref{sec:RCG} for details). An informal statement of the result is given below,
\begin{equation}\label{eq-informal}
    \left\|\cX^{(t)} - \cX^*\right\|_\tF^2 \lesssim \xi^2 + \left(1 - c\right)^{t} \left\|\cX^{(0)} - \cX^*\right\|_\tF^2 \quad \text{for all $t\geq 1$}
\end{equation}
with high probability. Here, we use $\xi^2$ to characterize the  statistical noise and its definition and interpretation will be given in section \ref{sec:theory}.
Then for specific statistical models, including sub-Gaussian tensor PCA, tensor regression, Poisson tensor PCA, and binomial tensor PCA, based on the general result \eqref{eq-informal}, we prove that the proposed algorithm achieves the minimax optimal rate of convergence in estimation error. Specifically for the low-rank tensor regression problem, Table \ref{tab:reg_comparison} illustrates the advantage of our method through a comparison with existing ones.

\begin{table}[ht]
	\centering
	{\footnotesize
	\begin{tabular}{c|c|c|c}
		\hline
		\multirow{2}{1.6cm}{\centering Algorithm}  & \multirow{2}{2.7cm}{\centering Sample complexity$^*$} & \multirow{2}{2.6cm}{\centering Estimation Error Upper Bound} & \multirow{2}{2.8cm}{\centering Recovery (noiseless)}\\ & & \\ \hline
        \multirow{2}{2cm}{\centering Our Method} & \multirow{2}{2cm}{\centering $p^{3/2}r$} & \multirow{2}{2cm}{\centering $\sigma^2 pr/n$} & \multirow{2}{1.8cm}{\centering Exact} \\ & & \\ \hline
        \multirow{2}{4cm}{\centering Tucker-Reg. \\ {\tiny\citep{zhou2013tensor}}}  & \multirow{2}{2cm}{\centering N.A.} & \multirow{2}{2cm}{\centering N.A.} & \multirow{2}{1.8cm}{\centering Exact}  \\ & & \\ \hline
        \multirow{2}{4cm}{\centering Nonconvex-PGD. \\ {\tiny\citep{chen2019non}}}  & \multirow{2}{2cm}{\centering $p^2r$} & \multirow{2}{2cm}{\centering $\sigma^2 p^2r/n$} & \multirow{2}{1.8cm}{\centering Exact} \\ & & \\\hline
        \multirow{2}{4cm}{\centering Nuclear Norm Min. \\ {\tiny\citep{raskutti2019convex}}}  & \multirow{2}{2cm}{\centering N.A.} & \multirow{2}{2cm}{\centering $\sigma^2 pr^2/n$} & \multirow{2}{1.8cm}{\centering Exact} \\ & & \\ \hline
        \multirow{2}{4cm}{\centering Schatten-1 Norm Min. \\ {\tiny\citep{tomioka2013convex}}}  & \multirow{2}{2cm}{\centering $p^2r$} & \multirow{2}{2cm}{\centering $\sigma^2 p^2r/n$} & \multirow{2}{1.8cm}{\centering Exact} \\& & \\ \hline
        \multirow{2}{4cm}{\centering ISLET  {\tiny\citep{zhang2018ISLET}}}  & \multirow{2}{2cm}{\centering $p^{3/2}r$} & \multirow{2}{2cm}{\centering $\sigma^2pr/n$} & \multirow{2}{1.8cm}{\centering Inexact} \\ & & \\ \hline
        \multirow{2}{4cm}{\centering Iterative Hard Thresholding\protect\footnotemark \\ {\tiny\citep{rauhut2017low}}}   & \multirow{2}{2cm}{\centering $pr$} & \multirow{2}{2cm}{\centering $\sigma^2$} & \multirow{2}{1.8cm}{\centering Exact}  \\ & &\\ \hline
	\end{tabular}
	}
	\caption{Comparison of different tensor regression methods when the rank is known. For simplicity, we assume $r_1=r_2=r_3=r$, $p_1=p_2=p_3=p$ and $\sigma^2 \ll \|\cX^*\|_\tF^2$. Here, the sample complexity$^*$ is the minimal sample size required to achieve the corresponding estimation error. 
	}
	\label{tab:reg_comparison}
\end{table}

\footnotetext[1]{The analysis in \cite{rauhut2017low} relies on an assumption that the projection on low-rank tensor manifold can be approximately done by High-Order SVD. It is, however, unclear whether this assumption holds in general.}

Finally, we apply the proposed framework to synthetic and real data examples, including photon-limited 4D-STEM (scanning transmission electron microscopy) imaging data and click-through e-commerce data. The comparison of performance with existing methods illustrates the merit of our proposed procedure.

\subsection{Related Literature}\label{sec:literature}

This work is related to a broad range of literature on tensor analysis. For example, tensor decomposition/SVD/PCA focuses on the extraction of low-rank structures from noisy tensor observations \citep{richard2014statistical,anandkumar2014tensor,hopkins2015tensor,montanari2017limitation,lesieur2017statistical,johndrow2017tensor,chen2019phase}. Correspondingly, a number of methods have been proposed and analyzed under either deterministic or random Gaussian noise, such as the maximum likelihood estimation \citep{richard2014statistical}, (truncated) power iterations \citep{anandkumar2014tensor,sun2017provable}, higher-order SVD \citep{de2000multilinear}, higher-order orthogonal iteration (HOOI) \citep{de2000best,zhang2018tensor}, STAT-SVD \citep{zhang2017optimal-statsvd}.

Since non-Gaussian-valued tensor data also commonly appear in practice, \cite{signoretto2011tensor,chi2012tensors,hong2018generalized} considered the generalized tensor decomposition and introduced computational efficient algorithms. However, the theoretical guarantees for these procedures and the statistical performances of the generalized tensor decomposition still remain open. 

 Our proposed framework includes the topic of tensor recovery and tensor regression. Various methods, such as the convex regularization \citep{tomioka2013convex,raskutti2019convex}, alternating minimization \citep{zhou2013tensor}, hard thresholding iteration \citep{chen2019non,rauhut2017low,rauhut2015tensor}, importance-sketching \citep{zhang2018ISLET} were introduced and studied. A more detailed comparison of these methods is summarized in Table \ref{tab:reg_comparison}. 

In addition,  high-order interaction pursuits \citep{hao2018sparse}, tensor completion \citep{liu2013tensor,yuan2014tensor,montanari2018spectral,xia2017polynomial,xia2017statistically,zhang2019cross,cai2019nonconvex}, and tensor block models \citep{ChiGaiSunZhoYan2018,lei2019consistent,wang2019multiway} are important topics in tensor analysis that have attracted enormous attention recently. Departing from the existing results, this paper, to the best of our knowledge, 
is the first to give a unified treatment for a broad range of tensor estimation problems with both statistical optimality and computational efficiency.

This work is also related to a substantial body of literature on low-rank matrix recovery, where the goal is to estimate a low-rank matrix based on a limited number of observations. Specific examples of this topic include matrix completion \citep{candes2009exact,candes2010matrix}, phase retrieval \citep{candes2015phase,cai2016optimal}, blind deconvolution \citep{ahmed2013blind}, low-rank matrix trace regression \citep{keshavan2010matrix,koltchinskii2011nuclear,chen2018harnessing,fan2019generalized}, and many others. 
A common approach for low-rank matrix recovery is via explicit low-rank factorization: one can decompose the target $p_1$-by-$p_2$ rank-$r$ matrix $\X$ into $\X=\U\V^\top$, where $\U\in\bbR^{p_1 \times r}, \V \in \bbR^{p_2\times r}$, then minimize the loss function $L(\U\V^\top)$ with respect to both $\U$ and $\V$ \citep{wen2012solving}.
Previously, \cite{zhao2015nonconvex} considered the noiseless setting of trace regression and proved that under good initialization, the first order alternating optimization on $\U$ and $\V$ achieves exact recovery. \cite{tu2016low,park2018finding} established the local convergence of gradient descent for strongly convex and smooth loss function $L$. The readers are referred to a recent survey paper \citep{chi2019nonconvex} on the applications and optimization landmarks of the non-convex factorized optimization. 
Despite significant developments in low-rank matrix recovery and non-convex optimization, they cannot be directly generalized to tensor estimation problems for many reasons. First, many basic matrix concepts or methods cannot be directly generalized to high-order ones \citep{hillar2013most}. Naive generalization of matrix concepts (e.g., operator norm, singular values, eigenvalues) are possible but often computationally NP-hard. Second, tensors have more complicated algebraic structure than matrices. As what we will illustrate later, one has to simultaneously handle all arm matrices (i.e., $\U_1$, $\U_2$, and $\U_3$) and the core tensor (i.e., $\cS$) with distinct dimensions in the theoretical error contraction analysis. To this end, we develop new technical tools on tensor algebra and perturbation results (e.g., \ref{lm-equivalent-criteria}, Lemmas \ref{lm-X-decomposition} in the Appendix). More technical issues of generalized tensor estimation will be addressed in Section \ref{sec:proof_sketch}.

The projected gradient schemes, which apply gradient descent on the parameter tensor $\cX$ followed by the low-rank tensor retraction/projection operators, form another important class of methods in the literature \citep{rauhut2015tensor,rauhut2017low,chen2019non}:
\begin{equation*}\label{eq:low-rank-ortohgonal-projection}
    \cX^{(t+1)} = \mathcal P\left(\cX^{(t)} - \eta \nabla L(\cX^{(t)})\right).
\end{equation*}
Different from the low-rank projection for matrices, the exact low-rank tensor projection (i.e., the best rank-$(r_1,r_2,r_3)$ approximation: $\mathcal P (\cX) = \argmin_{\mathcal{T}\text{ is rank-$(r_1, r_2, r_3)$}}\|\mathcal{T} - \cX\|_\tF$) is NP-hard in general \citep{hillar2013most} and less practical. Several inexact but efficient projection methods were developed and studied to overcome this issue. In particular, \cite{rauhut2015tensor} proposed a polynomial-time computable projected gradient scheme that converges linearly to the true tensor parameter for the noiseless tensor completion problem, given the initialization is sufficiently close to the solution. On the other hand, it is not clear if such schemes with inexact projection operators can achieve optimal statistical rate in the noisy setting \citep{chen2019non}. In contrast, the proposed method in this paper is both computationally efficient and statistically optimal in a variety of settings with provable guarantees.

\subsection{Organization of the Paper}

The rest of the article is organized as follows. After a brief introduction of the notation and preliminaries in Section \ref{sec:notation}, we introduce the general problem formulation in Section \ref{sec:model-descent}. A deterministic error and local convergence analysis of the projected gradient descent algorithm for order-3 tensor estimation is discussed in Section \ref{sec:PRGD}. Then we apply the results on a variety of generalized tensor estimation problems in Section \ref{sec:instances}, including sub-Gaussian tensor PCA, tensor regression, Poisson tensor PCA, and binomial tensor PCA. We develop the upper and minimax matching lower bounds in each of these scenarios. In Section \ref{sec:rank-estimation}, we propose a data-driven rank selection method with theoretical guarantee. The extension to general order-$d$ tensor estimation is discussed in Section \ref{sec:extension}. Simulation and real data analysis are presented in Section \ref{sec:numerics}. All proofs of technical results and more implementation details of algorithms are collected in the supplementary materials. 

\section{Generalized Tensor Estimation Model}\label{sec:formulation}

\subsection{Notation and Preliminaries}\label{sec:notation}

The following notation and preliminaries are used throughout this paper. The lowercase letters, e.g., $x, y, u, v$, are used to denote scalars or vectors. For any $a, b \in \bbR$, let $a\wedge b$ and $a\vee b$ be the minimum and maximum of $a$ and $b$, respectively. We use $C, C_0, C_1, \ldots$ and $c, c_0, c_1,\ldots$ to represent generic large and small positive constants respectively. The actual values of these generic symbols may differ from line to line.

We use bold uppercase letters $\A$, $\B$ to denote matrices. Let $\bbO_{p,r}$ be the collection of all $p$-by-$r$ matrices with orthonormal columns: $\bbO_{p,r} = \{\U \in \bbR^{p\times r}: \U^\top \U = \I_r\}$, where $\I_r$ is the $r$-by-$r$ identity matrix. For any matrix $\A\in \mathbb{R}^{p_1\times p_2}$, let $\sigma_1(\A) \geq \cdots \geq \sigma_{p_1\wedge p_2}(\A) \ldots \geq 0$ be its singular values in descending order. We also define $\SVD_r(\A) \in \mathbb{O}_{p, r}$ to be the matrix comprised of the top $r$ left singular vectors of $\A$. For any matrix $\A$, let $\A_{ij}, \A_{i\cdot}$, and $\A_{\cdot j}$ be the entry on the $i$th row and $j$th column, the $i$th row, and the $j$th column of $\A$, respectively. The inner product of two matrices with the same dimension is defined as $\langle\A,\B\rangle = \tr(\A^\top\B)$, where $\tr(\cdot)$ is the trace operator. We use $\|\A\| = \sigma_1(\A)$ to denote the spectral norm of $\A$, use $\|\A\|_\tF = \sqrt{\sum_{i,j} \A_{ij}^2} = \sqrt{\sum_{k=1}^{p_1 \wedge p_2} \sigma_k^2}$ to denote the Frobenius norm of $\A$, and use $\|\A\|_* = \sum_{k=1}^{p_1\wedge p_2}\sigma_k$ to denote the nuclear norm of $\A$. The $l_{2,\infty}$ norm of $\A$ is defined as the largest row-wise $l_2$ norm of $\A$: $\left\|\A\right\|_{2,\infty} = \max_{i} \left\|\A_{i\cdot}\right\|_2$. For any matrix $\A=[a_1,\ldots, a_J]\in \mathbb{R}^{I\times J}$ and $\B\in \mathbb{R}^{K\times L}$, the \emph{Kronecker product} is defined as the $(IK)$-by-$(JL)$ matrix $\A\otimes\B = [a_1\otimes \B \cdots a_J\otimes \B]$. 

In addition, we use calligraphic letters, e.g., $\cS, \cX, \cY$, to denote higher-order tensors. To simplify the presentation, we mainly focuses on order-3 tensors in this paper while all results for higher-order tensors can be carried out similarly. For tensor $\cS\in \bbR^{r_1\times r_2 \times r_3}$ and matrix $\U_1 \in \bbR^{p_1 \times r_1}$, the mode-1 tensor-matrix product is defined as:
\begin{equation*}
\cS \times_1 \U_1 \in \bbR^{p_1 \times r_2 \times r_3},\quad (\cS \times \U_1)_{i_1i_2i_3} = \sum_{j=1}^{r_1} \cS_{j i_2 i_3} (\U_1)_{i_1 j}.
\end{equation*}
For any $\U_2\in \mathbb{R}^{p_2\times r_2}, \U_3\in \mathbb{R}^{p_3\times r_3}$, the tensor-matrix products $\cS\times_2\U_2$ and $\cS\times_3\U_3$ are defined in a similarly way. Importantly, multiplication along different directions is commutative invariant: $(\cS\times_{k_1}\U_{k_1}) \times_{k_2}\U_{k_2} = (\cS \times_{k_2}\U_{k_2}) \times_{k_1} \U_{k_1}$ for any $k_1\neq k_2$. We simply denote
$$\cS\times_1 \U_1\times_2 \U_2\times_3 \U_3 = \llbracket\cS; \U_1, \U_2, \U_3 \rrbracket,$$
as this formula commonly appears in the analysis. We also introduce the matricization operator that transforms tensors to matrices: for $\cX \in \bbR^{p_1 \times p_2\times p_3}$, define
\begin{equation*}
    \begin{split}
	    \cM_1 (\cX) \in \bbR^{p_1 \times p_2p_3},\quad \text{where}\quad [\cM_1(\cX)]_{i_1,i_2 + p_2 (i_3-1)} = \cX_{i_1i_2i_3}, \\
	    \cM_2 (\cX) \in \bbR^{p_2 \times p_1p_3},\quad \text{where}\quad [\cM_2(\cX)]_{i_2,i_3 + p_3 (i_1-1)} = \cX_{i_1i_2i_3}, \\
	    \cM_3 (\cX) \in \bbR^{p_3 \times p_1p_2},\quad \text{where}\quad [\cM_3(\cX)]_{i_3,i_1 + p_1 (i_2-1)} = \cX_{i_1i_2i_3},
	\end{split}
\end{equation*}
and $\cM_k^{-1}:\mathbb{R}^{p_k \times(p_{k+1}p_{k+2})}\to \mathbb{R}^{p_1\times p_2\times p_3}$ as the inverse operator of $\cM_k$ where $k+1$ and $k+2$ are computed modulo 3. Essentially, $\cM_k$ ``flattens'' all but the $k$th directions of any tensor. The following identity that relates the matrix-tensor product and matricization plays an important role in our analysis:
\begin{equation*}
    \cM_k(\cS \times_1 \U_1 \times_2 \U_2 \times_3 \U_3) = \U_k \cM_k(\cS) \left(\U_{k+2} \otimes \U_{k+1}\right)^\top, \quad k=1,2,3.
\end{equation*}
Here again, $k+1$ and $k+2$ are computed modulo 3. The inner product of two tensors with the same dimension is defined as $\langle \cX, \cY \rangle = \sum_{ijk}\cX_{ijk}\cY_{ijk}$. The Frobenius norm of a tensor $\cX$ is defined as $\|\cX\|_\tF = \sqrt{\sum_{i,j,k} \cX_{ijk}^2}$. For any smooth tensor-variate function $f: \mathbb{R}^{p_1\times p_2\times p_3}\to \mathbb{R}$, let $\nabla f:\mathbb{R}^{p_1\times p_2\times p_3} \to \mathbb{R}^{p_1\times p_2\times p_3}$ be the gradient function such that $\left(\nabla f(\cX)\right)_{ijk} = \frac{\partial f}{\partial \cX_{ijk}}$. We simply write this as $\nabla f$ when there is no confusion. Finally, the readers are also referred to \cite{kolda2009tensor} for a comprehensive discussions on tensor algebra. The focus of this paper is on the following low-Tucker-rank tensors:
\begin{Definition}[Low Tucker Rank]\label{def:lowrank}
	We say $\cX^* \in \mathbb{R}^{p_1\times p_2\times p_3}$ is Tucker rank-$(r_1,r_2,r_3)$ if and only if $\cX^*$ can be decomposed as 
	$$\cX^* = \cS^* \times_1 \U_1^* \times_2 \U_2^* \times_3 \U_3^* =: \llbracket\cS^*; \U_1^*, \U_2^*, \U_3^*\rrbracket$$
	for some $\cS^* \in \bbR^{r_1 \times r_2 \times r_3}$ and $\U_k^* \in \bbR^{p_k \times r_k}, k=1,2,3$.
\end{Definition}
In addition, $\cX^*$ is Tucker rank-($r_1, r_2, r_3$) if and only if $\rank(\mathcal{M}_k(\cX^*))\leq r_k$ for $k=1,2,3$. For convenience of presentation, we denote $\bar p = \max\{p_1,p_2,p_3\}$, $\bar r = \max\{r_1, r_2, r_3\}$, $\underline p = \min\{p_1,p_2,p_3\}$, $\underline r = \min\{r_1,r_2,r_3\}$, $p_{-k} = p_1p_2p_3/p_k$ and $r_{-k}=r_1r_2r_3/r_{-k}$. 

\subsection{Generalized Tensor Estimation}\label{sec:model-descent}

Suppose we observe a dataset $D$ associated with an unknown parameter $\cX^*$. Here, $\cX^*$ is a $p_1$-by-$p_2$-by-$p_3$ rank-$(r_1,r_2,r_3)$ tensor and $r_k \ll p_k$. For example, $D$ can be a random sample drawn from some distribution parametrized by $\cX^*$ or a deterministic perturbation of $\cX^*$. The central goal of this paper is to have an efficient and accurate estimation of $\cX^*$. 

Let $L(\cX; D)$ be an empirical loss function known a priori, such as the negative log-likelihood function from the generating distribution or more general objective function. Then the following rank constrained optimization provides a straightforward way to estimate $\cX^*$:
\begin{equation}\label{eq:framework}
    \min_{\cX} L(\cX; D) \quad \text{subject to}\quad \rank(\cM_k(\cX)) \leq r_k, \quad  k=1,2,3.
\end{equation}
As mentioned earlier, this framework includes many instances arising from applications in various fields. Due to the connection between low Tucker rank and the decomposition discussed in Section \eqref{sec:notation}, it is natural to consider the following minimization problem
\begin{equation}\label{eq:rank-constraint-opt}
    \begin{split}
        \hat \cX &= \hat \cS \times_1 \hat \U_1 \times_2 \hat \U_2 \times_3 \hat \U_3,\\
        (\hat \cS, \hat \U_1, \hat \U_2 ,\hat \U_3) &= \argmin_{\cS, \U_1, \U_2, \U_3} L(\cS \times_1 \U_1 \times_2 \U_2 \times_3 \U_3; D),
    \end{split}
\end{equation}
and consider a gradient-based optimization algorithm to estimate $\cX^*$. Let $\nabla L:\mathbb{R}^{p_1\times p_2\times p_3}\to \mathbb{R}^{p_1\times p_2\times p_3}$ be the gradient of loss function. The following Lemma gives the partial gradients of $L$ on $\U_k$ and $\cS$. The proof is provided in the supplementary material (Appendix \ref{sec:proof-lemma1}). 
\begin{Lemma}[Partial Gradients of Loss]\label{lm:partial-gradient}
\begin{equation}\label{eq-four-part-gradient}
\begin{split}
	& \nabla_{\U_1} L\left(\llbracket\cS; \U_1, \U_2, \U_3\rrbracket\right) = \cM_1(\nabla L)(\U_3\otimes \U_2)\cM_1(\cS)^\top, \\
	& \nabla_{\U_2} L(\llbracket\cS; \U_1, \U_2, \U_3\rrbracket) = \cM_2(\nabla L)(\U_1\otimes \U_3)\cM_2(\cS)^\top, \\
	& \nabla_{\U_3} L(\llbracket\cS; \U_1, \U_2, \U_3\rrbracket) = \cM_3(\nabla L)(\U_2\otimes \U_1)\cM_3(\cS)^\top, \\
	& \nabla_{\cS} L\left(\llbracket\cS; \U_1, \U_2, \U_3\rrbracket\right) = \nabla L
	\times_1 \U_1^\top \times_2 \U_2^\top \times_3 \U_3^\top. 
\end{split}
\end{equation}
Here, $\nabla L$ is short for $\nabla L\left(\llbracket\cS; \U_1, \U_2, \U_3\rrbracket\right)$.
\end{Lemma}
As mentioned earlier, we consider optimizing the following objective function:
\begin{equation}\label{eq-loss-regularize}
    \begin{split}
        F(\cS&,\U_1 ,\U_2,\U_3) = L\left(\llbracket\cS; \U_1,\U_2,\U_3\rrbracket; D\right) + \frac{a}{2} \sum_{k=1}^3\left\|\U_k^\top \U_k - b^2\I_{r_k}\right\|_\tF^2,
    \end{split}
\end{equation}
where $a\geq 0, b>0$ are tuning parameters to be discussed later. By adding regularizers $\|\U_k^\top \U_k - b^2\I_{r_k}\|_\tF^2$, we can prevent $\U_k$ from being singular throughout gradient descent, while do not alter the minimizer. This can be summarized as the following proposition, whose proof is provided in Appendix \ref{sec:proof-prop1}.
\begin{Proposition}\label{prop:loss-equivalent}
     Suppose $(\hat \cS,\hat \U_1, \hat \U_2, \hat \U_2) = \argmin F(\cS,\U_1,\U_2,\U_3)$ for $F$ defined in \eqref{eq-loss-regularize}. Then
     $$\hat \cX = \llbracket \hat \cS,\hat \U_1, \hat \U_2, \hat \U_2 \rrbracket = \argmin\limits_{\cX:\rank(\cX)\leq (r_1,r_2,r_3)}L(\cX;D).$$
\end{Proposition}
Similar regularizers have been widely used on non-convex low-rank matrix optimization \citep{tu2016low,park2018finding} and more technical interpretations are provided in Section \ref{sec:proof_sketch}.

\section{Projected Gradient Descent}\label{sec:PRGD}

In this section, we study the local convergence of the projected gradient descent under a general deterministic framework.
\subsection{Restricted Correlated Gradient Condition}\label{sec:RCG}
We first introduce the regularity condition on the loss function $L$ and set $\cC$.
\begin{Definition}[Restricted Correlated Gradient (RCG)]\label{asmp-restricted-convexity}
Let $f$ be a real-valued function. We say $f$ satisfies $RCG\left(\alpha,\beta,\cC\right)$ condition for $\alpha, \beta > 0$ and the set $\cC$ if
\begin{equation}\label{eq-restricted-convex}
\left\langle \nabla f(x) - \nabla f(x^*), x - x^* \right\rangle \geq \alpha \left\|x - x^*\right\|_2^2 + \beta\left\|\nabla f(x) - \nabla f(x^*)\right\|_2^2
\end{equation}
for any $x \in \mathcal{C}$. Here, $x^*$ is some fixed target parameter.
\end{Definition}
Our later analysis will be based on the assumption that $L$ satisfies the RCG condition on to-be-specified sets of tensors with $x^* = \cX^*$ being the true parameter tensor.

\begin{Remark}[Interpretation of the RCG Condition]
The RCG condition is similar to the ``regularity condition'' appearing in recent nonconvex optimization literature \citep{chen2015solving,candes2015phase, chi2019nonconvex,yonel2020deterministic}:
\begin{equation}\label{eq-literature-RCG}
    \left\langle \nabla f(x), x - x^\# \right\rangle \geq \alpha \left\|x-x^\#\right\|_2^2 + \beta \left\|\nabla f(x)\right\|_2^2,
\end{equation}
where $f(\cdot)$ is the objective function in their context and $x^\#$ is the minimizer of $f(\cdot)$. The RCG condition can be seen as a generalization of \eqref{eq-literature-RCG}: in the deterministic case without statistical noise, the target $x^*$ usually becomes an exact stationery point of $f(\cdot)$ and \eqref{eq-restricted-convex} reduces to \eqref{eq-literature-RCG}. In addition, it is worthy noting that RCG condition does not require the function $f$ to be convex since $x^*$ is only a fixed target parameter in the requirement \eqref{eq-literature-RCG} (also see Figure 1 in \cite{chi2019nonconvex} for an example).

\end{Remark}

\subsection{Theoretical Analysis}\label{sec:theory}
We now consider a general setting that the loss function $L$ satisfies the RCG condition in a constrained domain:
\begin{equation}\label{eq:constraint-C}
\cC = \left\{\cX\in\mathbb{R}^{p_1\times p_2\times p_3}: \cX = \llbracket \cS; \U_1, \U_2, \U_3\rrbracket, \U_k \in \cC_k, \cS\in \cC_\cS\right\},
\end{equation}
where the true parameter tensor $\cX^*$ is feasible -- that is, $\cX^* \in \cC$. Here, $\cC_k$ and $\cC_\cS$ are some convex and rotation invariant sets: for any $\U_k \in \cC_k$, $\cS \in \cC_\cS$, we have $\U_k \R_k \in \cC_k$ and $\llbracket \cS; \R_1,\R_2, \R_3 \rrbracket \in \cC_\cS$ for arbitrary orthogonal matrices $\R_k \in \bbO_{r_k}$. Some specific problems of this general setting will be discussed in Section \ref{sec:instances}.

When $L$ and $\cX^*$ satisfy the condition above, we introduce the projected gradient descent in Algorithm \ref{alg:PGD}. In addition to the vanilla gradient descent, the proposed Algorithm \ref{alg:PGD} includes multiple projection steps to ensure that $\cX^*$ is in the regularized domain $\cC$ throughout the iterations.
\begin{algorithm}
\caption{Projected Gradient Descent}
\label{alg:PGD}
\begin{algorithmic}
\REQUIRE Initialization $\left(\cS^{(0)}, \U_1^{(0)}, \U_2^{(0)}, \U_3^{(0)}\right)$, constraint sets $\{\cC_k\}_{k=1}^3$, $\cC_\cS$,
tuning parameters $a,b>0$, step size $\eta$.
\FORALL {$t=0$ to $T-1$} 
	\FORALL {$k=1,2,3$}
	\STATE {$\tilde \U_k^{(t+1)} = \U_k^{(t)} - \eta\left(\nabla_{\U_k} L(\cS^{(t)}, \U_1^{(t)}, \U_2^{(t)}, \U_3^{(t)})+ a\U_k^{(t)}(\U_k^{(t)\top}\U_k^{(t)} - b^2 \I)\right)$}
	\STATE {$\U_k^{(t+1)} = \cP_{\cC_k}(\tilde \U_k^{(t+1)})$, where $\cP_{\cC_k}(\cdot)$ is the projection onto $\cC_k$.}
	\ENDFOR
	\STATE {$\tilde\cS^{(t+1)} = \cS^{(t)} - \eta \nabla_\cS L(\cS^{(t)}, \U_1^{(t)}, \U_2^{(t)}, \U_3^{(t)})$}
	\STATE {$\cS^{(t+1)} = \cP_{\cC_\cS}(\tilde \cS^{(t+1)})$, where $\cP_{\cC_\cS}(\cdot)$ is the projection onto $\cC_\cS$.}
\ENDFOR
\RETURN {$\cX^{(T)} = \cS^{(T)} \times_1 \U_1^{(T)} \times_2 \U_2^{(T)} \times_3 \U_3^{(T)}$}
\end{algorithmic}
\end{algorithm}

Suppose the true parameter $\cX^*$ is of Tucker rank-($r_1, r_2, r_3$). We also introduce the following value to quantify how different the $\cX^*$ is from being a stationary point of $L(\cX; D)$:
\begin{equation}\label{eq-noise-xi}
\begin{split}
\xi &:=  \sup_{\substack{\mathcal{T}\in \mathbb{R}^{p_1\times p_2\times p_3}, \|\mathcal{T}\|_\tF \leq 1 \\ \rank(\mathcal{T})\leq (r_1,r_2,r_3)}} \left|\left\langle \nabla L(\cX^*), \mathcal{T} \right\rangle\right|.
\end{split}
\end{equation}
Intuitively speaking, $\xi$ measures the amplitude of $\nabla L(\cX^*)$ projected onto the manifold of low-rank tensors. In many statistical models, $\xi$ essentially characterizes the amplitude of statistical noise. Specifically in the noiseless setting, $\cX^\ast$ is exactly a stationary point of $L$, then $\nabla L(\cX^\ast)=0$, $\xi = 0$. In various probabilistic settings, a suitable $L$ often satisfies $\bbE\nabla L(\cX^*) = 0$; then $\xi$ reflects the reduction of variance of $\nabla L(\cX^*)$ after projection onto the low-rank tensor manifold. We also define
\begin{equation*}
	\begin{split}
		\overline \lambda &:= \max \left\{\left\|\cM_1(\cX^*)\right\|, \left\|\cM_2(\cX^*)\right\|,\left\|\cM_3(\cX^*)\right\|\right\},\\
		\underline{\lambda} &:= \min \left\{\sigma_{r_1}\left(\cM_1(\cX^*)\right), \sigma_{r_2}\left(\cM_2(\cX^*)\right), \sigma_{r_3}\left(\cM_3(\cX^*)\right)\right\},
	\end{split}
\end{equation*}
and $\kappa=\overline{\lambda}/\underline{\lambda}$ can be regarded as a tensor condition number, as similarly defined for matrices. It is note worthy that the curvature of Tucker rank-($r_1,r_2,r_3$) tensor manifold on $\cX^*$ can be bounded by $\underline\lambda^{-1}$ \cite[Lemma 4.5]{lubich2013dynamical}.

We are now in position to establish a deterministic upper bound on the estimation error and a linear rate of convergence for the proposed Algorithm \ref{alg:PGD} when a warm initialization is provided. Specific initialization algorithms for different applications will be discussed in Section \ref{sec:instances}. 
\begin{Theorem}[Local Convergence]\label{thm:local-convergence-PGD}
Suppose $L$ satisfies $RCG(\alpha,\beta,\cC)$ for $\cC$ defined in \eqref{eq:constraint-C} and $b \asymp \overline \lambda^{1/4}$, $a = \frac{4\alpha b^4}{3\kappa^2}$ in \eqref{eq-loss-regularize}. Assume $\cX^* = \llbracket \cS^*; \U_1^*, \U_2^*, \U_3^*\rrbracket$ such that $\U_k^{*\top}\U_k^* = b^2\I_{r_k}, \U_k^* \in \cC_k, k=1,2,3$, and $\cS^* \in \cC_\cS$. Suppose the initialization $\cX^{(0)} = \llbracket \cS^{(0)};  \U_1^{(0)}, \U_2^{(0)}, \U_3^{(0)}\rrbracket$ satisfies $\left\|\cX^{(0)} - \cX^*\right\|_\tF^2 \leq c\frac{\alpha\beta\underline \lambda^2}{\kappa^{2}}$ 
for some small constant $c>0$, $\U_k^{(0)}\U_k^{(0)} = b^2\I_{r_k}, \U_k^{(0)} \in \cC_k, k=1,2,3$ and $\cS^{(0)} \in \cC_\cS$. Also, the signal-noise-ratio satisfies $\underline{\lambda}^2 \geq C_0\frac{\kappa^4}{\alpha^3\beta} \xi^2$ for some universal constant $C_0$.  Then there exists a constant $c>0$ such that if $\eta = \frac{\eta_0\beta}{b^6}$ for $\eta_0\leq c$, we have
\begin{equation*}
	\left\|\cX^{(t)} - \cX^*\right\|_\tF^2 \leq C\left(\frac{\kappa^4}{\alpha^2}\xi^2 + \kappa^2\left(1 - \frac{2\rho\alpha\beta\eta_0}{\kappa^2}\right)^t \left\|\cX^{(0)}-\cX^*\right\|_\tF^2\right).
\end{equation*}
\end{Theorem}
In addition, the following corollary provides a theoretical guarantee for the estimation loss of the proposed Algorithm \ref{alg:PGD} after a logarithmic number of iterations. 
\begin{Corollary}\label{coro:upper-bound}
Suppose the conditions of Theorem \ref{thm:local-convergence-PGD} hold and $\alpha,\beta,\kappa$ are constants. Then after at most $T = \Omega\left(\log\left(\left\|\cX^{(0)}-\cX^*\right\|_\tF/\xi\right) \right)$ iterations and for a constant $C$ that only relies on $\alpha, \beta, \kappa>0$, we have
\begin{equation*}
	\left\|\cX^{(T)} - \cX^*\right\|_{\rm F}^2 \leq C\xi^2. 
\end{equation*}
\end{Corollary}

\begin{Remark}\label{rm:Interpretation_xi}
When $\nabla L(\cX^*) = 0$, i.e., there is no statistical noise or perturbation, we have $\xi=0$. In this case, Theorem \ref{thm:local-convergence-PGD} and Corollary \ref{coro:upper-bound} imply that the proposed algorithm converges to the true target parameter $\cX^*$ at a linear rate:
\begin{equation*}
	\left\|\cX^{(t)}-\cX^*\right\|_\tF^2 \leq C\left(1 - \frac{2\rho\alpha  \beta\eta_0}{\kappa^2}\right)^t \left\|\cX^{(0)}-\cX^*\right\|_\tF^2.
\end{equation*}
When $\nabla L(\cX^*)\neq 0$, we have $\xi > 0$ and $\cX^*$ is not an exact stationary point of the loss function $L$. Then the estimation error $\left\|\cX^{(t)}-\cX^*\right\|_\tF^2$ is naturally not expected to go to zero, which matches the upper bounds of Theorem \ref{thm:local-convergence-PGD} and Corollary \ref{coro:upper-bound}. 
In a statistical model where noise or perturbation is in presence, the upper bound on the estimation error can be determined by evaluating $\xi$ under the specific random environment and these bounds are often minimax-optimal. See Section \ref{sec:instances} for more detail. 
\end{Remark}
\begin{Remark}
    If $\cC_\cS$ and $\cC_k$ are unbounded domains, then $\cC$ is the set of all rank-$(r_1,r_2,r_3)$ tensors, $\cP_{\cC_\cS}$, $\cP_{\cC_k}$ are identity operators, and the proposed Algorithm \ref{alg:PGD} essentially becomes the vanilla gradient descent. When $\cC_\cS$ and $\cC_k$ are non-trivial convex subsets, the projection steps ensure that $\cX^{(t)} = \llbracket \cS^{(t)};\U_1^{(t)},\U_2^{(t)},\U_3^{(t)} \rrbracket \in \cC$ and the RCG condition can be applied throughout the iterations. In fact, we found that the projection steps can be omitted in many numerical cases even if $L$ does not satisfy the RCG condition for the full set of low-rank tensors, such as the forthcoming Poisson and binomial tensor PCA. See Sections \ref{sec:instances} and \ref{sec:numerics} for more discussions. 
\end{Remark}

\subsection{Proof Sketch of Main Results}\label{sec:proof_sketch}

We briefly discuss the idea for the proof of Theorem \ref{thm:local-convergence-PGD} here. The complete proof is provided in Appendix \ref{sec:proof-main}. A key step in our analysis is to establish an error contraction inequality to characterize the estimation error of $\cX^{(t+1)}$ based on the one of $\cX^{(t)}$. Since the proposed non-convex gradient descent is performed on $\cS^{(t)}, \U_1^{(t)}, \U_2^{(t)}, \U_3^{(t)}$ jointly in lieu of $\cX^{(t)}$ directly, it becomes technically difficult to develop a direct link between $\left\|\cX^{(t+1)} - \cX^*\right\|_\tF^2$ and $\left\|\cX^{(t)} - \cX^*\right\|_\tF^2$. To overcome this difficulty, a ``lifting'' scheme was proposed and widely used in the recent literature on low-rank asymmetric matrix optimization \citep{tu2016low, zhu2017global, park2018finding}: one can factorize any rank-$r$ matrix estimator $\A^{(t)}$ and the target matrix parameter $\A^*$ into $\A^{(t)} = \U^{(t)}(\V^{(t)})^\top, \A^* = \U^*(\V^*)^\top$, where $\U^{(t)}, \V^{(t)}$ (or $\U^*, \V^*$) both have $r$ columns and share the same singular values. Then, one can stack them into one matrix

\begin{equation*}
    \W^{(t)} = \begin{bmatrix}
        \U^{(t)} \\
        \V^{(t)} \end{bmatrix},\quad \W^{\ast} = \begin{bmatrix}
        \U^{\ast} \\
        \V^{\ast}
        \end{bmatrix}.
\end{equation*}
By establishing the equivalence between $\min_{\R\in \bbO_r}\left\|\W^{(t)} - \W^*\R\right\|_\tF^2$ and $\left\|\A^{(t)} - \A^*\right\|_\tF^2$, and analyzing on $\min_{\R\in \bbO_r}\left\|\W^{(t)} - \W^*\R\right\|_\tF^2$, a local convergence of $\A^{(t)}$ to $\A^*$ can be established. However, the ``lifting'' scheme is not applicable to the tensor problem here since $\U_1,\U_2,\U_3, \cS$ have distinct shapes and cannot be simply stacked together. To overcome this technical issue in the generalized tensor estimation problems, we propose to assess the following criterion: 
\small
\begin{equation}\label{eq:E^t}
\begin{split}
E^{(t)} = \min_{\substack{\R_k \in \bbO_{p_k,r_k}\\k = 1,2,3}} \left\{\sum_{k=1}^3\left \| \U_k^{(t)} - \U_k^* \R_k\right\|_\tF^2 + \left\| \cS^{(t)} -  \llbracket\cS^*;\R_1^\top,\R_2^\top,\R_3^\top \rrbracket\right\|_\tF^2\right\}.
\end{split}
\end{equation}
\normalsize
Intuitively, $E^{(t)}$ measures the difference between a pair of tensor components $(\cS^{*},\U_1,\U_2,\U_3)$ and $(\cS^{(t)},\U_1^{(t)},\U_2^{(t)},\U_3^{(t)})$ under rotation.
The introduction of $E^{(t)}$ enables a convenient error contraction analysis as being an additive form of tensor components. In particular, the following lemma exhibits that $E^{(t)}$ is equivalent to the estimation error $\|\cX^{(t)}-\cX^*\|_\tF^2$ under regularity conditions.
\begin{Lemma}[An informal version of Lemma \ref{lm-equivalent-criteria}]\label{lm:simple-equivalence}
    \begin{equation*}
        c E^{(t)} \leq \|\cX^{(t)} - \cX^*\|_\tF^2 + C\sum_{k=1}^3\left\|(\U_k^{(t)})^\top\U_k^{(t)} - \U_k^{*\top}\U_k^*\right\|_\tF^2 \leq CE^{(t)}
    \end{equation*}
	under the regularity conditions to be specified in Lemma \ref{lm-equivalent-criteria}.
\end{Lemma}
Note that there is no equivalence between $E^{(t)}$ and $\left\|\cX^{(t)}-\cX^*\right\|_\tF^2$ unless we force $\U_k$ and $\U_k^*$ have similar singular structures, and this is the reason why we introduce the regularizer term in \eqref{eq-loss-regularize} to keep $\U_k^{(t)}$ from being singular.

Based on Lemma \ref{lm:simple-equivalence}, the proof of Theorem \ref{thm:local-convergence-PGD} reduces to establishing an error contraction inequality between $E^{(t)}$ and $E^{(t+1)}$:
\begin{equation}\label{eq-sketch-E-improve}
E^{(t+1)} \leq (1-\gamma)E^{(t)} + C\xi^2
\end{equation}
for constants $0<\gamma<1$ and $C>0$. Define the best rotation matrices
\small
\begin{equation*}
	\begin{split}
	(\R_1^{(t)},\R_2^{(t)},\R_3^{(t)}) = \argmin_{\substack{\R_k \in \bbO_{p_k, r_k}\\k = 1,2,3}} \left\{\sum_{k=1}^3\left \| \U_k^{(t)} - \U_k^* \R_k\right\|_\tF^2 + \left\|\cS^{(t)} - \llbracket \cS^*; \R_1^\top, \R_2^\top, \R_3^\top\rrbracket\right\|_\tF^2\right\}.
	\end{split}
\end{equation*}
\normalsize
By plugging in the gradient of $L(\cX)$ and $\cX = \cS\times_1 \U_1\times_2\U_2\times_3\U_3$, we can show
\small
\begin{equation}\label{eq-sketch-term-improve-U}
    \begin{split}
        &\left\|\U_k^{(t+1)} - \U_k^*\R_{k}^{(t+1)}\right\|_\tF^2 - \left\|\U_k^{(t)}-\U_k^*\R_{k}^{(t)}\right\|_\tF^2 \\
        \approx &  - 2\eta  \left\langle \cX^{(t)}-\cX_k^{(t)}, \nabla L(\cX^{(t)}) \right\rangle - \frac{a\eta}{2}\left\|\U_k^{(t)}\U_k^{(t)}-b^2 \I_{r_k}\right\|_\tF^2,
    \end{split}
\end{equation}
\normalsize
\small
\begin{equation}\label{eq-sketch-term-improve-S}
    \begin{split}
        & \left\|\cS^{(t+1)} - \llbracket\cS^*; \R_{1}^{(t+1)\top}, \R_2^{(t+1)\top}, \R_3^{(t+1)\top}\rrbracket\right\|_\tF^2 - \left\|\cS^{(t)} - \llbracket\cS^*; \R_{1}^{(t)\top}, \R_2^{(t)\top}, \R_3^{(t)\top}\rrbracket\right\|_\tF^2\\
        \approx &  - 2\eta\left\langle \cX^{(t)}-\cX_\cS^{(t)}, \nabla L(\cX^{(t)})  \right\rangle,
    \end{split}
\end{equation}
\normalsize
where 
$$\cX_k^{(t)} := \cS^{(t)}\times_k (\U_k^{*}\R_k^{(t)}) \times_{k+1} \U_{k+1}^{(t)} \times_{k+2} \U_{k+2}^{(t)}, \quad k=1,2,3;$$ 
$$\text{and} \quad \cX_\cS^{(t)}= \left\llbracket \cS^*;  \U_1^{(t)}\R_1^{(t)^\top}, \U_2^{(t)} \R_2^{(t)\top},  \U_3^{(t)}\R_3^{(t)\top}\right\rrbracket.$$
Note $E^{(t+1)}-E^{(t)}$ corresponds to the summation of \eqref{eq-sketch-term-improve-U} and \eqref{eq-sketch-term-improve-S}, whose right hand sides are dominated by the inner product between $\nabla L(\cX^{(t)})$ and $3\cX^{(t)}-\sum_{k=1}^3\cX_k^{(t)} - \cX_\cS^{(t)}$. We develop a new tensor perturbation Lemma to characterize $3\cX^{(t)}-\sum_{k=1}^3\cX_k^{(t)} - \cX_\cS^{(t)}$.
\begin{Lemma}[An informal version of Lemma \ref{lm-X-decomposition}]\label{lm-simple-decomposition}
    Under regularity conditions to be specified in Lemma \ref{lm-X-decomposition}, we have
    \begin{equation}\label{eq-sketch-decomposition}
        \cX^{(t)} - \cX^* = \left(\cX^{(t)} - \cX_\cS^{(t)}\right) + \sum_{k=1}^3 \left(\cX^{(t)} - \cX_k^{(t)}\right) + \cH_\varepsilon,
    \end{equation}
    where $\cH_\varepsilon$ is some low-rank residual tensor with $\|\cH_\varepsilon\|_\tF^2 = o\left((E^{(t)})^2\right)$. 
\end{Lemma}
Combining \eqref{eq-sketch-term-improve-U}\eqref{eq-sketch-term-improve-S} and Lemma \ref{lm-simple-decomposition}, we can connect $E^{(t+1)}$ and $E^{(t)}$ as
\small
\begin{equation}\label{eq-sketch-E-improve-1}
    \begin{split}
        E^{(t+1)} & \approx  E^{(t)} - 2\eta \left\langle \cX^{(t)} - \cX^* - \cH_\varepsilon, \nabla L(\cX^{(t)}) \right\rangle  - \frac{a\eta}{2}\sum_{k=1}^3 \left\|\U_k^{(t)\top}\U_k^{(t)} - \U_k^{*\top}\U_k^*\right\|_\tF^2.
    \end{split}
\end{equation}
\normalsize
Then, we introduce another decomposition
\small
\begin{equation}\label{eq-sketch-bound-inner-decompose}
\begin{split}
\MoveEqLeft\Big\langle  \cX^{(t)}  - \cX^* -  \cH_\varepsilon, \nabla L(\cX^{(t)}) \Big\rangle  = \underbrace{\left\langle \cX^{(t)} - \cX^*, \nabla L(\cX^{(t)}) - \nabla L(\cX^*) \right\rangle}_{A_1} \\
    &  - \underbrace{\left\langle \cH_\varepsilon, \nabla L(\cX^{(t)}) - \nabla L(\cX^*) \right\rangle}_{A_2} + \underbrace{\left\langle \cX^{(t)} - \cX^* +  \cH_\varepsilon, \nabla L(\cX^{*}) \right\rangle}_{A_3}
    \end{split}
\end{equation}
\normalsize
The three terms can be bounded separately:
\begin{equation}\label{ineq:sketch-separate-bound}
\begin{split}
    A_1 & \geq \alpha\|\cX^{(t)}-\cX^*\|_\tF^2 + \beta\left\|\nabla L(\cX^{(t)})-\nabla L(\cX^*)\right\|_\tF^2,\\
    |A_2| & \leq \frac{\beta}{2}\left\|\nabla L(\cX^{(t)})-\nabla L(\cX^*)\right\|_\tF^2 + \frac{1}{2\beta}\|\cH_\varepsilon\|_\tF^2, \\
    |A_3| & \leq \sqrt{E^{(t)}}\cdot \xi \leq cE^{(t)} + C\xi^2.
\end{split}
\end{equation}
Here the first inequality comes from RCG condition; the second inequality comes from Cauchy-Schwarz inequality and the fact $ab \leq \frac{1}{2}(a^2+b^2)$; and the last inequality utilizes the definition of $\xi$, as well as Lemma \ref{lm-simple-decomposition}. Combining \eqref{eq-sketch-bound-inner-decompose} and \eqref{ineq:sketch-separate-bound}, we obtain
\small
\begin{equation}\label{eq-sketch-bound-inner-product}
\begin{split}
\MoveEqLeft \Big\langle  \cX^{(t)}  - \cX^* -  \cH_\varepsilon, \nabla L(\cX^{(t)}) \Big\rangle + \left\langle \cX^{(t)} - \cX^* +  \cH_\varepsilon, \nabla L(\cX^{*}) \right\rangle \\
        \geq & \, \left(\alpha\|\cX^{(t)}-\cX^*\|_\tF^2 + \beta\left\|\nabla L(\cX^{(t)})-\nabla L(\cX^*)\right\|_\tF^2\right) \\
        & \qquad - \frac{\beta}{2}\left\|\nabla L(\cX^{(t)})-\nabla L(\cX^*)\right\|_\tF^2 - \left(cE^{(t)} + C\xi^2\right).
    \end{split}
\end{equation}
\normalsize
Then by choosing a suitable step size $\eta$ and applying \eqref{eq-sketch-E-improve-1} together with \eqref{eq-sketch-bound-inner-product}, one obtains
\begin{equation*}
    \begin{split}
        E^{(t+1)} & \leq E^{(t)}  + c E^{(t)} + C\xi^2  - c_2\left\|\cX^{(t)}-\cX^*\right\|_\tF^2 - c_3\sum_{k=1}^3 \left\|\U_k^{(t)\top}\U_k^{(t)} - \U_k^{*\top}\U_k^*\right\|_\tF^2 .
    \end{split}
\end{equation*}
Applying the equivalence between $\left\|\cX^{(t)} - \cX^*\right\|_\tF^2+C\sum_k\|(\U_k^{(t)})^\top\U_k^{(t)}-\U_k^{*\top}\U_k^*\|_\tF^2$ and $E^{(t)}$ (Lemma \ref{lm:simple-equivalence}), we can obtain \eqref{eq-sketch-E-improve} and finish the proof of Theorem \ref{thm:local-convergence-PGD}.

\section{Applications of Generalized Tensor Estimation}
\label{sec:instances}

Next, we apply the deterministic result to a number of generalized tensor estimation problems, including  sub-Gaussian tensor PCA, tensor regression, Poisson tensor PCA, and binomial tensor PCA to obtain the estimation error bound of (projected) gradient descent. In each case, Algorithm \ref{alg:PGD} is used with different initialization schemes specified by the problem settings. All the proofs are provided in Appendix \ref{sec:proof-others}. In addition, the generalized tensor estimation framework covers many other problems. A non-exhaustive list is provided in the introduction. See Section \ref{sec:discuss} for more discussions. 

\subsection{Sub-Gaussian Tensor PCA}\label{sec:sub-Gaussian}

Suppose we observe $\cY\in\mathbb{R}^{p_1\times p_2\times p_3}$, where $\mathbb{E}\cY = \cX^*$, $\cX^*$ is Tucker low-rank, and $\{\cY_{ijk}-\cX^*_{ijk}\}_{ijk}$ are independent and sub-Gaussian distributed. In literature, much attention has been focused on various setups related to this model, e.g., $\cY_{ijk}-\cX^*_{ijk}$ are i.i.d. Gaussian, $\cX^*$ is sparse, symmetric, rank-1, or CP-low-rank, etc \citep{richard2014statistical,sun2015guaranteed,perry2016statistical,montanari2017limitation,lesieur2017statistical,zhang2017optimal-statsvd,chen2019phase}. Particularly when $\{\cY_{ijk}-\cX^*_{ijk}\}_{ijk}$ are i.i.d. Gaussian distributed, it has been shown that the higher-order orthogonal iteration (HOOI) \citep{de2000best} achieves the optimal statistical performance on the estimation of $\cX^*$ \citep{zhang2018tensor}. It is however unclear whether HOOI works in the more general heteroskedastic setting, where the entries of $\cY$ have different variances.

Departing from the existing methods, we consider the estimation of $\cX^*$ via minimizing the quadratic loss function $L(\cX)=\frac{1}{2}\left\|\cX - \cY\right\|_\tF^2$ using  gradient descent.
It is easy to check that $L$ satisfies $RCG(\frac{1}{2},\frac{1}{2}, \bbR^{p_1\times p_2 \times p_3})$, so the projection steps in Algorithm \ref{alg:PGD} can be skipped throughout the iterations. To accommodate possible heteroskedastic noise, we apply HeteroPCA \citep{zhang2018heteroskedastic}, an iterative algorithm for PCA when heteroskedastic noise appears instead of the regular PCA for initialization. (The implementation of HeteroPCA in Algorithm \ref{alg:TD_Gaussian} is provided in Appendix \ref{sec:implementaionts}).

\begin{algorithm}
\caption{Initialization for Sub-Gaussian Tensor PCA}
\label{alg:TD_Gaussian}
\begin{algorithmic}
\REQUIRE $\cY \in \bbR^{p_1\times p_2 \times p_3}$, Tucker rank $(r_1,r_2,r_3)$, scaling parameter $b$.
\STATE{$\tilde{\cX} = \cY, \quad\tilde \U_k = \HeteroPCA_{r_k}\left(\cM_k(\tilde\cX)\cM_k(\tilde\cX)^\top\right)$\quad for $k=1,2,3$}
\STATE{$\tilde{\cS} = \llbracket\tilde \cX; \tilde \U_1^\top, \tilde \U_2^\top, \tilde \U_3 ^\top\rrbracket$}
\STATE{$\cS^{(0)} = \tilde{\cS} / b^3$ \quad $\U_k^{(0)} = b \tilde\U_k^{(0)}$\quad for $k=1,2,3$}
\RETURN {$(\cS^{(0)}, \U_1^{(0)}, \U_2^{(0)}, \U_3^{(0)})$}
\end{algorithmic}
\end{algorithm}
Now we can establish the theoretical guarantee of Algorithms \ref{alg:PGD} and \ref{alg:TD_Gaussian}  for sub-Gaussian tensor PCA based on the deterministic result in Theorem \ref{thm:local-convergence-PGD}.
\begin{Theorem}\label{thm:Gaussian}
Suppose we observe $\cY \in\mathbb{R}^{p_1\times p_2\times p_3}$, where $\mathbb{E}\cY = \cX^* = \llbracket \cS; \U_1, \U_2, \U_3\rrbracket$ is Tucker rank-$(r_1, r_2, r_3)$, $\U_k \in \bbO_{p_k,r_k}$ and $\left\|\U_k\right\|_{2,\infty} \leq c$ for some constant $c>0$. Suppose all entries of $\cY-\cX^*$ are independent mean-zero sub-Gaussian random variables such that
\begin{equation*}
    \|\cY_{ijk} - \cX_{ijk}\|_{\psi_2} = \sup_{q\geq 1} \mathbb{E}\left(|\cZ_{ijk}|^q\right)^{1/q}/q^{1/2} \leq \sigma.
\end{equation*}
Assume $\underline \lambda/\sigma \geq C_1 \overline p^{3/4}\overline r^{1/4}$. Then with probability at least $1-\exp(c\overline p)$, Algorithms \ref{alg:PGD} and \ref{alg:TD_Gaussian} yield
\begin{equation}\label{eq-GTD-upper-bound}
	\left\|\hat{\cX} - \cX^*\right\|_\tF^2 \leq C_2\sigma^2\left(r_1r_2r_3 + \sum_{k=1}^3 p_kr_k\right),
\end{equation}
where $C_1$, $C_2$ are constants that do not depend on $p_k$ and $r_k$. 
\end{Theorem}
\begin{Remark}
The proposed method turns out to achieve the minimax optimal rate of estimation error in a general class of sub-Gaussian tensor PCA settings since the order of upper bound \eqref{eq-GTD-upper-bound} matches the lower bound in literature \citep[Theorem 3]{zhang2018tensor}. Moreover, the condition $\underline\lambda /\sigma \geq C_1\overline p^{3/4}\overline r^{1/4}$ is optimal w.r.t. $\bar p$ in the sense that all the polynomial-time feasible algorithms cannot achieve consistent estimation when $\underline\lambda/\sigma < \bar p^{3/4-\varepsilon}$ for any $\varepsilon >0$ \citep[Theorem 4]{zhang2018tensor}.
\end{Remark}

\subsection{Low-rank Tensor Regression}\label{sec:tensor_reg}
Motivated by applications of neuroimaging analysis \citep{zhou2013tensor,li2017parsimonious,guhaniyogi2017bayesian}, spatio-temporal forecasting \citep{bahadori2014fast}, high-order interaction pursuit \citep{hao2018sparse}, longitudinal relational data analysis \citep{hoff2015multilinear}, 3D imaging processing \citep{guo2012tensor}, among many others, we consider the low-rank tensor regression next. Suppose we observe a collection of data $D=\{y_i, \cA_i\}_{i=1}^n$ that are associated through the following equation:
\begin{equation}\label{eq:regression-model}
    y_i = \left \langle \cA_i, \cX^\ast \right\rangle + \varepsilon_i, \quad \varepsilon_i\overset{iid}{\sim} N(0,\sigma^2), \quad i=1,\ldots,n.
\end{equation}
By exploiting the negative log-likelihood, it is natural to set $L$ to be the squared loss function
\begin{equation*}
    L\left(\cX;D\right) = \sum_{i=1}^n \left(\left\langle \cA_i, \cX \right\rangle - y_i\right)^2.
\end{equation*}
To estimate $\cX^*$, we first perform spectral method (Algorithm \ref{alg:regression}) to obtain initializer $\cX^{(0)} = \llbracket \cS^{(0)}, \U_1^{(0)}, \U_2^{(0)}, \U_3^{(0)}\rrbracket$, then perform the gradient descent (Algorithm \ref{alg:PGD}) without the projection steps to obtain the final estimator $\hat\cX$. A key step of Algorithm \ref{alg:regression} is HOSVD or HOOI, which are described in detail in Appendix \ref{sec:implementaionts}.

\begin{algorithm}
\caption{Initialization of Low-rank Tensor Regression}
\label{alg:regression}
\begin{algorithmic}
\REQUIRE $\left\{\cA_i,y_i\right\}$, $i=1,\ldots,n$, rank $(r_1,r_2,r_3)$, scaling parameter $b$.
\STATE{$\tilde \cX = \frac{1}{n}\sum y_i \cA_i$}
\STATE{$(\tilde \cS, \tilde \U_1, \tilde \U_2, \tilde \U_3) = \text{HOSVD}\left(\tilde \cX\right)$ or $(\tilde \cS, \tilde \U_1, \tilde \U_2, \tilde \U_3) = \text{HOOI}\left(\tilde \cX\right)$}
\STATE{$ \U_k^{(0)} = b \tilde\U_k$,\quad for $k=1,2,3$}
\STATE{$\cS^{(0)} = \tilde{\cS} / b^3$}
\RETURN {$(\cS^{(0)}, \U_1^{(0)}, \U_2^{(0)}, \U_3^{(0)})$}
\end{algorithmic}
\end{algorithm}

For technical convenience, we assume the covariates $\{\cA_i\}_{i=1}^n$ are randomly designed that all entries of $\cA_i$ are i.i.d. drawn from sub-Gaussian distribution with mean 0 and variance 1. The following theorem gives an estimation error upper bound for Algorithms \ref{alg:PGD} and \ref{alg:regression}.
\begin{Theorem}\label{thm:regression-sharp}
    Consider the low-rank tensor regression model \eqref{eq:regression-model}. Suppose $\sigma^2 \leq C_1\left\|\cX^*\right\|_\tF^2$, $\underline\lambda \geq C_2$, and the sample size $n \geq C_3\max\{\overline p^{3/2}\overline r, \overline p \cdot \overline r^{2}, \overline r^4\}$ for constants $C_1, C_2, C_3 > 0$. Then with probability at least $1-\exp\left\{-c(r_1r_2r_3+\sum_{k=1}^3 p_kr_k)\right\}$, the output of Algorithms \ref{alg:PGD} and \ref{alg:regression} satisfies
    \begin{equation*}
        \left\|\hat \cX - \cX^*\right\|_\tF^2 \leq C\sigma^2\left(r_1r_2r_3 + \sum_{k=1}^3 p_kr_k\right)/n,
    \end{equation*}
where $C_1,C_2,C_3,C$ are constants depending only on $\kappa$ and $c>0$ is a universal constant.
\end{Theorem}
Theorem \ref{thm:regression-sharp} together with the lower bound in \cite[Theorem 5]{zhang2018ISLET} shows that the proposed procedure achieves the minimax optimal rate of estimation error in the class of all $p_1$-by-$p_2$-by-$p_3$ tensors with rank-$(r_1, r_2, r_3)$ for tensor regression.

\begin{Remark}
    The assumption on the covariates $\left\{\cA_i\right\}_{i=1}^n$ in Theorem \ref{thm:regression-sharp} ensures that $\tilde \cX$ in Algorithm \ref{alg:regression} is an unbiased estimator of $\cX^*$. Such a setting has been considered as a benchmark setting in the high-dimensional statistical inference literature (see, e.g., \cite{candes2011tight,chen2019non,javanmard2018debiasing}). When $\left\{\cA_i\right\}_{i=1}^n$ are heteroskedastic, the spectral initialization may fail and an alternative idea is the following unfolded nuclear norm minimization:
    \begin{equation}\label{eq:SDP}
        \tilde \cX' = \argmin_{\cX \in \bbR^{p_1\times p_2 \times p_3}} \sum_{k=1}^3 \left\|\cM_k(\cX)\right\|_*,\quad s.t.~~\cA(\cX) = y.
    \end{equation}
    \eqref{eq:SDP} is equivalent to a semidefinite programming and can be solved by the interior-point method \citep{gandy2011tensor}.
\end{Remark}

\begin{Remark}
    There is a significant gap between the required sample size in Theorem \ref{thm:regression-sharp} ($O(p^{3/2}r)$) and the possible sample size lower bound, i.e., the degree of freedom of all rank-$r$ dimension-$p$ tensor $(O\left(r^3+pr\right))$. The existing algorithms achieving the sample size lower bound are often NP-hard to compute and thus intractable in practice. 
    We also note that the existence of a tractable algorithm for tensor completion that provably works with less than $p^{3/2-\varepsilon}$ measurements would disapprove an open conjecture in theoretical computer science on strongly random 3-SAT \citep[Corollary 16]{barak2016noisy}. Since tensor completion can be seen as a special case of tensor recovery, this suggests that it may be impossible to substantially improve the sample complexity required in Theorem \ref{thm:regression-sharp} using a polynomial-time algorithm. Therefore, our procedure can be taken as the first computationally efficient algorithm to achieve minimax optimal rate of convergence and exact recovery in the noiseless setting as illustrated in Table \ref{tab:reg_comparison}. 
\end{Remark}

\subsection{Poisson Tensor PCA}\label{sec:poisson}

Tensor data with count values commonly arise from various scientific applications, such as the photon-limited imaging \citep{timmerman1999multiscale,willett2007multiscale,salmon2014poisson,yankovich2016non}, online click-through data analysis \citep{shan2016predicting,sun2016sparse}, and metagenomic sequencing \citep{flores2014temporal}. In this section, we consider the Poisson tensor PCA model: assume we observe $\cY \in \mathbb{N}^{p_1\times p_2 \times p_3}$ that satisfies
\begin{equation}\label{eq:Poisson-independently}
\cY_{ijk} \sim \text{Poisson}(I\exp(\cX^*_{ijk}))\quad \text{independently},
\end{equation}
where $\cX^*$ is the low-rank tensor parameter and $I>0$ is the intensity parameter. When $\cX^*$ is entry-wise bounded (Assumption \ref{asmp-incoherence}), one can set $I$ as the average intensity of all entries of $\cY$ so that $I$ essentially quantifies the signal-to-noise ratio. Rather than estimating $I\exp(\cX^*)$, we focus on estimating $\cX^*$, the key tensor that captures the salient geometry or structure of the data.

 Then, the following negative log-likelihood is a natural choice of the loss function for estimating $\cX^*$,
\begin{equation}\label{eq:L-Poisson-likelihood}
L(\cX) = \sum_{i=1}^{p_1}\sum_{j=1}^{p_2}\sum_{k=1}^{p_3}\left(-\cY_{ijk}\cX_{ijk} + I\exp(\cX_{ijk})\right).
\end{equation} 
Unfortunately, $L(\cX)$ defined in \eqref{eq:L-Poisson-likelihood} satisfies $RCG(\alpha,\beta, \cC)$ only for a bounded set $\cC$ since the Poisson likelihood function is not strongly convex and smooth in the unbounded domain. We thus introduce the following assumption on $\cX^*$ to ensure that $\cX^*$ is in a bounded set $\cC$.
\begin{Assumption}\label{asmp-incoherence}
\allowbreak
Suppose $\cX^* = \llbracket\cS^*; \U_1^*, \U_2^*, \U_3^*\rrbracket$, where $\U_k^* \in \bbO_{p_k,r_k}$ is a $p_k$-by-$r_k$ orthogonal matrix for $k=1,2,3$. There exist some constants $\{\mu_k\}_{k=1}^3, B$ such that $\left\|\U_k^*\right\|_{2,\infty}^2 \leq \frac{\mu_k r_k}{p_k}$ for $k=1,2,3$ and $\overline \lambda \leq B\sqrt{\frac{\Pi_{k=1}^3p_k}{\Pi_{k=1}^3\mu_kr_k}}$ where $\overline \lambda:= \max_k \left\|\cM_k(\cS^*)\right\|$. Here, $\|\U_k^*\|_{2,\infty} = \max_i \|(\U_k^*)_{i\cdot}\|_2$ is the largest row-wise $\ell_2$ norm of $\U_k^*$.
\end{Assumption}
Assumption \ref{asmp-incoherence} requires that the loading $\U_k$ satisfies the incoherence condition, i.e., the amplitude of the tensor is ``balanced'' in all parts. Previously, the incoherence condition and its variations were commonly used in the matrix estimation literature \citep{candes2009exact,ma2017exploration} and  Poisson-type inverse problems (e.g., Poisson sparse regression \cite[Assumption 2.1]{jiang2015minimax}, Poisson matrix completion \cite[Equation (10)]{cao2015poisson}, compositional matrix estimation \cite[Equation (7)]{cao2017microbial}, Poisson auto-regressive models \citep{hall2016inference}). Assumption \ref{asmp-incoherence} also requires an upper bound on the spectral norm of each matricization of the core tensor $\cS^*$. Together with the incoherence condition on $\U_k^*$, this condition guarantees that 
$\cX^*$ is entry-wise upper bounded by $B$. In fact, the entry-wise bounded assumption is also widely used in high-dimensional matrix/tensor generalized linear models since it guarantees the local strong convexity and smoothness of the negative log-likelihood function \citep{ma2017exploration, wang2018learning,xu2019generalized}.

Next, we set $(\{\cC_k\}_{k=1}^3, \cC_\cS)$ as follows:
\begin{equation}\label{eq-specify-sets}
    \begin{split}
        \cC_k &= \left\{\U_k \in \bbR^{p_k\times r_k}: \left\|\U_k\right\|_{2,\infty} \leq b\sqrt{\frac{\mu_k r_k}{p_k}}\right\}, \\
        \cC_\cS &= \left\{\cS \in \bbR^{r_1\times r_2 \times r_3}: \max_k \left\|\cM_k(\cS)\right\| \leq b^{-3}B\sqrt{\frac{\Pi_{k=1}^3p_k}{\Pi_{k=1}^3\mu_kr_k}}\right\}.
    \end{split}
\end{equation}
Specifically for the Poisson tensor PCA, we can prove that if Assumption \ref{asmp-incoherence} holds, the loss function \eqref{eq:L-Poisson-likelihood} satisfies $RCG(\alpha,\beta,\left\{\cC_k\right\}_{k=1}^3, \cC_\cS)$ for constants $\alpha, \beta$ that only depend on $I$ and $B$ (see the proof of Theorem \ref{thm:Poisson} for details). We can also show that the following Algorithm \ref{alg:TD_Poisson} provides a sufficiently good initialization with high probability. 
\begin{algorithm}
\caption{Initialization for Poisson Tensor PCA}
\label{alg:TD_Poisson}
\begin{algorithmic}
\REQUIRE Initialization observation tensor $\cY \in \bbN^{p_1\times p_2 \times p_3}$, Tucker rank $(r_1,r_2,r_3)$, scaling parameter $b$, intensity parameter $I$.
\STATE{$\tilde \cX = \log\left((\cY_{jkl}+\frac{1}{2})/I\right)$}
\STATE{$(\tilde \cS, \tilde \U_1, \tilde \U_2, \tilde \U_3) = \text{HOSVD}\left(\tilde \cX\right)$ or $(\tilde \cS, \tilde \U_1, \tilde \U_2, \tilde \U_3) = \text{HOOI}\left(\tilde \cX\right)$}
\STATE{$ \U_k^{(0)} = b \tilde\U_k$,\quad for $k=1,2,3$}
\STATE{$\cS^{(0)} = \tilde{\cS} / b^3$}
\RETURN {$(\cS^{(0)}, \U_1^{(0)}, \U_2^{(0)}, \U_3^{(0)})$}
\end{algorithmic}
\end{algorithm}

Now we establish the estimation error upper bound for Algorithms \ref{alg:PGD} and \ref{alg:TD_Poisson}.
\begin{Theorem}\label{thm:Poisson}
\allowbreak
Suppose Assumption \ref{asmp-incoherence} holds and $I > C_1 \max\{\bar p, \underline{\lambda}^{-2} \sum_{k=1}^3(p_{-k}r_k + p_kr_k)\}$, where $p_{-k}:=p_1p_2p_3/p_k$. Then with probability at least $1-c/(p_1p_2p_3)$, the output of Algorithms \ref{alg:PGD} and \ref{alg:TD_Poisson} yields 
\begin{equation*}
	\left\|\hat \cX-\cX^*\right\|_\tF^2 \leq C_2I^{-1}\left(r_1r_2r_3 + \sum_{k=1}^3 p_kr_k\right).
\end{equation*}
Here $C_1,C_2$ are constants that do not depend on $p_k$ or $r_k$.
\end{Theorem}

We further consider the following class of low-rank tensors $\cF_{\p, \r}$, where the restrictions in $\cF_{\p, \r}$ correspond to the conditions in Theorem \ref{thm:Poisson}:
\small
\begin{equation}\label{eq:F_p,r}
\begin{split}
& \cF_{\p,\r} =
\left\{\cX = \llbracket \cS; \U_1,\U_2,\U_3\rrbracket: \begin{array}{l}
\U_k \in \bbO_{p_k,r_k},\quad \left\|\U_k\right\|_{2,\infty}^2 \leq \frac{\mu_k r_k}{p_k},\\ \max \limits_{k}\left\|\cM_k(\cS)\right\| \leq B\sqrt{\frac{\Pi_{k=1}^3p_k}{\Pi_{k=1}^3\mu_kr_k}}
\end{array}\right\}.
\end{split}
\end{equation}
\normalsize

With some technical conditions on tensor rank and the intensity parameter, we can develop the following lower bound in estimation error for Poisson PCA.
\begin{Theorem}[Lower Bound for Poisson tensor PCA]\label{thm:lower-bound-poisson}
	Assume $\bar{r} \leq C_1\underline{p}^{1/2}$, $\underline r> C_2$ and $\min_k\mu_k \geq C_3$ for constants $C_1,C_2,C_3>1$. Suppose one observes $\cY\in \mathbb{R}^{p_1\times p_2\times p_3}$,  where $\cY_{jkl} \sim \text{Poisson}\left(I\exp(\cX_{jkl})\right)$ independently, $\cX \in \cF_{\p,\r}$, and $I\geq c_0$. There exists a uniform constant $c$ that does not depend on $p_k$ or $r_k$, such that
	\begin{equation*}
	\begin{split}
	& \inf_{\hat{\cX}} \sup_{\cX \in \mathcal{F}_{\p,\r}} \bbE\left\|\hat{\cX} - \cX \right\|_\tF^2 \geq cI^{-1}\left( r_1r_2r_3 + \sum_{k=1}^3p_kr_k\right).
	\end{split}
	\end{equation*}
\end{Theorem}
Theorems \ref{thm:Poisson} and \ref{thm:lower-bound-poisson} together yield the optimal rate of estimation error for Poisson tensor PCA problem over the class of $\mathcal{F}_{\p, \r}$:
\begin{equation*}
    \inf_{\hat\cX}\sup_{\cX\in \mathcal{F}_{\p, \r}}\bbE \left\|\hat \cX - \cX\right\|_\tF^2 \asymp I^{-1}\left(r_1r_2r_3+\sum_{k=1}^3 p_kr_k\right).
\end{equation*}

\subsection{Binomial Tensor PCA}\label{sec:binomial}

The binomial tensor data commonly arise in the analysis of proportion when raw counts are available. For example, in the Human Mortality Database \citep{wilmoth2016human}, the number of deaths and the total number of population are summarized into a three-way tensor, where the $x$-, $y$-, $z$-coordinates are counties, ages, and years, respectively. Given the sufficiently large number of population in each country, one can generally assume that each entry of this data tensor satisfies the binomial distribution independently. 

Suppose we observe a count tensor $\cY \in \bbN^{p_1\times p_2 \times p_3}$ and a total population tensor $\cN \in \bbN^{p_1\times p_2 \times p_3}$ such that $\cY_{jkl} \sim \text{Binomial}(\cN_{jkl}, \cP_{jkl}^*)$ independently. Here, $\cP^* \in [0,1]^{p_1 \times p_2 \times p_3}$ is a probability tensor linked to an underlying latent parameter $\cX^* \in \bbR^{p_1 \times p_2 \times p_3}$ through $\cP^*_{jkl} = s(\cX_{jkl}^*)$, where $s(x)= 1/(1+e^{-x})$ is the sigmoid function. Our goal is to estimate $\cX^*$. 
To this end, we consider to minimize the following loss function: 
\begin{equation*}
L(\cX) = -\sum_{jkl}\left(\hat\cP_{jkl}\cX_{jkl} + \log\left(1-\sigma(\cX_{jkl})\right)\right),
\end{equation*}
where $\hat \cP_{jkl} := \cY_{jkl}/\cN_{jkl}$.

We assume $\cX^*$ satisfies Assumption \ref{asmp-incoherence} for the same reasons as in Poisson tensor PCA. We propose to estimate $\cP^*$ by applying Algorithm \ref{alg:TD_Binomial} (initialization) and Algorithm \ref{alg:PGD} (projected 
gradient descent) with the following constraint sets:
\begin{equation*}
\begin{split}
\cC_k &= \left\{\U_k \in \bbR^{p_k\times r_k}: \left\|\U_k\right\|_{2,\infty} \leq b\sqrt{\frac{\mu_k r_k}{p_k}}\right\}, \\
\cC_\cS &= \left\{\cS \in \bbR^{r_1\times r_2 \times r_3}: \max_k \left\|\cM_k(\cS)\right\| \leq b^{-3}B\sqrt{\frac{\Pi_{k=1}^3p_k}{\Pi_{k=1}^3\mu_kr_k}}\right\}.
\end{split}
\end{equation*} 
\begin{algorithm}
\caption{Initialization for Binomial Tensor PCA}
\label{alg:TD_Binomial}
\begin{algorithmic}
\REQUIRE $\cY,\cN \in \mathbb{N}^{p_1\times p_2 \times p_3}$, Tucker rank $(r_1,r_2,r_3)$, scaling parameter $b$
\STATE{$\tilde \cX_{jkl} = \log\left(\frac{\cY_{jkl} + 1/2}{\cN_{jkl} - \cY_{jkl} + 1/2}\right),~\forall j,k,l$}
\STATE{$(\tilde \cS, \tilde \U_1, \tilde \U_2, \tilde \U_3) = \text{HOSVD}\left(\tilde \cX\right)$ or $(\tilde \cS, \tilde \U_1, \tilde \U_2, \tilde \U_3) = \text{HOOI}\left(\tilde \cX\right)$}
\STATE{$ \U_k^{(0)} = b \tilde\U_k$,\quad for $k=1,2,3$}
\STATE{$\cS^{(0)} = \tilde{\cS} / b^3$}
\RETURN {$(\cS^{(0)}, \U_1^{(0)}, \U_2^{(0)}, \U_3^{(0)})$}
\end{algorithmic}
\end{algorithm}

We have the following theoretical guarantee for the estimator obtained by Algorithms \ref{alg:PGD} and \ref{alg:TD_Binomial} in binomial tensor PCA.
\begin{Theorem}[Upper Bound for Binomial Tensor PCA]\label{thm:Binomial}
Suppose Assumption \ref{asmp-incoherence} is satisfied and $N = \min_{jkl}\cN_{jkl}$ satisfies $N \geq C_1\max\left\{\overline p, \underline{\lambda}^{-2} \sum_k \left(p_{-k}r_k + p_kr_k\right)\right\}$. Then with probability at least $1-c/(p_1p_2p_3)$, we have the following estimation upper bound for the output of Algorithms \ref{alg:PGD} and \ref{alg:TD_Binomial}:
\begin{equation*}
	\left\|\hat \cX-\cX^*\right\|_\tF^2 \leq C_2N^{-1}\left(r_1r_2r_3+\sum_{k=1}^3p_kr_k\right).
\end{equation*}
Here, $C_1,C_2$ are some absolute constants that do not depend on $p_k$ or $r_k$.
\end{Theorem}
\begin{Remark}
We assume $N \geq C_1\max\left\{\overline p, \underline{\lambda}^{-2} \sum_k \left(p_{-k}r_k + p_kr_k\right)\right\}$ in Theorem \ref{thm:Binomial} as a technical condition to prove the estimation error upper bound of $\hat\cX$. $N$ here essentially characterizes the signal-noise ratio of the binomial tensor PCA problem.
\end{Remark}
Let $\mathcal{F}_{\p, \r}$ be the class of low-rank tensors defined in \eqref{eq:F_p,r}. We can prove the following lower bound result, which establishes the minimax optimality of the proposed procedure over the class of $\cF_{\p, \r}$ in binomial tensor PCA. 

\begin{Theorem}[Lower Bound for Binomial Tensor PCA]\label{thm:lower-bound-binomial}
Denote $N = \min_{jkl}\cN_{jkl}$. Assume $\bar{r} \leq C_1\underline{p}^{1/2}$, $\underline r> C_2$ and $\min_k\mu_k \geq C_3$ for some constants $C_1,C_2,C_3>1$. Suppose one observes $\cY\in \mathbb{N}^{p_1\times p_2\times p_3}$,  where $\cY_{jkl} \sim \text{Binomial}\left(\cN_{jkl}, \sigma(\cX_{jkl})\right)$ independently, $\cX^* \in \cF_{\p,\r}$, and $\max_{jkl} \cN_{jkl} \leq C\min_{jkl} \cN_{jkl}$. There exists constant $c$ that does not depend on $p_k$ or $r_k$, such that
\begin{equation*}
	\begin{split}
	&\inf_{\hat{\cX}} \sup_{\cX \in \mathcal{F}_{\p, \r}} \left\|\hat{\cX} - \cX \right\|_\tF^2 \geq c N^{-1} \left(r_1r_2r_3+\sum_{k=1}^3p_kr_k\right).
	\end{split}
\end{equation*}
\end{Theorem}

\section{Rank Selection}\label{sec:rank-estimation}

The tensor rank $(r_1,r_2,r_3)$ is required as an input to Algorithm \ref{alg:PGD} and plays a crucial role in the proposed non-convex optimization framework. While various empirical methods have been proposed for rank selection in specific applications of low-rank tensor estimation (e.g., \cite{yokota2016robust}), there is a paucity of theoretical guarantees in the literature. In this section, we  provide a rank estimation procedure with provable guarantees. Recall that in each application in Section \ref{sec:instances}, we first specify the initialization $\cX^{(0)} = \llbracket \cS^{(0)}; \U_1^{(0)}, \U_2^{(0)}, \U_3^{(0)}\rrbracket$ based on a spectral algorithm on some preliminary tensor $\tilde \cX$ (see their definitions in Algorithms \ref{alg:TD_Gaussian}-\ref{alg:TD_Binomial}). Since $\tilde \cX$ reflects the target tensor $\cX^*$ and $\sigma_s(\mathcal{M}_k(\cX^*))=0$ for $s\geq r_k+1$, we consider the following rank selection method by exploiting the singular values of $\tilde{\cX}$:
\begin{equation}\label{eq:rank-selection-rk}
    \hat r_k = \max\left\{r: \sigma_r\left(\cM_k(\tilde \cX)\right) \geq t_k \right\},\qquad k=1, 2, 3.
\end{equation}
Here, $t_k>0$ is the thresholding level whose value depends on specific problem settings. 
Next, we specifically consider the sub-Gaussian tensor PCA and tensor regression. 
\begin{Proposition}\label{prop:rank-selection}
    Suppose $\rank(\cX^*) = (r_1,r_2,r_3)$, $r_k = o(p_k)$, $p_k = o(p_{-k})$. Define $\hat r_k$ as \eqref{eq:rank-selection-rk} and
\begin{equation}\label{eq:rank-select-deltak}
    \delta_k := \text{Median}\left\{\sigma_1(\cM_k(\tilde\cX)),\ldots, \sigma_{p_k}\left(\cM_k(\tilde\cX)\right)\right\}.
\end{equation}
    \begin{enumerate}
        \item[(a)] In sub-Gaussian tensor PCA (Section \ref{sec:sub-Gaussian}), we set $\tilde \cX = \cY, t_k = 1.5\delta_k$. Suppose $\underline \lambda \geq C\overline p \sigma$. 
        Then, we have $\hat r_k = r_k, \forall k \in [d]$ with probability at least $1-Ce^{-c\bar p}$. 
        \item[(b)] In low-rank tensor regression (Section \ref{sec:tensor_reg}), we set $\tilde \cX = \frac{1}{n}\sum_{i=1}^n y_i\cA_i, t_k = 1.5\delta_k$. Suppose $n \geq C\kappa \bar p^2\bar r^{3/2}\log \bar r$ and $\underline \lambda \geq C\frac{\bar p \sigma}{\sqrt{n}}$. Then, $\hat r_k = r_k, \forall k \in [d]$ with probability at least $1-Ce^{-c\bar p}$. 
    \end{enumerate}
\end{Proposition}

In practice, we can also apply a simple criterion of the \emph{cumulative percentage of total variation} \cite[Chapter 6.1.1]{jolliffe1986principal} originating from  principle component analysis:
\begin{equation}\label{eq:rk-select-pca}
    \hat r_k = \argmin\left\{r: \sum_{i=1}^r \sigma_i^2(\cM_k(\tilde \cX)) \Big/  \sum_{i=1}^{p_k}\sigma_i^2(\cM_k(\tilde \cX)) \geq \rho\right\}
\end{equation}
Here, $\rho \in (0,1)$ is some empirical thresholding level. We will illustrate this principle on real data analysis in Section \ref{sec:real-data}. Under the general deterministic setting, the accurate (or optimal) estimation of tensor rank may be much more challenging and we leave it as future work.

\section{Extensions to General Order-$d$ Tensors}\label{sec:extension}

While our previous sections mainly focus on order-3 tensor estimation, our results can be generalized to the order-$d$ low-rank tensor estimation with the key ideas outline in this section. First, the constraint set $\cC$, RCG condition and noise quantity $\xi$ can be defined similarly by replacing the order-$3$ tensor with the general low-rank order-$d$ tensors; second, the local convergence analysis can be similarly conducted as Theorem \ref{thm:local-convergence-PGD}. Define
\begin{equation*}
E^{(t)} := \min_{\substack{\R_k \in \bbO_{p_k,r_k}\\k = 1,\ldots,d}} \left\{\sum_{k=1}^d\left \| \U_k^{(t)} - \U_k^* \R_k\right\|_\tF^2 + \left\| \cS^{(t)} -  \llbracket \cS^*; \R_1^\top,\ldots, \R_d^\top\rrbracket\right\|_\tF^2\right\}.
\end{equation*}
We can build the equivalence between $E^{(t)}$ and
\begin{equation*}
    \left\|\cX^{(t)} - \cX^*\right\|_\tF^2 + \frac{a}{2}\sum_{k=1}^d \left\|\U_k^{(t)\top}\U_k^{(t)} - b^2\I_{r_k}\right\|_{\tF}^2
\end{equation*}
by setting $b \asymp \overline \lambda^{1/(d+1)}$ and $a \asymp \overline \lambda$. Then, we can establish the following theoretical guarantee under a good initialization:
\begin{Theorem}[Informal]\label{thm:high-order-local-convergence}
Suppose $L$ satisfies RCG$(\alpha,\beta,\mathcal C)$ and assume $\kappa,\alpha,\beta$ are constants. Assume $\cX^* = \llbracket \cS^*;\U_1^*,\ldots,\U_d^*\rrbracket$ and the initialization $\cX^{(0)} = \llbracket \cS^{(0)};\U_1^{(0)},\ldots,U_d^{(0)}\rrbracket$ satisfy
\begin{equation}\label{eq:high-order-incoherence}
    \begin{split}
        \U_k^{*\top}\U_k^* = \U_k^{(0)\top}\U_k^{(0)} = b^2 \I_{r_k},~ \U_k^*, \U_k^{(0)} \in \cC_k, k = 1,\ldots, d;~~ \cS^*, \cS^{(k)} \in \cC_\cS.
    \end{split}
\end{equation}
Suppose the initialization error satisfies $\left\|\cX^{(0)} - \cX^*\right\| \leq c_d\underline \lambda^2$ and the signal-noise-ratio satisfies $\underline\lambda^2 \geq C_d\xi^2$. Then, by taking step size $\eta \leq c_d b^{-2d}$, the output of Algorithm \ref{thm:local-convergence-PGD} satisfies
\begin{equation*}
    \left\|\cX^{(t)} - \cX^*\right\|_\tF^2 \leq C_d'\left(\xi^2 + \left(1-c'_d\eta\right)^t \left\|\cX^{(0)}-\cX^*\right\|_\tF^2\right).
\end{equation*}
Here, $C_d, c_d, C_d', c_d'$ are constants only depending on $d$. 
\end{Theorem}
It is worth mentioning that Theorem \ref{thm:high-order-local-convergence} applies to low-rank matrix estimation (i.e. $d=2$): suppose $\X^* = \U_1^*\S^*\U_2^{*\top}$ and $\X^{(0)} = \U_1^{(0)}\S^{(0)}\U_2^{(0)\top}$  for $\U_k,\U_k^{(0)} \in \bbR^{p_k\times r}$, $\S^*,\S^{(0)} \in \bbR^{r\times r}$, we have
\begin{equation*}
    \left\|\X^{(t)} - \X^*\right\|_\tF^2 \leq C\left(\xi^2 +  (1-c\eta)^t\left\|\X^{(0)}-\X^*\right\|_\tF^2\right).
\end{equation*}
While this framework is more complicated than necessary since one can always decompose a low-rank matrix as the product of two factor matrices $\X = \U_1\U_2^\top$ without explicitly introducing the ``core matrix'' $\S \in \bbR^{r \times r}$ (see our previous discussions in Section \ref{sec:proof_sketch}).

Based on Theorem \ref{thm:high-order-local-convergence}, we can further extend the minimax optimal bounds for the proposed procedure in each application of Section \ref{sec:instances}, i.e., Theorems \ref{thm:Gaussian}--\ref{thm:lower-bound-binomial} to high-order scenarios. We summarize the results to Table \ref{tab:applications-oder-d}. 
\begin{table}
	\centering
	{\footnotesize
	\begin{tabular}{c|c|c|c}
		\hline
		{\centering Application}  & {\centering SNR condition} & {\centering Estimation error} & {\centering Lower bound}\\  \hline
        {\centering sub-Gaussian tensor-PCA} & {\centering $\underline \lambda/\sigma \gtrsim p^{d/4}r^{1/(d+1)}$} & {\centering $\sigma^2\left(pr + r^{d}\right)$} & {\centering $\sigma^2\left(pr + r^{d}\right)$} \\ \hline
        {\centering Tensor regression}  & {\centering $n \gtrsim p^{d/2}r$} & {\centering $n^{-1}\sigma^2\left(pr + r^{d}\right)$} & {\centering $n^{-1}\sigma^2\left(pr + r^{d}\right)$}  \\  \hline
        {\centering Poisson tensor-PCA}  & {\centering $\sqrt{I}\underline\lambda \gtrsim p^{(d-1)/2}r^{1/2}$} & {\centering $I^{-1}\left(pr + r^{d}\right)$} & {\centering $I^{-1}\left(pr + r^{d}\right)$} \\\hline
        {\centering Binomial tensor-PCA}  & {\centering $\sqrt{N}\underline\lambda \gtrsim p^{(d-1)/2}r^{1/2}$} & {\centering $N^{-1}\left(pr + r^{d}\right)$} & {\centering $N^{-1}\left(pr + r^{d}\right)$} \\\hline
	\end{tabular}
	}
	\caption{Minimax optimal estimation error bounds for order-$d$ low-rank tensor estimation in specific applications ($d\geq 2$). Here, for simplicity, $r_1=r_2=r_3=r$, $p_1=p_2=p_3=p$ and $r \leq p^{1/2}$.}
	\label{tab:applications-oder-d}
\end{table}

\section{Numerical Studies}\label{sec:numerics}

\subsection{Synthetic Data Analysis}\label{sec:synthetic-data}

In this section, we investigate the numerical performance of the proposed methods on the problems discussed in Section \ref{sec:instances} with simulated data. We assume the true rank $(r_1,r_2,r_3)$ is known to us and the algorithm only involves two tuning parameters: $a$ and $b$. According to Theorem \ref{thm:local-convergence-PGD}, a proper choice of $a$ and $b$ primarily depends on the unknown value $\overline \lambda$. In practice, we propose to use the initial estimate $\cX^{(0)}$ as an approximation of $\cX^*$, use $\bar\lambda^{(0)} = \max_k \left\|\cM_k(\cX^{(0)})\right\|$ as a plug-in estimate of $\overline \lambda$, then choose $a = \overline{\lambda}^{(0)}, b = (\overline{\lambda}^{(0)})^{1/4}.$ 
We consider the following root mean squared error (RMSE) to assess the estimation accuracy in all settings: 
\begin{equation}\label{eq-def-loss}
    \text{Loss}(\hat\cX,\cX^*) = (p_1p_2p_3)^{-1/2}\|\hat\cX - \cX^*\|_\tF.
\end{equation}
Average loss over 100 repetitions are reported in following different scenarios.

\ \par
\noindent{\bf Tensor Regression.} We investigate the numerical performance of the proposed procedure in low-rank tensor regression discussed in Section \ref{sec:tensor_reg}. For all simulation settings, we first generate an $r_1$-by-$r_2$-by-$r_3$ core tensor $\bar\cS$ with i.i.d. standard Gaussian entries and rescale it as $\cS = \bar \cS \cdot \lambda/\min_{k=1}^3\sigma_{r}\left(\cM_k(\bar\cS)\right)$. Here, $\lambda$ quantifies the signal level and will be specified later. Then we generate $\U_k$ uniformly at random from the Stiefel manifold $\bbO_{p_k, r_k}$ and calculate the true parameter as $\cX^* = \llbracket\cS; \U_1, \U_2, \U_3\rrbracket$. The rescaling procedure here ensures that $\min_{k} \sigma_{r_k}\left(\cM_k(\cX^*)\right) \geq \lambda$. Now, we draw a random sample based on the regression model \eqref{eq:regression-model}. 

\begin{figure} \centering
	\subfigure{
		\includegraphics[width =0.45\linewidth,height=2.0in]{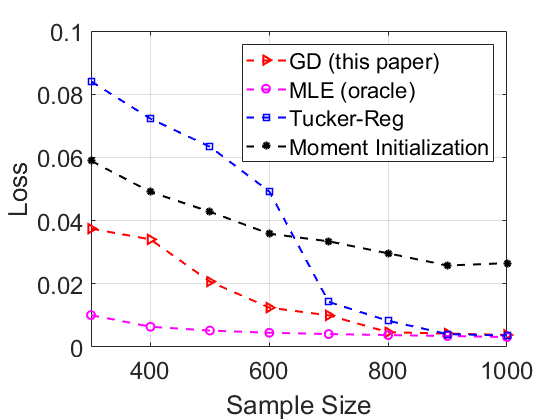}
	}
	\subfigure{
		\includegraphics[width =0.45\linewidth,height=2.0in]{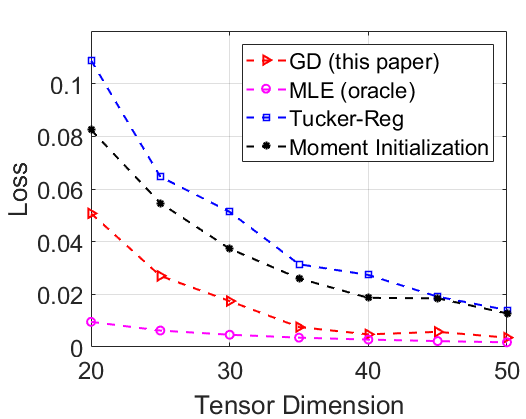}
	}
	\caption{Average estimation errors in low-rank tensor regression. The MLE (oracle) is approximated by running gradient descent with the initialization chosen at $\cX^*$, which would not be useable in practice. Here, $r = 3$, $\lambda = 2$, $\sigma = 1$. Left panel: $p = 30$, $n \in [300,1000]$. Right panel: $p \in [25, 50]$, $n = 1.2 p^{3/2} r$.}
	\label{fig:regression-error}
\end{figure}

We aim to compare the proposed method (Algorithms \ref{alg:PGD} and \ref{alg:regression}) with the initialization estimator (Algorithm \ref{alg:regression} solely), Tucker-Regression method\footnote{The implementation is based on \citep{zhou2013tensor,zhou2017matlab}.}, and MLE. Since the MLE corresponds to the global minimum of the rank-constrained optimization \eqref{eq:rank-constraint-opt} and is often computationally intractable, we instead consider a warm-start gradient descent estimator, i.e., performing Algorithm \ref{alg:PGD} starting from the true parameter $\cX^*$. We expect that the output of this procedure can well approximate MLE. We implement all four procedures under two settings: (a) $p_1=p_2=p_3 = p = 30$, $r_1=r_2=r_3=r=5$, $\lambda = 2$, $\sigma = 1$, $n$ varies from $300$ to $1000$; and (b) $p$ varies from 20 to 50, $r=3, \lambda=2,\sigma=1$, $n = 1.2p^{3/2}r$. The results are collected in Figure \ref{fig:regression-error}. We see from the left panel that for small sample size ($n\leq600$), the proposed gradient descent method significantly outperforms the Tucker-Regression and initialization estimator while has larger estimation errors than MLE. When the sample size increases ($n\geq700$), the performance of the proposed gradient descent and Tucker regression algorithms tend to be as good as MLE. Compared to the initialization,   gradient descent achieves a great improvement on the estimation accuracy. The right panel of Figure \ref{fig:regression-error} shows that the gradient descent performs as good as the warm-start gradient descent asymptotically and is significantly better than the initialization and Tucker-Regression estimators.

\ \par

\noindent{\bf Poisson Tensor PCA.} Next, we study the numerical performance of the proposed procedure on Poisson tensor PCA. As mentioned earlier in Section \ref{sec:poisson}, we found that the projection steps in Algorithm \ref{alg:PGD} are not essential to the numerical performance, thus we apply Algorithm \ref{alg:PGD} without the projection steps there in all numerical experiments for Poisson tensor PCA. 

For each experiment, we first generate a random core tensor $\cS \in \bbR^{r_1\times r_2 \times r_3}$ with i.i.d. standard normal entries and random orthogonal matrices $\U_k$ uniformly on Stiefel manifold $\bbO^{p_k \times r_k}$. Then we calculate $\bar\cX = \cS \times_1 \U_1 \times_2 \U_2 \times_3 \U_3$ and rescale it as $\cX^* = \bar\cX \cdot B/\left\|\bar \cX\right\|_\infty$ to ensure that each entry of $\cX^*$ is bounded by $B$. Now, we generate a count tensor $\cY \in \bbN^{p_1\times p_2\times p_3}$: $\cY_{jkl} \sim \text{Poisson}(I\exp(\cX^*_{jkl}))$ independently, and aim to estimate the low-rank tensor $\cX^*$ based on $\cY$. In addition to the proposed method, we also consider the baseline methods of Poisson-HOSVD and Poisson-HOOI that perform HOSVD and HOOI on $\log\left((\cY+1/2)/I\right)$ (i.e., Algorithm \ref{alg:TD_Poisson}).

\begin{figure} \centering
	\subfigure{
		\includegraphics[width =0.45\linewidth,height=2.0in]{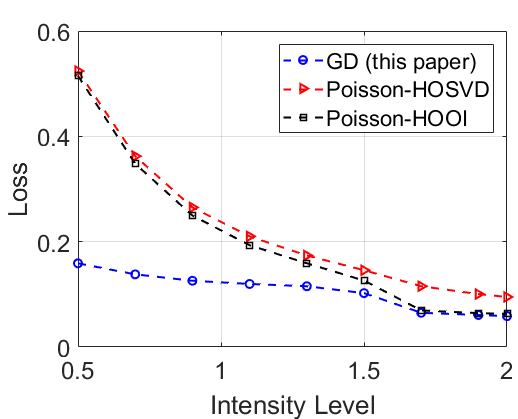}
	}
	\subfigure{
		\includegraphics[width =0.45\linewidth,height=2.0in]{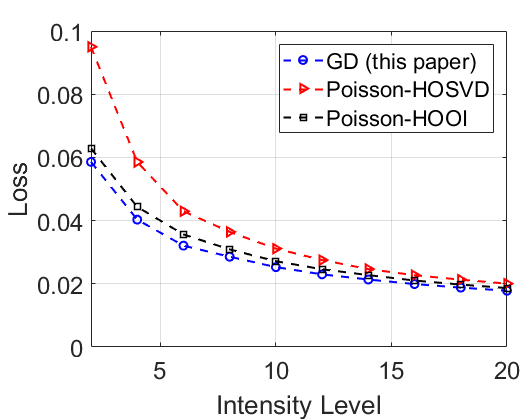}
	}
	\caption{Average estimation error of Poisson tensor PCA. Here $p=50$, $r=5$, $B = 2$. Left panel: intensity parameter $I \in [.5, 2]$. Right panel: $I\in [2, 20]$.}
	\label{fig:poisson-intensity}
	\subfigure{
		\includegraphics[width =0.45\linewidth,height=2.0in]{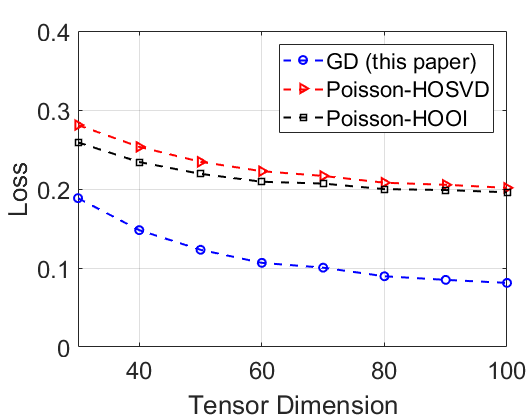}
	}
	\subfigure{
		\includegraphics[width =0.45\linewidth,height=2.0in]{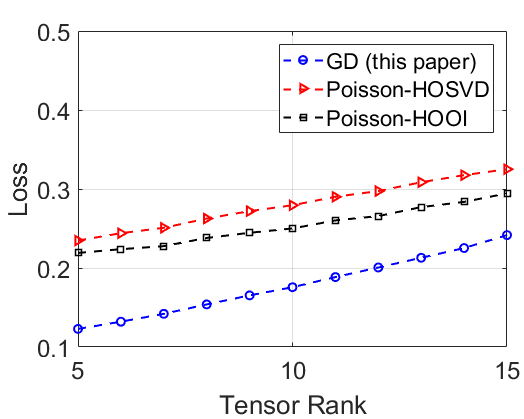}
	}
	\caption{Average estimation error of Poisson tensor PCA with different dimensions and ranks. Here, $B = 2$, $I=1$. Left panel: $r=5$, $p\in [30, 300]$. Right panel: $p=50$, $r\in [5,15]$.}
	\label{fig:poisson-dimension}
\end{figure}

First, we fix $p=50, r=5$, vary the intensity value $I$, and study the effect of $I$ to the numerical performance. As we can see from Figure \ref{fig:poisson-intensity}, for low intensity, the gradient method is significantly better than two baselines (left panel); for high intensity, three methods are comparable while the Poisson gradient descent is the best (right panel). Next, we study their performance for different tensor dimensions and ranks. In the left panel of Figure \ref{fig:poisson-dimension}, we set $r=5$ and vary $p$ from $30$ to $100$; in the right panel of Figure \ref{fig:poisson-dimension}, we fix $p=50$ and vary $r$ from $5$ to 15. As one can see, our method significantly outperforms the baselines in all settings. All these simulation results illustrate the benefits of applying gradient descent on the Poisson likelihood function. 

\ \par

\begin{figure} \centering
	\subfigure{
		\includegraphics[width =0.45\linewidth,height=2.0in]{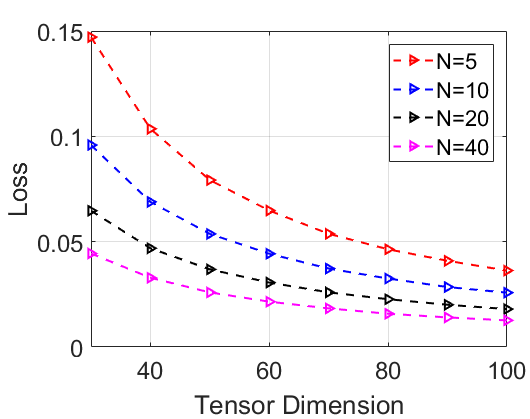}
	}
	\subfigure{
		\includegraphics[width =0.45\linewidth,height=2.0in]{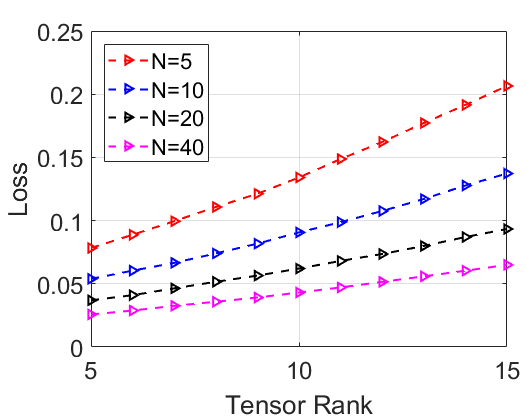}
	}
	\caption{Average estimation error of binomial tensor PCA. Here, $B = 2$. Left panel: $r=5$, $p\in[30, 100]$. Right panel: $p=50, r \in [5:15]$.}
	\label{fig:binomial-dimension}
\end{figure}

\noindent{\bf Binomial Tensor PCA.} We generate $\cX^*$ in the same way as the Poisson tensor PCA settings. Suppose we observe $\cY \in \bbN^{p_1\times p_2\times p_3}$ generated from
\begin{equation*}
    \cY_{jkl} \sim \text{Binomial}\left(\cN_{jkl}, s(\cX_{jkl}^*)\right),\quad \text{independently}.
\end{equation*}
We take all entries with the same population size (i.e., $\cN_{ijk} = N$) for simplification. 
We can see from the simulation results in Figure \ref{fig:binomial-dimension} that a larger population size $N$ yields smaller estimation error. In addition, according to Theorem \ref{thm:local-convergence-PGD}, the estimation error in theory is of order $O(p^{-1/2}r^{1/2})$, 
which matches the trend of estimation error curves in Figure \ref{fig:binomial-dimension}.

\ \par

\subsection{Real Data Analysis}\label{sec:real-data}
In this section, we apply the proposed framework to real data applications in 4D-STEM image denoising. An additional real data example on click-through prediction is postponed to Appendix \ref{sec:click-through} in the supplementary materials.

\begin{figure} \centering
	\subfigure{
		\includegraphics[width =1\linewidth,height=5in]{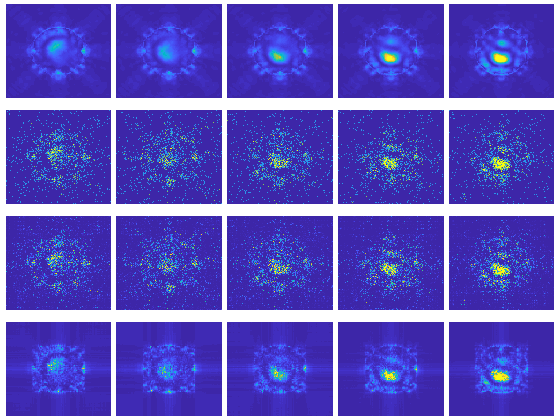}
	}
	\caption{Recovery results for the first five images of 4D-STEM data. First row: original images; second row: photon-limited observations; third row: denoised images by matrix method; forth row: denoised images by the proposed tensor method}
	\label{fig:4DSTEM-example}
\end{figure}

The 4D-scanning transmission electron microscopy (4D-STEM) is an important technique in modern material science that has been used to detect local material composition of structures such as films, defects and nanostructures \citep{krivanek1999towards,yankovich2016non}. In 4D-STEM imaging technology, a focused probe is usually rastered across part of the specimen and an X-ray and/or electron energy loss spectrum is recorded at each probe position, generating a series of photon-limited images. The data generated from 4D-STEM technique are typically order-4 tensors with approximate periodic structures, as a focused probe is located on a 2-D grid and one 2-D image is generated for each probe position (see \cite{yankovich2016non} for more details). Due to the physical conditions, the observable images are often photon-limited, highly noisy, and in the form of count matrices (see the second row of Figure \ref{fig:4DSTEM-example} for an example). A sufficient imaging denoising is often a crucial first step before the subsequent procedures.

We aim to illustrate the merit of the proposed method through denoising of data in 4D-STEM experiments. Specifically, we collect $160$ images generated from a row of electron probe positions\footnote{Simultaneously denoising the order-4 image data requires extremely large memory and computation source. Thus we focus on one row of images. It is also common to perform row-wise image denoising in 4D-STEM imaging analysis \citep{yankovich2016non}.}. Since the resolution of each image is $183\times 183$, the data images can be stacked into a non-negative tensor of size $160\times 183\times 183$. 
We assume the observational images $\cY$ are generated from Poisson distribution $\cY_{ijk} \overset{iid}{\sim}\Poisson(\exp(\cX_{ijk}^*))$. Our goal is to recover the original images based on the photon-limited observation $\cY$. Since $\cY$ is sparse ($\approx 88\%$ pixels are zero), we take the pre-initializer $\tilde\cX = \log(\cY+1/30)$ and take the input rank according to the \eqref{eq:rk-select-pca} with $\rho = 0.98$.
We apply the proposed gradient descent (Algorithms \ref{alg:PGD} and \ref{alg:TD_Poisson}) with the rank estimation $(\hat r_1, \hat r_2, \hat r_3) = (44,36,34)$ to obtain the estimator $\hat \cX$, then calculate $\exp(\hat \cX)$ as the collection of denoised images. 
We also denoise these images one by one via the matrix Procrustes flow \citep{park2018finding}, a variant of the matrix-version gradient descent method.\footnote{This algorithm also requires the specification of matrix rank. We empirically choose $\hat r = \argmin\left\{r: \sum_{i=1}^r \sigma_i^2(\X)\big/\sum_{i=1}^{p_k}\sigma_i^2(\X) \geq 0.98\right\}$ for each slice of the pre-initializer $\X = \tilde \cX_{i::}$.} The original, observational, and recovered images are provided in Figure \ref{fig:4DSTEM-example}. In addition, we calculate the recovery loss for each of the 160 images, i.e., $\|\exp(\cX^*_{i::}) - \exp(\hat \cX_{i::})\|_\tF/\|\exp(\cX^*_{i::})\|_\tF$, and the averaged recovery loss (and standard error) of matrix and tensor methods are 0.861 (0.183) and 0.303 (0.054), respectively. One can clearly see the advantage of the proposed tensor method that utilizes the tensor structure of the whole set of images. 

\ \par

\section{Discussions}\label{sec:discuss}

In this paper, we introduce a non-convex optimization framework for the generalized tensor estimation. Compared to the convex relaxation methods in the literature, the proposed scheme is computationally efficient and achieves desirable statistical error rate under suitable initialization and signal-to-noise ratio conditions. We apply the proposed framework on several problems, including sub-Gaussian denoising, tensor regression, Poisson and binomial tensor PCA. We can show that the proposed gradient descent procedure achieves the minimax optimal rate of estimation error under these statistical models.

In addition to the above-mentioned problems, the proposed framework can incorporate a broader range of settings. For example, the developed result is applicable to solve the {\it noisy tensor completion} problem \citep{xia2017statistically,montanari2018spectral,shah2019iterative,cai2019nonconvex}, which aims to recover the low-rank tensor $\cX^*$ based on a number of noisy observable entries, say $\{\cY_{ijk} = \cX^*_{ijk} + \cZ_{ijk}\}_{(i,j,k)\in\Omega}$, where $\Omega$ is a subset of indices.

Another example is {\it binary tensor PCA} \citep{wang2018learning}, where the central goal is to factorize the tensor from 0-1 valued observations. Suppose one observes $\cY_{ijk} \sim \text{Bernoulli}(\cP_{ijk})$ independently, where $\cP_{ijk} = s(\cX^*_{ijk})$, $\cX^*$ is low-rank, and $s(\cdot)$ is some link function. Then the proposed projected gradient descent method can be applied to estimate $\cX^*$ with provable guarantees. 

Community detection in social network has attracted enormous recent attention. Although most of the existing results focused on a single-layer of network, the \emph{multilayer network}, i.e., the connections between different nodes are reflected in multiple modalities, also commonly appear in practice \citep{han2015consistent,lei2019consistent,pensky2019spectral}. Consider a stack of multilayer network data with shared community structure. It is reasonable to assume that the adjacency tensor $\cA$ has a low-rank tensor structure: $\cA \sim \text{Bernoulli}(\cX^*)$ independently, where $\cX^* = \llbracket \cS^*, \Z^*, \Z^*, \T^* \rrbracket$, $\Z^*$ is the latent space of nodes features (or the indicator matrix for the community that each node belongs to), and $\T^*$ models the trend along the time. Then the community detection for multilayer networks essentially becomes the generalized tensor estimation problem.

In addition to the standard linear regression model discussed in Section \ref{sec:tensor_reg}, the proposed framework can be applied to a range of generalized tensor regression problems. Recall that the classical generalized linear model focuses on an exponential family, where the response $y_i$ satisfies the following density or probability mass function \citep{nelder1972generalized},
\begin{equation}\label{eq:y_i-density}
p(y_i|\theta_i, \phi) = \exp\left\{\frac{y\theta_i - b(\theta_i)}{a(\phi)} + c(y,\phi)\right\}.
\end{equation}
Here, $a, b, c$ are prespecified functions determined by the problem; $\theta_i$ and $\phi>0$ are natural and dispersion parameters, respectively. For the generalized tensor regression, it is natural to relate the tensor covariate and response \citep{zhou2013tensor} via 
\begin{equation}\label{eq:y_i-given-X}
\mu_i = \bbE(y_i|\cX^*)\quad g(\mu_i) = \langle\cA_i, \cX^*\rangle,
\end{equation}
where $g(\cdot)$ is a link function. To estimate $\cX^*$, we can apply the proposed Algorithm \ref{alg:regression} on the negative log-likelihood function
\begin{equation*}
    \begin{split}
        \sum_{i=1}^n \frac{y_i\theta_i - b(\theta_i)}{a(\phi)} + \sum_{i=1}^n c(y_i, \phi),
    \end{split}
\end{equation*}
where $\theta_i$ is determined by \eqref{eq:y_i-density} and \eqref{eq:y_i-given-X}.

Some other possible applications of the proposed framework include the \emph{high-order interaction pursuit} \citep{hao2018sparse}, \emph{generalized regression among multiple modes} \citep{xu2019generalized}, \emph{mixed-data-type tensor data analysis} \citep{baker2019feature}, etc. In all these problems, by exploring the log-likelihood of data and the domain $\cC$ that satisfies RCG condition, the proposed projected gradient descent can be applied and the theoretical guarantees can be developed based on the proposed framework.

\section*{Acknowledgement}

The authors thank Paul Voyles and Chenyu Zhang for providing the 4D-STEM dataset and for helpful discussions. The research of R. H. and A. R. Z. was supported in part by NSF DMS-1811868, NSF CAREER-1944904, and NIH R01-GM131399. The research of R. W. was supported in part by AFOSR FA9550-18-1-0166, DOE DE-AC02-06CH11357, NSF OAC-1934637, and NSF DMS-2023109. The research of R. H. was also supported in part by a RAship from Institute for Mathematics of Data Science at UW-Madison. 

\bibliographystyle{imsart-nameyear}
\bibliography{reference}

\begin{thebibliography}{114}

\bibitem[\protect\citeauthoryear{Ahmed, Recht and
  Romberg}{2013}]{ahmed2013blind}
\begin{barticle}[author]
\bauthor{\bsnm{Ahmed},~\bfnm{Ali}\binits{A.}},
  \bauthor{\bsnm{Recht},~\bfnm{Benjamin}\binits{B.}} \AND
  \bauthor{\bsnm{Romberg},~\bfnm{Justin}\binits{J.}}
(\byear{2013}).
\btitle{Blind deconvolution using convex programming}.
\bjournal{IEEE Transactions on Information Theory}
\bvolume{60}
\bpages{1711--1732}.
\end{barticle}
\endbibitem

\bibitem[\protect\citeauthoryear{Anandkumar, Hsu and
  Kakade}{2012}]{anandkumar2012method}
\begin{binproceedings}[author]
\bauthor{\bsnm{Anandkumar},~\bfnm{Animashree}\binits{A.}},
  \bauthor{\bsnm{Hsu},~\bfnm{Daniel}\binits{D.}} \AND
  \bauthor{\bsnm{Kakade},~\bfnm{Sham~M}\binits{S.~M.}}
(\byear{2012}).
\btitle{A method of moments for mixture models and hidden Markov models}.
In \bbooktitle{Conference on Learning Theory}
\bpages{33--1}.
\end{binproceedings}
\endbibitem

\bibitem[\protect\citeauthoryear{Anandkumar
  et~al.}{2014}]{anandkumar2014tensor}
\begin{barticle}[author]
\bauthor{\bsnm{Anandkumar},~\bfnm{Animashree}\binits{A.}},
  \bauthor{\bsnm{Ge},~\bfnm{Rong}\binits{R.}},
  \bauthor{\bsnm{Hsu},~\bfnm{Daniel}\binits{D.}},
  \bauthor{\bsnm{Kakade},~\bfnm{Sham~M}\binits{S.~M.}} \AND
  \bauthor{\bsnm{Telgarsky},~\bfnm{Matus}\binits{M.}}
(\byear{2014}).
\btitle{Tensor decompositions for learning latent variable models}.
\bjournal{The Journal of Machine Learning Research}
\bvolume{15}
\bpages{2773--2832}.
\end{barticle}
\endbibitem

\bibitem[\protect\citeauthoryear{Arroyo et~al.}{2019}]{arroyo2019inference}
\begin{barticle}[author]
\bauthor{\bsnm{Arroyo},~\bfnm{Jes{\'u}s}\binits{J.}},
  \bauthor{\bsnm{Athreya},~\bfnm{Avanti}\binits{A.}},
  \bauthor{\bsnm{Cape},~\bfnm{Joshua}\binits{J.}},
  \bauthor{\bsnm{Chen},~\bfnm{Guodong}\binits{G.}},
  \bauthor{\bsnm{Priebe},~\bfnm{Carey~E}\binits{C.~E.}} \AND
  \bauthor{\bsnm{Vogelstein},~\bfnm{Joshua~T}\binits{J.~T.}}
(\byear{2019}).
\btitle{Inference for multiple heterogeneous networks with a common invariant
  subspace}.
\bjournal{arXiv preprint arXiv:1906.10026}.
\end{barticle}
\endbibitem

\bibitem[\protect\citeauthoryear{Bahadori, Yu and Liu}{2014}]{bahadori2014fast}
\begin{binproceedings}[author]
\bauthor{\bsnm{Bahadori},~\bfnm{Mohammad~Taha}\binits{M.~T.}},
  \bauthor{\bsnm{Yu},~\bfnm{Qi~Rose}\binits{Q.~R.}} \AND
  \bauthor{\bsnm{Liu},~\bfnm{Yan}\binits{Y.}}
(\byear{2014}).
\btitle{Fast multivariate spatio-temporal analysis via low rank tensor
  learning}.
In \bbooktitle{Advances in neural information processing systems}
\bpages{3491--3499}.
\end{binproceedings}
\endbibitem

\bibitem[\protect\citeauthoryear{Baker, Tang and
  Allen}{2019}]{baker2019feature}
\begin{barticle}[author]
\bauthor{\bsnm{Baker},~\bfnm{Yulia}\binits{Y.}},
  \bauthor{\bsnm{Tang},~\bfnm{Tiffany~M}\binits{T.~M.}} \AND
  \bauthor{\bsnm{Allen},~\bfnm{Genevera~I}\binits{G.~I.}}
(\byear{2019}).
\btitle{Feature Selection for Data Integration with Mixed Multi-view Data}.
\bjournal{arXiv preprint arXiv:1903.11232}.
\end{barticle}
\endbibitem

\bibitem[\protect\citeauthoryear{Barak and Moitra}{2016}]{barak2016noisy}
\begin{binproceedings}[author]
\bauthor{\bsnm{Barak},~\bfnm{Boaz}\binits{B.}} \AND
  \bauthor{\bsnm{Moitra},~\bfnm{Ankur}\binits{A.}}
(\byear{2016}).
\btitle{Noisy tensor completion via the sum-of-squares hierarchy}.
In \bbooktitle{Conference on Learning Theory}
\bpages{417--445}.
\end{binproceedings}
\endbibitem

\bibitem[\protect\citeauthoryear{Bi, Qu and Shen}{2018}]{bi2018multilayer}
\begin{barticle}[author]
\bauthor{\bsnm{Bi},~\bfnm{Xuan}\binits{X.}},
  \bauthor{\bsnm{Qu},~\bfnm{Annie}\binits{A.}} \AND
  \bauthor{\bsnm{Shen},~\bfnm{Xiaotong}\binits{X.}}
(\byear{2018}).
\btitle{Multilayer tensor factorization with applications to recommender
  systems}.
\bjournal{The Annals of Statistics}
\bvolume{46}
\bpages{3308--3333}.
\end{barticle}
\endbibitem

\bibitem[\protect\citeauthoryear{Birg{\'e}}{2001}]{birge2001alternative}
\begin{barticle}[author]
\bauthor{\bsnm{Birg{\'e}},~\bfnm{Lucien}\binits{L.}}
(\byear{2001}).
\btitle{An alternative point of view on Lepski's method}.
\bjournal{Lecture Notes-Monograph Series}
\bpages{113--133}.
\end{barticle}
\endbibitem

\bibitem[\protect\citeauthoryear{Boucheron, Lugosi and
  Massart}{2013}]{boucheron2013concentration}
\begin{bbook}[author]
\bauthor{\bsnm{Boucheron},~\bfnm{St{\'e}phane}\binits{S.}},
  \bauthor{\bsnm{Lugosi},~\bfnm{G{\'a}bor}\binits{G.}} \AND
  \bauthor{\bsnm{Massart},~\bfnm{Pascal}\binits{P.}}
(\byear{2013}).
\btitle{Concentration inequalities: A nonasymptotic theory of independence}.
\bpublisher{Oxford university press}.
\end{bbook}
\endbibitem

\bibitem[\protect\citeauthoryear{Cai, Li and Ma}{2016}]{cai2016optimal}
\begin{barticle}[author]
\bauthor{\bsnm{Cai},~\bfnm{T~Tony}\binits{T.~T.}},
  \bauthor{\bsnm{Li},~\bfnm{Xiaodong}\binits{X.}} \AND
  \bauthor{\bsnm{Ma},~\bfnm{Zongming}\binits{Z.}}
(\byear{2016}).
\btitle{Optimal rates of convergence for noisy sparse phase retrieval via
  thresholded Wirtinger flow}.
\bjournal{The Annals of Statistics}
\bvolume{44}
\bpages{2221--2251}.
\end{barticle}
\endbibitem

\bibitem[\protect\citeauthoryear{Cai and Zhang}{2018}]{cai2018rate}
\begin{barticle}[author]
\bauthor{\bsnm{Cai},~\bfnm{T~Tony}\binits{T.~T.}} \AND
  \bauthor{\bsnm{Zhang},~\bfnm{Anru}\binits{A.}}
(\byear{2018}).
\btitle{Rate-optimal perturbation bounds for singular subspaces with
  applications to high-dimensional statistics}.
\bjournal{The Annals of Statistics}
\bvolume{46}
\bpages{60--89}.
\end{barticle}
\endbibitem

\bibitem[\protect\citeauthoryear{Cai et~al.}{2019}]{cai2019nonconvex}
\begin{binproceedings}[author]
\bauthor{\bsnm{Cai},~\bfnm{Changxiao}\binits{C.}},
  \bauthor{\bsnm{Li},~\bfnm{Gen}\binits{G.}},
  \bauthor{\bsnm{Poor},~\bfnm{H~Vincent}\binits{H.~V.}} \AND
  \bauthor{\bsnm{Chen},~\bfnm{Yuxin}\binits{Y.}}
(\byear{2019}).
\btitle{Nonconvex Low-Rank Tensor Completion from Noisy Data}.
In \bbooktitle{Advances in Neural Information Processing Systems}
\bpages{1861--1872}.
\end{binproceedings}
\endbibitem

\bibitem[\protect\citeauthoryear{Candes, Li and
  Soltanolkotabi}{2015}]{candes2015phase}
\begin{barticle}[author]
\bauthor{\bsnm{Candes},~\bfnm{Emmanuel~J}\binits{E.~J.}},
  \bauthor{\bsnm{Li},~\bfnm{Xiaodong}\binits{X.}} \AND
  \bauthor{\bsnm{Soltanolkotabi},~\bfnm{Mahdi}\binits{M.}}
(\byear{2015}).
\btitle{Phase retrieval via Wirtinger flow: Theory and algorithms}.
\bjournal{IEEE Transactions on Information Theory}
\bvolume{61}
\bpages{1985--2007}.
\end{barticle}
\endbibitem

\bibitem[\protect\citeauthoryear{Candes and Plan}{2010}]{candes2010matrix}
\begin{barticle}[author]
\bauthor{\bsnm{Candes},~\bfnm{Emmanuel~J}\binits{E.~J.}} \AND
  \bauthor{\bsnm{Plan},~\bfnm{Yaniv}\binits{Y.}}
(\byear{2010}).
\btitle{Matrix completion with noise}.
\bjournal{Proceedings of the IEEE}
\bvolume{98}
\bpages{925--936}.
\end{barticle}
\endbibitem

\bibitem[\protect\citeauthoryear{Candes and Plan}{2011}]{candes2011tight}
\begin{barticle}[author]
\bauthor{\bsnm{Candes},~\bfnm{Emmanuel~J}\binits{E.~J.}} \AND
  \bauthor{\bsnm{Plan},~\bfnm{Yaniv}\binits{Y.}}
(\byear{2011}).
\btitle{Tight oracle inequalities for low-rank matrix recovery from a minimal
  number of noisy random measurements}.
\bjournal{IEEE Transactions on Information Theory}
\bvolume{57}
\bpages{2342--2359}.
\end{barticle}
\endbibitem

\bibitem[\protect\citeauthoryear{Cand{\`e}s and Recht}{2009}]{candes2009exact}
\begin{barticle}[author]
\bauthor{\bsnm{Cand{\`e}s},~\bfnm{Emmanuel~J}\binits{E.~J.}} \AND
  \bauthor{\bsnm{Recht},~\bfnm{Benjamin}\binits{B.}}
(\byear{2009}).
\btitle{Exact matrix completion via convex optimization}.
\bjournal{Foundations of Computational mathematics}
\bvolume{9}
\bpages{717}.
\end{barticle}
\endbibitem

\bibitem[\protect\citeauthoryear{Cao and Xie}{2015}]{cao2015poisson}
\begin{barticle}[author]
\bauthor{\bsnm{Cao},~\bfnm{Yang}\binits{Y.}} \AND
  \bauthor{\bsnm{Xie},~\bfnm{Yao}\binits{Y.}}
(\byear{2015}).
\btitle{Poisson matrix recovery and completion}.
\bjournal{IEEE Transactions on Signal Processing}
\bvolume{64}
\bpages{1609--1620}.
\end{barticle}
\endbibitem

\bibitem[\protect\citeauthoryear{Cao, Zhang and Li}{2019}]{cao2017microbial}
\begin{barticle}[author]
\bauthor{\bsnm{Cao},~\bfnm{Yuanpei}\binits{Y.}},
  \bauthor{\bsnm{Zhang},~\bfnm{Anru}\binits{A.}} \AND
  \bauthor{\bsnm{Li},~\bfnm{Hongzhe}\binits{H.}}
(\byear{2019}).
\btitle{Multi-sample estimation of bacterial composition matrix in metagenomics
  data}.
\bjournal{Biometrika}.
\end{barticle}
\endbibitem

\bibitem[\protect\citeauthoryear{Chen}{2019}]{chen2019phase}
\begin{barticle}[author]
\bauthor{\bsnm{Chen},~\bfnm{Wei-Kuo}\binits{W.-K.}}
(\byear{2019}).
\btitle{Phase transition in the spiked random tensor with rademacher prior}.
\bjournal{The Annals of Statistics}
\bvolume{47}
\bpages{2734--2756}.
\end{barticle}
\endbibitem

\bibitem[\protect\citeauthoryear{Chen and Candes}{2015}]{chen2015solving}
\begin{binproceedings}[author]
\bauthor{\bsnm{Chen},~\bfnm{Yuxin}\binits{Y.}} \AND
  \bauthor{\bsnm{Candes},~\bfnm{Emmanuel}\binits{E.}}
(\byear{2015}).
\btitle{Solving random quadratic systems of equations is nearly as easy as
  solving linear systems}.
In \bbooktitle{Advances in Neural Information Processing Systems}
\bpages{739--747}.
\end{binproceedings}
\endbibitem

\bibitem[\protect\citeauthoryear{Chen and Chi}{2018}]{chen2018harnessing}
\begin{barticle}[author]
\bauthor{\bsnm{Chen},~\bfnm{Yudong}\binits{Y.}} \AND
  \bauthor{\bsnm{Chi},~\bfnm{Yuejie}\binits{Y.}}
(\byear{2018}).
\btitle{Harnessing structures in big data via guaranteed low-rank matrix
  estimation}.
\bjournal{arXiv preprint arXiv:1802.08397}.
\end{barticle}
\endbibitem

\bibitem[\protect\citeauthoryear{Chen, Raskutti and Yuan}{2019}]{chen2019non}
\begin{barticle}[author]
\bauthor{\bsnm{Chen},~\bfnm{Han}\binits{H.}},
  \bauthor{\bsnm{Raskutti},~\bfnm{Garvesh}\binits{G.}} \AND
  \bauthor{\bsnm{Yuan},~\bfnm{Ming}\binits{M.}}
(\byear{2019}).
\btitle{Non-convex projected gradient descent for generalized low-rank tensor
  regression}.
\bjournal{The Journal of Machine Learning Research}
\bvolume{20}
\bpages{172--208}.
\end{barticle}
\endbibitem

\bibitem[\protect\citeauthoryear{Chi and Kolda}{2012}]{chi2012tensors}
\begin{barticle}[author]
\bauthor{\bsnm{Chi},~\bfnm{Eric~C}\binits{E.~C.}} \AND
  \bauthor{\bsnm{Kolda},~\bfnm{Tamara~G}\binits{T.~G.}}
(\byear{2012}).
\btitle{On tensors, sparsity, and nonnegative factorizations}.
\bjournal{SIAM Journal on Matrix Analysis and Applications}
\bvolume{33}
\bpages{1272--1299}.
\end{barticle}
\endbibitem

\bibitem[\protect\citeauthoryear{Chi, Lu and Chen}{2019}]{chi2019nonconvex}
\begin{barticle}[author]
\bauthor{\bsnm{Chi},~\bfnm{Yuejie}\binits{Y.}},
  \bauthor{\bsnm{Lu},~\bfnm{Yue~M}\binits{Y.~M.}} \AND
  \bauthor{\bsnm{Chen},~\bfnm{Yuxin}\binits{Y.}}
(\byear{2019}).
\btitle{Nonconvex optimization meets low-rank matrix factorization: An
  overview}.
\bjournal{IEEE Transactions on Signal Processing}
\bvolume{67}
\bpages{5239--5269}.
\end{barticle}
\endbibitem

\bibitem[\protect\citeauthoryear{Chi et~al.}{2018}]{ChiGaiSunZhoYan2018}
\begin{bmisc}[author]
\bauthor{\bsnm{Chi},~\bfnm{Eric~C.}\binits{E.~C.}},
  \bauthor{\bsnm{Gaines},~\bfnm{Brian~R.}\binits{B.~R.}},
  \bauthor{\bsnm{Sun},~\bfnm{Will~Wei}\binits{W.~W.}},
  \bauthor{\bsnm{Zhou},~\bfnm{Hua}\binits{H.}} \AND
  \bauthor{\bsnm{Yang},~\bfnm{Jian}\binits{J.}}
(\byear{2018}).
\btitle{Provable Convex Co-clustering of Tensors}.
\bhowpublished{arXiv:1803.06518 [stat.ME]}.
\arxiv{https://arxiv.org/abs/1803.06518}
\end{bmisc}
\endbibitem

\bibitem[\protect\citeauthoryear{De~Lathauwer, De~Moor and
  Vandewalle}{2000a}]{de2000multilinear}
\begin{barticle}[author]
\bauthor{\bsnm{De~Lathauwer},~\bfnm{Lieven}\binits{L.}},
  \bauthor{\bsnm{De~Moor},~\bfnm{Bart}\binits{B.}} \AND
  \bauthor{\bsnm{Vandewalle},~\bfnm{Joos}\binits{J.}}
(\byear{2000}a).
\btitle{A multilinear singular value decomposition}.
\bjournal{SIAM journal on Matrix Analysis and Applications}
\bvolume{21}
\bpages{1253--1278}.
\end{barticle}
\endbibitem

\bibitem[\protect\citeauthoryear{De~Lathauwer, De~Moor and
  Vandewalle}{2000b}]{de2000best}
\begin{barticle}[author]
\bauthor{\bsnm{De~Lathauwer},~\bfnm{Lieven}\binits{L.}},
  \bauthor{\bsnm{De~Moor},~\bfnm{Bart}\binits{B.}} \AND
  \bauthor{\bsnm{Vandewalle},~\bfnm{Joos}\binits{J.}}
(\byear{2000}b).
\btitle{On the best rank-1 and rank-(r 1, r 2,..., rn) approximation of
  higher-order tensors}.
\bjournal{SIAM Journal on Matrix Analysis and Applications}
\bvolume{21}
\bpages{1324--1342}.
\end{barticle}
\endbibitem

\bibitem[\protect\citeauthoryear{Fan, Gong and Zhu}{2019}]{fan2019generalized}
\begin{barticle}[author]
\bauthor{\bsnm{Fan},~\bfnm{Jianqing}\binits{J.}},
  \bauthor{\bsnm{Gong},~\bfnm{Wenyan}\binits{W.}} \AND
  \bauthor{\bsnm{Zhu},~\bfnm{Ziwei}\binits{Z.}}
(\byear{2019}).
\btitle{Generalized high-dimensional trace regression via nuclear norm
  regularization}.
\bjournal{Journal of Econometrics}.
\end{barticle}
\endbibitem

\bibitem[\protect\citeauthoryear{Faust et~al.}{2012}]{faust2012microbial}
\begin{barticle}[author]
\bauthor{\bsnm{Faust},~\bfnm{Karoline}\binits{K.}},
  \bauthor{\bsnm{Sathirapongsasuti},~\bfnm{J~Fah}\binits{J.~F.}},
  \bauthor{\bsnm{Izard},~\bfnm{Jacques}\binits{J.}},
  \bauthor{\bsnm{Segata},~\bfnm{Nicola}\binits{N.}},
  \bauthor{\bsnm{Gevers},~\bfnm{Dirk}\binits{D.}},
  \bauthor{\bsnm{Raes},~\bfnm{Jeroen}\binits{J.}} \AND
  \bauthor{\bsnm{Huttenhower},~\bfnm{Curtis}\binits{C.}}
(\byear{2012}).
\btitle{Microbial co-occurrence relationships in the human microbiome}.
\bjournal{PLoS computational biology}
\bvolume{8}
\bpages{e1002606}.
\end{barticle}
\endbibitem

\bibitem[\protect\citeauthoryear{Fazel}{2002}]{fazel2002matrix}
\begin{bmisc}[author]
\bauthor{\bsnm{Fazel},~\bfnm{Maryam}\binits{M.}}
(\byear{2002}).
\btitle{Matrix rank minimization with applications}.
\end{bmisc}
\endbibitem

\bibitem[\protect\citeauthoryear{Flores et~al.}{2014}]{flores2014temporal}
\begin{barticle}[author]
\bauthor{\bsnm{Flores},~\bfnm{Gilberto~E}\binits{G.~E.}},
  \bauthor{\bsnm{Caporaso},~\bfnm{J~Gregory}\binits{J.~G.}},
  \bauthor{\bsnm{Henley},~\bfnm{Jessica~B}\binits{J.~B.}},
  \bauthor{\bsnm{Rideout},~\bfnm{Jai~Ram}\binits{J.~R.}},
  \bauthor{\bsnm{Domogala},~\bfnm{Daniel}\binits{D.}},
  \bauthor{\bsnm{Chase},~\bfnm{John}\binits{J.}},
  \bauthor{\bsnm{Leff},~\bfnm{Jonathan~W}\binits{J.~W.}},
  \bauthor{\bsnm{V{\'a}zquez-Baeza},~\bfnm{Yoshiki}\binits{Y.}},
  \bauthor{\bsnm{Gonzalez},~\bfnm{Antonio}\binits{A.}},
  \bauthor{\bsnm{Knight},~\bfnm{Rob}\binits{R.}} \betal{et~al.}
(\byear{2014}).
\btitle{Temporal variability is a personalized feature of the human
  microbiome}.
\bjournal{Genome biology}
\bvolume{15}
\bpages{531}.
\end{barticle}
\endbibitem

\bibitem[\protect\citeauthoryear{Friedland and
  Lim}{2018}]{friedland2018nuclear}
\begin{barticle}[author]
\bauthor{\bsnm{Friedland},~\bfnm{Shmuel}\binits{S.}} \AND
  \bauthor{\bsnm{Lim},~\bfnm{Lek-Heng}\binits{L.-H.}}
(\byear{2018}).
\btitle{Nuclear norm of higher-order tensors}.
\bjournal{Mathematics of Computation}
\bvolume{87}
\bpages{1255--1281}.
\end{barticle}
\endbibitem

\bibitem[\protect\citeauthoryear{Gandy, Recht and
  Yamada}{2011}]{gandy2011tensor}
\begin{barticle}[author]
\bauthor{\bsnm{Gandy},~\bfnm{Silvia}\binits{S.}},
  \bauthor{\bsnm{Recht},~\bfnm{Benjamin}\binits{B.}} \AND
  \bauthor{\bsnm{Yamada},~\bfnm{Isao}\binits{I.}}
(\byear{2011}).
\btitle{Tensor completion and low-n-rank tensor recovery via convex
  optimization}.
\bjournal{Inverse Problems}
\bvolume{27}
\bpages{025010}.
\end{barticle}
\endbibitem

\bibitem[\protect\citeauthoryear{Guhaniyogi, Qamar and
  Dunson}{2017}]{guhaniyogi2017bayesian}
\begin{barticle}[author]
\bauthor{\bsnm{Guhaniyogi},~\bfnm{Rajarshi}\binits{R.}},
  \bauthor{\bsnm{Qamar},~\bfnm{Shaan}\binits{S.}} \AND
  \bauthor{\bsnm{Dunson},~\bfnm{David~B}\binits{D.~B.}}
(\byear{2017}).
\btitle{Bayesian tensor regression}.
\bjournal{The Journal of Machine Learning Research}
\bvolume{18}
\bpages{2733--2763}.
\end{barticle}
\endbibitem

\bibitem[\protect\citeauthoryear{Guo, Kotsia and Patras}{2012}]{guo2012tensor}
\begin{barticle}[author]
\bauthor{\bsnm{Guo},~\bfnm{Weiwei}\binits{W.}},
  \bauthor{\bsnm{Kotsia},~\bfnm{Irene}\binits{I.}} \AND
  \bauthor{\bsnm{Patras},~\bfnm{Ioannis}\binits{I.}}
(\byear{2012}).
\btitle{Tensor learning for regression}.
\bjournal{IEEE Transactions on Image Processing}
\bvolume{21}
\bpages{816--827}.
\end{barticle}
\endbibitem

\bibitem[\protect\citeauthoryear{Hall, Raskutti and
  Willett}{2016}]{hall2016inference}
\begin{barticle}[author]
\bauthor{\bsnm{Hall},~\bfnm{Eric~C}\binits{E.~C.}},
  \bauthor{\bsnm{Raskutti},~\bfnm{Garvesh}\binits{G.}} \AND
  \bauthor{\bsnm{Willett},~\bfnm{Rebecca}\binits{R.}}
(\byear{2016}).
\btitle{Inference of high-dimensional autoregressive generalized linear
  models}.
\bjournal{arXiv preprint arXiv:1605.02693}.
\end{barticle}
\endbibitem

\bibitem[\protect\citeauthoryear{Han, Xu and Airoldi}{2015}]{han2015consistent}
\begin{binproceedings}[author]
\bauthor{\bsnm{Han},~\bfnm{Qiuyi}\binits{Q.}},
  \bauthor{\bsnm{Xu},~\bfnm{Kevin}\binits{K.}} \AND
  \bauthor{\bsnm{Airoldi},~\bfnm{Edoardo}\binits{E.}}
(\byear{2015}).
\btitle{Consistent estimation of dynamic and multi-layer block models}.
In \bbooktitle{International Conference on Machine Learning}
\bpages{1511--1520}.
\end{binproceedings}
\endbibitem

\bibitem[\protect\citeauthoryear{Hao, Zhang and Cheng}{2019}]{hao2018sparse}
\begin{barticle}[author]
\bauthor{\bsnm{Hao},~\bfnm{Botao}\binits{B.}},
  \bauthor{\bsnm{Zhang},~\bfnm{Anru}\binits{A.}} \AND
  \bauthor{\bsnm{Cheng},~\bfnm{Guang}\binits{G.}}
(\byear{2019}).
\btitle{Sparse and Low-rank Tensor Estimation via Cubic Sketchings}.
\bjournal{The Annals of Statistics}
\bvolume{revision under review}.
\end{barticle}
\endbibitem

\bibitem[\protect\citeauthoryear{Henriques and
  Madeira}{2019}]{henriques2019triclustering}
\begin{barticle}[author]
\bauthor{\bsnm{Henriques},~\bfnm{Rui}\binits{R.}} \AND
  \bauthor{\bsnm{Madeira},~\bfnm{Sara~C}\binits{S.~C.}}
(\byear{2019}).
\btitle{Triclustering algorithms for three-dimensional data analysis: A
  comprehensive survey}.
\bjournal{ACM Computing Surveys (CSUR)}
\bvolume{51}
\bpages{95}.
\end{barticle}
\endbibitem

\bibitem[\protect\citeauthoryear{Hillar and Lim}{2013}]{hillar2013most}
\begin{barticle}[author]
\bauthor{\bsnm{Hillar},~\bfnm{Christopher~J}\binits{C.~J.}} \AND
  \bauthor{\bsnm{Lim},~\bfnm{Lek-Heng}\binits{L.-H.}}
(\byear{2013}).
\btitle{Most tensor problems are NP-hard}.
\bjournal{Journal of the ACM (JACM)}
\bvolume{60}
\bpages{45}.
\end{barticle}
\endbibitem

\bibitem[\protect\citeauthoryear{Hoff}{2015}]{hoff2015multilinear}
\begin{barticle}[author]
\bauthor{\bsnm{Hoff},~\bfnm{Peter~D}\binits{P.~D.}}
(\byear{2015}).
\btitle{Multilinear tensor regression for longitudinal relational data}.
\bjournal{The annals of applied statistics}
\bvolume{9}
\bpages{1169}.
\end{barticle}
\endbibitem

\bibitem[\protect\citeauthoryear{Hong, Kolda and
  Duersch}{2018}]{hong2018generalized}
\begin{barticle}[author]
\bauthor{\bsnm{Hong},~\bfnm{David}\binits{D.}},
  \bauthor{\bsnm{Kolda},~\bfnm{Tamara~G}\binits{T.~G.}} \AND
  \bauthor{\bsnm{Duersch},~\bfnm{Jed~A}\binits{J.~A.}}
(\byear{2018}).
\btitle{Generalized canonical polyadic tensor decomposition}.
\bjournal{arXiv preprint arXiv:1808.07452}.
\end{barticle}
\endbibitem

\bibitem[\protect\citeauthoryear{Hopkins, Shi and
  Steurer}{2015}]{hopkins2015tensor}
\begin{binproceedings}[author]
\bauthor{\bsnm{Hopkins},~\bfnm{Samuel~B}\binits{S.~B.}},
  \bauthor{\bsnm{Shi},~\bfnm{Jonathan}\binits{J.}} \AND
  \bauthor{\bsnm{Steurer},~\bfnm{David}\binits{D.}}
(\byear{2015}).
\btitle{Tensor principal component analysis via sum-of-square proofs}.
In \bbooktitle{Proceedings of The 28th Conference on Learning Theory, COLT}
\bpages{3--6}.
\end{binproceedings}
\endbibitem

\bibitem[\protect\citeauthoryear{Javanmard
  et~al.}{2018}]{javanmard2018debiasing}
\begin{barticle}[author]
\bauthor{\bsnm{Javanmard},~\bfnm{Adel}\binits{A.}},
  \bauthor{\bsnm{Montanari},~\bfnm{Andrea}\binits{A.}} \betal{et~al.}
(\byear{2018}).
\btitle{Debiasing the lasso: Optimal sample size for gaussian designs}.
\bjournal{The Annals of Statistics}
\bvolume{46}
\bpages{2593--2622}.
\end{barticle}
\endbibitem

\bibitem[\protect\citeauthoryear{Jiang, Raskutti and
  Willett}{2015}]{jiang2015minimax}
\begin{barticle}[author]
\bauthor{\bsnm{Jiang},~\bfnm{Xin}\binits{X.}},
  \bauthor{\bsnm{Raskutti},~\bfnm{Garvesh}\binits{G.}} \AND
  \bauthor{\bsnm{Willett},~\bfnm{Rebecca}\binits{R.}}
(\byear{2015}).
\btitle{Minimax optimal rates for Poisson inverse problems with physical
  constraints}.
\bjournal{IEEE Transactions on Information Theory}
\bvolume{61}
\bpages{4458--4474}.
\end{barticle}
\endbibitem

\bibitem[\protect\citeauthoryear{Johndrow, Bhattacharya and
  Dunson}{2017}]{johndrow2017tensor}
\begin{barticle}[author]
\bauthor{\bsnm{Johndrow},~\bfnm{James~E}\binits{J.~E.}},
  \bauthor{\bsnm{Bhattacharya},~\bfnm{Anirban}\binits{A.}} \AND
  \bauthor{\bsnm{Dunson},~\bfnm{David~B}\binits{D.~B.}}
(\byear{2017}).
\btitle{Tensor decompositions and sparse log-linear models}.
\bjournal{The Annals of Statistics}
\bvolume{45}
\bpages{1--38}.
\end{barticle}
\endbibitem

\bibitem[\protect\citeauthoryear{Jolliffe}{1986}]{jolliffe1986principal}
\begin{bincollection}[author]
\bauthor{\bsnm{Jolliffe},~\bfnm{Ian~T}\binits{I.~T.}}
(\byear{1986}).
\btitle{Principal components in regression analysis}.
In \bbooktitle{Principal component analysis}
\bpages{129--155}.
\bpublisher{Springer}.
\end{bincollection}
\endbibitem

\bibitem[\protect\citeauthoryear{Keshavan, Montanari and
  Oh}{2010}]{keshavan2010matrix}
\begin{barticle}[author]
\bauthor{\bsnm{Keshavan},~\bfnm{Raghunandan~H}\binits{R.~H.}},
  \bauthor{\bsnm{Montanari},~\bfnm{Andrea}\binits{A.}} \AND
  \bauthor{\bsnm{Oh},~\bfnm{Sewoong}\binits{S.}}
(\byear{2010}).
\btitle{Matrix completion from noisy entries}.
\bjournal{Journal of Machine Learning Research}
\bvolume{11}
\bpages{2057--2078}.
\end{barticle}
\endbibitem

\bibitem[\protect\citeauthoryear{Kolda and Bader}{2009}]{kolda2009tensor}
\begin{barticle}[author]
\bauthor{\bsnm{Kolda},~\bfnm{Tamara~G}\binits{T.~G.}} \AND
  \bauthor{\bsnm{Bader},~\bfnm{Brett~W}\binits{B.~W.}}
(\byear{2009}).
\btitle{Tensor decompositions and applications}.
\bjournal{SIAM review}
\bvolume{51}
\bpages{455--500}.
\end{barticle}
\endbibitem

\bibitem[\protect\citeauthoryear{Koltchinskii, Lounici and
  Tsybakov}{2011}]{koltchinskii2011nuclear}
\begin{barticle}[author]
\bauthor{\bsnm{Koltchinskii},~\bfnm{Vladimir}\binits{V.}},
  \bauthor{\bsnm{Lounici},~\bfnm{Karim}\binits{K.}} \AND
  \bauthor{\bsnm{Tsybakov},~\bfnm{Alexandre~B}\binits{A.~B.}}
(\byear{2011}).
\btitle{Nuclear-norm penalization and optimal rates for noisy low-rank matrix
  completion}.
\bjournal{The Annals of Statistics}
\bvolume{39}
\bpages{2302--2329}.
\end{barticle}
\endbibitem

\bibitem[\protect\citeauthoryear{Krivanek, Dellby and
  Lupini}{1999}]{krivanek1999towards}
\begin{barticle}[author]
\bauthor{\bsnm{Krivanek},~\bfnm{OL}\binits{O.}},
  \bauthor{\bsnm{Dellby},~\bfnm{N}\binits{N.}} \AND
  \bauthor{\bsnm{Lupini},~\bfnm{AR}\binits{A.}}
(\byear{1999}).
\btitle{Towards sub-A electron beams}.
\bjournal{Ultramicroscopy}
\bvolume{78}
\bpages{1--11}.
\end{barticle}
\endbibitem

\bibitem[\protect\citeauthoryear{Kroonenberg}{2008}]{kroonenberg2008applied}
\begin{bbook}[author]
\bauthor{\bsnm{Kroonenberg},~\bfnm{Pieter~M}\binits{P.~M.}}
(\byear{2008}).
\btitle{Applied multiway data analysis}
\bvolume{702}.
\bpublisher{John Wiley \& Sons}.
\end{bbook}
\endbibitem

\bibitem[\protect\citeauthoryear{Lei, Chen and Lynch}{2019}]{lei2019consistent}
\begin{barticle}[author]
\bauthor{\bsnm{Lei},~\bfnm{Jing}\binits{J.}},
  \bauthor{\bsnm{Chen},~\bfnm{Kehui}\binits{K.}} \AND
  \bauthor{\bsnm{Lynch},~\bfnm{Brian}\binits{B.}}
(\byear{2019}).
\btitle{Consistent community detection in multi-layer network data}.
\bjournal{Biometrika}.
\end{barticle}
\endbibitem

\bibitem[\protect\citeauthoryear{Lesieur et~al.}{2017}]{lesieur2017statistical}
\begin{binproceedings}[author]
\bauthor{\bsnm{Lesieur},~\bfnm{Thibault}\binits{T.}},
  \bauthor{\bsnm{Miolane},~\bfnm{L{\'e}o}\binits{L.}},
  \bauthor{\bsnm{Lelarge},~\bfnm{Marc}\binits{M.}},
  \bauthor{\bsnm{Krzakala},~\bfnm{Florent}\binits{F.}} \AND
  \bauthor{\bsnm{Zdeborov{\'a}},~\bfnm{Lenka}\binits{L.}}
(\byear{2017}).
\btitle{Statistical and computational phase transitions in spiked tensor
  estimation}.
In \bbooktitle{2017 IEEE International Symposium on Information Theory (ISIT)}
\bpages{511--515}.
\bpublisher{IEEE}.
\end{binproceedings}
\endbibitem

\bibitem[\protect\citeauthoryear{Li and Li}{2010}]{li2010tensor}
\begin{binproceedings}[author]
\bauthor{\bsnm{Li},~\bfnm{Nan}\binits{N.}} \AND
  \bauthor{\bsnm{Li},~\bfnm{Baoxin}\binits{B.}}
(\byear{2010}).
\btitle{Tensor completion for on-board compression of hyperspectral images}.
In \bbooktitle{2010 IEEE International Conference on Image Processing}
\bpages{517--520}.
\bpublisher{IEEE}.
\end{binproceedings}
\endbibitem

\bibitem[\protect\citeauthoryear{Li and Zhang}{2017}]{li2017parsimonious}
\begin{barticle}[author]
\bauthor{\bsnm{Li},~\bfnm{Lexin}\binits{L.}} \AND
  \bauthor{\bsnm{Zhang},~\bfnm{Xin}\binits{X.}}
(\byear{2017}).
\btitle{Parsimonious tensor response regression}.
\bjournal{Journal of the American Statistical Association}
\bpages{1--16}.
\end{barticle}
\endbibitem

\bibitem[\protect\citeauthoryear{Li et~al.}{2018}]{li2018tucker}
\begin{barticle}[author]
\bauthor{\bsnm{Li},~\bfnm{Xiaoshan}\binits{X.}},
  \bauthor{\bsnm{Xu},~\bfnm{Da}\binits{D.}},
  \bauthor{\bsnm{Zhou},~\bfnm{Hua}\binits{H.}} \AND
  \bauthor{\bsnm{Li},~\bfnm{Lexin}\binits{L.}}
(\byear{2018}).
\btitle{Tucker tensor regression and neuroimaging analysis}.
\bjournal{Statistics in Biosciences}
\bvolume{10}
\bpages{520--545}.
\end{barticle}
\endbibitem

\bibitem[\protect\citeauthoryear{Liu et~al.}{2013}]{liu2013tensor}
\begin{barticle}[author]
\bauthor{\bsnm{Liu},~\bfnm{Ji}\binits{J.}},
  \bauthor{\bsnm{Musialski},~\bfnm{Przemyslaw}\binits{P.}},
  \bauthor{\bsnm{Wonka},~\bfnm{Peter}\binits{P.}} \AND
  \bauthor{\bsnm{Ye},~\bfnm{Jieping}\binits{J.}}
(\byear{2013}).
\btitle{Tensor completion for estimating missing values in visual data}.
\bjournal{IEEE Transactions on Pattern Analysis and Machine Intelligence}
\bvolume{35}
\bpages{208--220}.
\end{barticle}
\endbibitem

\bibitem[\protect\citeauthoryear{Lubich et~al.}{2013}]{lubich2013dynamical}
\begin{barticle}[author]
\bauthor{\bsnm{Lubich},~\bfnm{Christian}\binits{C.}},
  \bauthor{\bsnm{Rohwedder},~\bfnm{Thorsten}\binits{T.}},
  \bauthor{\bsnm{Schneider},~\bfnm{Reinhold}\binits{R.}} \AND
  \bauthor{\bsnm{Vandereycken},~\bfnm{Bart}\binits{B.}}
(\byear{2013}).
\btitle{Dynamical approximation by hierarchical Tucker and tensor-train
  tensors}.
\bjournal{SIAM Journal on Matrix Analysis and Applications}
\bvolume{34}
\bpages{470--494}.
\end{barticle}
\endbibitem

\bibitem[\protect\citeauthoryear{Ma and Ma}{2017}]{ma2017exploration}
\begin{barticle}[author]
\bauthor{\bsnm{Ma},~\bfnm{Zhuang}\binits{Z.}} \AND
  \bauthor{\bsnm{Ma},~\bfnm{Zongming}\binits{Z.}}
(\byear{2017}).
\btitle{Exploration of Large Networks via Fast and Universal Latent Space Model
  Fitting}.
\bjournal{arXiv preprint arXiv:1705.02372}.
\end{barticle}
\endbibitem

\bibitem[\protect\citeauthoryear{McMahan et~al.}{2013}]{mcmahan2013ad}
\begin{binproceedings}[author]
\bauthor{\bsnm{McMahan},~\bfnm{H~Brendan}\binits{H.~B.}},
  \bauthor{\bsnm{Holt},~\bfnm{Gary}\binits{G.}},
  \bauthor{\bsnm{Sculley},~\bfnm{David}\binits{D.}},
  \bauthor{\bsnm{Young},~\bfnm{Michael}\binits{M.}},
  \bauthor{\bsnm{Ebner},~\bfnm{Dietmar}\binits{D.}},
  \bauthor{\bsnm{Grady},~\bfnm{Julian}\binits{J.}},
  \bauthor{\bsnm{Nie},~\bfnm{Lan}\binits{L.}},
  \bauthor{\bsnm{Phillips},~\bfnm{Todd}\binits{T.}},
  \bauthor{\bsnm{Davydov},~\bfnm{Eugene}\binits{E.}},
  \bauthor{\bsnm{Golovin},~\bfnm{Daniel}\binits{D.}} \betal{et~al.}
(\byear{2013}).
\btitle{Ad click prediction: a view from the trenches}.
In \bbooktitle{Proceedings of the 19th ACM SIGKDD international conference on
  Knowledge discovery and data mining}
\bpages{1222--1230}.
\bpublisher{ACM}.
\end{binproceedings}
\endbibitem

\bibitem[\protect\citeauthoryear{Montanari, Reichman and
  Zeitouni}{2017}]{montanari2017limitation}
\begin{barticle}[author]
\bauthor{\bsnm{Montanari},~\bfnm{Andrea}\binits{A.}},
  \bauthor{\bsnm{Reichman},~\bfnm{Daniel}\binits{D.}} \AND
  \bauthor{\bsnm{Zeitouni},~\bfnm{Ofer}\binits{O.}}
(\byear{2017}).
\btitle{On the Limitation of Spectral Methods: From the Gaussian Hidden Clique
  Problem to Rank One Perturbations of Gaussian Tensors}.
\bjournal{IEEE Transactions on Information Theory}
\bvolume{63}
\bpages{1572--1579}.
\end{barticle}
\endbibitem

\bibitem[\protect\citeauthoryear{Montanari and
  Sun}{2018}]{montanari2018spectral}
\begin{barticle}[author]
\bauthor{\bsnm{Montanari},~\bfnm{Andrea}\binits{A.}} \AND
  \bauthor{\bsnm{Sun},~\bfnm{Nike}\binits{N.}}
(\byear{2018}).
\btitle{Spectral algorithms for tensor completion}.
\bjournal{Communications on Pure and Applied Mathematics}
\bvolume{71}
\bpages{2381--2425}.
\end{barticle}
\endbibitem

\bibitem[\protect\citeauthoryear{Nelder and
  Wedderburn}{1972}]{nelder1972generalized}
\begin{barticle}[author]
\bauthor{\bsnm{Nelder},~\bfnm{John~Ashworth}\binits{J.~A.}} \AND
  \bauthor{\bsnm{Wedderburn},~\bfnm{Robert~WM}\binits{R.~W.}}
(\byear{1972}).
\btitle{Generalized linear models}.
\bjournal{Journal of the Royal Statistical Society: Series A (General)}
\bvolume{135}
\bpages{370--384}.
\end{barticle}
\endbibitem

\bibitem[\protect\citeauthoryear{Nesterov}{1998}]{nesterov1998introductory}
\begin{barticle}[author]
\bauthor{\bsnm{Nesterov},~\bfnm{Yurii}\binits{Y.}}
(\byear{1998}).
\btitle{Introductory lectures on convex programming volume i: Basic course}.
\bjournal{Lecture notes}
\bvolume{3}
\bpages{5}.
\end{barticle}
\endbibitem

\bibitem[\protect\citeauthoryear{Oymak et~al.}{2015}]{oymak2015simultaneously}
\begin{barticle}[author]
\bauthor{\bsnm{Oymak},~\bfnm{Samet}\binits{S.}},
  \bauthor{\bsnm{Jalali},~\bfnm{Amin}\binits{A.}},
  \bauthor{\bsnm{Fazel},~\bfnm{Maryam}\binits{M.}},
  \bauthor{\bsnm{Eldar},~\bfnm{Yonina~C}\binits{Y.~C.}} \AND
  \bauthor{\bsnm{Hassibi},~\bfnm{Babak}\binits{B.}}
(\byear{2015}).
\btitle{Simultaneously structured models with application to sparse and
  low-rank matrices}.
\bjournal{IEEE Transactions on Information Theory}
\bvolume{61}
\bpages{2886--2908}.
\end{barticle}
\endbibitem

\bibitem[\protect\citeauthoryear{Park et~al.}{2018}]{park2018finding}
\begin{barticle}[author]
\bauthor{\bsnm{Park},~\bfnm{Dohyung}\binits{D.}},
  \bauthor{\bsnm{Kyrillidis},~\bfnm{Anastasios}\binits{A.}},
  \bauthor{\bsnm{Caramanis},~\bfnm{Constantine}\binits{C.}} \AND
  \bauthor{\bsnm{Sanghavi},~\bfnm{Sujay}\binits{S.}}
(\byear{2018}).
\btitle{Finding low-rank solutions via nonconvex matrix factorization,
  efficiently and provably}.
\bjournal{SIAM Journal on Imaging Sciences}
\bvolume{11}
\bpages{2165--2204}.
\end{barticle}
\endbibitem

\bibitem[\protect\citeauthoryear{Pensky et~al.}{2019}]{pensky2019spectral}
\begin{barticle}[author]
\bauthor{\bsnm{Pensky},~\bfnm{Marianna}\binits{M.}},
  \bauthor{\bsnm{Zhang},~\bfnm{Teng}\binits{T.}} \betal{et~al.}
(\byear{2019}).
\btitle{Spectral clustering in the dynamic stochastic block model}.
\bjournal{Electronic Journal of Statistics}
\bvolume{13}
\bpages{678--709}.
\end{barticle}
\endbibitem

\bibitem[\protect\citeauthoryear{Perry, Wein and
  Bandeira}{2016}]{perry2016statistical}
\begin{barticle}[author]
\bauthor{\bsnm{Perry},~\bfnm{Amelia}\binits{A.}},
  \bauthor{\bsnm{Wein},~\bfnm{Alexander~S}\binits{A.~S.}} \AND
  \bauthor{\bsnm{Bandeira},~\bfnm{Afonso~S}\binits{A.~S.}}
(\byear{2016}).
\btitle{Statistical limits of spiked tensor models}.
\bjournal{arXiv preprint arXiv:1612.07728}.
\end{barticle}
\endbibitem

\bibitem[\protect\citeauthoryear{Raskutti et~al.}{2019}]{raskutti2019convex}
\begin{barticle}[author]
\bauthor{\bsnm{Raskutti},~\bfnm{Garvesh}\binits{G.}},
  \bauthor{\bsnm{Yuan},~\bfnm{Ming}\binits{M.}},
  \bauthor{\bsnm{Chen},~\bfnm{Han}\binits{H.}} \betal{et~al.}
(\byear{2019}).
\btitle{Convex regularization for high-dimensional multiresponse tensor
  regression}.
\bjournal{The Annals of Statistics}
\bvolume{47}
\bpages{1554--1584}.
\end{barticle}
\endbibitem

\bibitem[\protect\citeauthoryear{Rauhut, Schneider and
  Stojanac}{2015}]{rauhut2015tensor}
\begin{bincollection}[author]
\bauthor{\bsnm{Rauhut},~\bfnm{Holger}\binits{H.}},
  \bauthor{\bsnm{Schneider},~\bfnm{Reinhold}\binits{R.}} \AND
  \bauthor{\bsnm{Stojanac},~\bfnm{{$\v{Z}$}eljka}\binits{v.}}
(\byear{2015}).
\btitle{Tensor completion in hierarchical tensor representations}.
In \bbooktitle{Compressed sensing and its applications}
\bpages{419--450}.
\bpublisher{Springer}.
\end{bincollection}
\endbibitem

\bibitem[\protect\citeauthoryear{Rauhut, Schneider and
  Stojanac}{2017}]{rauhut2017low}
\begin{barticle}[author]
\bauthor{\bsnm{Rauhut},~\bfnm{Holger}\binits{H.}},
  \bauthor{\bsnm{Schneider},~\bfnm{Reinhold}\binits{R.}} \AND
  \bauthor{\bsnm{Stojanac},~\bfnm{Zeljka}\binits{Z.}}
(\byear{2017}).
\btitle{Low rank tensor recovery via iterative hard thresholding}.
\bjournal{Linear Algebra and its Applications}
\bvolume{523}
\bpages{220--262}.
\end{barticle}
\endbibitem

\bibitem[\protect\citeauthoryear{Recht, Fazel and
  Parrilo}{2010}]{recht2010guaranteed}
\begin{barticle}[author]
\bauthor{\bsnm{Recht},~\bfnm{Benjamin}\binits{B.}},
  \bauthor{\bsnm{Fazel},~\bfnm{Maryam}\binits{M.}} \AND
  \bauthor{\bsnm{Parrilo},~\bfnm{Pablo~A}\binits{P.~A.}}
(\byear{2010}).
\btitle{Guaranteed minimum-rank solutions of linear matrix equations via
  nuclear norm minimization}.
\bjournal{SIAM review}
\bvolume{52}
\bpages{471--501}.
\end{barticle}
\endbibitem

\bibitem[\protect\citeauthoryear{Richard and
  Montanari}{2014}]{richard2014statistical}
\begin{binproceedings}[author]
\bauthor{\bsnm{Richard},~\bfnm{Emile}\binits{E.}} \AND
  \bauthor{\bsnm{Montanari},~\bfnm{Andrea}\binits{A.}}
(\byear{2014}).
\btitle{A statistical model for tensor PCA}.
In \bbooktitle{Advances in Neural Information Processing Systems}
\bpages{2897--2905}.
\end{binproceedings}
\endbibitem

\bibitem[\protect\citeauthoryear{Salmon et~al.}{2014}]{salmon2014poisson}
\begin{barticle}[author]
\bauthor{\bsnm{Salmon},~\bfnm{Joseph}\binits{J.}},
  \bauthor{\bsnm{Harmany},~\bfnm{Zachary}\binits{Z.}},
  \bauthor{\bsnm{Deledalle},~\bfnm{Charles-Alban}\binits{C.-A.}} \AND
  \bauthor{\bsnm{Willett},~\bfnm{Rebecca}\binits{R.}}
(\byear{2014}).
\btitle{Poisson noise reduction with non-local PCA}.
\bjournal{Journal of mathematical imaging and vision}
\bvolume{48}
\bpages{279--294}.
\end{barticle}
\endbibitem

\bibitem[\protect\citeauthoryear{Sewell and Chen}{2015}]{sewell2015latent}
\begin{barticle}[author]
\bauthor{\bsnm{Sewell},~\bfnm{Daniel~K}\binits{D.~K.}} \AND
  \bauthor{\bsnm{Chen},~\bfnm{Yuguo}\binits{Y.}}
(\byear{2015}).
\btitle{Latent space models for dynamic networks}.
\bjournal{Journal of the American Statistical Association}
\bvolume{110}
\bpages{1646--1657}.
\end{barticle}
\endbibitem

\bibitem[\protect\citeauthoryear{Shah and Yu}{2019}]{shah2019iterative}
\begin{barticle}[author]
\bauthor{\bsnm{Shah},~\bfnm{Devavrat}\binits{D.}} \AND
  \bauthor{\bsnm{Yu},~\bfnm{Christina~Lee}\binits{C.~L.}}
(\byear{2019}).
\btitle{Iterative Collaborative Filtering for Sparse Noisy Tensor Estimation}.
\bjournal{arXiv preprint arXiv:1908.01241}.
\end{barticle}
\endbibitem

\bibitem[\protect\citeauthoryear{Shan et~al.}{2016}]{shan2016predicting}
\begin{barticle}[author]
\bauthor{\bsnm{Shan},~\bfnm{Lili}\binits{L.}},
  \bauthor{\bsnm{Lin},~\bfnm{Lei}\binits{L.}},
  \bauthor{\bsnm{Sun},~\bfnm{Chengjie}\binits{C.}} \AND
  \bauthor{\bsnm{Wang},~\bfnm{Xiaolong}\binits{X.}}
(\byear{2016}).
\btitle{Predicting ad click-through rates via feature-based fully coupled
  interaction tensor factorization}.
\bjournal{Electronic Commerce Research and Applications}
\bvolume{16}
\bpages{30--42}.
\end{barticle}
\endbibitem

\bibitem[\protect\citeauthoryear{Shi, Zhou and Zhang}{2018}]{shi2018high}
\begin{barticle}[author]
\bauthor{\bsnm{Shi},~\bfnm{Pixu}\binits{P.}},
  \bauthor{\bsnm{Zhou},~\bfnm{Yuchen}\binits{Y.}} \AND
  \bauthor{\bsnm{Zhang},~\bfnm{Anru}\binits{A.}}
(\byear{2018}).
\btitle{High-dimensional Log-Error-in-Variable Regression with Applications to
  Microbial Compositional Data Analysis}.
\bjournal{arXiv preprint arXiv:1811.11709}.
\end{barticle}
\endbibitem

\bibitem[\protect\citeauthoryear{Signoretto
  et~al.}{2011}]{signoretto2011tensor}
\begin{barticle}[author]
\bauthor{\bsnm{Signoretto},~\bfnm{Marco}\binits{M.}},
  \bauthor{\bparticle{Van~de} \bsnm{Plas},~\bfnm{Raf}\binits{R.}},
  \bauthor{\bsnm{De~Moor},~\bfnm{Bart}\binits{B.}} \AND
  \bauthor{\bsnm{Suykens},~\bfnm{Johan~AK}\binits{J.~A.}}
(\byear{2011}).
\btitle{Tensor versus matrix completion: A comparison with application to
  spectral data}.
\bjournal{IEEE Signal Processing Letters}
\bvolume{18}
\bpages{403--406}.
\end{barticle}
\endbibitem

\bibitem[\protect\citeauthoryear{Sun and Li}{2016}]{sun2016sparse}
\begin{barticle}[author]
\bauthor{\bsnm{Sun},~\bfnm{Will~Wei}\binits{W.~W.}} \AND
  \bauthor{\bsnm{Li},~\bfnm{Lexin}\binits{L.}}
(\byear{2016}).
\btitle{Sparse Low-rank Tensor Response Regression}.
\bjournal{arXiv preprint arXiv:1609.04523}.
\end{barticle}
\endbibitem

\bibitem[\protect\citeauthoryear{Sun and Luo}{2015}]{sun2015guaranteed}
\begin{binproceedings}[author]
\bauthor{\bsnm{Sun},~\bfnm{Ruoyu}\binits{R.}} \AND
  \bauthor{\bsnm{Luo},~\bfnm{Zhi-Quan}\binits{Z.-Q.}}
(\byear{2015}).
\btitle{Guaranteed matrix completion via nonconvex factorization}.
In \bbooktitle{Foundations of Computer Science (FOCS), 2015 IEEE 56th Annual
  Symposium on}
\bpages{270--289}.
\bpublisher{IEEE}.
\end{binproceedings}
\endbibitem

\bibitem[\protect\citeauthoryear{Sun et~al.}{2017}]{sun2017provable}
\begin{barticle}[author]
\bauthor{\bsnm{Sun},~\bfnm{Will~Wei}\binits{W.~W.}},
  \bauthor{\bsnm{Lu},~\bfnm{Junwei}\binits{J.}},
  \bauthor{\bsnm{Liu},~\bfnm{Han}\binits{H.}} \AND
  \bauthor{\bsnm{Cheng},~\bfnm{Guang}\binits{G.}}
(\byear{2017}).
\btitle{Provable sparse tensor decomposition}.
\bjournal{Journal of the Royal Statistical Society: Series B (Statistical
  Methodology)}
\bvolume{79}
\bpages{899--916}.
\end{barticle}
\endbibitem

\bibitem[\protect\citeauthoryear{Tibshirani}{1996}]{tibshirani1996regression}
\begin{barticle}[author]
\bauthor{\bsnm{Tibshirani},~\bfnm{Robert}\binits{R.}}
(\byear{1996}).
\btitle{Regression shrinkage and selection via the lasso}.
\bjournal{Journal of the Royal Statistical Society. Series B (Methodological)}
\bpages{267--288}.
\end{barticle}
\endbibitem

\bibitem[\protect\citeauthoryear{Timmerman and
  Nowak}{1999}]{timmerman1999multiscale}
\begin{barticle}[author]
\bauthor{\bsnm{Timmerman},~\bfnm{Klaus}\binits{K.}} \AND
  \bauthor{\bsnm{Nowak},~\bfnm{Robert~David}\binits{R.~D.}}
(\byear{1999}).
\btitle{Multiscale modeling and estimation of Poisson processes with
  application to photon-limited imaging}.
\bjournal{IEEE Transactions on Information Theory}
\bvolume{45}
\bpages{846--842}.
\end{barticle}
\endbibitem

\bibitem[\protect\citeauthoryear{Tomioka and Suzuki}{2013}]{tomioka2013convex}
\begin{binproceedings}[author]
\bauthor{\bsnm{Tomioka},~\bfnm{Ryota}\binits{R.}} \AND
  \bauthor{\bsnm{Suzuki},~\bfnm{Taiji}\binits{T.}}
(\byear{2013}).
\btitle{Convex tensor decomposition via structured Schatten norm
  regularization}.
In \bbooktitle{Advances in neural information processing systems}
\bpages{1331--1339}.
\end{binproceedings}
\endbibitem

\bibitem[\protect\citeauthoryear{Tomioka et~al.}{2011}]{tomioka2011statistical}
\begin{binproceedings}[author]
\bauthor{\bsnm{Tomioka},~\bfnm{Ryota}\binits{R.}},
  \bauthor{\bsnm{Suzuki},~\bfnm{Taiji}\binits{T.}},
  \bauthor{\bsnm{Hayashi},~\bfnm{Kohei}\binits{K.}} \AND
  \bauthor{\bsnm{Kashima},~\bfnm{Hisashi}\binits{H.}}
(\byear{2011}).
\btitle{Statistical performance of convex tensor decomposition}.
In \bbooktitle{Advances in Neural Information Processing Systems}
\bpages{972--980}.
\end{binproceedings}
\endbibitem

\bibitem[\protect\citeauthoryear{Tu et~al.}{2016}]{tu2016low}
\begin{binproceedings}[author]
\bauthor{\bsnm{Tu},~\bfnm{Stephen}\binits{S.}},
  \bauthor{\bsnm{Boczar},~\bfnm{Ross}\binits{R.}},
  \bauthor{\bsnm{Simchowitz},~\bfnm{Max}\binits{M.}},
  \bauthor{\bsnm{Soltanolkotabi},~\bfnm{Mahdi}\binits{M.}} \AND
  \bauthor{\bsnm{Recht},~\bfnm{Ben}\binits{B.}}
(\byear{2016}).
\btitle{Low-rank Solutions of Linear Matrix Equations via Procrustes Flow}.
In \bbooktitle{International Conference on Machine Learning}
\bpages{964--973}.
\end{binproceedings}
\endbibitem

\bibitem[\protect\citeauthoryear{Vershynin}{2010}]{vershynin2010introduction}
\begin{barticle}[author]
\bauthor{\bsnm{Vershynin},~\bfnm{Roman}\binits{R.}}
(\byear{2010}).
\btitle{Introduction to the non-asymptotic analysis of random matrices}.
\bjournal{arXiv preprint arXiv:1011.3027}.
\end{barticle}
\endbibitem

\bibitem[\protect\citeauthoryear{Wang, Fischer and Song}{2017}]{wang2017three}
\begin{barticle}[author]
\bauthor{\bsnm{Wang},~\bfnm{Miaoyan}\binits{M.}},
  \bauthor{\bsnm{Fischer},~\bfnm{Jonathan}\binits{J.}} \AND
  \bauthor{\bsnm{Song},~\bfnm{Yun~S}\binits{Y.~S.}}
(\byear{2017}).
\btitle{Three-way clustering of multi-tissue multi-individual gene expression
  data using constrained tensor decomposition}.
\bjournal{bioRxiv}
\bpages{229245}.
\end{barticle}
\endbibitem

\bibitem[\protect\citeauthoryear{Wang and Li}{2018}]{wang2018learning}
\begin{barticle}[author]
\bauthor{\bsnm{Wang},~\bfnm{Miaoyan}\binits{M.}} \AND
  \bauthor{\bsnm{Li},~\bfnm{Lexin}\binits{L.}}
(\byear{2018}).
\btitle{Learning from Binary Multiway Data: Probabilistic Tensor Decomposition
  and its Statistical Optimality}.
\bjournal{arXiv preprint arXiv:1811.05076}.
\end{barticle}
\endbibitem

\bibitem[\protect\citeauthoryear{Wang and Zeng}{2019}]{wang2019multiway}
\begin{binproceedings}[author]
\bauthor{\bsnm{Wang},~\bfnm{Miaoyan}\binits{M.}} \AND
  \bauthor{\bsnm{Zeng},~\bfnm{Yuchen}\binits{Y.}}
(\byear{2019}).
\btitle{Multiway clustering via tensor block models}.
In \bbooktitle{Advances in Neural Information Processing Systems}
\bpages{713--723}.
\end{binproceedings}
\endbibitem

\bibitem[\protect\citeauthoryear{Wen, Yin and Zhang}{2012}]{wen2012solving}
\begin{barticle}[author]
\bauthor{\bsnm{Wen},~\bfnm{Zaiwen}\binits{Z.}},
  \bauthor{\bsnm{Yin},~\bfnm{Wotao}\binits{W.}} \AND
  \bauthor{\bsnm{Zhang},~\bfnm{Yin}\binits{Y.}}
(\byear{2012}).
\btitle{Solving a low-rank factorization model for matrix completion by a
  nonlinear successive over-relaxation algorithm}.
\bjournal{Mathematical Programming Computation}
\bvolume{4}
\bpages{333--361}.
\end{barticle}
\endbibitem

\bibitem[\protect\citeauthoryear{Willett and
  Nowak}{2007}]{willett2007multiscale}
\begin{barticle}[author]
\bauthor{\bsnm{Willett},~\bfnm{Rebecca~M}\binits{R.~M.}} \AND
  \bauthor{\bsnm{Nowak},~\bfnm{Robert~D}\binits{R.~D.}}
(\byear{2007}).
\btitle{Multiscale Poisson intensity and density estimation}.
\bjournal{IEEE Transactions on Information Theory}
\bvolume{53}
\bpages{3171--3187}.
\end{barticle}
\endbibitem

\bibitem[\protect\citeauthoryear{Wilmoth and
  Shkolnikov}{2006}]{wilmoth2016human}
\begin{bmisc}[author]
\bauthor{\bsnm{Wilmoth},~\bfnm{J.~R.}\binits{J.~R.}} \AND
  \bauthor{\bsnm{Shkolnikov},~\bfnm{V.}\binits{V.}}
(\byear{2006}).
\btitle{Human mortality database, available at: http://www.mortality.org}.
\end{bmisc}
\endbibitem

\bibitem[\protect\citeauthoryear{Xia and Yuan}{2017}]{xia2017polynomial}
\begin{barticle}[author]
\bauthor{\bsnm{Xia},~\bfnm{Dong}\binits{D.}} \AND
  \bauthor{\bsnm{Yuan},~\bfnm{Ming}\binits{M.}}
(\byear{2017}).
\btitle{On Polynomial Time Methods for Exact Low Rank Tensor Completion}.
\bjournal{arXiv preprint arXiv:1702.06980}.
\end{barticle}
\endbibitem

\bibitem[\protect\citeauthoryear{Xia, Yuan and
  Zhang}{2017}]{xia2017statistically}
\begin{barticle}[author]
\bauthor{\bsnm{Xia},~\bfnm{Dong}\binits{D.}},
  \bauthor{\bsnm{Yuan},~\bfnm{Ming}\binits{M.}} \AND
  \bauthor{\bsnm{Zhang},~\bfnm{Cun-Hui}\binits{C.-H.}}
(\byear{2017}).
\btitle{Statistically Optimal and Computationally Efficient Low Rank Tensor
  Completion from Noisy Entries}.
\bjournal{arXiv preprint arXiv:1711.04934}.
\end{barticle}
\endbibitem

\bibitem[\protect\citeauthoryear{Xu, Hu and Wang}{2019}]{xu2019generalized}
\begin{barticle}[author]
\bauthor{\bsnm{Xu},~\bfnm{Zhuoyan}\binits{Z.}},
  \bauthor{\bsnm{Hu},~\bfnm{Jiaxin}\binits{J.}} \AND
  \bauthor{\bsnm{Wang},~\bfnm{Miaoyan}\binits{M.}}
(\byear{2019}).
\btitle{Generalized tensor regression with covariates on multiple modes}.
\bjournal{arXiv preprint arXiv:1910.09499}.
\end{barticle}
\endbibitem

\bibitem[\protect\citeauthoryear{Yang and Barron}{1999}]{yang1999information}
\begin{barticle}[author]
\bauthor{\bsnm{Yang},~\bfnm{Yuhong}\binits{Y.}} \AND
  \bauthor{\bsnm{Barron},~\bfnm{Andrew}\binits{A.}}
(\byear{1999}).
\btitle{Information-theoretic determination of minimax rates of convergence}.
\bjournal{Annals of Statistics}
\bpages{1564--1599}.
\end{barticle}
\endbibitem

\bibitem[\protect\citeauthoryear{Yankovich et~al.}{2016}]{yankovich2016non}
\begin{barticle}[author]
\bauthor{\bsnm{Yankovich},~\bfnm{Andrew~B}\binits{A.~B.}},
  \bauthor{\bsnm{Zhang},~\bfnm{Chenyu}\binits{C.}},
  \bauthor{\bsnm{Oh},~\bfnm{Albert}\binits{A.}},
  \bauthor{\bsnm{Slater},~\bfnm{Thomas~JA}\binits{T.~J.}},
  \bauthor{\bsnm{Azough},~\bfnm{Feridoon}\binits{F.}},
  \bauthor{\bsnm{Freer},~\bfnm{Robert}\binits{R.}},
  \bauthor{\bsnm{Haigh},~\bfnm{Sarah~J}\binits{S.~J.}},
  \bauthor{\bsnm{Willett},~\bfnm{Rebecca}\binits{R.}} \AND
  \bauthor{\bsnm{Voyles},~\bfnm{Paul~M}\binits{P.~M.}}
(\byear{2016}).
\btitle{Non-rigid registration and non-local principle component analysis to
  improve electron microscopy spectrum images}.
\bjournal{Nanotechnology}
\bvolume{27}
\bpages{364001}.
\end{barticle}
\endbibitem

\bibitem[\protect\citeauthoryear{Yokota, Lee and
  Cichocki}{2016}]{yokota2016robust}
\begin{barticle}[author]
\bauthor{\bsnm{Yokota},~\bfnm{Tatsuya}\binits{T.}},
  \bauthor{\bsnm{Lee},~\bfnm{Namgil}\binits{N.}} \AND
  \bauthor{\bsnm{Cichocki},~\bfnm{Andrzej}\binits{A.}}
(\byear{2016}).
\btitle{Robust multilinear tensor rank estimation using higher order singular
  value decomposition and information criteria}.
\bjournal{IEEE Transactions on Signal Processing}
\bvolume{65}
\bpages{1196--1206}.
\end{barticle}
\endbibitem

\bibitem[\protect\citeauthoryear{Yonel and
  Yazici}{2020}]{yonel2020deterministic}
\begin{barticle}[author]
\bauthor{\bsnm{Yonel},~\bfnm{Bariscan}\binits{B.}} \AND
  \bauthor{\bsnm{Yazici},~\bfnm{Birsen}\binits{B.}}
(\byear{2020}).
\btitle{A Deterministic Convergence Framework for Exact Non-Convex Phase
  Retrieval}.
\bjournal{arXiv preprint arXiv:2001.02855}.
\end{barticle}
\endbibitem

\bibitem[\protect\citeauthoryear{Yu et~al.}{2018}]{yu2018recovery}
\begin{barticle}[author]
\bauthor{\bsnm{Yu},~\bfnm{Ming}\binits{M.}},
  \bauthor{\bsnm{Wang},~\bfnm{Zhaoran}\binits{Z.}},
  \bauthor{\bsnm{Gupta},~\bfnm{Varun}\binits{V.}} \AND
  \bauthor{\bsnm{Kolar},~\bfnm{Mladen}\binits{M.}}
(\byear{2018}).
\btitle{Recovery of simultaneous low rank and two-way sparse coefficient
  matrices, a nonconvex approach}.
\bjournal{arXiv preprint arXiv:1802.06967}.
\end{barticle}
\endbibitem

\bibitem[\protect\citeauthoryear{Yuan and Zhang}{2014}]{yuan2014tensor}
\begin{barticle}[author]
\bauthor{\bsnm{Yuan},~\bfnm{Ming}\binits{M.}} \AND
  \bauthor{\bsnm{Zhang},~\bfnm{Cun-Hui}\binits{C.-H.}}
(\byear{2014}).
\btitle{On tensor completion via nuclear norm minimization}.
\bjournal{Foundations of Computational Mathematics}
\bpages{1--38}.
\end{barticle}
\endbibitem

\bibitem[\protect\citeauthoryear{Zhang}{2019}]{zhang2019cross}
\begin{barticle}[author]
\bauthor{\bsnm{Zhang},~\bfnm{Anru}\binits{A.}}
(\byear{2019}).
\btitle{Cross: Efficient low-rank tensor completion}.
\bjournal{The Annals of Statistics}
\bvolume{47}
\bpages{936--964}.
\end{barticle}
\endbibitem

\bibitem[\protect\citeauthoryear{Zhang, Cai and
  Wu}{2018}]{zhang2018heteroskedastic}
\begin{barticle}[author]
\bauthor{\bsnm{Zhang},~\bfnm{Anru}\binits{A.}},
  \bauthor{\bsnm{Cai},~\bfnm{T~Tony}\binits{T.~T.}} \AND
  \bauthor{\bsnm{Wu},~\bfnm{Yihong}\binits{Y.}}
(\byear{2018}).
\btitle{Heteroskedastic PCA: Algorithm, optimality, and applications}.
\bjournal{arXiv preprint arXiv:1810.08316}.
\end{barticle}
\endbibitem

\bibitem[\protect\citeauthoryear{Zhang and
  Han}{2018}]{zhang2017optimal-statsvd}
\begin{barticle}[author]
\bauthor{\bsnm{Zhang},~\bfnm{Anru}\binits{A.}} \AND
  \bauthor{\bsnm{Han},~\bfnm{Rungang}\binits{R.}}
(\byear{2018}).
\btitle{Optimal sparse singular value decomposition for high-dimensional
  high-order data}.
\bjournal{Journal of the American Statistical Association}
\bpages{to appear}.
\end{barticle}
\endbibitem

\bibitem[\protect\citeauthoryear{Zhang and Xia}{2018}]{zhang2018tensor}
\begin{barticle}[author]
\bauthor{\bsnm{Zhang},~\bfnm{Anru}\binits{A.}} \AND
  \bauthor{\bsnm{Xia},~\bfnm{Dong}\binits{D.}}
(\byear{2018}).
\btitle{Tensor {SVD}: Statistical and Computational Limits}.
\bjournal{IEEE Transactions on Information Theory}
\bvolume{64}
\bpages{7311-7338}.
\end{barticle}
\endbibitem

\bibitem[\protect\citeauthoryear{Zhang et~al.}{2019}]{zhang2018ISLET}
\begin{barticle}[author]
\bauthor{\bsnm{Zhang},~\bfnm{Anru}\binits{A.}},
  \bauthor{\bsnm{Luo},~\bfnm{Yuetian}\binits{Y.}},
  \bauthor{\bsnm{Raskutti},~\bfnm{Garvesh}\binits{G.}} \AND
  \bauthor{\bsnm{Yuan},~\bfnm{Ming}\binits{M.}}
(\byear{2019}).
\btitle{{ISLET}: Fast and Optimal Low-rank Tensor Regression via Importance
  Sketching}.
\bjournal{arXiv preprint arXiv:1911.03804}.
\end{barticle}
\endbibitem

\bibitem[\protect\citeauthoryear{Zhao, Wang and Liu}{2015}]{zhao2015nonconvex}
\begin{binproceedings}[author]
\bauthor{\bsnm{Zhao},~\bfnm{Tuo}\binits{T.}},
  \bauthor{\bsnm{Wang},~\bfnm{Zhaoran}\binits{Z.}} \AND
  \bauthor{\bsnm{Liu},~\bfnm{Han}\binits{H.}}
(\byear{2015}).
\btitle{A nonconvex optimization framework for low rank matrix estimation}.
In \bbooktitle{Advances in Neural Information Processing Systems}
\bpages{559--567}.
\end{binproceedings}
\endbibitem

\bibitem[\protect\citeauthoryear{Zhou}{2017}]{zhou2017matlab}
\begin{bmisc}[author]
\bauthor{\bsnm{Zhou},~\bfnm{Hua}\binits{H.}}
(\byear{2017}).
\btitle{Matlab TensorReg Toolbox Version 1.0}.
\bnote{Available online at https://hua-zhou.github.io/TensorReg/}.
\end{bmisc}
\endbibitem

\bibitem[\protect\citeauthoryear{Zhou, Li and Zhu}{2013}]{zhou2013tensor}
\begin{barticle}[author]
\bauthor{\bsnm{Zhou},~\bfnm{Hua}\binits{H.}},
  \bauthor{\bsnm{Li},~\bfnm{Lexin}\binits{L.}} \AND
  \bauthor{\bsnm{Zhu},~\bfnm{Hongtu}\binits{H.}}
(\byear{2013}).
\btitle{Tensor regression with applications in neuroimaging data analysis}.
\bjournal{Journal of the American Statistical Association}
\bvolume{108}
\bpages{540--552}.
\end{barticle}
\endbibitem

\bibitem[\protect\citeauthoryear{Zhu et~al.}{2017}]{zhu2017global}
\begin{barticle}[author]
\bauthor{\bsnm{Zhu},~\bfnm{Zhihui}\binits{Z.}},
  \bauthor{\bsnm{Li},~\bfnm{Qiuwei}\binits{Q.}},
  \bauthor{\bsnm{Tang},~\bfnm{Gongguo}\binits{G.}} \AND
  \bauthor{\bsnm{Wakin},~\bfnm{Michael~B}\binits{M.~B.}}
(\byear{2017}).
\btitle{The global optimization geometry of nonsymmetric matrix factorization
  and sensing}.
\bjournal{arXiv preprint arXiv:1703.01256}.
\end{barticle}
\endbibitem

\end{thebibliography}

\newpage
\appendix

\setcounter{page}{1}

\begin{center}
    	{\LARGE Supplement to ``An Optimal Statistical and Computational Framework for Generalized Tensor Estimation"	
	}

	\bigskip\medskip
	{Rungang Han, ~ Rebecca Willett, ~ and ~ Anru R. Zhang}
\end{center}

In this supplement, we provide the implementation details of HeteroPCA, higher-order orthogonal iteration (HOOI), higher-order SVD (HOSVD), additional real data example, proofs of all the technical results, and the key technical lemmas.

\section{Implementation of Additional Algorithms}\label{sec:implementaionts}

We collect the implementations of HeteroPCA \citep{zhang2018heteroskedastic}, higher-order SVD (HOSVD), and higher-order orthogonal iteration (HOOI) \citep{de2000best, anandkumar2012method} in this section. For any square matrix $\A$, let $\Delta(\A)$ be $\A$ with all diagonal entries set to zero and $D(\A)$ be $\A$ with all off-diagonal entries set to zero.
\begin{algorithm}
	\caption{Heteroskedastic PCA (HeteroPCA)}
	\begin{algorithmic}
	\REQUIRE{symmetric matrix $\mathbf{\hat \Sigma}$, rank $r$, max iteration time $t_{max}$}
    \STATE{Set $\N^{(0)} = \Delta(\mathbf{\hat \Sigma})$.}
    \FORALL{$t = 1,\ldots, t_{max}$}
	\STATE{Calculate SVD: $\N^{(t)} = \sum_i \lambda_i u_i^{(t)} (v_i^{(t)})^\top$, where $\lambda_1 \geq \lambda_2 \geq \ldots \geq 0$.}
	\STATE{Let $\tilde \N^{(t)} = \sum_{i=1}^r \lambda_i^{(t)}u_i^{(t)}(v_i^{(t)})^\top$.}
	\STATE{Update diagonal entries $\N^{(t+1)} = D(\tilde \N^{(t)}) + \Delta(\N^{(t)})$.}
	\ENDFOR
    \RETURN{$\U = \left(u_1^{(t_{max})}, \ldots, u_r^{(t_{max})}\right)$}
	\end{algorithmic}\label{alg:HeteroPCA}
\end{algorithm}
\begin{algorithm}
	\caption{Higher-order Singular Value Decomposition (HOSVD)}
	\begin{algorithmic}
	\REQUIRE{$\cY \in \bbR^{p_1\times p_2 \times p_3}$, Tucker rank $(r_1,r_2,r_3)$}
    \STATE{$\U_k = \SVD_{r_k}\left(\cM_k(\cY)\right)$,\quad for $k=1,2,3$}
    \STATE{$\cS = \llbracket \cY; \U_1^\top, \U_2^\top, \U_3^\top \rrbracket$}
    \RETURN{$\left(\cS, \U_1, \U_2, \U_3\right)$}
	\end{algorithmic}\label{alg:HOSVD}
\end{algorithm}
\begin{algorithm}
	\caption{Higher-order Orthogonal Iteration (HOOI)}
	\begin{algorithmic}
	\REQUIRE{$\cY \in \bbR^{p_1\times p_2 \times p_3}$, Tucker rank $(r_1,r_2,r_3)$, max iteration $t_{max}$}
    \STATE{Initialize $\U_k^{(0)} = \SVD_{r_k}\left(\cM_k(\cY)\right)$,\quad for $k=1,2,3$}
    \FORALL{$t = 1,\ldots t_{max}$}
	\FORALL{$k = 1,2, 3$}
	\STATE{$ \A_k^{(t)} = \mathcal{M}_k\left(\cY\times_{k+1}  \U_{k+1}^{(t-1)\top} \times_{k+2}  \U_{k+2}^{(t-1)\top}\right) $}
	\STATE{$\U_k^{(t)} = \SVD_{r_k}\left(\A_k^{(t)}\right)$}
	\ENDFOR
	\STATE $t = t+1$
	\ENDFOR
    \RETURN{$\left( \cS^{(t_{max})}, \U_1^{(t_{max})}, \U_2^{(t_{max})}, \U_3^{(t_{max})}\right)$}
	\end{algorithmic}\label{alg:HOOI}
\end{algorithm}

\newpage
\section{Additional Real Data Example}\label{sec:click-through}
We study the prediction of users' online click-through behavior on \emph{Taobao.com}, one of the most popular online shopping website in China. The data\footnote{Available at:  \url{https://tianchi.aliyun.com/dataset/dataDetail?dataId=649}} are collected from Nov 25 to Dec 02, 2017 and arranged into 32 periods as each day is divided into four periods: 00:00 -- 06:00, 06:00 -- 12:00, 12:00 -- 18:00, and 18:00 -- 24:00. By this means, the dataset is in the form of a count-valued third-order tensor, where the $(i, j, k)$th entry represents the total count of clicks by the $i$th user on $j$th item-category in the $k$th period along eight days. Due to the high-dimensionality of the original dataset ($\approx 10^6$ users and $\approx 10^4$ item categories), we only focus on the most active 100 users and the most popular 50 categories. 

To investigate the predictive power of the proposed procedure, we consider all $\binom{8}{4}=70$ even partitions of the eight days: $\pi\cup\pi^c=\{1,\ldots, 8\}$, $\pi\cap \pi=\emptyset$. For each $\pi$, we aggregate the dataset into two tensors $\cY_1^\pi, \cY_2^\pi \in \bbN^{100 \times 50 \times 4}$, where $(\cY_1^\pi)_{ijk}$ and $(\cY_2^\pi)_{ijk}$ are the sums of count clicks made by the $i$th customer on the $j$th item-category in the $k$th daily time interval on the days of $\pi$ and $\pi^c$, respectively. 
We apply the proposed method (Algorithms \ref{alg:PGD} and \ref{alg:TD_Poisson} in Poisson tensor PCA), HOSVD \citep{de2000best}, and HOOI \citep{de2000multilinear} respectively on $\cY_1^\pi$ to obtain the tensor estimator $\hat{\cY}_1^\pi$. 
Then, we evaluate both the training error $\|\cY_1^\pi - \hat \cY_1^\pi\|_\tF / \sqrt{100\cdot 50\cdot 4}$ and the prediction error $\|\cY_2^\pi - \hat \cY_1^\pi\|_\tF / \sqrt{100\cdot 50\cdot 4}$ for each $\pi$ and provide both the average error and standard deviation in Table \ref{tab:CTC}. As we can see, although the proposed method yields a larger training error, there is a significant advantage in the prediction error to the classic Gaussian-likelihood-based methods (HOSVD or HOOI). This data analysis also illustrates the merit of the generalized tensor estimation framework for handling non-Gaussian tensor data.
\begin{table}
	\centering
	\begin{tabular}{l|l|l}
		\hline
		Methods       & Training Error & Testing Error 
		\\ \hline
        HOSVD  & 4.43(0.18) & 5.34(0.24) 
        \\ \hline
        HOOI  & \textbf{4.31}(0.15) & 5.33(0.25)
        \\ \hline
        Poisson-PCA  & 4.91(0.17) & \textbf{5.27}(0.22)
        \\ \hline 
	\end{tabular}
	\caption{Average training and prediction errors of click-through counts. The standard error of prediction errors is provided in the parentheses. After numerical explorations, we set the Tucker-rank to $(5,5,2)$ in these experiments.}
	\label{tab:CTC}
\end{table}

\newpage

\section{Proof of Theorem \ref{thm:local-convergence-PGD}}\label{sec:proof-main}
In this section, we provide the proof of technical results on error contraction. We divide the proof into five steps. In Step 1, we introduce the notations and conditions that are used to develop the theory. Then in Steps 2, 3, 4, we prove the one-step error contraction and provide the convergence analysis. Finally in Step 5, we verify the conditions imposed in Step 1. 
\begin{enumerate}[leftmargin=*]
	\item[Step 1] (Notations and Conditions) To simplify the rest of the proof, we assume $b = \overline{\lambda}^{1/4}$ holds, since the following argument  also holds by changing the absolute constants when we are under the condition that $c\overline{\lambda}^{1/4} \leq b \leq C\overline{\lambda}^{1/4}$ for fixed $c$ and $C$. We first introduce or rephrase the following list of notations. Let $\cX^*$ be the target low-rank tensor that satisfies $\cX^* = \llbracket\cS^*; \U_1^*, \U_2^*, \U_3^*\rrbracket$ such that $\U_k^{*\top}\U_k^* = b^2\I_{r_k}$, $\U_k^* \in \cC_k$ for $k=1,2,3$. For each step $t=0,1,\ldots, t_{max}$, we define
	\begin{enumerate}[label=(\subscript{N}{\arabic*})]
		\item Error measurement
			\begin{equation*}
				\begin{split}
					E^{(t)} = \min_{\substack{\R_k \in \bbO_{p_k,r_k}\\ k = 1,2,3}} \left\{\sum_{k=1}^3\left \| \U_k^{(t)} - \U_k^* \R_k\right\|_\tF^2 + \left\| \cS^{(t)} - \llbracket\cS^*; \R_1^\top, \R_2^\top, \R_3^\top\rrbracket\right\|_\tF^2\right\}
				\end{split}
			\end{equation*}
			\begin{equation*}
				\begin{split}
				& (\R_1^{(t)},\R_2^{(t)},\R_3^{(t)}) = \argmin_{\substack{\R_k \in \bbO_{p_k, r_k}\\ k = 1,2,3}} \left\{\sum_{k=1}^3\left \| \U_k^{(t)} - \U_k^* \R_k\right\|_\tF^2 + \left\|\cS^{(t)} - \llbracket\cS^*; \R_1^\top, \R_2^\top, \R_3^\top\rrbracket\right\|_\tF^2\right\} \\
				& (\tilde\R_1^{(t)},\tilde\R_2^{(t)},\tilde\R_3^{(t)}) = \argmin_{\substack{\R_k \in \bbO_{p_k, r_k}\\ k = 1,2,3}} \left\{\sum_{k=1}^3\left \| \tilde\U_k^{(t)} - \U_k^* \R_k\right\|_\tF^2 + \left\|\tilde\cS^{(t)} - \llbracket\cS^*; \R_1^\top, \R_2^\top, \R_3^\top\rrbracket\right\|_\tF^2\right\}
				\end{split}
			\end{equation*}
		\item Dual loadings
			\begin{equation}\label{eq:def-V}
				\begin{split}
					 \V_k^{(t)} = ( \U_{k+2}^{(t)} \otimes  \U_{k+1}^{(t)}) \left(\cM_k( \cS^{(t)})\right)^\top , \quad k=1,2,3
				\end{split}
			\end{equation}
		\item Signal tensors
			\begin{equation}\label{eq:def-signal-tensor}
				\begin{split}
					&\cX^{(t)} = \llbracket\cS^{(t)};  \U_1^{(t)}, \U_2^{(t)}, \U_3^{(t)}\rrbracket\\
					&\cX_\cS^{(t)} = \llbracket\cS^*; \U_1^{(t)}\R_1^{(t)\top}, \U_2^{(t)} \R_2^{(t)\top}, \U_3^{(t)}\R_3^{(t)\top}\rrbracket\\
					& \cX_k^{(t)} = \cS^{(t)} \times_k \U_k^*\R_k^{(t)} \times_{k+1} \U_{k+1}^{(t)} \times_{k+2} \U_{k+2}^{(t)}, \quad k=1,2,3
				\end{split}
			\end{equation}
	\end{enumerate}
	We also assume the following conditions hold: (They are verified in Step 5)
	\begin{enumerate}[label=(\subscript{A}{\arabic*})]
			\item For any $t=0,1,\ldots,t_{max}$, we have
				\begin{equation}\label{eq-Assumption1}
					\begin{split}
						&\left\| \U_k^{(t)} \right\| \leq 1.01b, \quad  \left\|\cM_k\left( \cS^{(t)}\right)\right\| \leq \frac{1.01\overline\lambda}{b^3},\quad  k=1,2,3.
					\end{split}
				\end{equation}
			\item For $c_0 = \frac{1}{20000}$ and any $t=0,1,\ldots,t_{max}$, we have	\begin{equation}\label{eq-Assumption2}
					E^{(t)} \leq c_0\frac{\alpha\beta \sqrt{\underline{\lambda}}}{\kappa^{3/2}}.
				\end{equation}
			Note that by RCG-condition \eqref{eq-restricted-convex}, one has
	    \begin{equation*}
	        \begin{split}
	            & \alpha \|\cX-\cX^*\|_\tF^2 + \frac{1}{4\alpha} \|\nabla f(\cX) - \nabla f(\cX^*)\|_\tF^2   \\
	            & \qquad \geq \|\cX-\cX^*\|_\tF \cdot \|\nabla f(\cX) - f(\cX^*)\|_\tF \geq \alpha\|\cX-\cX^*\|_\tF^2 + \beta\|\nabla f(\cX) - \nabla f(\cX^*)\|_\tF^2
	        \end{split}
	    \end{equation*}
	    and thus $\alpha\beta \leq 1/4$. Therefore, \eqref{eq-Assumption2} further implies 
	    \begin{equation*}\label{eq-Assumption2-imply}
	        E^{(t)}\leq \frac{c_0\sqrt{\underline \lambda}}{4\kappa^{3/2}}.
	    \end{equation*}
	\end{enumerate}
	\item[Step 2] (Descent of $E^{(t)}$) In this step, we show that under conditions \eqref{eq-Assumption1} and \eqref{eq-Assumption2},
	\begin{equation}\label{eq-Et-descent}
	    E^{(t+1)} \leq E^{(t)} - 2\eta\left(Q_{\cS,1}+\sum_{k=1}^3 Q_{k,1}\right) + \eta^2\left(Q_{\cS,2} + \sum_{k=1}^3 Q_{k,2}\right),
	\end{equation}
where
\begin{equation}\label{eq:D_k1,D_k2}
	\begin{split}
		& Q_{k,1} = \left\langle \cX^{(t)} - \cX_k^{(t)}, \nabla L(\cX^{(t)}) \right\rangle + \frac{a}{4} \left\| \U_k^{(t)\top} \U_k^{(t)} - b^2\I_{r_k}\right\|_\tF^2 - \frac{c_0a\underline{\lambda}^{1/2}}{16}\left\| \U_k^{(t)} - \U_k^* \R_k^{(t)}\right\|_\tF^2, \\
		& Q_{k,2} = 8\overline\lambda^{2}b^{-2}\left(\xi^2 + \left\|\nabla L(\cX^{(t)}) - \nabla L(\cX^*)\right\|_\tF^2  \right) + \frac{5}{2}a^2b^2 \left\| \U_k^{(t)\top} \U_k^{(t)}-b^2 \I_{r_k}\right\|_\tF^2, \\
		& Q_{\cS,1} = \left\langle \cX^{(t)} - \cX_\cS^{(t)}, \nabla L(\cX^{(t)}) \right\rangle, \\
		& Q_{\cS,2} = 4b^6\left(\xi^2 + \left\|\nabla L(\cX^{(t)}) - \nabla L (\cX^*)\right\|_\tF^2\right).
	\end{split}
\end{equation}
By definition of $E^{(t+1)}$, we first have
\begin{equation}\label{eq-link-step}
    \begin{split}
        E^{(t+1)} &= \sum_{k=1}^3\left \| \U_k^{(t+1)} - \U_k^* \R_k^{(t+1)}\right\|_\tF^2 + \left\|\cS^{(t+1)} - \llbracket\cS^*; \R_1^{(t+1)\top}, \R_2^{(t+1)\top}, \R_3^{(t+1)\top}\rrbracket\right\|_\tF^2 \\
        & \overset{(a)}{\leq} \sum_{k=1}^3\left \| \U_k^{(t+1)} - \U_k^* \tilde\R_k^{(t+1)}\right\|_\tF^2 + \left\|\cS^{(t+1)} - \llbracket\cS^*; \tilde\R_1^{(t+1)\top}, \tilde\R_2^{(t+1)\top}, \tilde\R_3^{(t+1)\top}\rrbracket\right\|_\tF^2 \\
        & \overset{(b)}{\leq} \sum_{k=1}^3\left \| \tilde\U_k^{(t+1)} - \U_k^* \tilde\R_k^{(t+1)}\right\|_\tF^2 + \left\|\tilde\cS^{(t+1)} - \llbracket\cS^*; \tilde\R_1^{(t+1)\top}, \tilde\R_2^{(t+1)\top}, \tilde\R_3^{(t+1)\top}\rrbracket\right\|_\tF^2 \\
        & \overset{(c)}{\leq} \sum_{k=1}^3\left \| \tilde\U_k^{(t+1)} - \U_k^* \R_k^{(t)}\right\|_\tF^2 + \left\|\tilde\cS^{(t+1)} - \llbracket \cS^*; \R_1^{(t)\top}, \R_2^{(t)\top}, \R_3^{(t)\top}\rrbracket\right\|_\tF^2.
    \end{split}
\end{equation}
Here, (a) and (c) comes from the definitions of $(\R_1^{(t+1)},\R_2^{(t+1)},\R_3^{(t+1)})$ and $(\tilde\R_1^{(t+1)},\tilde\R_2^{(t+1)},\tilde\R_3^{(t+1)})$, while (b) comes from the projection step $\U_k^{(t+1)} = \cP_{\cC_k}(\tilde \U_k^{(t+1)}), \cS^{(t+1)} = \cP_{\cC_\cS}(\tilde \cS^{(t+1)})$ as $\cC_k$ and $\cC_\cS$ are convex and rotation invariant sets.

Now we analyze the error on loadings and core tensor separately. Specifically, we will show:
\begin{equation}\label{eq-decompose-U-error-t+1-result}
	\begin{split}
		& \left\|\tilde\U_k^{(t+1)} - \U_k^*\R_k^{(t)}\right\|_\tF^2 \leq \left\|\U_k^{(t)} - \U_k^*\R_k^{(t)}\right\|_\tF^2 - 2\eta Q_{k,1} + \eta^2 Q_{k,2},\quad k=1,2,3.
	\end{split}
\end{equation}
\begin{equation}\label{eq-decompose-core-error-t+1-result}
	\begin{split}
		& \left\| \cS^{(t+1)} - \cS^* \times_1 \R_1^{(t)\top} \times_2 \R_2^{(t)\top} \times_3 \R_3^{(t)\top}\right\|_\tF^2 \\
		\leq & \left\|\cS^{(t)} - \cS^* \times_1 \R_1^{(t)\top} \times_2 \R_2^{(t)\top} \times_3 \R_3^{(t)\top}\right\|_\tF^2 - 2\eta Q_{\cS,1} + \eta^2 Q_{\cS,2}.
	\end{split}
\end{equation}
Then \eqref{eq-Et-descent} can be obtained by combining \eqref{eq-link-step}, \eqref{eq-decompose-U-error-t+1} and \eqref{eq-decompose-core-error-t+1-result}.
\begin{itemize}[leftmargin=*]
    \item To show \eqref{eq-decompose-U-error-t+1-result}, we first focus on $\left \| \tilde\U_1^{(t+1)} - \U_1^* \R_1^{(t)}\right\|_\tF^2$. By plugging in the gradient (Lemma \ref{lm:partial-gradient}), we have the following decomposition.
\begin{equation}\label{eq-decompose-U-error-t+1}
	\begin{split}
		& \left\|\tilde{\U}_1^{(t+1)} - \U_1^*\R_1^{(t)}\right\|_\tF^2\\ 
		= & \left\|{\U}_1^{(t)} - \U_1^*\R_1^{(t)} - \eta\left(\cM_1(\nabla L(\cX^{(t)})) \V_1^{(t)} + a \U_1^{(t)}( \U_1^{(t)\top} \U_1^{(t)}-b^2 \I)\right)\right\|_\tF^2 \\
		= & \left\|{\U}_1^{(t)} - \U_1^*\R_1^{(t)}\right\|_\tF^2 + \eta^2 \left\|\cM_1(\nabla L(\cX^{(t)})) \V_1^{(t)} + a \U_1^{(t)}( \U_1^{(t)\top} \U_1^{(t)}-b^2 \I)\right\|_\tF^2 \\
		& - 2\eta \left\langle  \U_1^{(t)} - \U_1^* \R_1^{(t)}, \cM_1(\nabla L(\cX^{(t)})) \V_1^{(t)}\right\rangle \\
		& - 2\eta a \left\langle \U_1^{(t)} - \U_1^* \R_1^{(t)}, \U_1^{(t)}( \U_1^{(t)\top} \U_1^{(t)}-b^2\I)\right\rangle
	\end{split}
\end{equation}
We bound the last three terms separately. First, we have
		\begin{equation*}
	\begin{split}
		&\left\|\cM_1(\nabla L(\cX^{(t)})) \V_1^{(t)} + a \U_1^{(t)}( \U_1^{(t)\top} \U_1^{(t)}-b^2 \I)\right\|_\tF^2  \\
		\leq & 2 \left(\left\|\cM_1(\nabla L(\cX^{(t)})) \V_1^{(t)}\right\|_\tF^2 +  a^2 \left\|\U_1^{(t)}( \U_1^{(t)\top} \U_1^{(t)}-b^2 \I)\right\|_\tF^2\right).
	\end{split}
\end{equation*}
The first term can be bounded as
\small
\begin{equation*}
	\begin{split}
		& \left\|\cM_1(\nabla L(\cX^{(t)})) \V_1^{(t)}\right\|_\tF^2 \\
		\leq & 2\left(\left\|\cM_1(\nabla L(\cX^*))\V_1^{(t)}\right\|_\tF^2 + \left\|\cM_1\left(\nabla L(\cX^{(t)}) - \nabla L(\cX^*)\right)\V_1^{(t)}\right\|_\tF^2\right) \\
		= & 2\left\|\cM_1\left(\nabla L(\cX^*)\right)\left(\U_3^{(t)} \otimes  \U_2^{(t)}\right)\cM_1( \cS^{(t)})^\top \right\|_\tF^2 + 2\left\|\cM_1\left(\nabla L(\cX^{(t)}) - \nabla L(\cX^*)\right) \V_1^{(t)}\right\|_\tF^2. \\
	\end{split}
\end{equation*}
\normalsize
Notice that by the duality of Frobenius norm, we have
\begin{equation*}
    \begin{split}
        & \left\|\cM_1\left(\nabla L(\cX^*)\right)\left(\U_3^{(t)} \otimes  \U_2^{(t)}\right)\cM_1( \cS^{(t)})^\top\right\|_\tF \\
        & = \sup_{\W_1 \in \bbR^{p_1\times r_1}, \|\W_1\|_\tF \leq 1} \left\langle \cM_1\left(\nabla L(\cX^*)\right)\left(\U_3^{(t)} \otimes  \U_2^{(t)}\right)\cM_1( \cS^{(t)})^\top, \W_1 \right\rangle \\
        & = \sup_{\W_1 \in \bbR^{p_1\times r_1}, \|\W_1\|_\tF \leq 1} \left\langle \cM_1\left(\nabla L(\cX^*)\right), \W_1 \cM_1(\cS^{(t)}) \left(\U_3^{(t)}\otimes \U_2^{(t)}\right)^\top \right\rangle \\
        & = \sup_{\W_1 \in \bbR^{p_1\times r_1}, \|\W_1\|_\tF \leq 1} \left\langle \nabla L(\cX^*), \cS^{(t)} \times_1 \W_1 \times_2 \U_2^{(t)} \times_3 \U_3^{(t)} \right\rangle \\
        & \leq \left\|\cM_1(\cS^{(t)})\right\|\cdot \left\|\U_3^{(t)} \otimes \U_2^{(t)}\right\| \cdot \xi,
    \end{split}
\end{equation*}
then it follows that
\small
\begin{equation*}
    \begin{split}
		\left\|\cM_1(\nabla L(\cX^{(t)})) \V_1^{(t)}\right\|_\tF^2 \leq & 2\left\|\cM_1(\cS^{(t)})\right\|^2\cdot \left\|\U_3^{(t)} \otimes \U_2^{(t)}\right\|^2 \cdot \xi^2 \\
		& + 2\left\|\nabla L(\cX^{(t)}) - \nabla L(\cX^*)\right\|_\tF^2 \left\| \U_3^{(t)} \otimes\U_2^{(t)}\right\|^2 \left\|\cM_1( \cS^{(t)})\right\|^2 \\
		\overset{\eqref{eq-Assumption1}}{\leq}  2\cdot (1.01b)^4 &\cdot \frac{(1.01\overline \lambda)^2}{b^6} \xi^2 + 2\cdot (1.01b)^4 \cdot \frac{(1.01\overline \lambda)^2}{b^6} \left\|\nabla L(\cX^{(t)}) - \nabla L(\cX^*)\right\|_\tF^2 \\
		\leq  4\overline \lambda^{2}b^{-2}&\left(\xi^2 + \left\|\nabla L(\cX^{(t)}) - \nabla L(\cX^*)\right\|_\tF^2  \right).
	\end{split}
\end{equation*}
\normalsize
In addition,
\begin{equation*}
	\begin{split}
		& \left\| \U_1^{(t)}\left( \U_1^{(t)\top} \U_1^{(t)}-b^2 \I_{r_1}\right)\right\|_\tF^2 \leq \left\| \U_1^{(t)}\right\|^2 \cdot \left\| \U_1^{(t)\top} \U_1^{(t)}-b^2 \I_{r_1}\right\|_\tF^2\\
		\overset{\eqref{eq-Assumption1}}{\leq} & (1.01b)^2 \left\| \U_1^{(t)\top} \U_1^{(t)}-b^2\I_{r_1}\right\|_\tF^2 \leq \frac{5}{4}b^2\left\| \U_1^{(t)\top} \U_1^{(t)}-b^2\I_{r_1}\right\|_\tF^2.
	\end{split}
\end{equation*}
Combining the two inequalities above, we have
\begin{equation}\label{eq-decompose-U-quadratic}
	\begin{split}
		&\left\|\cM_1(\nabla L(\cX^{(t)})) \V_1^{(t)} + a \U_1^{(t)}(\U_1^{(t)\top} \U_1^{(t)}-b^2 \I_{r_1})\right\|_\tF^2 \\
		\leq & 8 \overline\lambda^{2}b^{-2}\left(\xi_1^2 + \left\|\nabla L(\cX^{(t)}) - \nabla L(\cX^*)\right\|_\tF^2 \right) + \frac{5}{2}a^2b^2 \left\| \U_1^{(t)\top} \U_1^{(t)}-b^2 \I_{r_1}\right\|_\tF^2 = Q_{1,2}.
	\end{split}
\end{equation}
 
For the third term on the right hand side of \eqref{eq-decompose-U-error-t+1}, we have
		    \begin{equation}\label{eq-decompose-U-linear-1}
	            \begin{split}
		            & \left\langle  \U_1^{(t)} - \U_1^* \R_1^{(t)}, \cM_1\left(\nabla L(\cX^{(t)})\right) \V_1^{(t)}\right\rangle \\
		            = & \left\langle  \U_1^{(t)}  \V_1^{(t)\top}  - \U_1^* \R_1^{(t)}  \V_1^{(t)\top}, \cM_1\left(\nabla L(\cX^{(t)})\right)\right\rangle \\
		            \overset{\eqref{eq:def-V}}{=} & \left\langle \cX^{(t)} - \cS^{(t)} \times_1 \U_1^*  \R_1^{(t)} \times_2  \U_2^{(t)} \times_3  \U_3^{(t)}, \nabla L(\cX^{(t)}) \right\rangle \\     \overset{\eqref{eq:def-signal-tensor}}{=} &   \left\langle  \cX^{(t)} - \cX_1^{(t)}, \nabla L(\cX^{(t)}) \right\rangle.
	            \end{split}
            \end{equation}
For the last term on the right hand side of \eqref{eq-decompose-U-error-t+1}, we have
		\begin{equation}
	\begin{split}
		& \left\langle\U_1^{(t)}-\U_1^* \R_1^{(t)}, \U_1^{(t)}( \U_1^{(t)\top} \U_1^{(t)}-b^2\I_{r_1})\right\rangle \\
		= & \left\langle \U_1^{(t)\top} \U_1^{(t)} -  \U_1^{(t)\top} \U_1^* \R_1^{(t)},  \U_1^{(t)\top} \U_1^{(t)}-b^2\I_{r_1}\right\rangle \\
		= & \frac{1}{2}\left\langle  \U_1^{(t)\top} \U_1^{(t)}-\U_1^{*\top}\U_1^*,  \U_1^{(t)\top} \U_1^{(t)}-b^2\I_{r_1}\right\rangle \\
		&+ \frac{1}{2}\left\langle \U_1^{*\top}\U_1^* - 2 \U_1^{(t)\top} \U_1^* \R_1^{(t)} +   \U_1^{(t)\top} \U_1^{(t)},  \U_1^{(t)\top} \U_1^{(t)}-b^2\I_{r_1}\right\rangle \\
		\overset{(a)}{=} & \frac{1}{2}\left\| \U_1^{(t)\top} \U_1^{(t)} - b^2 \I_{r_1}\right\|_\tF^2 + \frac{1}{2}\left\langle  \U_1^{(t)\top}\left( \U_1^{(t)} - \U_1^* \R_1^{(t)}\right),  \U_1^{(t)\top} \U_1^{(t)}-b^2 \I_{r_1}\right\rangle \\
		&+ \frac{1}{2}\left\langle \U_1^{*\top}\U_1^* -  \U_1^{(t)\top} \U_1^*\R_1^{(t)},  \U_1^{(t)\top} \U_1^{(t)}-b^2\I_{r_1}\right\rangle.
	\end{split}
\end{equation}
Here, (a) is due to the assumption that $\U_1^{*\top}\U_1^* = b^2\I_{r_1}$. Since
\begin{equation*}
	\begin{split}
		& \left\langle \U_1^{*\top}\U_1^* -  \U_1^{(t)\top} \U_1^* \R_1^{(t)},  \U_1^{(t)\top} \U_1^{(t)}-b^2\I_{r_1}\right\rangle  \\
		\overset{(b)}{=} & \left\langle \U_1^{*\top}\U_1^* - \R_1^{(t)\top}\U_1^{*\top} \U_1^{(t)}, \U_1^{(t)\top} \U_1^{(t)}-b^2\I_{r_1}\right\rangle \\
		\overset{(c)}{=} & \left\langle \R_1^{(t)\top}\U_1^{*\top}\U_1^* \R_1^{(t)} - \R_1^{(t)\top}\U_1^{*\top} \U_1^{(t)},  \U_1^{(t)\top} \U_1^{(t)}-b^2\I_{r_1}\right\rangle \\
		= & \left\langle (\U_1^*  \R_1^{(t)})^\top\left(\U_1^* \R_1^{(t)} -  \U_1^{(t)}\right),   \U_1^{(t)\top} \U_1^{(t)}-b^2\I_{r_1}\right\rangle,
	\end{split}
\end{equation*}
where (b) is due to the fact that $\langle \A, \B\rangle = \langle \A^\top, \B\rangle$ for symmetric matrix $\B$ and (c) holds because $\U_1^{*\top}\U_1^* = b^2\I_{r_1}$ and $\R_1^\top \R_1 = \I_{r_1}$, we further have
\begin{equation}\label{eq-decompose-U-linear-2}
	\begin{split}
		& \left\langle \U_1^{(t)} - \U_1^* \R_1^{(t)}, \U_1^{(t)}( \U_1^{(t)\top} \U_1^{(t)}-b^2\I_{r_1})\right\rangle \\
		= & \frac{1}{2}\left\| \U_1^{(t)\top} \U_1^{(t)} - b^2 \I_{r_1}\right\|_\tF^2 \\
		& + \frac{1}{2}\left\langle \left(\U_1^*\R_1^{(t)} - \U_1^{(t)}\right)^\top \left(\U_1^* \R_1^{(t)} - \U_1^{(t)}\right), \U_1^{(t)\top} \U_1^{(t)}-b^2\I_{r_1}  \right\rangle \\
		\geq & \frac{1}{2}\left\| \U_1^{(t)\top} \U_1^{(t)} - b^2 \I_{r_1}\right\|_\tF^2 - \frac{1}{2}\left\|\U_1^* \R_1^{(t)} - \U_1^{(t)}\right\|_\tF^2\cdot \left\| \U_1^{(t)\top} \U_1^{(t)} - b^2 \I_{r_1}\right\|_\tF \\
		\geq & \frac{1}{2} \left\| \U_1^{(t)\top} \U_1^{(t)} - b^2 \I_{r_1}\right\|_\tF^2 - \left(\frac{1}{4} \left\|\U_1^{(t)\top}\U_1^{(t)} - b^2\I_{r_1}\right\|_\tF^2 + \frac{1}{4}\left\| \U_1^{(t)} - \U_1^*\R_1^{(t)}\right\|_\tF^4\right) \\
		\geq & \frac{1}{4} \left\| \U_1^{(t)\top} \U_1^{(t)} - b^2 \I_{r_1}\right\|_\tF^2 - \frac{1}{4} E^{(t)} \left\|\U_1^{(t)} - \U_1^*\R_1^{(t)}\right\|_\tF^2 \\
		\geq & \frac{1}{4} \left\| \U_1^{(t)\top} \U_1^{(t)} - b^2\I_{r_1}\right\|_\tF^2 - \frac{c_0}{16}\sqrt{\underline{\lambda}}\left\| \U_1^{(t)} - \U_1^* \R_1^{(t)}\right\|_\tF^2.
	\end{split}
\end{equation}
Here, the last inequality of \eqref{eq-decompose-U-linear-2} comes from the definition of $E^{(t)}$, \eqref{eq-Assumption2}, $\kappa\geq1$, and $\alpha\beta \leq 1/4$: $\left\| \U_1^{(t)} - \U_1^* \R_1^{(t)}\right\|_\tF^2 \leq E^{(t)} \leq c_0\frac{\alpha\beta\sqrt{\underline{\lambda}}}{\kappa^{3/2}} \leq \frac{c_0}{4}\sqrt{\underline{\lambda}}$.

Combining \eqref{eq-decompose-U-error-t+1}, \eqref{eq-decompose-U-quadratic}, \eqref{eq-decompose-U-linear-1}, and \eqref{eq-decompose-U-linear-2}, we obtain
\begin{equation*}
	\begin{split}
		& \left\|\tilde\U_1^{(t+1)} - \U_1^*\R_1^{(t)}\right\|_\tF^2 \leq \left\|\U_1^{(t)} - \U_1^*\R_1^{(t)}\right\|_\tF^2 - 2\eta Q_{1,1} + \eta^2 Q_{1,2}.
	\end{split}
\end{equation*}
Here, $Q_{1,1}, Q_{1, 2}$ are defined in \eqref{eq:D_k1,D_k2}. Then more generally,
\begin{equation*}
	\begin{split}
		& \left\|\tilde\U_k^{(t+1)} - \U_k^*\R_k^{(t)}\right\|_\tF^2 \leq \left\|\U_k^{(t)} - \U_k^*\R_k^{(t)}\right\|_\tF^2 - 2\eta Q_{k,1} + \eta^2 Q_{k,2},\quad k=1,2,3.
	\end{split}
\end{equation*}
This finishes the proof of \eqref{eq-decompose-U-error-t+1-result}.
\item Next we prove \eqref{eq-decompose-core-error-t+1-result}. Specifically, we have the following decomposition,
\begin{equation*}
	\begin{split}
		&\left\|\tilde\cS^{(t+1)} - \llbracket\cS^*; \R_1^{(t)\top}, \R_2^{(t)\top}, \R_3^{(t)\top}\rrbracket\right\|_\tF^2 \\
		= &\left\|\cS^{(t)} - \llbracket\cS^*; \R_1^{(t)\top}, \R_2^{(t)\top}, \R_3^{(t)\top}\rrbracket - \eta \llbracket\nabla L(\cX^{(t)}); \U_1^{(t)\top}, \U_2^{(t)\top}, \U_3^{(t)\top}\rrbracket\right\|_\tF^2 \\
		= & \left\|{\cS}^{(t)} - \llbracket\cS^*; \R_1^{(t)\top},  \R_2^{(t)\top}, \R_3^{(t)\top}\rrbracket\right\|_\tF^2 + \eta^2 \left\| \llbracket\nabla L(\cX^{(t)}); \U_1^{(t)\top}, \U_2^{(t)\top}, \U_3^{(t)\top}\rrbracket\right\|_\tF^2\\
		& - 2\eta \left\langle \cS^{(t)} - \llbracket\cS^*; \R_1^{(t)\top}, \R_2^{(t)\top}, \R_3^{(t)\top}\rrbracket, \llbracket\nabla L(\cX^{(t)}); \U_1^{(t)\top},  \U_2^{(t)\top}, \U_3^{(t)\top}\rrbracket \right\rangle.
	\end{split}
\end{equation*}
On the one hand, we have
\begin{equation*}
	\begin{split}
		& \left\langle \cS^{(t)} - \llbracket\cS^*; \R_1^{(t)\top}, \R_2^{(t)\top}, \R_3^{(t)\top}\rrbracket, \llbracket\nabla L(\cX^{(t)});  \U_1^{(t)\top}, \U_2^{(t)\top}, \U_3^{(t)\top}\rrbracket \right\rangle \\
		= & \left\langle \cM_1(\cS^{(t)}) - \R_1^{(t)\top}\cM_1(\cS^*) \left(\R_3^{(t)} \otimes \R_2^{(t)}\right), \U_1^{(t)\top}\cM_1\left(\nabla L(\cX^{(t)})\right)\left(\U_3^{(t)} \otimes \U_2^{(t)}\right) \right\rangle \\
		= & \Big\langle \U_1^{(t)}\cM_1(\cS^{(t)})\left( \U_3^{(t)} \otimes \U_2^{(t)}\right)^\top - \U_1^{(t)}\R_1^{(t)\top}\cM_1(\cS^*) \left(\R_3^{(t)} \U_3^{(t)\top} \otimes \R_2^{(t)} \U_2^{(t)\top}\right),\\ 
		& \qquad \cM_1\left(\nabla L(\cX^{(t)})\right)  \Big\rangle\\
		= & \left\langle \cX^{(t)} - \cS^* \times_1 \U_1^{(t)}\R_1^{(t)\top} \times_2 \U_2^{(t)}\R_2^{(t)\top} \times_3 \U_3^{(t)}\R_3^{(t)\top}, \nabla L(\cX^{(t)})\right\rangle \\
		= & \left\langle \cX^{(t)} - \cX_\cS^{(t)}, \nabla L(\cX^{(t)})\right\rangle = Q_{\cS,1}.
	\end{split}
\end{equation*}
On the other hand, by noticing that
\begin{equation*}
    \begin{split}
        & \left\|\llbracket\nabla L(\cX^{*}); \U_1^{(t)\top}, \U_2^{(t)\top}, \U_3^{(t)\top}\rrbracket\right\|_\tF = \sup_{\substack{\cS \in \bbR^{r_1\times r_2 \times r_3} \\ \left\|\cS\right\|_\tF \leq 1}} \left\langle \llbracket\nabla L(\cX^{*}); \U_1^{(t)\top}, \U_2^{(t)\top}, \U_3^{(t)\top}\rrbracket, \cS \right\rangle \\
        & = \sup_{\substack{\cS \in \bbR^{r_1\times r_2 \times r_3} \\ \left\|\cS\right\|_\tF \leq 1}} \left\langle \nabla L(\cX^*), \llbracket \cS; \U_1^{(t)}, \U_2^{(t)}, \U_3^{(t)}\rrbracket \right\rangle \\
        & \leq \left\|\U_1^{(t)}\right\| \cdot \left\|\U_2^{(t)}\right\| \cdot \left\|\U_3^{(t)}\right\| \cdot \xi,
    \end{split}
\end{equation*}
we also have
\small
\begin{equation*}
	\begin{split}
		&\left\|\llbracket\nabla L(\cX^{(t)}); \U_1^{(t)\top}, \U_2^{(t)\top}, \U_3^{(t)\top}\rrbracket\right\|_\tF^2 \\
		\leq & 2\left\|\llbracket\nabla L(\cX^{*}); \U_1^{(t)\top}, \U_2^{(t)\top}, \U_3^{(t)\top}\rrbracket\right\|_\tF^2 \\ 
		& \qquad + 2\left\|\llbracket\nabla L(\cX^{(t)})-\nabla L(\cX^{*}); \U_1^{(t)\top}, \U_2^{(t)\top}, \U_3^{(t)\top}\rrbracket\right\|_\tF^2  \\
		\leq & 2\left\|\U_1^{(t)}\right\|^2 \cdot \left\|\U_2^{(t)}\right\|^2 \cdot \left\|\U_3^{(t)}\right\|^2 \cdot \xi^2 + \\
		& \qquad +2\left\|\U_1^{(t)}\right\|^2 \cdot \left\|\U_2^{(t)}\right\|^2 \cdot \left\|\U_3^{(t)}\right\|^2 \cdot \left\|\nabla L(\cX^{(t)}) - \nabla L(\cX^{*})\right\|_\tF^2 \\
		\overset{\eqref{eq-Assumption1}}{\leq} & 4 b^6\left(\xi^2 + \left\|\nabla L(\cX^{(t)}) - \nabla L(\cX^{*})\right\|_\tF^2\right) = Q_{\cS, 2}.
	\end{split}
\end{equation*}
\normalsize
Therefore,
\begin{equation*}
	\begin{split}
		& \left\| \cS^{(t+1)} - \cS^* \times_1 \R_1^{(t)\top} \times_2 \R_2^{(t)\top} \times_3 \R_3^{(t)\top}\right\|_\tF^2 \\
		\leq & \left\|\cS^{(t)} - \cS^* \times_1 \R_1^{(t)\top} \times_2 \R_2^{(t)\top} \times_3 \R_3^{(t)\top}\right\|_\tF^2 - 2\eta Q_{\cS,1} + \eta^2 Q_{\cS,2},
	\end{split}
\end{equation*}
which proves \eqref{eq-decompose-core-error-t+1-result}.
\end{itemize}

\item[Step 3] In this Step 3, we aim to develop a lower bound for $Q_{\cS,1}+\sum_{k=1}^3Q_{k,1}$.
By definitions of $Q_{\cS,1}$ and $Q_{k,1}$, we have
\small
\begin{equation}\label{eq-linear-all}
	\begin{split}
		& Q_{\cS,1}+\sum_{k=1}^3 Q_{k,1}  \\
		= &\left\langle 4 \cX^{(t)} - \cX_S^{(t)} - \sum_{k=1}^3 \cX_k^{(t)}, \nabla L(\cX^{(t)})\right\rangle  \\
		& + a \sum_{k=1}^3\left(\frac{1}{4}\left\|\U_k^{(t)\top}\U_k^{(t)} - b^2 \I_{r_k}\right\|_\tF^2 - \frac{c_0}{16}\sqrt{\underline{\lambda}}\left\| \U_k^{(t)} - \U_k^*\R_k^{(t)}\right\|_\tF^2\right) \\
		= & \left\langle 4  \cX^{(t)} - \left(3 \cX^{(t)} + \cX^* - \cH_\varepsilon\right), \nabla L(\cX^{(t)})\right\rangle \\ 
		& + a\sum_{k=1}^3\left(\frac{1}{4}\left\| \U_k^{(t)\top} \U_k^{(t)} - b^2\I_{r_k}\right\|_\tF^2 - \frac{c_0}{16}\sqrt{\underline{\lambda}}\left\| \U_k^{(t)} - \U_k^* \R_k^{(t)}\right\|_\tF^2\right) \\
		= & \left\langle  \cX^{(t)} - \cX^* + \cH_\varepsilon, \nabla L(\cX^{(t)})\right\rangle \\
		& + a\sum_{k=1}^3\left(\frac{1}{4}\left\|\U_k^{(t)\top} \U_k^{(t)} - b^2\I_{r_k}\right\|_\tF^2 - \frac{c_0}{16}\sqrt{\underline{\lambda}}\left\| \U_k^{(t)} - \U_k^* \R_k^{(t)}\right\|_\tF^2\right).
	\end{split}
\end{equation}
\normalsize
Here, $\cH_\varepsilon$ is a tensor of small amplitude that is obtained by applying Lemma \ref{lm-X-decomposition}. By \eqref{eq-Assumption1}, the quantities $B_1, B_2, B_3$ in the context of Lemma \ref{lm-X-decomposition} satisfy $B_1 \leq 1.01b$, $B_2 \leq \frac{1.01\overline\lambda}{b^3}$ and $B_3 \leq E^{(t)}$. Thus, Lemma \ref{lm-X-decomposition} implies
\begin{equation}\label{ineq:H_epsilon}
	\begin{split}
		\left\|\cH_\varepsilon\right\|_\tF \leq & B_3^{3/2}B_2 + 3B_1B_2B_3 + 3B_3B_1^2\\
		\leq & 1.01\overline\lambda b^{-3}(E^{(t)})^{3/2} + 3(1.01b)^2E^{(t)} + 3(1.01b)^2E^{(t)} \\
		\overset{\eqref{eq-Assumption2}}{\leq} & \left(1.01 c_0^{1/2}\overline\lambda \underline{\lambda}^{1/4} b^{-3} + 6 \cdot 1.01^2 b^2\right)E^{(t)} \\
		= &  \left(1.01 c_0^{1/2}\kappa^{1/4} + 6 \cdot 1.01^2 \kappa^{1/2}\right)\underline{\lambda}^{1/2}E^{(t)} \\
		\leq & 6.5\kappa^{1/2}\underline{\lambda}^{1/2} E^{(t)}.
	\end{split}
\end{equation}

The first term on the right hand side of \eqref{eq-linear-all} can be further bounded as
\small
\begin{equation*}
	\begin{split}
		& \left\langle \cX^{(t)} - \cX^* + \cH_\varepsilon, \nabla L(\cX^{(t)}) \right\rangle = \left\langle \cX^{(t)} - \cX^*, \nabla L(\cX^{(t)}) - \nabla L(\cX^*) \right\rangle \\
		& + \left\langle \cX^{(t)} - \cX^* + \cH_\varepsilon, \nabla L(\cX^{*})\right\rangle + \left\langle \cH_\varepsilon, \nabla L(\cX^{(t)}) - \nabla L(\cX^*)\right\rangle \\
		\overset{(e)}{\geq} & \alpha\left\| \cX^{(t)} - \cX^*\right\|_\tF^2 + \beta \left\| \nabla L(\cX^{(t)}) - \nabla L(\cX^*)\right\|_\tF^2  - \left\|\cH_\varepsilon \right\|_\tF\left\|\nabla L(\cX^{(t)}) - \nabla L(\cX^*)\right\|_\tF \\
		& - \left| \left\langle \cX^{(t)} - \cX^* + \cH_\varepsilon, \nabla L(\cX^*)\right\rangle\right|. \\
	\end{split}
\end{equation*}
\normalsize
Here ($e$) is due to RCG condition and Cauchy-Schwarz inequality. Now we have
\small
\begin{equation*}
	\begin{split}
		& \left\|\cH_\varepsilon \right\|_\tF\left\|\nabla L(\cX^{(t)}) - \nabla L(\cX^*)\right\|_\tF \leq \frac{\beta}{2}\left\| \nabla L(\cX^{(t)}) - \nabla L(\cX^*)\right\|_\tF^2 + \frac{1}{2\beta} \left\|\cH_\varepsilon\right\|_\tF^2 \\
		\overset{\eqref{ineq:H_epsilon}}{\leq} & \frac{\beta}{2}\left\| \nabla L(\cX^{(t)}) - \nabla L(\cX^*)\right\|_\tF^2 + \frac{1}{2\beta}\left(6.5\kappa^{1/2}\underline{\lambda}^{1/2} E^{(t)}\right)^2 \\
		\overset{\eqref{eq-Assumption2}}{\leq} & \frac{\beta}{2}\left\| \nabla L(\cX^{(t)}) - \nabla L(\cX^*)\right\|_\tF^2 + 24c_0\frac{ \alpha\underline{\lambda}^{3/2}}{\kappa^{1/2}}E^{(t)},
	\end{split}
\end{equation*}
\normalsize
and
\small
\begin{equation}\label{eq-noise-inner-product}
	\begin{split}
		& \left| \left\langle \cX^{(t)} - \cX^* + \cH_\varepsilon, \nabla L(\cX^*)\right\rangle\right|  \\
		& =  \left| \left\langle 4 \cX^{(t)} - \cX_\cS^{(t)} - \sum_{k=1}^3 \cX_k^{(t)}, \nabla L(\cX^*)\right\rangle\right| \\
		\leq & \left| \left\langle \cX^{(t)} - \cX_\cS^{(t)}, \nabla L(\cX^*)\right\rangle\right| + \sum_{k=1}^3 \left| \left\langle  \cX^{(t)} - \cX_k^{(t)}, \nabla L(\cX^*)\right\rangle\right| \\
		= & \left| \left\langle \left\llbracket\cS^{(t)} - \llbracket \cS^*;\R_1^{(t)\top},\R_2^{(t)\top},\R_3^{(t)\top} \rrbracket; \U_1^{(t)}, \U_2^{(t)}, \U_3^{(t)}\right\rrbracket, \nabla L(\cX^*)\right\rangle\right| \\
		& + \sum_{k=1}^3 \left| \left\langle \cS^{(t)} \times_k (\U_k^{(t)} - \U_k^*\R_k^{(t)}) \times_{k+1}  \U_{k+1}^{(t)} \times_{k+2} \U_{k+2}^{(t)} , \nabla L(\cX^*)\right\rangle\right| \\
		\leq & \xi \left\|\cS^{(t)} - \llbracket \cS^*;\R_1^{(t)\top},\R_2^{(t)\top},\R_3^{(t)\top} \rrbracket\right\|_\tF \cdot \left\| \U_1^{(t)}\right\| \cdot \left\| \U_2^{(t)}\right\| \cdot \left\| \U_3^{(t)}\right\| \\
		& + \xi \cdot \sum_{k=1}^3 \left\|\cM_k(\cS^{(t)}) \right\| \cdot \left\|\U_k^{(t)} - \U_k^*\R_k^{(t)}\right\|_\tF \cdot \left\| \U_{k+2}^{(t)}\otimes \U_{k+1}^{(t)}\right\|  \\
		\leq & \left(\|\U_1^{(t)}\|\|\U_2^{(t)}\|\|\U_3^{(t)}\|+\sum_{k=1}^3\left\|\cM_k(\cS^{(t)})\right\|\|\U_{k+2}^{(t)}\|\|\U_{k+1}^{(t)}\|\right)\xi (E^{(t)})^{1/2} \\
		& \overset{\eqref{eq-Assumption1}}{\leq} \left(1.01^3b^3  + 3\cdot 1.01^3\overline\lambda b^{-3} b^2 \right)\xi(E^{(t)})^{1/2} \\
		\leq & \frac{9}{2} \left(b^3 + \overline\lambda b^{-1}\right)\xi\sqrt{E^{(t)}} \leq \frac{c_1}{4} \frac{\alpha (b^3 + \overline\lambda b^{-1})^2}{\kappa^{2}}E^{(t)} + \frac{81\kappa^2}{4c_1\alpha}\xi^2 \\
		= & c_1 \frac{\alpha\underline{\lambda}^{3/2}}{\kappa^{1/2}}E^{(t)} + \frac{81\kappa^2}{4c_1\alpha}\xi^2.
	\end{split}
\end{equation}
\normalsize
Here, $c_1$ is some small constant that can be specified as $c_1 = 1/3000$. Then it follows from the previous three inequalities that
\small
\begin{equation}\label{eq-decrease-lower-bound}
	\begin{split}
		& \left\langle \cX^{(t)} - \cX^* + \cH_\varepsilon, \nabla L(\cX^{(t)})\right\rangle  \\
		\geq &\alpha \left\| \cX^{(t)} - \cX^*\right\|_\tF^2 + \frac{\beta}{2}\left\|\nabla L(\cX^{(t)}) - \nabla L(\cX^*)\right\|_\tF^2 \\
		&  - (24c_0+c_1)\frac{\alpha\underline{\lambda}^{3/2}}{\kappa^{1/2}}E^{(t)} - \frac{81\kappa^2}{4c_1\alpha}\xi^2. \\
	\end{split}
\end{equation}
\normalsize
Now applying Lemma \ref{lm-equivalent-criteria} with $c_d = 0.01$ and $b = \overline{\lambda}$, we have
\small
\begin{equation}\label{eq-apply-equivalance-lemma}
	\begin{split}
		E^{(t)} &\leq 480 \kappa^2 b^{-6}\left\|\cX^{(t)} - \cX^*\right\|_\tF^2 + 80 b^{-2}\sum_{k=1}^3\left\| \U_k^{(t)\top} \U_k^{(t)} - b^2\I_{r_k}\right\|_\tF^2.
	\end{split}
\end{equation}
\normalsize
By setting $a = \frac{4\alpha b^{4}}{3\kappa^2} = \frac{4\alpha\underline{\lambda}}{3\kappa}$ and applying \eqref{eq-linear-all}, \eqref{eq-decrease-lower-bound} and \eqref{eq-apply-equivalance-lemma}, we have
\small
\begin{equation}\label{ineq:convergece-aggregate}
	\begin{split}
		& Q_{\cS,1}+\sum_{k=1}^3 Q_{k,1} \geq \alpha \left(\left\|\cX^{(t)} - \cX^*\right\|_\tF^2 + \frac{\kappa^{-2}b^4}{6}\sum_{k=1}^3\left\| \U_k^{(t)\top} \U_k^{(t)} - b^2\I_{r_k}\right\|_\tF^2\right) \\
		&\qquad +\frac{\beta}{2}\left\|\nabla L(\cX^{(t)}) - \nabla L(\cX^*)\right\|_\tF^2 + \frac{a}{8}\sum_{k=1}^3\left\| \U_k^{(t)\top} \U_k^{(t)} - b^2\I_{r_k}\right\|_\tF^2 \\
		&\qquad - \frac{c_0}{12}\frac{\alpha{\underline{\lambda}}^{3/2}}{\kappa} \sum_{k=1}^3 \left\| \U_k^{(t)} - \U_k^* \R_k^{(t)}\right\|_\tF^2 - (24c_0 + c_1)\frac{\alpha\underline{\lambda}^{3/2}}{\kappa^{1/2}}E^{(t)} - \frac{81\kappa^2}{4c_1\alpha}\xi^2\\
		\overset{\eqref{eq-apply-equivalance-lemma}}{\geq} &  \frac{\alpha\kappa^{3/2}\underline{\lambda}^{3/2}}{480\kappa^2}  E^{(t)} +\frac{\beta}{2}\left\|\nabla L(\cX^{(t)}) - \nabla L(\cX^*)\right\|_\tF^2 + \frac{a}{8}\sum_{k=1}^3\left\| \U_k^{(t)\top} \U_k^{(t)} - b^2\I_{r_k}\right\|_\tF^2 \\
		& -\left(24c_0 + c_1\right)\frac{\alpha{\underline{\lambda}}^{3/2}}{\kappa^{1/2}}E^{(t)} - \frac{81\kappa^2}{4c_1\alpha}\xi^2\\
		\geq &  \frac{\alpha\underline{\lambda}^{3/2}}{\kappa^{1/2}} \left( \frac{1}{480} - 24c_0 - c_1 \right)  E^{(t)} + \frac{\beta}{2}\left\|\nabla L(\cX^{(t)}) - \nabla L(\cX^*)\right\|_\tF^2  - \frac{81\kappa^2}{4c_1\alpha}\xi^2\\
		& + \frac{a}{8}\sum_{k=1}^3\left\|\U_k^{(t)\top} \U_k^{(t)} - b^2\I_{r_k}\right\|_\tF^2.
	\end{split}
\end{equation}
\normalsize
Finally we can plug in $c_0 = \frac{1}{20000}$ and $c_1 = \frac{1}{3000}$ and establish a lower bound for $Q_{\cS,1}+\sum_{k=1}^3Q_{k,1}$:
\small
\begin{equation}\label{eq-lower-bound-linear-term}
	\begin{split}
		Q_{\cS,1}+\sum_{k=1}^3 Q_{k,1} & \geq  \rho\frac{\alpha\underline{\lambda}^{3/2}}{\kappa^{1/2}} E^{(t)} + \frac{\beta}{2}\left\|\nabla L(\cX^{(t)}) - \nabla L(\cX^*)\right\|_\tF^2 \\
		& \qquad + \frac{a}{8}\sum_{k=1}^3\left\|\U_k^{(t)\top} \U_k^{(t)} - b^2\I_{r_k}\right\|_\tF^2 - C_1\frac{\kappa^2}{\alpha}\xi^2,
	\end{split}
\end{equation}
\normalsize
where $\rho = 1/2000$ and $C_1$ is some universal constant.
\item[Step 4] (Convergence Analysis) In this step we combine all results in previous steps to establish the error bound for $E^{(t)}$ and $\|\cX^{(t)}- \cX^*\|_\tF$. By plugging in $a=\frac{4\alpha b^{4}}{3\kappa^2}$ to the definitions of $Q_{\cS,2}$ and $Q_{k,2}$ in \eqref{eq:D_k1,D_k2}, 
we have
\begin{equation}\label{eq-upper-bound-quadratic-term}
	\begin{split}
		 & Q_{\cS,2} + \sum_{k=1}^3 Q_{k,2} \\
		 \leq & 28 b^6\left\| \nabla L(\cX^{(t)}) - \nabla L(\cX^*)\right\|_\tF^2 + \frac{5\alpha^2b^{10}}{\kappa^4} \sum_{k=1}^3\left\| \U_k^{(t)\top} \U_k^{(t)} - b^2\I_{r_k}\right\|_\tF^2 + 28 b^6\xi^2.
	\end{split}
\end{equation}
Combining \eqref{eq-link-step}, \eqref{eq-lower-bound-linear-term}, and \eqref{eq-upper-bound-quadratic-term} and setting $\eta = \eta_0\beta/b^6$, for any positive constant $\eta_0 \leq \frac{1}{28}$, we have:
\begin{equation*}
	\begin{split}
		E^{(t+1)} \leq & E^{(t)} - \frac{2\rho\alpha\beta\eta_0}{\kappa^2} E^{(t)} - \left(\frac{\beta^2}{b^6}\eta_0 - \frac{28\beta^2}{b^6}\eta_0^2\right)\left\| \nabla L(\cX^{(t)}) - \nabla L(\cX^*)\right\|_\tF^2\\
		& - \left(\frac{\alpha\beta b^{-2}}{3\kappa^2}\eta_0 - \frac{5\alpha^2\beta^2b^{-2}}{\kappa^4}\eta_0^2\right)\sum_{k=1}^3\left\| \U_k^{(t)\top} \U_k^{(t)} - b^2\I_{r_k}\right\|_\tF^2 \\
		& + \frac{2C_1\beta b^{-6}\kappa^2\xi^2}{\alpha}\eta_0  + 282\beta^2 b^{-6}\xi^2 \eta_0^2,
	\end{split}
\end{equation*}
\begin{equation*}
	\begin{split}
		&\frac{\beta^2}{b^6}\eta_0 - \frac{28\beta^2}{b^6}\eta_0^2 = \frac{\beta^2\eta_0}{b^6}(1 - 28\eta_0) > 0,\\
		&\frac{\alpha\beta b^{-2}}{3\kappa^2}\eta_0 - \frac{5\alpha^2\beta^2b^{-2}}{\kappa^4}\eta_0^2 = \frac{\alpha \beta\eta_0}{\kappa^2b^2}\left(\frac{1}{3}-\frac{5\alpha\beta}{\kappa^2}\eta_0\right)\overset{(f)}{>}0.
	\end{split}
\end{equation*}
($f$) is also due to $\beta \leq \frac{1}{4\alpha}$ and $\kappa \geq 1$. Then it follows that
\begin{equation}\label{eq-convergence-E-next-step}
	\begin{split}
		E^{(t+1)} & \leq \left(1 - \frac{2\rho\alpha\beta\eta_0}{\kappa^2}\right)E^{(t)} + \frac{2C_1\beta b^{-6}\kappa^2\xi^2}{\alpha}\eta_0  + 28\beta^2 b^{-6}\xi^2\eta_0^2 \\
		& \leq \left(1 - \frac{2\rho\alpha\beta\eta_0}{\kappa^2}\right)E^{(t)} + \frac{2C_1\beta b^{-6}\kappa^2\xi^2}{\alpha}\eta_0  + \frac{3\beta b^{-6}\xi^2}{\alpha}\eta_0^2 \\
		& \leq \left(1 - \frac{2\rho\alpha\beta\eta_0}{\kappa^2}\right)E^{(t)} + \frac{(2C_1+1)\beta b^{-6}\kappa^2\xi^2}{\alpha}\eta_0
	\end{split}
\end{equation}
for $t=0,1,\ldots, t_{max}$. Now we use induction to show
\begin{equation}\label{eq-upper-bound-Et}
    \begin{split}
        E^{(t)} \leq \frac{(2C_1+1)\kappa^4}{2\rho\alpha^2b^6}\xi^2 + \left(1 - \frac{2\rho\alpha\beta \eta_0}{\kappa^2}\right)^t E^{(0)}.
    \end{split}
\end{equation}
When $t=0$, \eqref{eq-upper-bound-Et} clearly holds. We assume it holds at $t=t_0$, then for $t = t_0 + 1$, we have
\begin{equation*}
    \begin{split}
        E^{(t_0+1)}  \overset{\eqref{eq-convergence-E-next-step}}{\leq} & \left(1 - \frac{2\rho\alpha\beta \eta_0}{\kappa^2}\right) E^{(t_0)} + \frac{(2C_1+1)\beta b^{-6}\kappa^2\xi^2}{\alpha}\eta_0\\
        \overset{\eqref{eq-upper-bound-Et}}{\leq} & \left(1 - \frac{2\rho\alpha\beta \eta_0}{\kappa^2}\right)^{t_0+1} E^{(0)} + \left(1 - \frac{2\rho\alpha\beta \eta_0}{\kappa^2}\right)\frac{(2C_1+1)\kappa^4}{2\rho\alpha^2b^6}\xi^2 \\
        & \qquad +  \frac{(2C_1+1)\beta b^{-6}\kappa^2\xi^2}{\alpha}\eta_0 \\
         = & \left(1 - \frac{2\rho\alpha\beta \eta_0}{\kappa^2}\right)^{t_0+1} E^{(0)} + \frac{(2C_1+1)\kappa^4}{2\rho\alpha^2b^6}\xi^2,
    \end{split}
\end{equation*}
thus \eqref{eq-upper-bound-Et} also holds at step $t_0+1$. By induction, \eqref{eq-upper-bound-Et} holds for any $t=0,1,\ldots, t_{max}$.
Then we can apply Lemma \ref{lm-equivalent-criteria} again and obtain
\begin{equation*}
    \begin{split}
        \left\|\cX^{(t)} - \cX^*\right\|_\tF^2 & \leq 42b^6 E^{(t)} \leq \frac{42(2C_1+1)\kappa^4}{\rho\alpha^2}\xi^2 + 42b^6\left(1 - \frac{\rho\alpha\beta\eta_0}{\kappa^2}\right)^t E^{(0)} \\
        & \leq C\left(\frac{\kappa^4}{\alpha^2}\xi^2 + \kappa^2\left(1 - \frac{2\rho\alpha\beta\eta_0}{\kappa^2}\right)^t \left\|\cX^{(0)} - \cX^*\right\|_\tF^2\right),
    \end{split}
\end{equation*}
where $C$ is some universal constant.
	\item[Step 5] In this step, we show that conditions \eqref{eq-Assumption1} and \eqref{eq-Assumption2} hold. We first apply induction on \eqref{eq-convergence-E-next-step} to prove \eqref{eq-Assumption2}. Since $\|\U_k^{*\top}\U_k^* - b^2 \I_{r_k}\|_\tF^2 = 0$, by Lemma \ref{lm-equivalent-criteria} and the initialization error bound, we have
	\begin{equation*}
	    E^{(0)} \leq  480b^{-6}\kappa^2 \left\|\cX^{(0)} - \cX^*\right\|_\tF^2 \leq c_0\frac{\alpha\beta \underline\lambda^{1/2}}{\kappa^{3/2}}.
	\end{equation*}
	Now assume \eqref{eq-Assumption2} holds at step $t$, then for step $t+1$, we have
\small
\begin{equation*}
	\begin{split}
		E^{(t+1)} & \overset{\eqref{eq-convergence-E-next-step}}{\leq} \left(1 - \frac{\rho\alpha\beta\eta_0}{\kappa^2}\right)E^{(t)} + \frac{(2C_1+1)\beta b^{-6}\kappa^2\xi^2}{\alpha}\eta_0 \\
		& \leq \left(1 - \frac{\rho\alpha\beta\eta_0}{\kappa^2}\right)\frac{c_0\alpha\beta \underline\lambda^{1/2}}{\kappa^{3/2}} + \frac{(2C_1+1)\beta b^{-6}\kappa^2\xi^2}{\alpha}\eta_0 \\
		& \leq c_0\frac{\alpha\beta \underline\lambda^{1/2}}{\kappa^{3/2}} - \eta_0\left(\rho c_0\frac{\alpha^2\beta^2 \underline\lambda^{1/2}}{\kappa^{7/2}} - \frac{(2C_1+1)\beta b^{-6}\kappa^2\xi^2}{\alpha}\right).
	\end{split}
\end{equation*}
\normalsize
By the signal-noise-ratio assumption that $\underline{\lambda}^2 \geq C\frac{\kappa^4}{\alpha^3\beta}\xi^2$ for some universal big constant $C$, we know 
\begin{equation*}
    \rho c_0\frac{\alpha^2\beta^2 \underline\lambda^{1/2}}{\kappa^{7/2}} - \frac{(2C_1+1)\beta b^{-6}\kappa^2\xi^2}{\alpha} \geq 0
\end{equation*}
and then $E^{(t+1)} \leq c_0\frac{\alpha\beta \underline\lambda^{1/2}}{\kappa^{3/2}}$. By induction, \eqref{eq-Assumption2} holds for any $t=0,1,\ldots, t_{max}$.

Then we use \eqref{eq-Assumption2} to prove \eqref{eq-Assumption1}. Note that in Step $t$, we have $E^{(t)} \leq c_0\frac{\alpha\beta\underline\lambda^{1/2}}{\kappa^{3/2}} = c_0\frac{\alpha\beta\underline b^2}{\kappa^{2}} \leq \frac{c_0}{4}b^2$, which implies that for any $k=1,2,3$,
\small
\begin{equation*}	
    \begin{split}
        \left\|\U_k^{(t)}\right\| &\leq \left\|\U_k^* \R_k^{(t)}\right\| + \left\| \U_k^{(t)} - \U_k^* \R_k^{(t)}\right\|\\
        &= b + \left\|\U_k^{(t)} - \U_k^* \R_k^{(t)}\right\|_\tF  \leq (1+\frac{\sqrt{c_0}}{2})b \leq 1.01b,
    \end{split}
\end{equation*}
	\begin{equation*}
		\begin{split}			\left\|\cM_k(\cS^{(t)})\right\| \leq &  \left\|\R_1^{(t)\top}\cM_k(\cS^*)\left(\R_3\otimes \R_2\right)\right\| \\
		&  + \left\|\R_k^{(t)\top}\cM_k(\cS^*)\left(\R_{k+1}^{(t)}\otimes \R_{k+2}^{(t)}\right) - \cM_k(\cS^{(t)})\right\| \\
			\leq & \overline\lambda/b^3 + \left\|\cS^{(t)} - \llbracket\cS^*; \R_1^{(t)}, \R_2^{(t)}, \R_3^{(t)}\rrbracket\right\|_\tF \\
			\leq & \overline\lambda/b^3 +\frac{\sqrt{c_0}}{2}b \leq 1.01\overline\lambda/b^3.
		\end{split}
	\end{equation*} 
\normalsize
Now we have finished the proof of Theorem \ref{thm:local-convergence-PGD}.\quad $\square$
\end{enumerate}

\section{Proofs of Other Theorems}\label{sec:proof-others}
We collect the proofs of all the other theorems in this section.
\subsection{Proof of Theorem \ref{thm:Gaussian}}
Without loss of generality, we assume $\sigma=1$. Let $\cX^* = \llbracket \cS, \U_1, \U_2, \U_3 \rrbracket$, where $\U_1,\U_2,\U_3$ are orthogonal matrices, and let $\U_k^* = b\U_k$ for $k=1,2,3$ and $\cS^* = \cS/b^3$. Also recall that $\tilde \U_k = \SVD_{r_k}(\cM_k(\cY))$, and $\U_k^{(0)} = b \tilde \U_k^{(0)}$. We would like to apply Theorem \ref{thm:local-convergence-PGD} to obtain the result. In the context of sub-Gaussian tensor PCA, we have
\begin{equation*}
    \nabla L(\cX^*) = \cX - \cY =: \cZ,
\end{equation*}
where all entries of $\Z$ are independent mean-zero sub-Gaussian random variables such that
\begin{equation*}
    \|\cZ_{ijk}\|_{\psi_2} = \sup_{q\geq 1} \mathbb{E}\left(|\cZ_{ijk}|^q\right)^{1/q}/q^{1/2} \leq 1.
\end{equation*}
We claim that with probability at least $1-C\exp(-c\underline p)$, the following conditions hold:
\begin{equation}\label{eq-Gaussian-hpc}
    \begin{split}
        &\left\|\sin \Theta (\tilde \U_k, \U_k)\right\| \leq \frac{\sqrt{p_k}\underline\lambda + \sqrt{p_1p_2p_3}}{\underline\lambda^2},\quad k=1,2,3, \\
        & \xi = \sup_{\substack{\mathcal{T}\in \mathbb{R}^{p_1\times p_2\times p_3}, \\ \rank(\mathcal{T})\leq (r_1,r_2,r_3), \\ \|\mathcal{T}\|_\tF \leq 1}} \left\langle \cZ, \mathcal{T} \right\rangle \leq C\left(\sqrt{r_1r_2r_3} + \sum_{k=1}^3 \sqrt{p_kr_k}\right).
    \end{split}
\end{equation}
Here, the first inequality of \eqref{eq-Gaussian-hpc} holds with probability at least $1-C\exp(-c\overline{p})$ ($\overline p = \max\{p_1,p_2,p_3\}$) by the proof of \citep[Theorem 4]{zhang2018heteroskedastic}; while the second inequality holds with probability at least $1-\exp(\sum p_kr_k)$ by Lemma \ref{lm-concentration-Gaussian-xi}. Applying union bound proves the above claim.

Now we start to bound $\left\|\cX^{(0)} - \cX^*\right\|_\tF$ under the condition \eqref{eq-Gaussian-hpc}. By the initialization procedure, we have 
\begin{equation*}
	\cX^{(0)} = \cY \times_1 \Proj_{\tilde{\U}_1} \times_2 \Proj_{\tilde{\U}_2} \times_3 \Proj_{\tilde{\U}_3}.
\end{equation*}
Then it follows that
\begin{equation*}
	\begin{split}
		& \left\|\cX^{(0)} - \cX^* \right\|_\tF \leq \left\|\cX^* - \cX^* \times_1 \Proj_{\tilde{\U}_1} \times_2 \Proj_{\tilde{\U}_2} \times_3 \Proj_{\tilde{\U}_3}\right\|_\tF + \left\|\cZ \times_1 \Proj_{\tilde{\U}_1} \times_2 \Proj_{\tilde{\U}_2} \times_3 \Proj_{\tilde{\U}_3}\right\|_\tF. \\
	\end{split}
\end{equation*}
For the first term, we have
\begin{equation*}
	\begin{split}
		& \left\|\cX^* - \cX^* \times_1 \Proj_{\tilde{\U}_1} \times_2 \Proj_{\tilde{\U}_2} \times_3 \Proj_{\tilde{\U}_3}\right\|_\tF \\
		& = \left\|\cX^* \times_1 \Proj_{\tilde{\U}_{1\perp}} + \cX^* \times_1 \Proj_{\tilde{\U}_1} \times_2 \Proj_{\tilde{\U}_{2\perp}} + \cX^* \times_1 \Proj_{\tilde{\U}_1} \times_2 \Proj_{\tilde{\U}_{2}} \times_3 \Proj_{\tilde{\U}_{3\perp}}\right\|_\tF \\
		& \leq \left\|\cX^* \times_1 \Proj_{\tilde{\U}_{1\perp}}\right\|_\tF + \left\|\cX^* \times_2 \Proj_{\tilde{\U}_{2\perp}}\right\|_\tF + \left\|\cX^* \times_3 \Proj_{\tilde{\U}_{3\perp}}\right\|_\tF \\
		& = \sum_{k=1}^3 \left\|\tilde{\U}_{k\perp}^\top \cM_k(\cX^*)\right\|_\tF
	\end{split}
\end{equation*}
and for any $k=1,2,3$,
\begin{equation}\label{eq-Gaussian-eq1}
	\begin{split}
		\left\|\tilde{\U}_{k\perp}^\top \cM_k(\cX^*)\right\|_\tF  &= \left\|\tilde{\U}_{k\perp}^\top \U_k\U_k^\top \cM_k(\cX^*)\right\|_\tF \\
		& \leq \left\|\tilde{\U}_{k\perp} \U_k\right\|_\tF \left\|\U_k^\top \cM_k(\cX^*)\right\| \\
		& = \left\|\sin \Theta (\tilde{\U}_k, \U_k)\right\|_\tF \cdot \left\| \cM_k(\cX^*)\right\| \\
		& \overset{\eqref{eq-Gaussian-hpc}}{\leq} \frac{\underline\lambda\sqrt{p_kr_k} + \sqrt{p_1p_2p_3r_k}}{\underline\lambda^2} \overline{\lambda}.
	\end{split}
\end{equation}
In the mean time, we also have
\begin{equation}\label{eq-Gaussian-eq2}
    \begin{split}
        & \left\|\cZ \times_1 \Proj_{\tilde{\U}_1} \times_2 \Proj_{\tilde{\U}_2} \times_3 \Proj_{\tilde{\U}_3}\right\|_\tF = \sup_{\substack{\mathcal{T}\in \mathbb{R}^{p_1\times p_2\times p_3}, \\ \|\cT\|_\tF \leq 1}} \left\langle \cZ, \llbracket \cT; \tilde{\U}_1^\top, \tilde\U_2^\top, \tilde \U_3^\top  \rrbracket \right\rangle \\ 
        &  \overset{\eqref{eq-Gaussian-hpc}}{\leq} C\left(\sqrt{r_1r_2r_3} + \sum_{k=1}^3 \sqrt{p_kr_k}\right).
    \end{split}
\end{equation}
Combining \eqref{eq-Gaussian-eq1} and \eqref{eq-Gaussian-eq2}, we obtain 
\begin{equation*}
	\begin{split}
		\left\|\cX^{(0)} - \cX^*\right\|_\tF^2 \leq C\left(\kappa^2\frac{p_1p_2p_3 \overline r}{\underline\lambda^2} + \kappa^2\sum_{k=1}^3 p_kr_k + r_1r_2r_3\right),
	\end{split}
\end{equation*}
where $\overline r = \max\{r_1,r_2,r_3\}$. 

Now we start to apply Theorem \ref{thm:local-convergence-PGD}. One can verify that the quadratic loss $L(\cX) = \frac{1}{2}\left\|\cX-\cY\right\|_\tF^2$ satisfies RCG$(\alpha,\beta,\bbR^{p_1\times p_2 \times p_3})$ with $\alpha = \beta = 1/2$. 
Given $\underline \lambda \geq C_0 \overline p^{3/4}\overline r^{1/4}$ for some sufficiently large $C_0$ which only depends on $\kappa$, we then have $\left\|\cX^{(0)} - \cX^*\right\|_\tF^2 \leq c\underline\lambda^2/\kappa^2$, which is the initialization condition required by Theorem \ref{thm:local-convergence-PGD} in the context of sub-Gaussian tensor PCA. Also, by the second inequality in \eqref{eq-Gaussian-hpc}, we have
\begin{equation*}
    \underline{\lambda}^2 \geq C\xi^2,
\end{equation*}
and the signal-noise-ratio condition in Theorem \ref{thm:local-convergence-PGD} is also satisfied. In conclusion, we see that, with probability at least $1-C\exp(c\underline p)$, the conditions in Theorem \ref{thm:local-convergence-PGD} are all satisfied, and vanilla gradient descent (Algorithm \ref{alg:PGD}) achieves the following statistical error bound after sufficient iterations:
\begin{equation*}
    \left\|\cX^{(T)} - \cX^*\right\|_\tF \leq C\xi \leq C\left(\sqrt{r_1r_2r_3} + \sum_{k=1}^3 \sqrt{p_kr_k}\right).
\end{equation*}
Now the proof is finished.\qquad $\square$

\subsection{Proof of Theorem \ref{thm:regression-sharp}}
By rescaling the overall model by a factor of $1/\sqrt{n}$, we can assume without loss of generality each entry of the design tensor $\cA_i$ comes from i.i.d. $N(0, 1/n)$ and the noise $\varepsilon_i \sim N(0, \sigma^2/n)$. Recall the loss function is $L(\cX) = \frac{1}{2}\sum_{i=1}^n\left(\langle \cA_i, \cX\rangle - y_i\right)^2$. Let $\cA:\bbR^{p_1\times p_2 \times p_3} \rightarrow \bbR^n$ be the linear operator such that $\left[\cA(\cX)\right]_i = \langle \cA_i, \cX\rangle$ and $\cA^*$ be the adjoint operator of $\cA$:
\begin{equation*}
	\cA^*(x) = \sum_{j=1}^n x_j \cA_j,\qquad x \in \bbR^n.
\end{equation*}
Then we can rewrite the model as
\begin{equation*}
    y = \cA(\cX^*) + \varepsilon \in \bbR^n,\quad \varepsilon \sim N(0, (\sigma^2/n) \I_n),
\end{equation*}
and the loss function can be written as $L(\cX) = \frac{1}{2}\left\|\cA(\cX) - y\right\|_2^2$.

The proof idea is the same as the proof of Theorem \ref{thm:local-convergence-PGD}. However, in tensor regression, we have an analytical form of the gradient: $\nabla L(\cX) = \cA^*(\cA(\cX)-y)$, which helps us to build some tighter results. To make the proof clear and comparable with the proof of Theorem \ref{thm:local-convergence-PGD}, we also separate the proof into several steps.  
\begin{enumerate}[leftmargin=*]
	\item[Step 1] (Notations and Conditions) We follow the same notations as we defined in step 1 in the proof of Theorem \ref{thm:local-convergence-PGD}. Similarly, we assume the following conditions:
\begin{equation}\label{eq-regression-asmp-1}
    \begin{split}
        &\left\|\U_k^{(t)}\right\| \leq 1.01b, \quad \left\|\cM_k(\cS^{(t)})\right\| \leq \frac{1.01\overline\lambda}{b^3}, \quad  k=1,2,3
    \end{split}
\end{equation}
\begin{equation}\label{eq-regression-asmp-2}
    E^{(t)} \leq c_0\frac{\underline\lambda^{1/2}}{\kappa^{3/2}},\qquad t=1,\ldots,t_{max}
\end{equation}
where $c_0$ is some absolute constant we will specify later. Besides, we impose the following high-probability conditions:
\begin{enumerate}[label=(\subscript{A}{\arabic*})]
		\item
		\begin{equation}\label{eq-regression-A-1}
			\begin{split}
				  \xi := \sup_{\substack{\mathcal{T}\in \mathbb{R}^{p_1\times p_2\times p_3}, \\ \rank(\mathcal{T})\leq (r_1,r_2,r_3), \\ \|\mathcal{T}\|_\tF \leq 1}} \left\langle \nabla L(\cX^*), \mathcal{T} \right\rangle = \sup_{\substack{\mathcal{T}\in \mathbb{R}^{p_1\times p_2\times p_3}, \\ \rank(\mathcal{T})\leq (r_1,r_2,r_3), \\ \|\mathcal{T}\|_\tF \leq 1}} \left\langle \cA^*(\varepsilon), \mathcal{T} \right\rangle \leq 2\sigma\sqrt{\frac{ df}{n}},
			\end{split}
		\end{equation}
		where $df := r_1r_2r_3 + \sum_{k=1}^3 p_kr_k$.
		\item For any tensor $\cX \in \bbR^{p_1 \times p_2 \times p_3}$ such that $\rank(\cX) \leq (3r_1,3r_2,3r_3)$,
		\begin{equation}\label{eq-regression-A-2} 
		    \begin{split}
		        &\frac{9}{10}\left\|\cX\right\|_\tF \leq \left\|\cA(\cX)\right\|_2 \leq \frac{11}{9} \left\|\cX\right\|_\tF^2, \\
		        & \left\|\cA^*(\cA(\cX))\right\|_2 \leq \frac{4}{3}\left(\sqrt{\frac{p_1p_2p_3}{n}} \vee 1\right) \left\|\cX\right\|_\tF.
		    \end{split}
		\end{equation}
	\end{enumerate}
We will verify the validity of these assumptions at the end of the proof. Also, we assume $b = \overline{\lambda}^{1/4}$ to simplify the proof.
\item[Step 2] In this step, we exactly follows the proof of \eqref{eq-Et-descent} in Theorem \ref{thm:local-convergence-PGD} and obtain
\begin{equation}\label{eq-regression-Et-descent}
    E^{(t+1)} \leq E^{(t)} - 2\eta\left(Q_{\cS,1} + \sum_{k=1}^3 Q_{k,1}\right) + \eta^2 \left(Q_{\cS,2} + \sum_{k=1}^3 Q_{k,2}\right), 
\end{equation}
where
\begin{equation*}
	\begin{split}
		& Q_{k,1} = \left\langle \cX^{(t)} - \cX_k^{(t)}, \nabla L(\cX^{(t)}) \right\rangle + \frac{a}{4} \left\| \U_k^{(t)\top} \U_k^{(t)} - b^2\I_{r_k}\right\|_\tF^2 - \frac{c_0}{4}a\sqrt{\underline{\lambda}}\left\| \U_k^{(t)} - \U_k^* \R_k^{(t)}\right\|_\tF^2, \\
		& Q_{k,2} = 8\overline\lambda^{2}b^{-2}\left(\xi^2 + \left\|\nabla L(\cX^{(t)}) - \nabla L(\cX^*)\right\|_\tF^2  \right) + \frac{5}{2}a^2b^2 \left\| \U_k^{(t)\top} \U_k^{(t)}-b^2 \I_{r_k}\right\|_\tF^2,\\
		& Q_{\cS,1} = \left\langle \cX^{(t)} - \cX_\cS^{(t)}, \nabla L(\cX^{(t)}) \right\rangle, \\
		& Q_{\cS,2} = 4b^6\left(\xi^2 + \left\|\nabla L(\cX^{(t)}) - \nabla L (\cX^*)\right\|_\tF^2\right).
	\end{split}
\end{equation*}
\item[Step 3] In this step, we provide a sharper lower bound for $Q_{\cS,1} + \sum_{k=1}^3 Q_{k,1}$:
\begin{equation}\label{eq-regression-lower-bound-linear-term}
    \begin{split}
		& Q_{\cS,1} + \sum_{k=1}^3 Q_{k,1} \\
		& \geq \rho\frac{\underline{\lambda}^{3/2}}{\kappa^{1/2}}E^{(t)} +  \frac{a}{8}\sum_{k=1}^3\left\|\U_k^{(t)\top}\U_k^{(t)}- b^2\I_{r_k}\right\|_\tF^2 + \frac{1}{20}\left\|\cX^{(t)} - \cX^*\right\|_\tF^2 - C_1\kappa^2\xi^2,
	\end{split}
\end{equation}
where $\rho$ is some universal small constant and $C_1$ is some universal big constant.
First of all, by the proof of Lemma \ref{lm-X-decomposition}, we have
\begin{equation*}
    4\cX^{(t)} - \cX_\cS^{(t)} - \sum_{k=1}^3 \cX_k^{(t)} = \cX^{(t)} - \cX^* + \cH_\varepsilon,
\end{equation*}
where
\begin{equation}\label{eq-regresion-Delta-definition}
    \begin{split}
        & \cH_\varepsilon = \llbracket\cS^*;\H_1,\H_2,\H_3\rrbracket + \sum_{k=1}^3 \cS^* \times_k \U_k\R_k^\top \times_{k+1} \H_{k+1} \times_{k+2} \H_{k+2} \\
		& \qquad + \sum_{k=1}^3 \cH_\cS \times_k \H_k \times_{k+1} \U_{k+1}\R_{k+1}^\top \times_{k+2} \U_{k+2}\R_{k+2}^\top,
    \end{split}
\end{equation}
\begin{equation*}
    \begin{split}
        & \H_k := \U_k^* - \U_k\R_k^\top, \quad k = 1,2,3,\\
        & \cH_\cS := \cS^* - \llbracket\cS;\R_1,\R_2,\R_3\rrbracket.
    \end{split}
\end{equation*}
Then it follows that 
\begin{equation}\label{eq-regression-linear-decomposition}
    \begin{split}
        & \left\langle 4\cX^{(t)} - \cX_\cS^{(t)} - \sum_{k=1}^3 \cX_k^{(t)}, \nabla L(\cX^{(t)}) \right\rangle = \left\langle \cX^{(t)} - \cX^*, \nabla L(\cX^{(t)}) - \nabla L(\cX^*)\right\rangle \\
        & + \left\langle \cH_\varepsilon, \nabla L(\cX^{(t)}) - \nabla L(\cX^*) \right\rangle + \left\langle 4\cX^{(t)} - \cX_\cS^{(t)} - \sum_{k=1}^3 \cX_k^{(t)}, \nabla L(\cX^*) \right\rangle.
    \end{split}
\end{equation}
Since $\nabla L(\cX^{(t)}) - \nabla L(\cX^*) = \cA^*\cA(\cX^{(t)} - \cX^*)$, we firstly have
\begin{equation}\label{eq-regression-signal-inner-product}
    \begin{split}
        & \left\langle \cX^{(t)} - \cX^*, \nabla L(\cX^{(t)}) - \nabla L(\cX^*)\right\rangle = \left\langle \cX^{(t)} - \cX^*, \cA^*\cA(\cX^{(t)}-\cX^*)\right\rangle \\
        & = \left\langle \cA(\cX^{(t)} - \cX^*), \cA(\cX^{(t)}-\cX^*)\right\rangle = \left\|\cA(\cX^{(t)} - \cX^*)\right\|_2^2 \overset{\eqref{eq-regression-A-2}}{\geq} \frac{4}{5}\left\|\cX^{(t)} - \cX^*\right\|_\tF^2.
    \end{split}
\end{equation}
Next we give an upper bound of $\left|\left\langle \cH_\varepsilon, \nabla L(\cX^{(t)}) - \nabla L(\cX^*) \right\rangle\right|$. To this end, we need to use the fact that $\cH_\varepsilon$ is a summation of rank-$(r_1,r_2,r_3)$ tensors. By Lemma \ref{lm-incoherent-inner-product-difference}, for any rank-$(r_1,r_2,r_3)$ tensor $\cX' \in \bbR^{p_1\times p_2 \times p_3}$, we have 
\begin{equation*}
    \begin{split}
        & \left|\left\langle \cX', \nabla L(\cX^{(t)}) - \nabla L(\cX^*)  \right\rangle\right| = \left|\left \langle \cA(\cX'), \cA(\cX^{(t)} - \cX^*) \right \rangle\right| \\
        & \overset{\eqref{eq-regression-A-2}}{\leq} \left|\left\langle \cX', \cX^{(t)} - \cX^* \right\rangle\right| + \frac{1}{2}\left\|\cX'\right\|_\tF \left\|\cX^{(t)} - \cX^*\right\|_\tF \leq \frac{3}{2}\left\|\cX'\right\|_\tF\left\|\cX^{(t)}-\cX^*\right\|_\tF.
    \end{split}
\end{equation*}
Thus by plugging each component in the right hand side of \eqref{eq-regresion-Delta-definition} in $\cX'$ and using triangle inequality, we obtain the following upper bound:
\begin{equation*}
    \left|\left\langle \cH_\varepsilon, \nabla L(\cX^{(t)}) - \nabla L(\cX^*)  \right\rangle\right|  \leq \frac{3}{2}C_\Delta\left\|\cX^{(t)} - \cX^*\right\|_\tF,
\end{equation*}
where $C_\Delta$ is defined as
\begin{equation*}
    \begin{split}
        & C_\Delta := \left\|\llbracket\cS^*;\H_1,\H_2, \H_3\rrbracket\right\|_\tF + \sum_{k=1}^3 \left\|\cS^* \times_k \U_k\R_k^\top \times_{k+1} \H_{k+1}\times_{k+2} \H_{k+2}\right\|_\tF \\
		& \qquad + \sum_{k=1}^3 \left\|\cH_\cS \times_k \H_k \times_{k+1} \U_{k+1}\R_{k+1}^\top \times_{k+2} \U_{k+2}\R_{k+2}^\top\right\|_\tF.
    \end{split}
\end{equation*}
By the proof of Lemma \ref{lm-X-decomposition}, we know that 
\begin{equation*}
	\begin{split}
		C_\Delta & \overset{\eqref{eq-regression-asmp-1}}{\leq} 1.01\overline\lambda b^{-3}(E^{(t)})^{3/2} + 3(1.01b)^2E^{(t)} + 3\cdot \frac{1.01\overline{\lambda}}{b^3}\cdot (1.01b)E^{(t)} \\
		\overset{\eqref{eq-regression-asmp-2}}{\leq} & \left(1.01 c_0^{1/2}\overline\lambda \underline{\lambda}^{1/4} b^{-3} + 6 \cdot 1.01^2 b^2\right)E^{(t)} \\
		\leq & 6.5\kappa^{1/2}\underline{\lambda}^{1/2} E^{(t)}. \\
	\end{split}
\end{equation*}
Note that the above inequality is similar to \eqref{ineq:H_epsilon} in the proof of Theorem \ref{thm:local-convergence-PGD}.
Combining all above, we obtain:
\begin{equation}\label{eq-regression-delta}
    \begin{split}
        & \left|\left\langle \cH_\varepsilon, \nabla L(\cX^{(t)}) - \nabla L(\cX^*) \right\rangle\right| \leq \frac{3}{2} C_\Delta\left\|\cX^{(t)}-\cX^*\right\|_\tF \\
        \leq & 10\left\|\cX^{(t)}-\cX^*\right\|_\tF \cdot \left(\kappa^{1/2}\underline\lambda^{1/2}E^{(t)}\right)\\
        \leq & \frac{1}{4}\left\|\cX^{(t)} - \cX^*\right\|_\tF^2 + 100\kappa\underline\lambda (E^{(t)})^2 \\
        \overset{\eqref{eq-regression-asmp-2}}{\leq} & \frac{1}{4}\left\|\cX^{(t)} - \cX^*\right\|_\tF^2 + \frac{100c_0\underline{\lambda}^{3/2}}{\kappa^{1/2}}E^{(t)}. 
    \end{split}
\end{equation}
In the meantime, by the same argument in \eqref{eq-noise-inner-product} in the proof of Theorem \ref{thm:local-convergence-PGD}, we can show that
\begin{equation}\label{eq-regression-noise-inner-product}
    \left|\left\langle 4\cX^{(t)} - \cX_\cS^{(t)} - \sum_{k=1}^3 \cX_k^{(t)}, \nabla L(\cX^*) \right\rangle\right| \leq c_1 \frac{\underline{\lambda}^{3/2}}{\kappa^{1/2}}E^{(t)} + \frac{81\kappa^2}{4c_1}\xi^2,
\end{equation}
and $c_1$ is some universal constant which will be specified later. Now
combining \eqref{eq-regression-linear-decomposition}, \eqref{eq-regression-signal-inner-product}, \eqref{eq-regression-delta} and \eqref{eq-regression-noise-inner-product} and specifying $a=\frac{2b^4}{3\kappa^2}$ (which is similar to \eqref{ineq:convergece-aggregate}), we obtain that 
\begin{equation*}
    \begin{split}
        & Q_{\cS,1} + \sum_{k=1}^3 Q_{k,1} \geq \frac{1}{2}\left(\left\|\cX^{(t)} - \cX^*\right\|_\tF^2 +  \frac{1}{6}\kappa^{-2}b^4 \sum_{k=1}^3 \left\|\U_k^{(t)\top}\U_k^{(t)} - b^2\I_{r_k}\right\|_\tF^2\right) \\
        & \qquad + \frac{a}{8}\sum_{k=1}^3\left\|\U_k^{(t)\top}\U_k^{(t)} - b^2\I_{r_k}\right\|_\tF^2- \frac{c_0\underline{\lambda}^{3/2}}{6\kappa^{1/2}}\sum_{k=1}^3\left\|\U_k^{(t)} - \U_k^*\R_k^{(t)}\right\|_\tF^2 \\
        & \qquad + \frac{1}{20}\left\|\cX^{(t)} - \cX^*\right\|_\tF^2 - \frac{(100c_0 + c_1)\underline{\lambda}^{3/2}}{\kappa^{1/2}}E^{(t)} - C_1\kappa^2 \xi^2 \\
        & \overset{\text{Lemma \ref{lm-equivalent-criteria}}}{\geq} \frac{\underline{\lambda}^{3/2}}{960\kappa^{1/2}}E^{(t)} + \frac{a}{8}\sum_{k=1}^3\left\|\U_k^{(t)\top}\U_k^{(t)} - b^2\I_{r_k}\right\|_\tF^2 + \frac{1}{20}\left\|\cX^{(t)} - \cX^*\right\|_\tF^2 \\
        & \qquad\qquad - \frac{(100c_0+c_2)\underline \lambda^{3/2}}{\kappa^{1/2}} E^{(t)} - C_1\kappa^2\xi^2\\
        & = \rho\frac{\underline{\lambda}^{3/2}}{\kappa^{1/2}}E^{(t)} +  \frac{a}{8}\sum_{k=1}^3\left\|\U_k^{(t)\top}\U_k^{(t)}- b^2\I_{r_k}\right\|_\tF^2 + \frac{1}{20}\left\|\cX^{(t)} - \cX^*\right\|_\tF^2 - C_1\kappa^2\xi^2,
    \end{split}
\end{equation*}
Here in the last step, we specify $c_0$ and $c_1$ to be small constants such that $\rho := \frac{1}{960} - 100c_0 - c_1 > 0$, and $C_1 = \frac{81}{4c_1}$, which gives \eqref{eq-regression-lower-bound-linear-term}. 
\item[Step 4] In this step, we provide the error contraction of $E^{(t)}$. By plugging in $a = \frac{2b^4}{3\kappa^2}$, we have
\begin{equation*}
	\begin{split}
		 & Q_{\cS,2} + \sum_{k=1}^3 Q_{k,2} \\
		 & \leq 12 b^6\left\| \nabla L(\cX^{(t)}) - \nabla L(\cX^*)\right\|_\tF^2 + \frac{2b^{10}}{\kappa^4} \sum_{k=1}^3\left\| \U_k^{(t)\top} \U_k^{(t)} - b^2\I_{r_k}\right\|_\tF^2 + 12 b^6 \xi^2,
	\end{split}
\end{equation*}
and
\begin{equation*}
    \left\|\nabla L(\cX^{(t)}) - \nabla L(\cX^*)\right\|_\tF^2 = \left\|\cA^*\cA(\cX^{(t)}-\cX^*)\right\|_\tF^2 \overset{\eqref{eq-regression-A-2}}{\leq }\frac{16}{9}\left(\frac{p_1p_2p_3}{n} \vee 1\right) \left\|\cX^{(t)} - \cX^*\right\|_\tF^2.
\end{equation*}
Now by taking $\eta=\eta_0/b^6$ in \eqref{eq-regression-Et-descent} and applying \eqref{eq-regression-lower-bound-linear-term}, we have
\begin{equation}
	\begin{split}
		E^{(t+1)} \leq & \left(1-\frac{2\rho\eta_0}{\kappa^2}\right)E^{(t)} - \left(\frac{1}{10b^6}\eta_0 - \frac{64}{3b^6}\left(\frac{p_1p_2p_3}{n}\vee 1\right)\eta_0^2\right)\left\| \cX^{(t)} - \cX^*\right\|_\tF^2\\
		& - \left(\frac{b^{-2}}{12\kappa^2}\eta_0 - \frac{2b^{-2}}{\kappa^4}\eta_0^2\right)\sum_{k=1}^3\left\| \U_k^{(t)\top} \U_k^{(t)} - b^2\I_{r_k}\right\|_\tF^2 + 2C_1b^{-6}\kappa^2\eta_0\xi^2 + 12b^{-6}\eta_0^2 \xi^2.
	\end{split}
\end{equation}
Then as long as $\eta_0 < \frac{3}{640}\left(\frac{n}{p_1p_2p_3} \wedge 1\right)$,
\begin{equation*}
	\begin{split}
		& \frac{1}{10b^6}\eta_0 - \frac{64}{3b^6}\left(\frac{p_1p_2p_3}{n}\vee 1\right)\eta_0^2 > 0,\\
		& \frac{b^{-2}}{12\kappa^2}\eta_0 - \frac{2b^{-2}}{\kappa^4}\eta_0^2 > 0,
	\end{split}
\end{equation*}
and it follows that
\begin{equation}\label{eq-regression-convergence-E-next-step}
	\begin{split}
		& E^{(t+1)} \leq \left(1-\frac{2\rho\eta_0}{\kappa^2}\right)E^{(t)} + (2C_1 + 1)\kappa^2b^{-6}\eta_0 \xi^2.
	\end{split}
\end{equation}
By induction (as we proved for \eqref{eq-upper-bound-Et}), we can then show that
\begin{equation}\label{eq-regression-upper-bound-Et}
    \begin{split}
        E^{(t)} & \leq \frac{(2C_1+1)\kappa^4\xi^2}{2\rho b^6} +  \left(1 - \frac{2\rho\eta_0}{\kappa^2}\right)^t E^{(0)} \\
        & \overset{\eqref{eq-regression-A-1}}{\leq} C\kappa^4\frac{\sigma^2 \cdot df}{n b^{6}}+  \left(1 - \frac{2\rho\eta_0}{\kappa^2}\right)^t E^{(0)}.
    \end{split}
\end{equation}
Applying Lemma \ref{lm-equivalent-criteria}, we then have
\begin{equation*}
    \begin{split}
        & \left\|\cX^{(t)} - \cX^*\right\|_\tF^2 \leq 42b^6 E^{(t)} \leq C\kappa^4\frac{\sigma^2 \cdot df}{n} + 42\left(1 - \frac{2\rho\eta_0}{\kappa^2}\right)^t E^{(0)}. \\
        & \leq C\left(\kappa^4\frac{\sigma^2 \cdot df}{n} +  \left(1 - \frac{2\rho\eta_0}{\kappa^2}\right)^t \kappa^2 \left\|\cX^{(0)}-\cX^*\right\|_\tF^2\right).
    \end{split}
\end{equation*}
Therefore, for sufficiently large $T$, one can see that we have $\left\|\cX^{(T)} - \cX^*\right\|_\tF^2 \leq C\kappa^4 \frac{\sigma^2 \cdot df}{n}$.
\item[Step 5] Finally, we  check the validity of conditions. We  first assume the conditions \eqref{eq-regression-A-1}, \eqref{eq-regression-A-2} hold, and verify \eqref{eq-regression-asmp-1} and \eqref{eq-regression-asmp-2}; then we show that \eqref{eq-regression-A-1} and \eqref{eq-regression-A-2} hold with high probability. We  start from \eqref{eq-regression-asmp-2}. By the proof of \cite[Theorem 4]{zhang2018ISLET}, we know as long as 
\begin{equation*}
    n \geq C\frac{\left(\left\|\cX^*\right\|_\tF^2 + \sigma^2\right)\left(\overline p^{3/2}+\kappa \overline p \cdot \overline{r}\right)}{\underline \lambda^2},
\end{equation*}
the initialization estimator has the following error bound:
\begin{equation*}
    \left\|\cX^{(0)} - \cX^*\right\|_\tF^2 \leq C'\frac{df(\sigma^2 + \left\|\cX^*\right\|_\tF^2)}{n} \leq C'\frac{df \cdot \left\|\cX^*\right\|_\tF^2}{n},
\end{equation*}
for some universal constants $C,C'$.
Then by Lemma \ref{lm-equivalent-criteria}, we have
\begin{equation*}
    \begin{split}
        & E^{(0)} \leq 11\kappa^2  b^{-6}\left\|\cX^{(0)} - \cX^*\right\|_\tF^2 \leq \frac{C\kappa^2 df \cdot \overline{r}\overline{\lambda}^2}{n\overline{\lambda}^{3/2}} \leq c_0\frac{\underline{\lambda}^{1/2}}{\kappa^{3/2}}.
    \end{split}
\end{equation*}
Here we use the assumption that
\begin{equation*}
    n \geq C\kappa^4df \cdot \overline{r}.
\end{equation*}
Thus we show that \eqref{eq-regression-asmp-2} holds at $t=0$. Now suppose \eqref{eq-regression-asmp-2} holds at $t$, we show it also holds at $t+1$. By \eqref{eq-regression-convergence-E-next-step},
\begin{equation*}
	\begin{split}
		E^{(t+1)} & \leq  \left(1-\frac{2\rho\eta_0}{\kappa^2}\right)E^{(t)} + (2C_1 + 1)\kappa^2b^{-6}\eta_0 \xi^2 \\
		& \leq \left(1 - \frac{2\rho\eta_0}{\kappa^2}\right)c_0\frac{ \underline\lambda^{1/2}}{\kappa^{3/2}} + (2C_1 + 1)\kappa^2b^{-6}\eta_0 \xi^2 \\
		& = c_0\frac{ \underline\lambda^{1/2}}{\kappa^{3/2}} - \left(\frac{2c_0\rho \underline\lambda^{1/2}}{\kappa^{3.5}} - (2C_1 + 1)\kappa^{1/2}\underline\lambda^{-3/2} \xi^2\right)\eta_0 \\
		& \leq c_0\frac{ \underline\lambda^{1/2}}{\kappa^{3/2}}.
	\end{split}
\end{equation*}  
Here in the last inequality we use the signal-noise-ratio assumption: $\underline{\lambda}^2 \geq C\kappa^4\frac{df \sigma^2}{n}$. Thus \eqref{eq-regression-asmp-2} is verified for any $t = 0,\ldots, t_{max}$. The verification of \eqref{eq-regression-asmp-1} is the same as we did for \eqref{eq-Assumption1}, and we omitted it here.
Now we start to show that \eqref{eq-regression-A-1}, \eqref{eq-regression-A-2} hold with high probability. Since $n \geq C df$, we know that \eqref{eq-regression-A-2} holds with probability at least $1 - C\exp\left(-c\cdot df\right)$ by Lemma \ref{lm-random-Gaussian-design}. 

On the other hand, by definition,
\begin{equation*}
    \begin{split}
		\xi := \sup_{\substack{\mathcal{\cS}\in \mathbb{R}^{r_1\times r_2\times r_3}, \left\|\cS\right\|_\tF \leq 1 \\ \W_k \in \bbR^{p_1\times r_1}, \left\|\W_k\right\| \leq 1, k=1,2,3 }} \left\langle \cA^*(\varepsilon), \llbracket \cS; \W_1, \W_2, \W_3\rrbracket \right\rangle.
	\end{split}
\end{equation*}
Notice that for any fixed $\cS,W_1,\W_2,\W_3$, one have
\begin{equation*}
    \begin{split}
        & \left\langle \cA^*(\varepsilon), \llbracket \cS; \W_1, \W_2, \W_3\rrbracket \right\rangle = \left\langle \varepsilon, \cA\left(\llbracket \cS; \W_1, \W_2, \W_3\rrbracket\right)\right\rangle,
    \end{split}
\end{equation*}
which has normal distribution with mean zero and variance $\frac{\tau^2\sigma^2}{n}$ when conditional on $\{\cA_i\}_{i=1}^n$ , with $\tau := \left\|\cA\left(\llbracket \cS; \W_1, \W_2, \W_3\rrbracket\right)\right\|_2$. Thus we have
\begin{equation*}
    \bbP\left( \left.\left\langle \cA^*(\varepsilon), \llbracket \cS; \W_1, \W_2, \W_3\rrbracket \right\rangle \geq t \right|\{\cA_i\}_{i=1}^n\right) \leq e^{-\frac{nt^2}{2\tau^2\sigma^2}}.
\end{equation*}
Under the event $A := \left\{\text{\eqref{eq-regression-A-2} holds}\right\}$, $\tau^2 \leq \frac{3}{2}\left\|\llbracket \cS; \W_1, \W_2, \W_3\rrbracket\right\|_\tF^2 \leq \frac{3}{2}$, then it follows that
\begin{equation}\label{eq-regression-N1-fix-bound}
    \bbP\left( \left.\left\langle \cA^*(\varepsilon), \llbracket \cS; \W_1, \W_2, \W_3\rrbracket\right\rangle \geq t \right|A\right) \leq e^{-\frac{nt^2}{3\sigma^2}}.
\end{equation}
Now for $k=1,2,3$, we can construct an $\varepsilon$-net $\left\{\V_k^{(1)}, \ldots, \V_k^{(N_k)}\right\}$ of Stiefel $\bbO_{p_k, r_k}$ with metric $d(\V_1,\V_2) = \left\|\V_1\V_1^\top - \V_2\V_2^\top\right\|$,  such that 
\begin{equation*}
	\begin{split}
		\sup_{\V_k \in \bbO_{p_k,r_k}}\min_{i \leq N_k} d(\V_k, \V_k^{(i)}) \leq \varepsilon
	\end{split}
\end{equation*}
and $N_k \leq \left(\frac{c_0}{\varepsilon}\right)^{p_k(p_k-r_k)}$ for some absolute constant $c_0$.
Also, we can construct $\varepsilon$-net $\{\cS^{(1)}, \ldots, \cS^{(N_\cS)}\}$ for the core tensors on $\mathbb B_{r_1,r_2,r_3}$, which is the $l_2$ unit ball for $\bbR^{r_1\times r_2 \times r_3}$, such that
\begin{equation*}
	\sup_{\cS \in \mathbb B_{r_1,r_2,r_3}}\min_{i \leq N_\cS} \left\|\cS - \cS^{(i)}\right\|_\tF \leq \varepsilon
\end{equation*} 
and $N_\cS \leq \left(\frac{c_1}{\varepsilon}\right)^{r_1r_2r_3}.$
Then by the similar $\varepsilon-$net argument in Lemma \ref{lm-concentration-Gaussian-xi}, we can show that 
\begin{equation*}
    \begin{split}
        \bbP\left(\left.\xi > 2\sqrt{\frac{\sigma^2 df}{n}} \right| A\right) \leq C\exp\left(-c \cdot df)\right),
    \end{split}
\end{equation*}
and it follows that
\begin{equation*}
    \begin{split}
        & \bbP\left(\xi < 2\sqrt{\frac{\sigma^2 df}{n}} \right) \geq \bbP\left(A\right) \cdot \bbP\left(\left.N_1 > C\sqrt{\frac{\sigma^2 df}{n}} \right| A\right)  \\
        & \geq \left(1-Ce^{-c \cdot df}\right)\left(1-Ce^{-c \cdot df}\right) \geq 1 - 2Ce^{-c \cdot df}.
    \end{split}
\end{equation*}
Thus we proved that with probability at least $1-Ce^{-c \cdot df}$, \eqref{eq-regression-A-2} holds, now the proof is finished.\qquad $\square$

\end{enumerate}

\subsection{Proof of Theorem \ref{thm:Poisson}}
We first introduce some notations to simplify the proof. Recall for each $(j,k,l) \in [p_1] \times [p_2] \times [p_3]$, $\cY_{jkl} \sim \text{Poisson}(\nu_{jkl})$, where $\nu_{jkl} := I\exp(\cX^*_{jkl})$. By Assumption \ref{asmp-incoherence}, we have $|\cX^*_{jkl}| \leq B$. We next define a new random tensor $\cY'$ such that $\cY'_{jkl} = \cY_{jkl}1_{\left\{\frac{\nu_{jkl}}{10} \leq \cY_{jkl} \leq 10\nu_{jkl}\right\}} + \nu_{jkl}1_{\left\{\cY_{jkl} \not \in [\frac{\nu_{jkl}}{10}, 10\nu_{jkl}]\right\}}$, and we further define $\cZ:= \cZ_1 + \cZ_2$ where
\begin{equation*}
	\begin{split}
		\cZ_1 &= \log\left(\cY'+1/2\right) - \bbE \log(\cY'+1/2),\\
		\cZ_2 &= \bbE\log(\cY' + 1/2) - \cX^*-\log I.
	\end{split}
\end{equation*} 
The following conditions are introduced for the proof:
\begin{enumerate}[label=(\subscript{A}{\arabic*})]
	\item
		\begin{equation}\label{eq-Poisson-hpc-1}
            \cY_{jkl} = \cY'_{jkl}, \quad\forall (j,k,l) \in [p_1]\times[p_2]\times [p_3]
        \end{equation}
		\item
		\begin{equation}\label{eq-Poisson-hpc-2}
            \left\|\cM_k(\cZ_1)\right\| \leq 2K_0\sqrt{\frac{e^B}{I}}\left(\sqrt{p_{-k}}+\sqrt{p_k}\right),\quad k = 1,2,3
        \end{equation}
		\item
		\begin{equation}\label{eq-Poisson-hpc-3}
            \sup_{\left\|\V_k\right\| \leq 1, k=1,2,3} \left\|\cZ_1 \times_1 \V_1^\top \times_2 \V_2^{\top} \times_3 \V_3^{\top}\right\|_\tF \leq C_1\frac{\sqrt{df}}{\sqrt{I/e^B}},
        \end{equation}
        \item
		\begin{equation}\label{eq-Poisson-hpc-4}
            \sup_{\substack{\mathcal{T}\in \mathbb{R}^{p_1\times p_2\times p_3}, \\ \rank(\mathcal{T})\leq (r_1,r_2,r_3), \\ \|\mathcal{T}\|_\tF \leq 1}} \left\langle \cY/I - \exp(\cX^*), \mathcal{T} \right\rangle \leq C_2\frac{\sqrt{df}}{\sqrt{I/e^B}}.
        \end{equation}
	\end{enumerate}
Here $K_0,C_1,C_2$ are some universal constants and $df:=r_1r_2r_3 + \sum_{k=1}^3 p_kr_k$. We now start the analysis assuming that the above conditions hold, and we will verify that they hold with high probability at last.\\
First of all, by Lemma \ref{lm-log-error}, if we set $I \geq C e^{B}$, then for each $(j,k,l) \in [p_1]\times [p_2] \times [p_3]$, we have
\begin{equation}\label{eq-TPoisson-bias-subgaussian}
	\begin{split}
		& \left| \bbE\log(\cY_{jkl}' + 1/2) - \cX_{jkl}^* - \log I \right| \leq 4\nu_{jkl}^{-1} \leq \frac{4e^B}{I},\\
		& \left\|\log(\cY'_{jkl} + 1/2) - \bbE\log(\cY'_{jkl}+1/2)\right\|_{\psi_2} \leq K_0 \nu_{jkl}^{-1/2} \leq K_0\sqrt{\frac{e^B}{I}},
	\end{split}
\end{equation}
where $K_0$ is the constant defined in Lemma \ref{lm-log-error}. Now we have
\begin{equation*}
	\begin{split}
		& \left\|\cM_k(\cZ_2)\right\| \leq \left\|\cM_k(\cZ_2)\right\|_\tF\\
		& = \sqrt{\sum_{jkl}\left|\bbE\log(\cY_{jkl}' + 1/2) - \cX_{jkl}^* - \log I\right|^2} \overset{\eqref{eq-TPoisson-bias-subgaussian}}{\leq} \sqrt{p_1p_2p_3 \left(\frac{4e^B}{I}\right)^2} = \frac{4e^B\sqrt{p_1p_2p_3}}{I},
	\end{split}
\end{equation*}
and then 
\begin{equation*}
	\begin{split}
		& \left\|\cM_k(\cZ)\right\| \leq \left\|\cM_k(\cZ_1)\right\| + \left\|\cM_k(\cZ_2)\right\| \\
		& \overset{\eqref{eq-Poisson-hpc-2}}{\leq}  2K_0\sqrt{\frac{e^B}{I}}\left(\sqrt{p_{-k}}+\sqrt{p_{k}}\right) + \frac{4e^B\sqrt{p_1p_2p_3}}{I} \leq C\sqrt{\frac{e^B}{I}}\left(\sqrt{p_{-k}}+\sqrt{p_{k}}\right),
	\end{split}
\end{equation*}
where the last inequality comes from the assumption that $I > C\overline{p}e^B$, and $\overline p = \max\{p_1,p_2,p_3\}$.
Since $\tilde \U_k^{(0)}$ is the $r_k$ leading singular vectors of $\cM_k(\tilde\cX^{(0)})$, and $\tilde \cX^{(0)} = \log\left((\cY+1/2)/I\right) \overset{\eqref{eq-Poisson-hpc-1}}{=} \log\left((\cY' + 1/2)/I\right) = \cX^* + \cZ$, it follows that
\begin{equation*}
	\begin{split}
		& \left\|\tilde\U_{k\perp}^{(0)\top}\cM_k(\cX^*)\right\|_\tF \overset{\text{Lemma \ref{lm-low-rank-matrix-perturbation}}}{\leq} 2\sqrt{r_k}\left\|\cM_k(\cZ)\right\| \leq C\sqrt{\frac{e^B}{I}}\left(\sqrt{p_{-k}r_k} + \sqrt{p_kr_k}\right).
	\end{split}
\end{equation*}
Here $\tilde \U_{k\perp}^{(0)} \in \bbO_{p_k,p_k-r_k}$ and is orthogonal to $\tilde\U_k$. Now we can obtain the upper bound of $\left\|\cX^{(0)} - \cX^* \right\|_\tF^2$. Since
\begin{equation*}
	\begin{split}
		& \left\|\cX^* \times_1 \Proj_{\tilde\U_1^{(0)}} \times_2 \Proj_{\tilde\U_2^{(0)}} \times_3 \Proj_{\tilde\U_3^{(0)}} - \cX^*\right\|_\tF \\
		& = \left\|\cX^* \times_1 \Proj_{\tilde\U_{1\perp}^{(0)}} + \cX^* \times_1 \Proj_{\tilde\U_{1}^{(0)}} \times_2 \Proj_{\tilde\U_{2\perp}^{(0)}} + \cX^* \times_1 \Proj_{\tilde\U_{1}^{(0)}} \times_2 \Proj_{\tilde\U_{2}^{(0)}} \times_3 \Proj_{\tilde\U_{3\perp}^{(0)}}\right\|_\tF \\
		& \leq  \sum_{k=1}^3 \left\|\cX^* \times_k P_{\tilde\U_{k\perp}^{(0)}}\right\|_\tF = \sum_{k=1}^3 \left\|\tilde\U_{k\perp}^{(0)\top}\cM_k(\cX^*) \right\|_\tF \leq C\sqrt{\frac{e^B}{I}}\sum_{k=1}^3 \left(\sqrt{p_{-k}r_k} + \sqrt{p_{k}r_k}\right),
	\end{split}
\end{equation*}
we have
\begin{equation*}
	\begin{split}
		& \left\|\cX^{(0)} - \cX^*\right\|_\tF = \left\|\tilde\cX^{(0)} \times_1 \Proj_{\tilde\U_1^{(0)}} \times_2 \Proj_{\tilde\U_2^{(0)}} \times_3 \Proj_{\tilde\U_3^{(0)}} - \cX^*\right\|_\tF \\
		& \leq \left\|\cZ_1 \times_1 \Proj_{\tilde\U_1^{(0)}} \times_2 \Proj_{\tilde\U_2^{(0)}} \times_3 \Proj_{\tilde\U_3^{(0)}}\right\|_\tF + \left\|\cZ_2 \times_1 \Proj_{\tilde\U_1^{(0)}} \times_2 \Proj_{\tilde\U_2^{(0)}} \times_3 \Proj_{\tilde\U_3^{(0)}}\right\|_\tF \\
		& + \left\|\cX^* \times_1 \Proj_{\tilde\U_1^{(0)}} \times_2 \Proj_{\tilde\U_2^{(0)}} \times_3 \Proj_{\tilde\U_3^{(0)}} - \cX^*\right\|_\tF \\
		& \overset{\eqref{eq-Poisson-hpc-3}}{\leq}C\sqrt{\frac{e^B}{I}}\left(\sqrt{r_1r_2r_3} + \sum_{k=1}^3 \sqrt{p_k r_k} + \sum_{k=1}^3 \sqrt{p_{-k}r_k}\right) + \left\|\cZ_2\right\|_\tF \\
		& \leq C\sqrt{\frac{e^B}{I}}\left(\sqrt{r_1r_2r_3} + \sum_{k=1}^3 \sqrt{p_k r_k} + \sum_{k=1}^3 \sqrt{p_{-k}r_k}\right) + \frac{4e^B \sqrt{p_1p_2p_3}}{I} \\
		& \leq C\sqrt{\frac{e^B}{I}}\sum_{k=1}^3 \left(\sqrt{p_k r_k} + \sqrt{p_{-k}r_k}\right),
	\end{split}
\end{equation*}
where the last inequality comes from the assumption that  $I>e^B \overline p$. 

Now we turn to apply Theorem \ref{thm:local-convergence-PGD}. We take the loss function to be the normalized Poisson negative log-likelihood
\begin{equation*}
    L(\cX) = \frac{1}{I}\sum_{j=1}^{p_1}\sum_{k=1}^{p_2}\sum_{l=1}^{p_3} \left(-\cY_{jkl}\cX_{jkl} + I\exp(\cX_{ijk})\right)
\end{equation*}
with $\nabla L(\cX) = -\cY/I + \exp(\cX^*)$. In the convex region $\{\cX \in \bbR^{p_1\times p_2 \times p_3}: \left\|\cX\right\|_\infty \leq B\}$, $L(\cX)$ is $e^B$-smooth and $e^{-B}$-strongly convex, and by Lemma \ref{lm-convex-smooth}, we know $L(\cX)$ satisfies RCG$(\alpha,\beta,\cC)$ with $\alpha = \beta = \frac{1}{e^B+e^{-B}}$ and $\cC$ is defined according to \eqref{eq:constraint-C} and \eqref{eq-specify-sets}. Then as long as $I\underline \lambda^2 \geq C\kappa^2e^{3B}\sum_{k=1}^3 \left(p_{-k}r_k + p_kr_k\right)$ for some universal constant $C$, the initialization condition $\left\|\cX^{(0)}- \cX^*\right\|_\tF^2 \leq c_0\alpha\beta\kappa^{-2}\underline\lambda^2$ in Theorem \ref{thm:local-convergence-PGD} is satisfied. In the meantime, one can calculate that
\begin{equation*}
    \xi  = \sup_{\substack{\mathcal{T}\in \mathbb{R}^{p_1\times p_2\times p_3}, \\ \rank(\mathcal{T})\leq (r_1,r_2,r_3), \\ \|\mathcal{T}\|_\tF \leq 1}} \left\langle \cY/I - \exp(\cX^*), \mathcal{T} \right\rangle \overset{\eqref{eq-Poisson-hpc-4}}{\leq} C\frac{\sqrt{df}}{\sqrt{I/e^B}}.
\end{equation*}
Then $\underline{\lambda}^2 \geq C\kappa^4 e^{5B} \left(r_1r_2r_3 + \sum_{k=1}^3 p_kr_k\right)$ implies the signal-noise-ratio condition $\underline{\lambda}^2 \geq C\frac{\kappa^4}{\alpha^3\beta}\xi^2$ in Theorem \ref{thm:local-convergence-PGD}, and we can obtain the statistical error rate of gradient descent after sufficient steps:
\begin{equation*}
	\left\|\cX^{(T)} - \cX^*\right\|_\tF^2 \leq C\frac{\kappa^4}{\alpha^2}\xi^2 \leq \frac{C\kappa^4 e^{3B}}{I} \cdot df.
\end{equation*}
Now it suffices to check \eqref{eq-Poisson-hpc-1}-\eqref{eq-Poisson-hpc-4} actually hold with high probability. First of all, let $A$ be the event that $\cY_{jkl}' = \cY_{jkl}, \forall j,k,l$, then we have 
\begin{equation*}
	\begin{split}
		 & \bbP(A) = 1 - \bbP\left(\exists (j,k,l), \cY_{jkl} \not \in [\frac{\nu_{jkl}}{10},10\nu_{jkl}]\right) \\
		 & \geq 1 - \sum_{j,k,l} \left(\bbP(\cY_{jkl} < \frac{\nu_{jkl}}{10}) + \bbP(\cY_{jkl} > 10\nu_{jkl})\right) \\
		 & \overset{(a)}{\geq} 1 - 2\sum_{j,k,l} \exp\left(-\frac{4}{5}\nu_{jkl}\right)  \geq 1 - 2\sum_{jkl} \exp\left({-\frac{2}{5}Ie^{-B}}\right) \\
		 & = 1 - 2\exp\left(\log(p_1p_2p_3) - \frac{2e^{-B}}{5}I\right)  \overset{(b)}{\geq} 1-\frac{C}{p_1p_2p_3}.
	\end{split}
\end{equation*}
Here (a) comes from applying Poisson's tail bound (Lemma \ref{lm-Poisson-tail-bound}) and (b) is true as long as $I > Ce^{B}\log (p_1p_2p_3)$, thus \eqref{eq-Poisson-hpc-1} holds with probability at least $1-\frac{C}{p_1p_2p_3}$. By Lemma \ref{lm-log-error}, $\log\left(\cM_k(\cY')+1/2\right) - \bbE\log(\cM_k(\cY')+1/2)$ has independent sub-Gaussian entries with $\psi_2$ norm bounded by $K_0\frac{e^B}{I}$, then by \cite{vershynin2010introduction}, we have 
\begin{equation*}
	\begin{split}
		\left\|\cM_k(\cZ_1)\right\| \leq K_0\sqrt{\frac{e^B}{I}}\left(\sqrt{p_{-k}}+\sqrt{p_k} + t\right)
	\end{split}
\end{equation*}
hold with probability at least $1-e^{-t^2}$. By setting $t=\sqrt{p_k}$, we have \eqref{eq-Poisson-hpc-2} hold with probability at least $1-e^{-\underline{p}}$. Next, for \eqref{eq-Poisson-hpc-3}, note that
\begin{equation*}
	\begin{split}
		\sup_{\left\|\V_k\right\| \leq 1, k=1,2,3}\left\|\cZ_1 \times_1 \V_1^\top \times_2 \V_2^{\top} \times_3 \V_3^{\top}\right\|_\tF = \sup_{\substack{\left\|\V_k\right\| \leq 1, k=1,2,3\\ \cS \in\bbR^{r_1\times r_2 \times r_3}, \|\cS\|_\tF \leq 1}} \left\langle \cZ_1, \llbracket \cS;\V_1,\V_2,\V_3 \rrbracket\right\rangle.
	\end{split}
\end{equation*}
Since each entry of $\cZ_1$ has independent sub-Gaussian entry with $\psi_2$ norm bounded by $C_0\frac{e^B}{I}$, by Lemma \ref{lm-concentration-Gaussian-xi}, we have:
\begin{equation*}
	\begin{split}
		\bbP\left(\sqrt{\frac{e^B}{I}}\sup_{\left\|\V_k\right\| \leq 1, k=1,2,3}\left\|\cZ_1 \times_1 \V_1^\top \times_2 \V_2^{\top} \times_3 \V_3^{\top}\right\|_\tF > C_1\sqrt{df}\right) \leq 2\exp(-c \cdot df),
	\end{split}
\end{equation*}
which gives \eqref{eq-Poisson-hpc-3}. Finally, by Lemma \ref{lm-xi-Poisson}, we know that \eqref{eq-Poisson-hpc-4} holds with probability at least $1-c/(p_1p_2p_3)$. Thus, applying union bounds on the above probabilistic events, we know that \eqref{eq-Poisson-hpc-1}-\eqref{eq-Poisson-hpc-4} hold with probability at least $1-c/(p_1p_2p_3)$ and the proof is finished. $\quad\quad \square$

\subsection{Proof of Theorem \ref{thm:Binomial}}
The proof is very similar to that of Theorem \ref{thm:Poisson}. Let $g(\cP;\cN)$ be a mapping from $\bbR^{p_1\times p_2\times p_3}$ to $\bbR^{p_1\times p_2\times p_3}$, such that
\begin{equation*}
    [g(\cP;\cN)]_{jkl} = \log\left(\frac{\cP_{jkl} + 1/(2\cN_{jkl})}{1-\cP_{jkl} + 1/(2\cN_{jkl})}\right),
\end{equation*}
then we have $\tilde \cX^{(0)}_{jkl} = g(\hat \cP, \cN)$ with $\hat \cP_{jkl} = \cY_{jkl} / \cN_{jkl}$ for any $j,k,l$. Now we define 
\begin{equation*}
    \cP'_{jkl} = \hat\cP_{jkl}1_{\left\{|\hat\cP_{jkl} - \cP_{jkl} | \leq \cP_{jkl}/2\right\}} + \cP_{jkl}1_{\left\{|\hat\cP_{jkl} - \cP_{jkl} | > \cP_{jkl}/2\right\}},
\end{equation*}
where $\cP_{jkl}:= s(\cX^*_{jkl})$, and further denote $\cZ = \cZ_1 + \cZ_2$ with
\begin{equation*}
    \cZ_1 = g(\cP';\cN) - \bbE g(\cP';\cN),\quad \cZ_2 = \bbE g(\cP'; \cN) - \cX^*.
\end{equation*}
Again, we first impose the following conditions for the deterministic analysis, and then verify they actually hold with high probability.
\begin{enumerate}[label=(\subscript{A}{\arabic*})]
	\item
		\begin{equation}\label{eq-Binomial-hpc-1}
            \cP_{jkl}' = \hat \cP_{jkl}, \quad \forall (j,k,l) \in [p_1]\times[p_2]\times [p_3]
        \end{equation}
		\item
		\begin{equation}\label{eq-Binomial-hpc-2}
            \left\|\cM_k(\cZ_1)\right\| \leq 2K_0\sqrt{\frac{e^{2B}}{N}}\left(\sqrt{p_{-k}}+\sqrt{p_k} \right), \quad k = 1,2,3
        \end{equation}
		\item
		\begin{equation}\label{eq-Binomial-hpc-3}
            \sup_{\left\|\V_k\right\| \leq 1, k=1,2,3} \left\|\cZ_1 \times_1 \V_1^\top \times_2 \V_2^{\top} \times_3 \V_3^{\top}\right\|_\tF \leq C_1\frac{\sqrt{df}}{\sqrt{N/e^{2B}}},
        \end{equation}
        \item
		\begin{equation}\label{eq-Binomial-hpc-4}
            \sup_{\substack{\mathcal{T}\in \mathbb{R}^{p_1\times p_2\times p_3}, \\ \rank(\mathcal{T})\leq (r_1,r_2,r_3), \\ \|\mathcal{T}\|_\tF \leq 1}} \left\langle -\hat \cP + s(\cX^*), \mathcal{T} \right\rangle \leq C_2\frac{\sqrt{df}}{\sqrt{N}}.
        \end{equation}
	\end{enumerate}
By Lemma \ref{lm-binomial-bias-variance}, we know that given $N>C_0Be^{3B}$ for some universal constant $C_0$, 
\begin{equation*}
    \begin{split}
        & \left| \bbE [g(\cP'; \cN)]_{jkl} - \cX_{jkl}^* \right| \leq C\frac{e^{2B}}{N}, \\
		& \left\|g(\cP';\cN) - \bbE g(\cP';\cN)\right\|_{\psi_2} \leq C\sqrt{\frac{e^{2B}}{N}}.
    \end{split}
\end{equation*}
Then one can combine this with conditions \eqref{eq-Binomial-hpc-1}, \eqref{eq-Binomial-hpc-2} and \eqref{eq-Binomial-hpc-3} to obtain the initialization error by the same proof technique we used in Theorem \ref{thm:Poisson}:
\begin{equation*}
    \left\|\cX^{(0)} - \cX^*\right\|_\tF \leq C\sqrt{\frac{e^{2B}}{N}} \sum_{k=1}^3 \left(\sqrt{p_kr_k} + \sqrt{p_{-k}r_k}\right),
\end{equation*}
as long as $N \geq C\overline p e^{3B}$.

Now we start to check the conditions in Theorem \ref{thm:local-convergence-PGD}. Notice that the loss function $L(\cX) = -\sum_{jkl}\left(\cP_{jkl}\cX_{jkl} + \log(1-s(\cX_{jkl}))\right)$ is of $\frac{1}{e^B+3}$-smoothness and $\frac{1}{4}$-strong convexity, then by Lemma \ref{lm-convex-smooth}, we can set
\begin{equation*}
    \alpha = \frac{1}{2(e^B+3)} < \frac{\frac{1}{4(e^B+3)}}{\frac{1}{e^B+3}+1/4},\quad \beta = \frac{1}{2} \leq \frac{1}{\frac{1}{e^B+3}+1/4}
\end{equation*}
and $L(\cX)$ satisfies RCG$(\alpha,\beta,\cC)$ with $\cC$ defined according to \eqref{eq:constraint-C} and \eqref{eq-specify-sets}.
In addition, we can also evaluate the scale of $\xi$:
\begin{equation*}
    \xi := \sup_{\substack{\mathcal{T}\in \mathbb{R}^{p_1\times p_2\times p_3}, \\ \rank(\mathcal{T})\leq (r_1,r_2,r_3), \\ \|\mathcal{T}\|_\tF \leq 1}} \left\langle -\hat \cP + s(\cX^*), \mathcal{T} \right\rangle \overset{\eqref{eq-Binomial-hpc-4}}{\leq} C_2\frac{\sqrt{df}}{\sqrt{N}}.
\end{equation*}
Then as we have $N\underline{\lambda}^2 \geq C\kappa^2 e^{3B}\sum_{k=1}^3 \left(p_kr_k + p_{-k}r_k\right)$, it follows that
\begin{equation*}
    \left\|\cX^{(0)} - \cX^*\right\|_\tF^2 \leq C\frac{e^{2B}}{N}\sum_{k=1}^3\left(p_kr_k + p_{-k}r_k\right) \leq \frac{c\alpha\beta \underline{\lambda}^2}{\kappa^2}.
\end{equation*}
Thus the initialization condition is meet. In the meantime, one can also check that the signal-noise-ratio condition $\underline\lambda^2 \geq C\frac{\kappa^4}{\alpha^3\beta}\xi^2$ is also satisfied when $\kappa$ and $B$ are treated as constants. Thus by applying Theorem \ref{thm:local-convergence-PGD}, we obtain 
\begin{equation*}
    \left\|\cX^{T} - \cX^*\right\|_\tF^2 \leq \frac{\kappa^2}{\alpha^2}\xi^2 = C\kappa^2e^{2B}\frac{df}{N}.
\end{equation*}
Now we are going to show that conditions \eqref{eq-Binomial-hpc-1}-\eqref{eq-Binomial-hpc-4} hold with high probability. First of all,
\begin{equation*}
    \begin{split}
         \bbP\left(\cP = \cP'\right) & = 1 - \bbP\left(\exists(j,k,l), \cP_{jkl} \neq \cP_{jkl}'\right)  \\
         & \geq 1 - \sum_{j,k,l}\left(\bbP\left(\hat \cP_{jkl} \leq \frac{\cP_{jkl}}{2}\right) + \bbP\left(\hat\cP_{jkl}\geq \frac{3}{4}\cP_{jkl}\right)\right) \\
         & \geq 1 - 2\sum_{j,k,l}\exp\left(-\frac{1}{2}\cN_{jkl}\cP_{jkl}^2\right) \\
         & \geq 1 - 2p_1p_2p_3 \exp\left(-cNe^{-2B}\right) \\
         & \geq  1 - \frac{1}{p_1p_2p_3},
    \end{split}
\end{equation*}
where the last but one inequality comes from the assumption that $\|\cX^*\|_\infty < B$ and $\cP = s(\cX^*)$, and the last inequality comes from the condition that $N \geq Ce^{2B}\log(\overline p)$, so \eqref{eq-Binomial-hpc-1} holds with probability at least $1-\frac{1}{p_1p_2p_3}$. We can also prove that \eqref{eq-Binomial-hpc-2}, \eqref{eq-Binomial-hpc-3}, and \eqref{eq-Binomial-hpc-4} hold with probability at least $1-1/p_1p_2p_3$ as we did for the proof of \eqref{eq-Poisson-hpc-2}, \eqref{eq-Poisson-hpc-3} and \eqref{eq-Poisson-hpc-4} in Theorem \ref{thm:Poisson}, and we omitted them here. 
 \qquad $\square$

\subsection{Proof of Theorem \ref{thm:lower-bound-poisson}}
Since we assume $\max_k r_k \leq \min_k \sqrt{p_k}$, it suffices to show the following inequalities:
\begin{equation}\label{eq-lb-loadings}
	\inf_{\hat{\cX}} \sup_{\cX \in \mathcal{F}_{\p, \r}} \bbE\left\|\hat{\cX} - \cX \right\|_\tF^2 \geq cI^{-1}p_kr_k , \quad \forall k=1,2,3
\end{equation}
We only need to prove \eqref{eq-lb-loadings} for $k=1$. First of all, we let $\V \in \{-1,1\}^{\lfloor \frac{p_k}{2} \rfloor \times r_k}$ such that each entry are i.i.d. Rademacher random variables (i.e., $\bbP(\V_{ij} = 1) = \bbP(\V_{ij} = -1) = 1/2,~\forall i,j$). Then by \cite{vershynin2010introduction},
	\begin{equation*}
	    \bbP\left(\sqrt{\lfloor \frac{p_1}{2}}\rfloor - \sqrt{r_1} - t \leq \sigma_{r_1}(\V) \leq \left\|\V\right\| \leq \sqrt{\lfloor \frac{p_1}{2}}\rfloor + \sqrt{r_1} + t \right) > 1 - e^{-\frac{t^2}{2}}.
	\end{equation*}
	Taking $t = c_0\sqrt{r_1}$ for some sufficiently small constant $c_0$, since the above probability is positive, we know that there exists a $\V_0 \in \{-1,1\}^{\lfloor \frac{p_1}{2} \rfloor \times r_1}$, such that 
	\begin{equation*}
	    \left\|\V_0\right\| \leq \sqrt{\lfloor \frac{p_1}{2}}\rfloor + (1+c_0)\sqrt{r_1} \leq \sqrt{p_1r_1},
	\end{equation*}
	and 
	\begin{equation*}   
	    \begin{split}
	        \sigma_{r_1}(\V_0) & \geq \sqrt{\lfloor \frac{p_1}{2}\rfloor} - (1+c_0)\sqrt{r_1} \geq \sqrt{\frac{p_1-1}{2}} - (1+c_0)C_2^{-1/2}r_1  \\
	        & \geq \sqrt{\frac{p_1-1}{2}} - (1+c_0)C_1C_2^{-1/2}p_1^{1/2} \\
	        & \geq \sqrt{\frac{p_1-1}{2}} - c_0p_1^{1/2} \geq \sqrt{\frac{p_1}{3}}.
	    \end{split}
	\end{equation*}
	Here we obtain the last but one inequality by taking $C_2 > C_1^2/c_0^2$.
	
	In addition, for $i=1,2,\ldots, N$, let $\U^{(i)} \in \{-1, 1\}^{\lceil \frac{p_1}{2}\rceil \times r_1}$ be the i.i.d. copy of $\V_0$, and we denote
	\begin{equation*}
		\Omega^{(ij)} = \left\{1\leq a \leq \lceil \frac{p_1}{2} \rceil, 1\leq b \leq r_1: \U^{(i)}_{ab} \neq \U^{(j)}_{ab}\right\}.
	\end{equation*}	
	Since $\left|\Omega^{(ij)}\right| \sim \text{Binomial}\left(\lceil\frac{p_1}{2}\rceil r_1,\frac{1}{2}\right)$, we have
	\begin{equation}\label{ineq-lowerbound-1}
		\begin{split}
			& \bbP\left(\frac{\lceil\frac{p_1}{2}\rceil r_1}{4} \leq \left|\Omega^{(ij)}\right| \leq \frac{3\lceil\frac{p_1}{2}\rceil r_1}{4}, \forall i\neq j \right) \\
			& = 1 - \bbP\left(\exists i\neq j, \left|\Omega^{(ij)}\right| < \frac{\lceil\frac{p_1}{2}\rceil r_1}{4} ~\text{or}~ \left|\Omega^{(ij)}\right| > \frac{3\lceil\frac{p_1}{2}\rceil r_1}{4}, \forall i\neq j \right) \\
			& \geq  1 - N(N-1)/2\left(\bbP\left(\left|\Omega^{(ij)}\right| < \frac{\lceil\frac{p_1}{2}\rceil r_1}{4}\right) + \bbP\left(\left|\Omega^{(ij)}\right| > \frac{3\lceil\frac{p_1}{2}\rceil r_1}{4}\right)\right) \\
			& = 1 - \frac{N(N-1)}{2}\bbP\left(\left|\frac{\left|\Omega^{(ij)}\right|}{\lceil\frac{p_1}{2}\rceil r_1} - \frac{1}{2}\right| > \frac{1}{4}\right)  \\
			& \geq 1 - N(N-1)e^{-\lceil\frac{p_1}{2}\rceil r_1/32} \\
			& \geq 1 - N(N-1)e^{-p_1r_1/64}.
 		\end{split}
	\end{equation}
	In the meantime, by \cite{vershynin2010introduction}, we also have
	\begin{equation*}
	    \begin{split}
	        &\bbP\left(\left\|\U_1^{(i)}\right\| \leq \sqrt{\lceil \frac{p_1}{2} \rceil} + \sqrt{r_1} + t, ~\forall 1\leq i \leq N\right) \\
	        & = 1 - \bbP\left(\exists i\in [N], \left\|\U_1^{(i)}\right\| > \sqrt{\lceil \frac{p_1}{2} \rceil} + \sqrt{r_1} + t\right) \\
	        & \geq 1 - Ne^{-\frac{t^2}{2}}.
	    \end{split}
	\end{equation*}
	Taking $t = \sqrt{p_1r_1}$, we see that
	\begin{equation}\label{ineq-lowerbound-2}
	    \bbP\left(\left\|\U_1^{(i)}\right\| \leq 2\sqrt{p_1r_1}, ~\forall 1\leq i \leq N\right) \geq 1 - Ne^{-p_1r_1/2}.
	\end{equation}
	Combining \eqref{ineq-lowerbound-1} and \eqref{ineq-lowerbound-2}, we know that with probability at least $1 - N(N-1)e^{-p_1r_1/64} - Ne^{-p_1r_1/2}$, 
	\begin{equation}\label{ineq-U-property-1}
	    \begin{split}
	        \lceil\frac{p_1}{2}\rceil r_1 \leq \left\|\U^{(i)} - \U^{(j)}\right\|_\tF^2 \leq 3\lceil\frac{p_1}{2}\rceil r_1,\quad 1\leq i < j \leq N,
	    \end{split}
	\end{equation}
	\begin{equation}\label{ineq-U-property-2}
	    \begin{split}
	        \left\|\U_1^{(i)}\right\| \leq 2\sqrt{p_1r_1}, \quad 1\leq i \leq N.
	    \end{split}
	\end{equation}
	By taking $N = e^{cp_1r_1}$ for some small constant $c$, the probability defined in \eqref{ineq-lowerbound-2} is positive and thus there exists $\left\{\U_1^{(i)},\ldots,\U_1^{(N)}\right\} \subset \{-1, 1\}^{\lceil \frac{p_1}{2}\rceil \times r_1}$ such that \eqref{ineq-U-property-1} and \eqref{ineq-U-property-2} hold.
	
	Now we let $\cS \in \bbR^{r_1 \times r_2 \times r_3}$ be a fixed core tensor such that 
	\begin{equation}\label{ineq-lowerbound-core}
	   \frac{B}{C\sqrt{\mu_1\mu_2\mu_3}} \leq  \min\limits_k \sigma_{r_k}\left(\cM_1(\cS)\right) \leq \max\limits_k \left\|\cM_1(\cS)\right\| \leq \frac{B}{3\sqrt{\mu_1\mu_2\mu_3}}
	\end{equation}
	for some $C > 3$, and let $\U_k \in \bbO_{p_k,r_k}$ be the orthogonal matrix such that $\left\|\U_k\right\|_{2,\infty}^2 \leq \frac{\mu_k r_k}{p_k}$ for $k=2,3$. Let $\cX^{(i)} = \delta \cdot \llbracket \cS; \tilde\U_1^{(i)}, \U_2, \U_3 \rrbracket$ for some $0<\delta<1$ with
	\begin{equation*}
	    \tilde \U_1^{(i)} = \begin{bmatrix}
	        \U_1^{(i)} \\ \V_0
	    \end{bmatrix}.
	\end{equation*}
	Let $\tilde{\U}_1^{(i)} = \hat{\U}_1^{(i)}\R^{(i)}$ be the QR-decomposition of $\tilde\U_1^{(i)}$ where $\hat\U_1^{(i)} \in \bbO_{p_k,r_k}$ is an orthogonal matrix. Thus we can rewrite $\cX^{(i)} = \llbracket \delta\cS\times_1 \R^{(i)}; \hat \U_1^{(i)},\U_2,\U_3\rrbracket$. Note that by construction, $\sigma_{r_1}\left(\tilde\U_1^{(i)}\right) \geq \sigma_{r_1}\left(\V_0\right) \geq \sqrt{\frac{p_1}{3}}$, and we have
	\begin{equation*}
	    \begin{split}
	        \left\|\hat\U_1^{(i)}\right\|_{2,\infty}^2 &= \max_i \left\|e_i^\top \hat\U_1^{(i)}\right\|_2^2 = \max_i \left\|e_i^\top \tilde\U_1^{(i)}\left(\R^{(i)}\right)^{-1}\right\|_2^2 \\
	        & \leq \max_i \frac{\left\|e_i^\top \tilde\U_1^{(i)}\right\|_2^2}{\sigma_{r_1}^2\left(\R^{(i)}\right)}  \leq \frac{\mu_kr_k}{p_k}.
	    \end{split}
	\end{equation*}
	Here the last inequality comes from the facts that $\left\|e_i^\top \tilde \U_1^{(i)}\right\|_2^2 = r_k$ for each $i$ and $\sigma_{r_1}\left(\R^{(i)}\right) = \sigma_{r_1}\left(\tilde \U_1^{(i)}\right)$.
	
	In the meantime, we also have
	\begin{equation*}
	    \begin{split}
	        & \max_k \left\|\cM_k(\delta \cS \times_1 \R^{(i)})\right\| \leq  \max_k \left\|\cM_k(\cS )\right\| \cdot \left\|\R^{(i)}\right\| \\
	        & \qquad \overset{\eqref{ineq-lowerbound-core}}{\leq} \frac{B}{3\sqrt{\mu_1\mu_2\mu_3}} \left\|\tilde\U_1^{(i)}\right\| \leq \frac{B}{3\sqrt{\mu_1\mu_2\mu_3}}\left(\left\|\U_1^{(i)}\right\| + \left\|\V_0\right\|\right) \\
	        & \qquad \overset{\eqref{ineq-U-property-2}}{\leq}  \frac{B}{\sqrt{\mu_1\mu_2\mu_3}}\sqrt{p_1r_1} \leq B\sqrt{\frac{\Pi_{k=1}^3p_k}{\Pi_{k=1}^3\mu_kr_k}},
	    \end{split}
	\end{equation*}
	where the last inequality comes from the assumption that $\max_k r_k \leq \min_k \sqrt{p_k}$. Thus we have $\cX^{(i)} \in \cF_{\p,\r}$ for each $i \in [N]$.
	
	Now we provide the KL-divergence between $\cY^{(i_1)} \sim \Poisson(I\exp(\cX^{(i_1)}))$ and $\cY^{(i_2)} \sim \Poisson(I\exp(\cX^{(i_2)}))$ for $i_1 \neq i_2$. For any two Poisson distribution $P \sim \Poisson(\lambda_1)$ and $Q \sim \Poisson(\lambda_2)$, the KL-divergence between $P$ and $Q$ is
	\begin{equation*}
	    D_{KL}\left(P||Q\right) = \lambda_1\log\left(\frac{\lambda_1}{\lambda_2}\right) + \lambda_2 - \lambda_1.
	\end{equation*}
	Thus we have
	\begin{equation}\label{ineq-KL-bound}
	    \begin{split}
	        & D_{KL}\left(\cY^{(i_1)}||\cY^{(i_2)}\right) = \sum_{j,k,l} D_{KL}\left(\cY^{(i_1)}_{jkl}||\cY^{(i_2)}_{jkl}\right) \\
	        = & I \sum_{j,k,l} \left(\exp(\cX_{jkl}^{(i_1)})\left(\cX^{(i_1)}_{jkl} - \cX^{(i_2)}_{jkl}\right) + \exp(\cX^{(i_2)}_{jkl}) - \exp(\cX^{(i_1)}_{jkl})\right) \\
	        \overset{(a)}{=} & I \sum_{j,k,l} \frac{1}{2} \exp\left(\xi^{(i_1,i_2)}_{jkl}\right)\left(\cX^{(i_1)}_{jkl} - \cX^{(i_2)}_{jkl}\right)^2 \\
	        \overset{(b)}{\leq} & \frac{1}{2}Ie^B\left\|\cX^{(i)} -\cX^{(j)}\right\|_\tF^2  = \frac{1}{2}Ie^B\delta^2 \left\|\llbracket \cS; \tilde \U_1^{(i_1)} - \tilde \U_2^{(i_2)}, \U_2, \U_3\rrbracket \right\|_\tF^2 \\
	        \leq & \frac{1}{2}Ie^B \delta^2\left\|\cM_1(\cS)\right\|^2 \left\|\tilde \U_1^{(i_1)} - \tilde \U_1^{(i_2)}\right\|_\tF^2 \\
	         \overset{\eqref{ineq-lowerbound-core}}{\leq}& I\frac{\delta^2B^2e^B}{18\mu_1\mu_2\mu_3} \left\|\tilde{\U}_1^{(i_1)} - \tilde \U_1^{(i_2)}\right\|_{\tF}^2  \overset{\eqref{ineq-U-property-1}}{\leq} I\frac{\delta^2 B^2e^B}{12\mu_1\mu_2\mu_3} p_1r_1.
	    \end{split}
	\end{equation}
	Here (a) is obtained by applying second order Taylor's expansion on exponential function, and $\xi_{jkl}^{(i_1,i_2)}$ is some real number between $\cX_{jkl}^{(i_1)}$ and $\cX_{jkl}^{(i_2)}$; (b) comes from the fact that $\cX_{jkl}^{(i_1)} \vee \cX_{jkl}^{(i_2)} \leq B$ by construction.
	
	Besides, $\forall i_1 \neq i_2$, we also have
	\begin{equation}\label{ineq-separate-bound}
	    \begin{split}
	        \left\|\cX^{(i_1)} - \cX^{(i_2)}\right\|_\tF^2 &= \delta^2\left\|\llbracket \cS; \tilde \U_1^{(i_1)} - \tilde \U_2^{(i_2)}, \U_2, \U_3\rrbracket \right\|_\tF^2 \\
		    & \geq \delta^2 \sigma_{r_1}\left(\cM_1(\cS)\right)^2  \left\|\tilde \U_1^{(i_1)} - \tilde \U_2^{(i_2)}\right\|_\tF^2 \\
		    & \overset{\eqref{ineq-lowerbound-core}}{\geq} \frac{\delta^2 B^2}{C^2\mu_1\mu_2\mu_3}\left\|\tilde \U_1^{(i_1)} - \tilde \U_2^{(i_2)}\right\|_\tF^2 \\
		    &\overset{\eqref{ineq-U-property-1}}{\geq} \frac{\delta^2 B^2}{4C^2\mu_1\mu_2\mu_3}p_1r_1.
	    \end{split}
	\end{equation}
	Then by generalized Fano Lemma \citep{yang1999information},
	\begin{equation*}
		\inf_{\hat \cX} \sup_{\cX \in \{\cX^{(i)}\}_{i=1}^N} \bbE\left\|\hat \cX - \cX\right\|_\tF^2 \geq c\frac{\delta^2 B^2}{\mu_1\mu_2\mu_3}p_1r_1 \left(1 - \frac{IB^2e^{B}\delta^2 p_1r_1/(12\mu_1\mu_2\mu_3) + \log 2}{cp_1r_1 - 1}\right).
	\end{equation*}
	By setting $\delta^2 = \frac{c\mu_1\mu_2\mu_3}{100IB^2e^B}$, we obtain \eqref{ineq-lowerbound-1}.\qquad $\square$

\subsection{Proof of Theorem \ref{thm:lower-bound-binomial}}
The proof of Theorem \ref{thm:lower-bound-binomial} is similar to that of Theorem \ref{thm:lower-bound-poisson}. We construct the same $\left\{\cX^{(1)}, \ldots \cX^{(N)}\right\}_{i=1}^N$ as we did in the proof of Theorem \ref{thm:lower-bound-poisson}. 
Notice that for any two binomial distributions $P = \text{Binomial}(n,p)$, $Q = \text{Binomial}(n,q)$, the KL-divergence between $P$ and $Q$ is
\begin{equation*}
    D_{KL}\left(P||Q\right) = \log\left(\frac{p}{q}\right)np + \log\left(\frac{1-p}{1-q}\right)n(1-p).
\end{equation*}
Now let $P_1 = \text{Binomial}(n, s(x_1))$ and $P_2 = \text{Binomial}(n, s(x_2))$, for some $-B \leq x_1 \leq x_2 \leq B$ where $s(x) = 1/(1+e^{-x})$, then the KL-divergence of $P_1$ and $P_2$ is
\begin{equation*}
    \begin{split}
        &D_{KL}(P_1||P_2) = n\left( \log\left(\frac{1+e^{-x_2}}{1+e^{-x_1}}\right)\frac{1}{1+e^{-x_1}} + \log\left(\frac{1+e^{-x_2}}{1+e^{-x_1}}e^{x_2-x_1}\right)\frac{e^{-x_1}}{1+e^{-x_1}}\right) \\
        & = n\left(\log\left(\frac{1+e^{-x_2}}{1+e^{-x_1}}\right) + \frac{e^{-x_1}}{1+e^{-x_1}}(x_2-x_1)\right) \\
        & = n\left(\log\left(1+e^{-x_2}\right) - \log\left(1+e^{-x_1}\right) + \frac{e^{-x_1}}{1+e^{-x_1}}(x_2-x_1)\right) \\
        & \overset{(a)}{=} \frac{ne^{-\xi}}{(1+e^{-\xi})^2} = \frac{n}{2 + e^{\xi} + e^{-\xi}} \leq \frac{n}{4}.
        \end{split}
\end{equation*}
Here $(a)$ comes from the second order Taylor's expansion of $f(x) = \log\left(1+e^{-x}\right)$.
Then we obtain
\begin{equation*}
    D_{KL}\left(\cY^{(i_1)}||\cY^{(i_2)}\right) \leq c\delta^2\max_{jkl}\cN_{jkl}p_1r_1 \leq c\delta^2 \min_{jkl}\cN_{jkl}p_1r_1,
\end{equation*}
given $\cY_{jkl}^{(i)} \sim \text{Binomial}(\cN_{jkl},s(\cX_{jkl}))$ independently. By the similar argument in the proof of Theorem \ref{thm:lower-bound-poisson}, we can show 
\begin{equation*}
    \inf_{\hat{\cX}} \sup_{\cX \in \mathcal{F}_{\p, \r}} \left\|\hat{\cX} - \cX \right\|_\tF^2 \geq c N^{-1} \left(r_1r_2r_3+\sum_{k=1}^3p_kr_k\right),
\end{equation*}
and the proof is finished. \qquad $\square$

\subsection{Proof of Proposition \ref{prop:rank-selection}}
Assume $\sigma=1$ without loss of generality. We only prove for $k=1$, while the proof for other modes follows in the same way. We prove the two scenarios separately.
\begin{enumerate}[leftmargin=*]
    \item[(a)] \emph{Sub-Gaussian tensor PCA.} Let $\U_1 \in \bbO_{p_1,r_1}$ be the left singular subspace of $\cM_1(\cX^*)$ and $\U_{1\perp} \in \bbO_{p_1,p_1-r_1}$ be the orthogonal subspace of $\U_1$. Let $ \cZ:=\tilde\cX - \cX^*$ and define the event:
    \begin{equation*}
        A = \left\{\frac{9}{10}\sqrt{p_{-1}} \leq \sigma_{p_1}\left(\cM_1(\cZ)\right) \leq \|\cM_1(\cZ)\| \leq \frac{11}{10}\sqrt{p_{-1}}\right\}.
    \end{equation*}
    Since $\cM_1(\cZ)$ has independent mean-zero unit-variance sub-Gaussian random variables, by the concentration of singular values of random matrix \citep[Corollary 5.35]{vershynin2010introduction}, we have 
    \begin{equation*}
        \bbP\left(\sqrt{p_{-1}}-\sqrt{p_{1}}-t\leq \sigma_{p_1}\left(\cM_1(\cZ)\right) \leq \|\cM_1(\cZ)\| \leq \sqrt{p_{-1}} + \sqrt{p_1} + t\right) \geq 1-Ce^{-ct^2}.
    \end{equation*}
    Since $p_1 \leq c p_{-1}$, taking $t = c\sqrt{p_1}$ in the above inequality yields that $\bbP(A) \geq 1-2e^{-cp_{-1}}$. 
    Now we provide a lower bound for $\sigma_{p_1}\left(\cM_1(\tilde\cX)\right)$ assuming $A$ holds. By definition,
    \begin{equation*}
        \begin{split}
            \sigma_{p_1}^2\left(\cM_1(\tilde \cX)\right) &= \inf_{\substack{u \in \bbR^{p_1} \\ \|u\|=1}} \left\|u^\top \cM_1(\tilde \cX)\right\|_2^2 \\
            & = \inf_{\substack{u \in \bbR^{p_1} \\ \|u\|=1}} \left(\left\|u^\top  \U_1\U_1^\top \cM_1(\tilde \cX)\right\|_2^2 + \left\|u^\top  \U_{1\perp}\U_{1\perp}^\top \cM_1(\tilde \cX)\right\|_2^2\right).
        \end{split}
    \end{equation*}
    Since
    \begin{equation*}
        \begin{split}
            \left\|u^\top\U_1\U_1^\top \cM_1(\tilde \cX)\right\|_2^2 & \geq \|u^\top \U_1\|_2^2 \cdot \sigma_{r_1}^2\left(\U_1^\top\cM_1(\tilde \cX)\right) \\
            & \geq \|u^\top \U_1\|_2^2 \cdot \left(\sigma_{r_1}\left(\U_1^\top\cM_1(\cX^*)\right)-\|\cM_1(\cZ)\|\right)^2 \\
            & \geq \|u^\top\U_1\|_2^2 \cdot (\underline \lambda - 1.1\sqrt{p_{-1}})^2 \geq \frac{1}{2}\underline \lambda^2 \|u^\top\U_1\|_2^2
        \end{split}
    \end{equation*}
    and
    \begin{equation*}
        \begin{split}
            \left\|u^\top  \U_{1\perp}\U_{1\perp}^\top \cM_1(\tilde \cX)\right\|_2^2 & = \left\|u^\top  \U_{1\perp}\U_{1\perp}^\top \cM_1(\cZ)\right\|_2^2 \\
            &\geq \|u^\top\U_{1\perp}\U_{1\perp}^\top\|\cdot \sigma_{p_1}^2\left(\cM_1(\cZ)\right) \geq \frac{4p_{-1}}{5}\|u^\top\U_{1\perp}\|_2^2,
        \end{split}
    \end{equation*}
    we further have
    \begin{equation}\label{ineq:rk-select-1}
        \begin{split}
            \sigma_{p_1}^2\left(\cM_1(\tilde \cX)\right) \geq \inf_{\substack{u \in \bbR^{p_1} \\ \|u\|=1}} \left(\frac{\underline \lambda^2}{2}\|u^\top \U_1^\top\| + \frac{4p_{-1}}{5}\|u^\top\U_{1\perp}\|_2^2\right) \geq \frac{4p_{-1}}{5}.
        \end{split}
    \end{equation}
    Recall $\delta_1$ is the median of non-zero singular values of $\cM_1(\tilde \cX)$, we then have
    \begin{equation*}
        \begin{split}
            \delta_1 \geq \sigma_{p_1}(\cM_1(\tilde{ \cX})) \geq \frac{9\sqrt{p_{-1}}}{10}.
        \end{split}
    \end{equation*}
    On the other hand, since $r_1 = o(p_1)$, 
    \begin{equation*}
        \delta_1 \leq \sigma_{r_1+1}\left(\cM_k(\tilde\cX)\right) \leq \|\cM_k(\cZ)\| \leq \frac{11}{10}\sqrt{p_{-1}}.
    \end{equation*}
    Thus,
    \begin{equation}\label{ineq:rk-select-2}
        \begin{split}
        \sigma_{r_1+1}(\cM_1(\tilde \cX)) &\leq 1.5\delta_1, \\
        \sigma_{r_1}(\cM_1(\tilde{\cX})) & \geq \underline \lambda - \frac{11\sqrt{p_{-1}}}{10} \\
        & \geq 2\sqrt{p_{-1}} \geq 1.5\delta_1.\\    
        \end{split}
    \end{equation}
    \eqref{ineq:rk-select-2} implies $\hat r_1 = r_1$ and the proof is finished.
    
    \item[(b)] \emph{Tensor regression.} Following the proof of Theorem \ref{thm:regression-sharp}, we assume each entry of $\cA_i$ comes from i.i.d. $N(0,1/n)$ and $\varepsilon \sim N(0,\frac{1}{n}\I_n)$. Then we have $\tilde{\cX} = \sum_{i=1}^n y_i\cA_i = \cA^*(y)$. Define $\cZ := \tilde\cX - \cA^*\cA(\cX^*) = \cA^*(\varepsilon)$ and define the linear operator $\cA_1: \bbR^{p_1\times p_{-1}} \rightarrow \bbR^n$:
    \begin{equation*}
        \left(\cA_1(\X)\right)_i = \left\langle \cM_1(\cA_i), \X \right\rangle.
    \end{equation*}
    Let $\cA_1^*$ be the adjoint operator of $\cA_1$. We consider the following events:
    \begin{equation*}
        \begin{split}
            A_1 & := \left\{\frac{9}{10}\sqrt{\frac{p_{-1}}{n}} \leq \sigma_{p_1}(\cM_1(\cZ))\leq \|\cM_1(\cZ)\| \leq \frac{11}{10}\sqrt{\frac{p_{-1}}{n}} \right\},\\
            A_2 & := \left\{(1-\frac{1}{2\kappa\sqrt{r_1}}) \leq \frac{\|\cA_1(\X)\|_2^2}{\|\X\|_\tF^2} \leq (1+\frac{1}{2\kappa\sqrt{r_1}})~~\text{for any rank-$2r_1$ matrix $\X \in \bbR^{p_1\times p_{-1}}$}\right\}, \\
            A_3 & := \left\{ \sigma_{r_1+1}\left(\cA_1^*\cA_1\left(\cM_1(\cX^*)\right)\right)  \leq \frac{11}{10}\sqrt{\frac{p_{-1}}{n}}\|\cX^*\|_\tF \right\},\\
            A_4 & := \left\{ \sigma_{p_1-r_1}\left(\cA_1^*\cA_1\left(\cM_1(\cX^*)\right)\right)  \geq \frac{9}{10}\sqrt{\frac{p_{-1}}{n}}\|\cX^*\|_\tF \right\}.
         \end{split}
    \end{equation*}
    We first show that each of the above events hold with probability at least $1-Ce^{-cp_{-1}}$. 
    \begin{itemize}
        \item $A_1$: Note that conditional on $\varepsilon$, $\cM_1(\cZ) = \cA_1^*(\varepsilon)$ has i.i.d. mean-zero normal entries with variance $\|\varepsilon\|_2^2$. Therefore, by random matrix theory, we have
    \begin{equation}\label{ineq:rank-regression-0}
        \begin{split}
            \bbP\left(\sigma_{p_1}(\cA_1^*(\varepsilon)) \geq \frac{\|\varepsilon\|_2}{\sqrt{n}}\left(\sqrt{p_{-1}}-\sqrt{p_1}-t\right)\Big| \varepsilon\right) \geq 1-2e^{-t^2/2}, \\
            \bbP\left(\left\|\cA_1^*(\varepsilon))\right\| \leq \frac{\|\varepsilon\|_2}{\sqrt{n}}\left(\sqrt{p_{-1}}+\sqrt{p_1}+t\right)\Big| \varepsilon\right) \geq 1-2e^{-t^2/2}, \\
        \end{split}
    \end{equation}
    In addition, by the concentration of Chi-square random variable \cite[Lemma 8.1]{birge2001alternative}, we have
    \begin{equation*}
        \begin{split}
            \bbP\left( \|\varepsilon\|_2^2 \geq 1+2\sqrt{x/n}+2x/n\right) \leq e^{-x}, \\
            \bbP\left( \|\varepsilon\|_2^2 \leq 1-2\sqrt{x/n}\right) \leq e^{-x}.
        \end{split}
    \end{equation*}
    Taking $x=cn$, we have with probability at least $1-2\exp(-cp_{-1})$ that $\|\varepsilon\|_2 \geq 9/10$. Now we specify $t = c\sqrt{p_{-1}}$ in \eqref{ineq:rank-regression-0} and obtain
    \begin{equation}\label{ineq:rank-regression-0-1}
        \begin{split}
            & \bbP \left(\sigma_{p_1}(\cA_1^*(\varepsilon))  \geq \frac{4}{5}\sqrt{\frac{p_{-1}}{n}}\right) \\
            \geq & \bbP\left(\|\varepsilon\|\geq \frac{4}{5}\right) \cdot \bbP\left(\sigma_{p_1}(\cA_1^*(\varepsilon)) \geq \frac{4}{5}\sqrt{\frac{p_{-1}}{n}}\Big| \|\varepsilon\| \geq 1/2\right) \\
            \geq & (1-e^{-cp_{-1}}) \bbP\left(\sigma_{p_1}(\cA_1^*(\varepsilon)) \geq \|\varepsilon\|\cdot\frac{8}{9}\sqrt{\frac{p_{-1}}{n}}\Big| \|\varepsilon\| \geq 1/2\right) \\
            \geq & (1-e^{-cp_{-1}})(1-2e^{-cp_{-1}}) \geq 1-Ce^{-cp_{-1}}.
        \end{split}
    \end{equation}
    Similarly one can prove $\bbP\left(\left\|\cA_1^*(\varepsilon)\right\| \leq \frac{6}{5}\sqrt{\frac{p_{-1}}{n}}\right) \geq 1-Ce^{-cp_{-1}}$. Thus, we have proved $\bbP(A_1) \geq 1-Ce^{-cp_{-1}}$. 
        \item $A_2$: By \cite[Lemma 4.3]{recht2010guaranteed}, when $n \geq C\kappa p_{-1}r_{1}^{3/2}\log r_1$, one has $\bbP(A_2) \geq 1-\exp(-cp_{-1}r_1)$.
        \item $A_3,A_4$: Since $\cM_1(\cA_i)$ has i.i.d. Gaussian entries and $\rank(\cM_1(\cX^*))=r_1$, we can assume that the last $p_1-r_1$ rows of $\cM_1(\cX^*)$ are zeros without loss of generality. Therefore, we can write
        \begin{equation*}
            \cM_1(\cX^*) = \begin{bmatrix}
                \X_1 \\ \mathbf O
            \end{bmatrix},\qquad \cM_1(\cA_i) = \begin{bmatrix}
                \A_{i1} \\ \A_{i2}
            \end{bmatrix}
        \end{equation*}
        for $\X_1,\A_{i1} \in \bbR^{r_1 \times p_{-1}}$, $\A_{i2} \in \bbR^{p_1-r_1 \times p_{-1}}$. Then we have
        \begin{equation}\label{ineq:rank-regression-1}
        \begin{split}
            & \sigma_{r_1+1}\left(\cA_1^*\cA_1\left(\cM_1(\cX^*)\right)\right)  = \sigma_{r_1+ 1}\left(\sum_{i=1}^n \langle \A_{i1}, \X_1 \rangle \begin{bmatrix}
            \A_{i1} \\ \A_{i2}
            \end{bmatrix}\right)\\
            & \qquad = \inf_{\U \in \bbO_{p_1,p_1-r_1}} \sup_{u \in \U} \left\|u^\top \sum_{i=1}^n \langle \A_{i1}, \X_1 \rangle \begin{bmatrix}
            \A_{i1} \\ \A_{i2}
            \end{bmatrix}\right\|_2 \leq   \left\|\sum_{i=1}^n \langle \A_{i1}, \X_1 \rangle \A_{i2} \right\|
        \end{split}
        \end{equation}
        and
        \begin{equation}\label{ineq:rank-regression-2}
        \begin{split}
            & \sigma_{p_1-r_1}\left(\cA_1^*\cA_1\left(\cM_1(\cX^*)\right) \right)  = \sup_{\U \in \bbO_{p_1,r_1}}\inf_{u \in \U} \left\|u^\top \sum_{i=1}^n \langle \A_{i1}, \X_1 \rangle \begin{bmatrix}
            \A_{i1} \\ \A_{i2}
            \end{bmatrix} \right\|_2\\
            & \qquad \geq \sup_{\U \in \bbO_{p_1-r_1,r_1}}\inf_{u \in \U} \left\|u^\top \sum_{i=1}^n \langle \A_{i1}, \X_1 \rangle \A_{i2} \right\|_2 \\
             & \qquad = \sigma_{p_1-r_1}\left(\sum_{i=1}^n \langle \A_{i1}, \X_1 \rangle \A_{i2}\right).
        \end{split}
        \end{equation}
        Recall that $\A_{i1}$ and $\A_{i2}$ are independent. Conditional on $\A_{i1}$, $\sum_{i=1}^n \langle \A_{i1}, \X_1 \rangle \A_{i2}$ is a $(p_1-r_1)$-by-$p_{-1}$ random matrix with i.i.d. $N\left(0,\tilde\sigma^2\right)$ entries where $\tilde\sigma^2 := \frac{1}{n}\sum_{i=1}^n\langle \A_{i_1},\X_1\rangle^2$. Therefore, by random matrix theory,
        \small
        \begin{equation*}
        \begin{split}
            & \bbP\left(\left\|\sum_{i=1}^n \langle \A_{i1}, \X_1 \rangle \A_{i2}\right\| \leq \tilde\sigma(\sqrt{p_{-1}}-\sqrt{p_1-r_1}+t)\Bigg | \{\A_{i1}\}_{i=1}^n\right) \leq 1-2e^{-t^2/2}, \\
            &\bbP\left(\sigma_{p_1-r_1}\left(\sum_{i=1}^n \langle \A_{i1}, \X_1 \rangle \A_{i2}\right) \geq \tilde\sigma(\sqrt{p_{-1}}-\sqrt{p_1-r_1}-t)\Bigg | \{\A_{i1}\}_{i=1}^n\right) \leq 1-2e^{-t^2/2}.
        \end{split}
        \end{equation*}
        \normalsize
        On the other hand, since $n^2\tilde\sigma^2/\|\X_1\|_\tF^2 \sim \chi_n^2$, we have
        \begin{equation*}
            \begin{split}
            \bbP\left(\tilde\sigma^2 \leq \frac{\|\X_1\|_\tF^2}{n}\left(1+2\sqrt{x/n}+2x/n\right)\right) &\geq 1-e^{-x},\\
            \bbP\left(\tilde\sigma^2 \geq \frac{\|\X_1\|_\tF^2}{n}\left(1-2\sqrt{x/n}\right)\right) &\geq 1-e^{-x}.
            \end{split}
        \end{equation*}
        Combining the concentrations inequalities above and following the same argument as \eqref{ineq:rank-regression-0-1}, we obtain
        \begin{equation}\label{ineq:rank-regression-3}
        \begin{split}
            \bbP\left(\sigma_{p_1-r_1}\left(\sum_{i=1}^n \langle \A_{i1}, \X_1 \rangle \A_{i2}\right) \geq \frac{9}{10}\|\X_1\|_\tF\sqrt{\frac{p_{-1}}{n}}\right) & \leq 1-Ce^{-cp_{-1}}, \\
            \bbP\left(\left\|\sum_{i=1}^n \langle \A_{i1}, \X_1 \rangle \A_{i2}\right\| \leq \frac{11}{10}\|\X_1\|_\tF\sqrt{\frac{p_{-1}}{n}}\right) & \leq 1-Ce^{-cp_{-1}}. \\
        \end{split}
        \end{equation}
        \eqref{ineq:rank-regression-1},\eqref{ineq:rank-regression-2} and \eqref{ineq:rank-regression-3} together imply that $\bbP(A_3) \geq 1-Ce^{-cp_{-1}}$, $\bbP(A_4) \geq 1-Ce^{-cp_{-1}}$.
        
    \end{itemize}
    By union bound, we know that $\bbP\left(A_1\cap A_2\cap A_3 \cap A_4\right) \geq 1-Ce^{-cp_{-1}}$. Now we show that $\hat r_1 = r_1$ under the event $A_1 \cap A_2 \cap A_3 \cap A_4$. 
    On the one hand,
    \begin{equation*}
        \begin{split}
            \sigma_{r_1+1}\left(\cM_1(\tilde \cX)\right) & = \sigma_{r_1+1}\left(\cA_1^*\cA_1(\cM_1(\cX^*))+\cA_1^*(\varepsilon))\right) \\
            & \leq \sigma_{r_1+1}\left(\cA_1^*\cA_1(\cM_1(\cX^*))\right) + \left\|\cA_1^*(\varepsilon)\right\| \\
            & \overset{A_1,A_3}{\leq} \frac{11}{10}\sqrt{\frac{p_{-1}}{n}}(\|\cX^*\|_\tF + 1) \leq \frac{6}{5}\sqrt{\frac{p_{-1}}{n}}\|\cX^*\|_\tF.
        \end{split}
    \end{equation*}
    Here in the last inequality we use the assumption that $1=\sigma \leq c\|\cX^*\|_\tF$. Meanwhile, 
    \begin{equation*}
        \begin{split}
            \sigma_{p_1-r_1}\left(\cM_1(\tilde\cX)\right) & \geq \sigma_{p_1-r_1}\left(\cA_1^*\cA_1\left(\cM_1(\cX^*)\right)\right) - \left\|\cA_1^*(\varepsilon)\right\| \\
            & \overset{A_1,A_3}{\geq} \frac{9}{10}\sqrt{\frac{p_{-1}}{n}}\left(\|\cX^*\|_\tF -1 \right) \geq \frac{4}{5}\sqrt{\frac{p_{-1}}{n}}\|\cX^*\|_\tF.
        \end{split}
    \end{equation*}
    Consequently, the median $\delta_1$ satisfies 
    \begin{equation}\label{ineq:rank-regression-4}
        \begin{split}
            1.5\delta_1 \geq 1.5\sigma_{p_1-r_1}\left(\cM_1(\tilde\cX)\right) \geq \frac{6}{5}\sqrt{\frac{p_{-1}}{n}}\|\cX^*\|_\tF \geq \sigma_{r_1+1}\left(\cM_1(\tilde \cX)\right).
        \end{split}
    \end{equation}
    On the other hand, 
    \begin{equation*}
        \begin{split}
            & \sigma_{r_1}\left(\cA_1^*\cA_1(\cM_1(\cX^*))\right) = \inf_{\U \in \bbO_{p_1,r_1-1}} \left\|(\I_{p_1}- \U\U^\top)\cA_1^*\cA_1(\cM_1(\cX^*))\right\| \\
            & \qquad = \inf_{\U \in \bbO_{p_1,r_1-1}} \sup_{\substack{u \in \bbR^{p_1}, \|u\|=1 \\ v \in \bbR^{ p_{-1}},\|v\|=1}}\left\langle (\I_{p_1}- \U\U^\top)\cA_1^*\cA_1(\cM_1(\cX^*)), uv^\top \right\rangle \\
            &  \qquad = \inf_{\U \in \bbO_{p_1,r_1-1}} \sup_{\substack{u \in \bbR^{p_1}, \|u\|=1 \\ v \in \bbR^{ p_{-1}},\|v\|=1}}\left\langle \cA_1(\cM_1(\cX^*)), \cA_1((\I_{p_1}- \U\U^\top)uv^\top) \right\rangle \\
            & \qquad \overset{\text{Lemma \ref{lm-incoherent-inner-product-difference}}}{\geq} \inf_{\U \in \bbO_{p_1,r_1-1}} \sup_{\substack{u \in \bbR^{p_1}, \|u\|=1 \\ v \in \bbR^{ p_{-1}},\|v\|=1}}\left\langle \cM_1(\cX^*), (\I_{p_1}- \U\U^\top)uv^\top \right\rangle - \frac{1}{2\kappa\sqrt{r_1}}\|\cM_1(\cX^*)\|_\tF \\
            & \qquad =\inf_{\U \in \bbO_{p_1,r_1-1}} \left\|(\I_{p_1}- \U\U^\top)\cM_1(\cX^*)\right\| - \frac{1}{2\kappa\sqrt{r_1}}\|\cM_1(\cX^*)\|_\tF \\
            & \qquad = \sigma_{r_1}\|\cM_1(\cX^*)\| - \frac{1}{2\kappa\sqrt{r_1}}\|\cM_1(\cX^*)\|_\tF \geq \lambda_1 - \frac{1}{2\kappa\sqrt{r_1}}\kappa\sqrt{\overline{r}} \lambda_1 \geq \frac{\lambda_1}{2}.
        \end{split}
    \end{equation*}
    Then, 
    \begin{equation}\label{ineq:rank-regression-5}
        \begin{split}
            \sigma_{r_1}\left(\cM_1(\tilde{\cX})\right) & \geq \sigma_{r_1}\left(\cA_1^*\cA_1(\cM_1(\cX^*))\right) - \left\|\cA_1^*(\varepsilon)\right\| \\
            & \geq \frac{\lambda_1}{2} - \sqrt{\frac{p_{-1}}{n}}\left\|\cX^*\right\|_\tF \\
            & \geq \left(\frac{1}{2\kappa\sqrt{r_1}}-\sqrt{\frac{p_{-1}}{n}}\right)\left\|\cX^*\right\|_\tF \\
            & \geq 2\sqrt{\frac{p_{-1}}{n}}\left\|\cX^*\right\|_\tF \geq 1.5\delta_1.
        \end{split}
    \end{equation}
    Here the last but two inequality comes from the assumption on $n$. Now $\hat r_1 = r_1$ is implied by \eqref{ineq:rank-regression-4} and \eqref{ineq:rank-regression-5} and the proof is finished.  \qquad\qquad $\square$
\end{enumerate}

\section{Technical Lemmas}\label{sec:lemmas}
We collect the technical lemmas with their proofs in this sections. They are widely used in the proof of the major theorems.

\subsection{Proof of Lemma \ref{lm:partial-gradient}}\label{sec:proof-lemma1}
Let $\cX := \llbracket \cS; \U_1,\U_2,\U_3 \rrbracket$, we first calculate the partial gradient of $L(\cS,\U_1,\U_2,\U_3)$ with respect to $\U_1$ by chain rule:
\begin{equation*}
    \begin{split}
        \frac{\partial L}{\partial \U_{1,ij}} & = \sum_{a=1}^{p_1} \sum_{b=1}^{p_2} \sum_{c=1}^{p_3} \frac{\partial L}{\partial \cX_{abc}} \cdot \frac{\partial \cX_{abc}}{\partial \U_{1,ij}} \\
        & = \sum_{a=1}^{p_1} \sum_{b=1}^{p_2} \sum_{c=1}^{p_3} \frac{\partial L}{\partial \cX_{abc}} \cdot \left(I_{\{a=i\}} \sum_{l_2=1}^{r_2}\sum_{l_3=1}^{r_3} \cS_{jl_2l_3}\U_{2,bl_2}\U_{3,cl_3} \right)\\
        & = \sum_{b=1}^{p_2} \sum_{c=1}^{p_3} \frac{\partial L}{\partial \cX_{ibc}} \left(\sum_{l_2=1}^{r_2}\sum_{l_3=1}^{r_3} \cS_{jl_2l_3}\U_{2,bl_2}\U_{3,cl_3}\right),
    \end{split}
\end{equation*}
where the second identity comes from the following fact:
$$\cX_{abc} = \sum_{l_1=1}^{r_1}\sum_{l_2=1}^{r_2}\sum_{l_3=1}^{r_3} \cS_{l_1l_2l_3}\U_{1,al_1}\U_{2,bl_2}\U_{3,cl_3}$$ In the meantime, one can verify that
\begin{equation*}
    \left(\cM_1(\nabla L(\cX))(\U_3\otimes \U_2)\cM_1(\cS)^\top\right)_{ij} = \sum_{k_1=1}^{p_2}\sum_{k_2=1}^{p_3}\sum_{k_3=1}^{r_2} \sum_{k_4=1}^{r_3} \frac{\partial L}{\partial \cX_{ik_1k_2}}\cdot \U_{2, k_1k_3}\U_{3, k_2k_4} \cS_{jk_3k_4},
\end{equation*}
which is exactly what we calculated for $\frac{\partial L}{\partial \U_{1,ij}}$ (by changing the order of summation). The partial gradient for $\U_2$ and $\U_3$ can be similarly calculated. For core tensor $\cS$, we have
\begin{equation*}
    \begin{split}
        \frac{\partial L}{\partial \cS_{ijk}} & = \sum_{a=1}^{p_1} \sum_{b=1}^{p_2} \sum_{c=1}^{p_3} \frac{\partial L}{\partial \cX_{abc}} \cdot \frac{\partial \cX_{abc}}{\partial \cS_{ijk}} \\
        & = \sum_{a=1}^{p_1} \sum_{b=1}^{p_2} \sum_{c=1}^{p_3} \frac{\partial L}{\partial \cX_{abc}} \cdot \U_{1,ai} \U_{2,bj} \U_{3,ck} \\
        & = \left(\llbracket \nabla L(\cX); \U_1^\top, \U_2^\top, \U_3^\top \rrbracket\right)_{ijk},
    \end{split}
\end{equation*}
which has finished the proof of this lemma. \quad\quad $\square$

\subsection{Proof of Proposition \ref{prop:loss-equivalent}}\label{sec:proof-prop1}
Let
\begin{equation*}
    \left( \hat \cS, \hat \U_1, \hat \U_2, \hat \U_3 \right) = \argmin\limits_{\cS,\U_1,\U_2,\U_3} L(\llbracket \cS; \U_1, \U_2, \U_3 \rrbracket;D) + \frac{a}{2}\sum_{k=1}^3\left\|\U_k^\top\U_k - b^2\I_{r_k}\right\|_\tF^2.
\end{equation*}
We claim $\hat \U_k^\top\hat \U_k = b^2\I_{r_k}$ for $k=1,2,3$. Otherwise, consider the QR-decomposition $\hat \U_k = \tilde\U_k \R_k$ where $\tilde \U_k \in \bbO_{p_k,r_k}$ and define $\bar \cS = b^{-3}\llbracket \hat\cS; \R_1, \R_2,\R_3 \rrbracket$, $\bar \U_k = b\tilde\U_k$. Then clearly $F(\hat \cS,\hat\U_1,\hat\U_2,\hat \U_3) > F(\bar \cS,\bar\U_1,\bar\U_2,\bar \U_3)$ and makes contradiction.

Similarly, For arbitrary rank-$(r_1,r_2,r_3)$ tensor $\cX = \llbracket \cS;\U_1,\U_2,\U_3\rrbracket$, we can redecompose it as $\cX = \llbracket \cS^*;\U_1^*,\U_2^*,\U_3^*\rrbracket$ such that $\U_k^{*\top}\U_k^* = b^2\I_{r_k}$. Then it follows that
\begin{equation*}
    L(\hat\cX;D) = F(\bar \cS,\bar \U_1,\bar \U_2,\bar \U_3) \leq F(\cS^*,\U_1^*,\U_2^*,\U_3^*) = L(\cX;D).
\end{equation*}
Thus, $\hat \cX = \llbracket \hat\cS; \hat\U_1, \hat\U_2, \hat\U_3 \rrbracket$ is a rank-constraint minimizer of $L(\cX;D)$. $~\square$

\subsection{Lemmas for Main Theory} 
We collect the technical lemmas that are used to develop Theorem \ref{thm:local-convergence-PGD} in this section.

The first lemma builds the relationship between the standard strongly convex and smooth assumption with RCG condition \ref{asmp-restricted-convexity}.
\begin{Lemma}\label{lm-convex-smooth}
Consider a continuously differentiable function $f$, given a convex domain $\mathbb{G}$, suppose for any $x, y\in \mathbb{G}$, we have
\begin{equation*}
\frac{m}{2} \left\|x-y\right\|^2 \leq f(y) - f(x) - \left\langle \nabla f(x), y - x \right\rangle \leq \frac{M}{2} \left\|x - y\right\|^2,
\end{equation*}
then
\begin{equation*}
	\left\langle \nabla f(x) - \nabla f(y), x - y \right\rangle \geq \frac{mM}{m+M}\left\|x - y\right\|^2 + \frac{1}{m+M}\left\|\nabla f(x)- \nabla f(y)\right\|^2.
\end{equation*}
\end{Lemma}
\textbf{Proof of Lemma \ref{lm-convex-smooth}.} See \cite{nesterov1998introductory}. \qquad $\square$

The next lemma establishes the relationship between $\left\|\cX^{(t)} - \cX^*\right\|_\tF^2$ and $E^{(t)}$ in the proof of Theorem \ref{thm:local-convergence-PGD}.
\begin{Lemma}\label{lm-equivalent-criteria}
Suppose $\cX^* = \llbracket\cS^*; \U_1^*, \U_2^*, \U_3^*\rrbracket$, $\U_k^{*\top}\U_k^* = b^2\I_{r_k}, k=1,2,3$, $\overline{\lambda} = \max_k \left\|\cM_k(\cX^*)\right\|$, and $\underline{\lambda} = \min_k \sigma_{r_k}(\cM_k(\cX^*))$. Let $\cX = \llbracket\cS; \U_1, \U_2, \U_3\rrbracket$ be another Tucker low-rank tensor with $\U_k \in \bbR^{p_k \times r_k}$, $\left\|\U_k\right\| \leq (1+c_d)b$, and $\max_k \left\|\cM_k(\cS)\right\| \leq (1+c_d)\frac{\overline \lambda}{b^3}$ for some constant $c_d>0$. Define 
\begin{equation*}
	\begin{split}
		E & :=  \min_{\R_k \in \bbO_{p_k, r_k}, k = 1,2,3} \left(\sum_{k=1}^3\left \| \U_k - \U_k^* \R_k\right\|_\tF^2 + \left\| \cS - \llbracket\cS^*; \R_1^\top, \R_2^\top, \R_3^\top\right\|_\tF^2\right).
	\end{split}
\end{equation*}
Then we have
\begin{equation*}
	\begin{split}
		& E \leq b^{-6}\left(7 + \frac{12b^8}{\underline{\lambda}^2}C_d\right)\left\|\cX - \cX^*\right\|_\tF^2 + 2b^{-2}C_d \sum_{k=1}^3\left\|\U_k^\top \U_k - b^2 \I\right\|_\tF^2,\\
		& \left\|\cX - \cX^*\right\|_\tF^2  \leq 4b^6\left( 1 + (3+2c_d)^2(1+c_d)^4 \overline \lambda^2 b^{-8} \right) E,
	\end{split}
\end{equation*}
where $C_d := 1 + 7\overline{\lambda}^2b^{-8}\left((1+c_d)^6 + (1+c_d)^4(2+c_d)^2\right)$.
\end{Lemma}
\textbf{Proof of Lemma \ref{lm-equivalent-criteria}.} First, note that 
\begin{equation*}
    \left\|\llbracket\cS; \R_1, \R_2, \R_3 \rrbracket - \cS^*\right\|_\tF = b^{-3}\left\| \llbracket\cS; \U_1^*\R_1, \U_2^*\R_2, \U_3^*\R_3 \rrbracket  - \cX^*\right\|_\tF.   
\end{equation*}
Then we can decompose $\llbracket\cS; \U_1^*\R_1, \U_2^*\R_2, \U_3^*\R_3 \rrbracket - \cX^*$ into seven terms,
\begin{equation*}
	\begin{split}
		&  \llbracket\cS;\U_1 + \U_1^*\R_1 - \U_1, \U_2 + \U_2^*\R_2 - \U_2,\U_3 + \U_3^*\R_3 - \U_3\rrbracket - \cX^*\\
		= & (\cX - \cX^*) + \sum_{k=1}^3\cS \times_k (\U_k^*\R_k - \U_k) \times_{k+1} \U_{k+1} \times_{k+2} \U_{k+2} \\
		&+ \sum_{k=1}^2 \cS \times_k (\U_k^*\R_k - \U_k) \times_{k+1} (\U_{k+1}^*\R_{k+1}-\U_{k+1}) \times_{k+2} \U_{k+2} \\
		&+ \cS \times_1 (\U_1^*\R_1 - \U_1) \times_{2} \U_{2}^*\R_{2} \times_{3} (\U_3^*\R_3 - \U_{3}).
	\end{split}
\end{equation*}
By the Inequality of arithmetic and geometric means, we further have 
\begin{equation*}
	\begin{split}
		& \left\|\cS - \llbracket\cS^*; \R_1^\top ,\R_2^\top,\R_3^\top\rrbracket\right\|_\tF^2 = \left\|\llbracket\cS;\R_1,\R_2,\R_3\rrbracket - \cS^*\right\|_\tF^2\\
		\leq & 7b^{-6}\|\cX - \cX^*\|_\tF^2 + 7 b^{-6} \sum_{k=1}^3\left(\left\|\cM_k(\cS)\right\|^2\left\|\U_k - \U_k^*\R_k\right\|_\tF^2 \left\|\U_{k+2}\otimes \U_{k+1}\right\|^2\right) \\
		&+ 7b^{-6} \sum_{k=1}^2\left(\left\|\cM_{k}(\cS)\right\|^2\left\|\U_k - \U_k^*\R_k\right\|_\tF^2 \left\|\U_{k+1} - \U_{k+1}^*\R_{k+1}\right\|^2\left\| \U_{k+2}\right\|^2\right)\\
		& + 7 b^{-6}\left\|\cM_2(\cS)\right\|^2\left\|\U_1^*\R_1 - \U_1\right\|_\tF^2 \left\|\U_3^*\R_3 - \U_3\right\|^2 \left\|\U_2^* \R_2\right\|^2 \\
		\leq & 7b^{-6}\|\cX - \cX^*\|_\tF^2 + 7 b^{-6}\left((1+c_d)^6\overline \lambda^2 b^{-2} + (1+c_d)^4(2+c_d)^2\overline \lambda^2 b^{-2} \right) \sum_{k=1}^3 \left\|\U_k - \U_k^* \R_k\right\|_\tF^2. \\
		= & 7b^{-6}\|\cX - \cX^*\|_\tF^2 + 7\overline{\lambda}^2b^{-8}\left((1+c_d)^6 + (1+c_d)^4(2+c_d)^2\right) \sum_{k=1}^3 \left\|\U_k - \U_k^* \R_k\right\|_\tF^2.
	\end{split}
\end{equation*}
Since the above inequality holds for any orthogonal matrices $\R_1,\R_2,\R_3$, it follows that
\begin{equation}\label{ineq:E-X-1}
	\begin{split}
		&E = \min_{\R_k \in \bbO_{p_k\times r_k}, k = 1,2,3} \left(\sum_{k=1}^3\left \|\U_k - \U_k^* \R_k\right\|_\tF^2 + \left\| \cS - \cS^* \times_1 \R_1^\top \times_2 \R_2^\top \times_3 \R_3^\top\right\|_\tF^2\right) \\
		& \leq 7b^{-6}\left\|\cX - \cX^*\right\|_\tF^2 + C_d \sum_{k=1}^3\min_{\R_k \in \bbO_{r_k}}\left\|\U_k - \U_k^* \R_k\right\|_\tF^2,
	\end{split}
\end{equation}
where $C_d = 1 + 7\overline{\lambda}^2b^{-8}\left((1+c_d)^6 + (1+c_d)^4(2+c_d)^2\right)$.\\
Now let $\U_k = \tilde \U \tilde{\mathbf{\Sigma}} \tilde{\V}^\top$ be the SVD of $\U_k$, where $\tilde \U \in \bbO_{p_k,r_k}$, $\tilde{\V} \in \bbO_{r_k}$, and $\tilde{\mathbf{\Sigma}}=\diag(\sigma_1,\ldots,\sigma_{r_k})$ is a diagonal matrix. Then we have
\begin{equation}\label{ineq:E-X-2}
	\begin{split}
		& \min_{\R_k \in \bbO_{r_k}}\left\|\U_k - \U_k^*\R_k\right\|_\tF^2 = \min_{\R_k \in \bbO_{r_k}}\left\|\U_k - b \tilde{\U} \tilde{\V}^\top + b\tilde{\U} \tilde{\V}^\top - \U_k^*\R_k\right\|_\tF^2 \\
		\leq & 2\left\|\tilde \U \tilde{\mathbf{\Sigma}} \tilde{\V}^\top - b \tilde{\U} \tilde{\V}^\top\right\|_\tF^2 + 2\min_{\R_k \in \bbO_{r_k}}\left\| b\tilde{\U} \tilde{\V}^\top - \U_k^*\R_k\right\|_\tF^2. \\
		= & 2\left\|\tilde{\mathbf{\Sigma}} - b\I\right\|_\tF^2 + 2\min_{\R_k \in \bbO_{r_k}}\left\|b\tilde{\U} -  \U_k^*\R_k\right\|_\tF^2.
	\end{split}
\end{equation}
On the one hand,
\begin{equation}\label{ineq:E-X-3}
	\begin{split}
		& \left\|\tilde{\mathbf{\Sigma}} - b\I\right\|_\tF^2 = b^2 \left\|\tilde{\mathbf{\Sigma}}/b - \I\right\|_\tF^2 \\
		= & b^2 \sum_{i=1}^{r_k} \left(\frac{\sigma_i}{b} - 1\right)^2 \leq b^2 \sum_{i=1}^{r_k} \left( \frac{\sigma_i^2}{b^2} - 1\right)^2 = b^{-2}\sum_{i=1}^{r_k}\left(\sigma_i^2 -b^2\right)^2 \\
		= & b^{-2} \left\|\tilde{\mathbf{\Sigma}}^2 - b^2\I\right\|_\tF^2 = b^{-2}\left\|\tilde{\V}\tilde{\mathbf{\Sigma}}^2\tilde{\V}^\top  - b^2\I\right\|_\tF^2 \\
		= & b^{-2}\left\|\U_k^\top \U_k - \U_k^{*\top}\U_k^*\right\|_\tF^2 = b^{-2}\left\|\U_k^\top \U_k - b^2\I\right\|_\tF^2.
	\end{split}
\end{equation}
Here we use the inequality $(x-1)^2 \leq (x^2-1)^2$ for any non-negative number $x$.
On the other hand, since $\tilde{\U}_k$ and $\U_k^*/b$ have orthonormal columns and span the left singular subspaces of $\cM_k(\cX)$ and $\cM_k(\cX^*)$ respectively, we have
\begin{equation*}
    \begin{split}
        \left\|\cX - \cX^*\right\|_\tF^2 & = \left\|\cM_k(\cX) - \cM_k(\cX^*)\right\|_\tF^2 \geq \left\|\tilde\U_{k\perp}^{\top} \left(\cM_k(\cX) - \cM_k(\cX^*)\right)\right\|_\tF^2 \\
        & = \left\|\tilde \U_{k\perp}^{\top} \cM_k(\cX^*)\right\|_\tF^2 = \left\|\tilde \U_{k\perp}^{\top} (\U_k^*/b)(\U_k^*/b)^\top\cM_k(\cX^*)\right\|_\tF^2 \\
        & \geq \sigma_{r_k}^2(\cM_k(\cX^*))\cdot \left\|\tilde{\U}_{k\perp}^{\top} \U_k^*/b\right\|_\tF^2 \geq \underline{\lambda}^2 \left\|\tilde{\U}_{k\perp}^{\top} \U_k^*/b\right\|_\tF^2,
    \end{split}
\end{equation*}
where $\tilde{\U}_{k\perp} \in \bbO_{p_k-r_k,r_k}$ is the perpendicular orthogonal matrix of $\tilde{\U}_{k}$. Then it follows by \cite[Lemma 1]{cai2018rate} that
\begin{equation}\label{ineq:E-X-4}
	\begin{split}
		\min_{\R_k \in \bbO_{r_k}}\left\|b\tilde{\U}_k - \U_k^*\R_k\right\|_\tF^2 \leq 2b^2 \left\|\tilde{\U}_{k\perp}^{\top} (\U_k^*/b)\right\|_\tF^2 \leq 2b^2 \frac{\|\cX - \cX^*\|_\tF^2}{\underline{\lambda}^2}.
	\end{split}
\end{equation}
Combing \eqref{ineq:E-X-2},\eqref{ineq:E-X-3} and \eqref{ineq:E-X-4}, we have:
\begin{equation*}
	\min_{\R_k \in \bbO_{r_k}} \left\|\U_k - \U_k^*\R_k\right\|_\tF^2 \leq 2b^{-2}\left\|\U_k^\top\U_k - b^2\I\right\|_\tF^2 + 4b^2\frac{\left\|\cX - \cX^*\right\|_\tF^2}{\underline{\lambda}^2}.
\end{equation*}
Then by \eqref{ineq:E-X-1}, we finally obtain
\begin{equation*}
	\begin{split}
		E &\leq  7b^{-6}\left\|\cX - \cX^*\right\|_\tF^2 + C_d \sum_{k=1}^3\min_{\R_k \in \bbO_{r_k}}\left\|\U_k - \U_k^* \R_k\right\|_\tF^2\\
		& \leq b^{-6}\left(7 + \frac{12b^8}{\underline{\lambda}^2}C_d\right)\left\|\cX - \cX^*\right\|_\tF^2 + 2b^{-2}C_d \sum_{k=1}^3\left\|\U_k^\top \U_k - b^2 \I\right\|_\tF^2.
	\end{split}
\end{equation*}

Next, we prove the second inequality. To this end, we now denote
\begin{equation*}
    (\R_1,\R_2,\R_3) = \argmin_{\substack{\R_k \in \bbO_{p_k, r_k}\\ k = 1,2,3}} \left\{\sum_{k=1}^3\left \| \U_k - \U_k^* \R_k\right\|_\tF^2 + \left\|\cS - \llbracket\cS^*; \R_1^\top, \R_2^\top, \R_3^\top\rrbracket\right\|_\tF^2\right\}.
\end{equation*}
Let $\cH_\cS = \cS^* - \llbracket\cS;\R_1,\R_2,\R_3\rrbracket$ and $\H_k = \U_k^* - \U_k\R_k^\top$. Then, we have
\begin{equation*}
	\cX^* = (\cH_\cS + \llbracket\cS;\R_1,\R_2, \R_3\rrbracket) \times_1 (\H_1 + \U_1\R_1^\top) \times_2 (\H_2 + \U_2\R_2^\top) \times_3 (\H_3 + \U_3\R_3^\top)
\end{equation*}
and it follows that
\begin{equation*}
	\begin{split}
		&\left\|\cX - \cX^*\right\|_\tF \leq \left\|\cH_\cS \times_1 \U_1^* \times_2 \U_2^* \times_3 \U_3^*\right\|_\tF + \sum_{k=1}^3 \left\|\cS \times_k \H_k \R_k \times_{k+1} \U_{k+1} \times_{k+2} \U_{k+2}\right\|_\tF \\
		& + \sum_{k=1}^2 \left\|\cS \times_k \H_k \R_k \times_{k+1} \H_{k+1}\R_{k+1} \times_{k+2} \U_{k+2}\right\|_\tF + \left\|\cS \times_1 \H_1\R_1 \times_{2} \U^*_{2}\R_{2} \times_{3} \H_3\R_{3}\right\|_\tF \\
		\leq & b^3 \left\|\cH_\cS\right\|_\tF + \sum_{k=1}^3 \left\|\cM_k(\cS)\right\| \left\|\U_{k+1}\otimes \U_{k+2}\right\| \left\|\H_k\right\|_\tF + \sum_{k=1}^2 \left\|\cM_k(\cS)\right\| \left\|\H_k\right\|_\tF \left\|\H_{k+1}\right\| \left\|\U_{k+2}\right\| \\
		& + \left\|\cM_3(\cS)\right\| \left\|\H_1\right\| \left\|\U_2^*\right\| \left\|\H_3\right\|_\tF \\
		\leq & b^3\left\|\cH_\cS\right\|_\tF + (1+c_d)^3\overline \lambda b^{-1} \sum_{k=1}^3 \left\|\H_k\right\|_\tF + (1+c_d)^2(2+c_d)\overline\lambda b^{-1}\sum_{k=1}^3 \left\|\H_k\right\|_\tF  \\
		= & b^3 \left\|\cH_\cS\right\|_\tF + (3+2c_d)(1+c_d)^2\overline \lambda b^{-1} \sum_{k=1}^3 \left\|\H_k\right\|_\tF.
	\end{split}
\end{equation*}
Thus,
\begin{equation*}
	\left\|\cX - \cX^*\right\|_\tF^2 \leq 4\left(b^6  \left\|\cH_\cS\right\|_\tF^2  + (3+2c_d)^2(1+c_d)^4 \overline \lambda^2 b^{-2} \sum_{k=1}^3\left\|\H_k\right\|_\tF^2\right).\quad\quad \square
\end{equation*}

Suppose $s, u_1, u_2, u_3$ are real values and $\varepsilon_s, \varepsilon_{u_1}, \varepsilon_{u_2}, \varepsilon_{u_3}$ are small perturbations, then by simple calculation we can see that
\begin{equation}\label{eq:expansion}
    (s+\varepsilon_s)(u_1+\varepsilon_{u_2})(u_1+\varepsilon_{u_2})(u_3+\varepsilon_{u_3}) = su_1u_2u_3 + \varepsilon_s u_1u_2u_3 + \sum_{k=1}^3 s\varepsilon_{u_k} u_{k+1}u_{k+2} + h_\varepsilon,
\end{equation}
where $h_\varepsilon$ is a high-order perturbation term. The following lemma can be seen as a tensor version of \eqref{eq:expansion}, which plays a key role in the proof of Theorem \ref{thm:local-convergence-PGD}.
\begin{Lemma}\label{lm-X-decomposition}
Suppose $\cX^* = \llbracket\cS^*; \U_1^*, \U_2^*, \U_3^*\rrbracket$, $\cX = \llbracket\cS; \U_1, \U_2, \U_3\rrbracket$ with $\cS,\cS^* \in \bbR^{r_1\times r_2 \times r_3}$, $\U_k, \U_k^* \in \bbR^{p_k \times r_k}, \R_k \in \bbO_{r_k}, k=1,2,3$. Let
\begin{equation*}
    \begin{split}
        & \cX_k = \cS \times_k \U_k^*\R_k \times_{k+1} \U_{k+1} \times_{k+2} \U_{k+2}, \quad \H_k = \U_k^* - \U_k\R_k^\top, \quad k = 1,2,3\\
        & \cX_\cS = \llbracket \cS^*; \U_1\R_1^{\top}, \U_2\R_2^{\top}, \U_3\R_3^{\top}\rrbracket, \quad \cH_\cS = \cS^* - \llbracket\cS; \R_1, \R_2, \R_3\rrbracket.
    \end{split}
\end{equation*}
Then we have the following decomposition of $\cX^*$:
\begin{equation*}
	\begin{split}
		& \cX^* = \cX_S + \sum_{k=1}^3 (\cX_k - \cX) + \cH_\varepsilon,	\quad \left\|\cH_\varepsilon\right\|_\tF \leq B_2B_3^{3/2}  + 3B_1B_2B_3 + 3B_1^2B_3.
	\end{split}
\end{equation*}
Here,
\begin{equation*}
	\begin{split}
  		&B_1 := \max_k \left\{\left\|\U_k\right\|,\left\|\U_k^*\right\|\right\},\quad B_2 := \max_k \left\{\left\|\cM_k(\cS)\right\|,\left\|\cM_k(\cS^*)\right\|\right\},\\
  		&B_3 := \max \left\{\left\|\cH_\cS\right\|_\tF^2, \left\|\H_1\right\|_\tF^2, \left\|\H_2\right\|_\tF^2, \left\|\H_3\right\|_\tF^2\right\}.
	\end{split}
\end{equation*}
\end{Lemma}
\textbf{Proof of Lemma \ref{lm-X-decomposition}.}
Since $\cS^* = \llbracket\cS; \R_1, \R_2, \R_3\rrbracket + \cH_\cS$, 
\begin{equation}\label{eq:X^*-decompose}
	\begin{split}
		& \cX^* =\llbracket\cS;\U_1^*\R_1,\U_2^* \R_2,\U_3^*\R_3\rrbracket + \llbracket\cH_\cS ;\U_1^*,\U_2^*,\U_3^*\rrbracket. \\
	\end{split}
\end{equation}
For the first term on the right hand side of \eqref{eq:X^*-decompose}, since $\U_k^* =  \U_k\R_k^\top + \H_k$, we have
\begin{equation*}
	\begin{split}
		& \llbracket\cS;\U_1^*\R_1,\U_2^* \R_2, \U_3^*\R_3\rrbracket = \llbracket\cS;\U_1 + \H_1\R_1,\U_2 + \H_2\R_2,\U_3 + \H_3\R_3\rrbracket \\
		 = & \cX + \sum_{k=1}^3 \cS \times_k \U_k \times_{k+1} \U_{k+1} \times_{k+2}\H_{k+2}\R_{k+2} + \cH_\varepsilon^{(1)}\\
		 = & \cX + \sum_{k=1}^3 \cS \times_k \U_k \times_{k+1} \U_{k+1} \times_{k+2}(\U_{k+2}^*\R_{k+2} - \U_{k+2}) + \cH_\varepsilon^{(1)} \\
		 = & \sum_{k=1}^3 \cX_k - 2\cX + \cH_\varepsilon^{(1)},
	\end{split}
\end{equation*}
where
\begin{equation*}
	\begin{split}
		\cH_\varepsilon^{(1)} &=  \llbracket\cS; \H_1\R_1 , \H_2\R_2, \H_3\R_3\rrbracket + \sum_{k=1}^3 \cS \times_k \U_k \times_{k+1} \H_{k+1}\R_{k+1}\times_{k+2} \H_{k+2}\R_{k+2} \\
		& = \llbracket\cS^* - \cH_\cS;\H_1, \H_2, \H_3\rrbracket + \sum_{k=1}^3 (\cS^* - \cH_\cS) \times_k \U_k\R_k^{\top} \times_{k+1} \H_{k+1} \times_{k+2}\H_{k+2}.
	\end{split}
\end{equation*}
For the second term on the right hand side of \eqref{eq:X^*-decompose}, we have
\begin{equation*}
	\begin{split}
		& \llbracket\cH_\cS; \U_1^*, \U_2^*, \U_3^* \rrbracket \\
		= & (\cS^* - \cS \times_1 \R_1 \times_2 \R_2 \times_3 \R_3) \times_1 (\H_1 + \U_1\R_1^\top) \times_2 (\H_2 + \U_2\R_2^\top) \times_3 (\H_3 + \U_3\R_3^\top) \\
		= & \cX_\cS - \cX + \cH_\varepsilon^{(2)},
	\end{split}
\end{equation*}
where
\begin{equation*}
	\begin{split}
		\cH_\varepsilon^{(2)} = & \llbracket\cH_\cS; \H_1,\H_2,\H_3\rrbracket + \sum_{k=1}^3 \cH_\cS \times_k \U_k\R_k^\top \times_{k+1} \H_{k+1}\times_{k+2} \H_{k+2} \\
		& + \sum_{k=1}^3 \cH_\cS \times_k \H_k \times_{k+1} \U_{k+1}\R_{k+1}^\top \times_{k+2} \U_{k+2}\R_{k+2}^\top.
	\end{split}
\end{equation*}
Then it follows that 
\begin{equation*}
	\begin{split}
		\cX^* = \cX_\cS + \sum_{k=1}^3 \cX_k - 3\cX + (\cH_\varepsilon^{(1)} + \cH_\varepsilon^{(2)}).
	\end{split}
\end{equation*}
Since
\begin{equation*}
	\begin{split}
		\cH_\varepsilon = \cH_\varepsilon^{(1)} + \cH_\varepsilon^{(2)} = & \llbracket\cS^*; \H_1 , \H_2 , \H_3\rrbracket + \sum_{k=1}^3 \cS^* \times_k \U_k\R_k^\top \times_{k+1} \H_{k+1} \times_{k+2} \H_{k+2} \\
		& + \sum_{k=1}^3 \cH_\cS \times_k \H_k \times_{k+1} \U_{k+1}\R_{k+1}^\top \times_{k+2} \U_{k+2}\R_{k+2}^\top,
	\end{split}
\end{equation*}
we have
\begin{equation*}
	\begin{split}
		\left\|\cH_{\varepsilon} \right\|_\tF \leq & \left\|\llbracket\cS^* ;\H_1 ,\H_2, \H_3\rrbracket\right\|_\tF + \sum_{k=1}^3 \left\|\cS^* \times_k \U_k\R_k^\top \times_{k+1} \H_{k+1}\times_{k+2} \H_{k+2}\right\|_\tF \\
		& + \sum_{k=1}^3 \left\|\cH_\cS \times_k \H_k \times_{k+1} \U_{k+1}\R_{k+1}^\top \times_{k+2} \U_{k+2}\R_{k+2}^\top\right\|_\tF \\
		= & \left\| \H_1 \cM_1 (\cS^*) (\H_3 \otimes \H_2)^\top\right\|_\tF + \sum_{k=1}^3 \left\|\U_k\R_k^\top \cM_k(\cS^*) (\H_{k+2} \otimes \H_{k+1})^\top\right\|_\tF \\
		& + \sum_{k=1}^3 \left\|\H_k \cM_k(\cH_\cS) (\U_{k+2}\R_{k+2}^\top \otimes \U_{k+1}\R_{k+1}^\top)^\top\right\|_\tF \\
		\leq & \left\|\H_1 \right\|_\tF \left\|\cM_1(\cS^*)\right\| \left\|\H_3 \otimes \H_2\right\|_\tF + \sum_{k=1}^3 \left\|\U_k \R_k^\top\cM_k(\cS^*)\right\| \left\|\H_{k+2}\otimes \H_{k+1}\right\|_\tF \\
		& + \sum_{k=1}^3 \left\|\H_k \cM_k(\cH_\cS)\right\|_\tF \left\|\U_{k+2}\R_{k+2}^\top \otimes \U_{k+1}\R_{k+1}^\top\right\|\\
		\leq & B_3^{3/2}B_2 + 3B_3 B_1B_2 + 3B_3B_1^2.
	\end{split}
\end{equation*}
Therefore, we have finished the proof. $\quad \square$

\subsection{Other Technical Lemmas}
We collect additional technical lemmas in this section.
\begin{Lemma}\label{lm-low-rank-matrix-perturbation}
Suppose $\X,\Z \in \bbR^{m\times n}$ and $\Y = \X + \Z$, $\rank(\cX)=r$. If the leading $r$ left and right singular vectors of $\Y$ are $\hat \U \in \bbO_{m,r}$ and $\hat \V \in \bbO_{n,r}$, then
\begin{equation*}
    \begin{split}
        & \max\left\{\left\|\hat \U_\perp^\top \X\right\|, \left\|\X \hat\V_\perp\right\|\right\} \leq 2\|\Z\|, \quad \max\left\{\left\|\hat \U_\perp^\top \X\right\|_\tF, \left\|\X \hat\V_\perp\right\|_\tF\right\} \leq \min\{2\sqrt{r}\|\Z\|, 2\|\Z\|_\tF\}.
    \end{split}
\end{equation*}
\end{Lemma}
\textbf{Proof of Lemma \ref{lm-low-rank-matrix-perturbation}.} See \cite[Lemma 6]{zhang2017optimal-statsvd}. \qquad $\square$

\begin{Lemma}\label{lm-concentration-Gaussian-xi}
	Assume all the entries of $\cZ\in \mathbb{R}^{p_1\times p_2\times p_3}$ are independent mean-zero random variables with bounded Orlicz-$\psi_2$ norm: 
	\begin{equation}\label{eq:Y-X-sub-gaussian}
	\|\cZ_{ijk}\|_{\psi_2} = \sup_{q\geq 1} \mathbb{E} \left(|\cZ_{ijk}|^q\right)^{1/q}/q^{1/2} \leq \sigma.
	\end{equation}
    Then there exist some universal constants $C, c$, such that
	\begin{equation*}
    \begin{split}
	    & \sup_{\substack{\mathcal{T}\in \mathbb{R}^{p_1\times p_2\times p_3}, \|\mathcal{T}\|_\tF \leq 1, \\ \rank(\mathcal{T})\leq (r_1,r_2,r_3)}} \left\langle\cZ, \mathcal{T} \right\rangle \leq C\sigma \left(r_1r_2r_3 + \sum_{k=1}^3 p_kr_k\right)^{1/2}
	\end{split}
	\end{equation*}
	with probability at least $1 - \exp\left(-c\sum_{k=1}^3 p_kr_k\right)$. 
\end{Lemma}
\textbf{Proof of Lemma \ref{lm-concentration-Gaussian-xi}.} First of all,
By \cite[Lemma 7]{zhang2018tensor}, for any $k=1,2,3$, we can construct $\varepsilon$-net $\left\{\V_k^{(1)}, \ldots, \V_k^{(N_k)}\right\}$ for $\left\{\V_k \in \bbR^{p_k \times r_k}: \left\|\V_k \right\|\leq 1\right\}$ such that
\begin{equation*}
	\sup_{\V_k: \left\|\V_k\right\|\leq 1}\min_{i \leq N_k} \left\|\V_k - \V_k^{(i)}\right\| \leq \varepsilon
\end{equation*}
with $N_k \leq \left(\frac{4+\varepsilon}{\varepsilon}\right)^{p_kr_k}$. Also, we can construct $\varepsilon$-net $\{\cS^{(1)}, \ldots, \cS^{(N_\cS)}\}$ for $\{\cS \in \bbR^{r_1\times r_2 \times r_3}: \left\|\cS\right\|_\tF \leq 1\}$ such that
\begin{equation*}
	\sup_{\cS: \left\|\cS\right\|_\tF \leq 1}\min_{i \leq N_\cS} \left\|\cS - \cS^{(i)}\right\| \leq \varepsilon
\end{equation*} 
with $N_\cS \leq ((2+\varepsilon)/\varepsilon)^{r_1r_2r_3}$. We use $\cI$ to denote the index set $[N_\cS] \times [N_1] \times [N_2] \times [N_3]$. Now suppose we have 
\begin{equation}
	\begin{split}
		\left(\cS^*, \V_1^*, \V_2^*, \V_3^*\right) = \argmax_{\substack{\cS \in \bbR^{r_1\times r_2 \times r_3}, \left\|\cS\right\|_\tF \leq 1 \\ \V_l \in \bbR^{p_l\times r_l}, \|\V_l\|\leq 1}} \left\langle \cZ, \llbracket \cS; \V_1, \V_2, \V_3\rrbracket \right\rangle, \\
	\end{split}
\end{equation}
and denote
\begin{equation*}
    T := \left\langle \cZ, \llbracket \cS^*; \V_1^*, \V_2^*, \V_3^*\rrbracket \right\rangle. 
\end{equation*}
Then we can find some index $i = (i_0, i_1,i_2,i_3) \in \cI$, such that
\begin{equation*}
	\begin{split}
	    \left\|\cS^* - \cS^{(i_0)}\right\|_\tF &\leq \varepsilon,\quad \left\|\V_k^* - \V_k^{(i_k)}\right\| \leq \varepsilon, ~k=1,2,3
	\end{split}
\end{equation*} 
and it follows that for any $j,k,l$, by taking $\varepsilon = 1/14$, we have
\begin{equation*}
	\begin{split}
		& T = \left\langle \cZ, \left\llbracket \cS^{(i_0)}; \V_1^{(i_1)}, \V_2^{(i_2)}, \V_3^{(i_3)}\right\rrbracket \right\rangle + T - \left\langle \cZ, \left\llbracket \cS^{(i_0)}; \V_1^{(i_1)}, \V_2^{(i_2)}, \V_3^{(i_3)}\right\rrbracket \right\rangle  \\
		& = \left\langle \cZ, \left\llbracket \cS^{(i_0)}; \V_1^{(i_1)}, \V_2^{(i_2)}, \V_3^{(i_3)}\right\rrbracket \right\rangle + \left\langle \cZ, \left\llbracket \cS^{*}; \V_1^{*}, \V_2^{*}, \V_3^{*}\right\rrbracket - \left\llbracket \cS^{(i_0)}; \V_1^{(i_1)}, \V_2^{(i_2)}, \V_3^{(i_3)}\right\rrbracket \right\rangle \\
		& \leq \left\langle \cZ, \left\llbracket \cS^{(i_0)}; \V_1^{(i_1)}, \V_2^{(i_2)}, \V_3^{(i_3)}\right\rrbracket \right\rangle + \left(3\varepsilon + 3\varepsilon^2 + \varepsilon^3\right)T  \\
		&  \leq \left\langle \cZ, \left\llbracket \cS^{(i_0)}; \V_1^{(i_1)}, \V_2^{(i_2)}, \V_3^{(i_3)}\right\rrbracket \right\rangle + T/2
	\end{split}
\end{equation*}
and $T \leq 2\left\langle \cZ, \llbracket \cS^{(i_0)}; \V_1^{(i_1)}, \V_2^{(i_2)}, \V_3^{(i_3)}\rrbracket \right\rangle$.\\
Notice that for any fixed $\cS$, $\left\{\V_k\right\}_{k=1}^3$ with $\|\cS\|_\tF \leq 1$ and $\left\|\V_k\right\| \leq 1$, we have $\left\|\llbracket \cS; \V_1,\V_2,\V_3 \rrbracket\right\|_\tF \leq 1$. Since $\cZ_{jkl}$ are independent sub-Gaussian random variable with $\left\|\cZ_{jkl}\right\|_{\psi_2} \leq \sigma$, by Hoeffding-type inequality for sub-Gaussian random variables' summation \citep{vershynin2010introduction}, we have 
\begin{equation*}
	\begin{split}
		\bbP\left(\left|\left\langle \cZ, \llbracket \cS^{(i_0)}; \V_1^{(i_1)}, \V_2^{(i_2)}, \V_3^{(i_3)}\rrbracket \right\rangle\right| \geq t\right) \leq \exp\left(1 - \frac{c t^2 }{\sigma^2}\right).
	\end{split}
\end{equation*}
Then it follows by union bound that
\begin{equation*}
	\begin{split}
		\bbP\left(T \geq t \right) & \leq \bbP\left(\max_{j_0,j_1,j_2,j_3} \left|\left\langle \cZ, \llbracket \cS^{(j_0)}; \V_1^{(j_1)}, \V_2^{(j_2)}, \V_3^{(j_3)}\rrbracket \right\rangle\right| \geq \frac{t}{2}\right) \\
		& \leq N_\cS N_1N_2N_3\exp\left(1-\frac{ct^2}{4\sigma^2}\right) \\
		& \leq \exp\left(1 - \frac{ct^2}{4\sigma^2} + C(\varepsilon)\left(r_1r_2r_3 + \sum_{k=1}^3 p_kr_k\right)\right).
	\end{split}
\end{equation*}
Now by taking $t = C\sigma \sqrt{r_1r_2r_3 + \sum_{k=1}^3 p_kr_k}$, we have with probability at least $1 - \exp\left(-c(r_1r_2r_3 + \sum_{k=1}^3 p_kr_k)\right)$,
\begin{equation*}
	T \leq C\sigma \left(r_1r_2r_3 + \sum_{k=1}^3 p_kr_k\right)^{1/2},
\end{equation*}
which has finished the proof of this lemma.$\quad \square$

\begin{Lemma}\label{lm-random-Gaussian-design}
Suppose $\cA_i$ are Gaussian ensembles with variance $1/n$, i.e., each entry of $\cA_i$ comes from $N(0,1/n)$ independent. Let $\cA$ be the affine map such that $[\cA(\cX)]_i = \langle \cA_i, \cX \rangle$, and let $\cA^*$ be the adjoint operator of $\cA$. Assume $n>C\left(r_1r_2r_3 + \sum_{k=1}^3 p_kr_k\right)$, then we have for any rank-$(r_1,r_2,r_3)$ tensor $\cX \in \bbR^{p_1 \times p_2 \times p_3}$,
\begin{equation}\label{eq-RIP-linear-map}
	\begin{split}
		\frac{9}{10} \left\|\cX\right\|_\tF \leq \left\|\cA(\cX)\right\| \leq \frac{11}{9} \left\|\cX\right\|_\tF, 
	\end{split}
\end{equation}
\begin{equation}\label{eq-RIP-adjoint-map}
	\begin{split}
		\left\|\cA^*\cA(\cX)\right\|_\tF \leq \frac{6}{5}\left(\sqrt{\frac{p_1p_2p_3}{n}} \vee 1\right) \left\|\cX\right\|_\tF 
	\end{split}
\end{equation}
hold with probability at least $1 - C\exp\left( -c\left(r_1r_2r_3 + \sum_{k=1}^3 p_kr_k\right)\right)$.
\end{Lemma}
\textbf{Proof of Lemma \ref{lm-random-Gaussian-design}.} 
By \cite[Theorem 2]{rauhut2017low}, for any $\delta,\varepsilon \in (0,1)$, we have 
\begin{equation*}
    (1-\delta)\|\cX\|_\tF^2 \leq \|\cA(\cX)\|_\tF^2 \leq (1+\delta)\|\cX\|_\tF^2
\end{equation*}
for arbitrary rank-$(r_1,r_2,r_3)$ tensors $\cX$ with probability at least $1-\varepsilon$ provided that
\begin{equation*}
    n \geq C\delta^{-2}\max\left\{(r_1r_2r_3 + \sum_{k=1}^3 p_k r_k), \log(\varepsilon^{-1})\right\}.
\end{equation*}
Setting $\delta = 1/10$ and $\varepsilon = \exp\left(-C(r_1r_2r_3 + \sum_{k=1}^3 p_kr_k)\right)$, \eqref{eq-RIP-linear-map} is established.
To prove \eqref{eq-RIP-adjoint-map}, we first note that $\|\cA^*\|$ is equivalent to the spectral norm of a (scaled) Gaussian random matrix of size $n \times (p_1p_2p_3)$. By random matrix theory \citep[Corollary 5.35]{vershynin2010introduction}, we have with probability at least $1-2\exp(-cp)$,
\begin{equation*}
	\left\|\sqrt{n}\cA^*\right\| \leq \sqrt{p_1p_2p_3} + \sqrt{n} + \frac{1}{11} (\sqrt{p_1p_2p_3} \vee \sqrt{n}) \leq \frac{12}{11}(\sqrt{p_1p_2p_3} \vee \sqrt{n}),
\end{equation*}
then \eqref{eq-RIP-adjoint-map} follows from \eqref{eq-RIP-linear-map} and the fact $\left\|\cA^* \cA (\cX)\right\| \leq \|\cA^*\| \left\|\cA (\cX)\right\|$.\quad  $\square$

\begin{Lemma}\label{lm-incoherent-inner-product-difference}
Let $\cA$ be a linear map from $\bbR^{p_1 \times p_2 \times p_3}$ to $\bbR^{n}$ such that for any rank-$(\r_1 + \r_2)$ tensor $\cX$,
\begin{equation*}
    (1-\delta)\left\|\cX\right\|_\tF^2 \leq \left\|\cA(\cX)\right\|_\tF^2 \leq (1+\delta)\left\|\cX\right\|_\tF^2,
\end{equation*}
Then, for all rank-$\r_1$ tensor $\cX$ and rank-$\r_2$ tensor $\cY$, we have
\begin{equation*}
    \left|\langle \cA(\cX), \cA(\cY)\rangle - \langle\cX , \cY \rangle\right| \leq \delta \left\|\cX\right\|_\tF \left\|\cY\right\|_\tF.
\end{equation*}
\end{Lemma}
\textbf{Proof of Lemma \ref{lm-incoherent-inner-product-difference}.} Without loss of generality we assume $\left\|\cX\right\|_\tF = \left\|\cY\right\|_\tF = 1$, then the proof simply follows from the application of the parallelogram identity:
\begin{equation*}
    \begin{split}
        & \left|\left \langle \cA(\cX), \cA(\cY) \right\rangle - \langle \cX, \cY \rangle\right| \\
        & = \frac{1}{4}\left|\left(\left\|\cA(\cX+\cY)\right\|_\tF^2 - \left\|\cA(\cX-\cY)\right\|_\tF^2\right) - \left(\left\|\cX+\cY\right\|_\tF^2 - \left\|\cX-\cY\right\|_\tF^2\right)\right| \\
        & \leq \frac{1}{4}\left|\left\|\cA(\cX+\cY)\right\|_\tF^2 - \left\|\cX+\cY\right\|_\tF^2\right| + \left|\left\|\cA(\cX-\cY)\right\|_\tF^2 - \left\|\cX-\cY\right\|_\tF^2\right| \\
        & \leq \frac{\delta}{4}\left(\left\|\cX + \cY\right\|_\tF^2 + \left\|\cX-\cY\right\|_\tF^2\right) \leq \delta.	\quad \quad \square
    \end{split}
\end{equation*}

\begin{Lemma}[Poisson Tail Bound]\label{lm-Poisson-tail-bound}
    Suppose $W\sim \Poisson(v)$, then for any $x\geq 0$,
    \begin{equation*}
        \bbP\left(W \geq v + x\right) \leq \exp\left(-\frac{x^2}{2v}\psi_{Benn}(x/v)\right).
    \end{equation*}
    For any $0 \leq x \leq v$,
    \begin{equation*}
        \bbP\left(W \leq v - x\right) \leq \exp\left(-\frac{x^2}{2v}\psi_{Benn}(-x/v)\right) \leq \exp\left(-\frac{x^2}{2v}\right),
    \end{equation*}
    where $\psi_{Benn}(t) = \frac{(1+t)\log(1+t)-t}{t^2/2}$ for $t>-1, t\neq 0$, and $\psi_{Benn}(0)=1$.\\
    Specifically, for any $0<x<2v$, we have
    \begin{equation*}
        \bbP(W \geq v+x) \leq \exp\left(-\frac{x^2}{4v}\right).
    \end{equation*}
\end{Lemma}
\textbf{Proof of Lemma \ref{lm-Poisson-tail-bound}.} See \cite{boucheron2013concentration}.\qquad $\square$

\begin{Lemma}[Bias and Subgaussian of Truncated Poisson Distribution]\label{lm-log-error}
Suppose $W \sim \text{Poisson}(\nu)$ and $W' = W1_{\left\{\frac{1}{10}\nu \leq W \leq 10\nu\right\}} + \nu 1_{\left\{W \not \in [\frac{1}{10}\nu, 10\nu]\right\}}$, then for any $\epsilon>0$, there exists $C_\epsilon > 0$ that only depends on $\epsilon$, such that for all $\nu \geq C_\epsilon$, we have 
\begin{equation*}
	\left|\bbE \log(W' + \frac{1}{2}) - \log \nu\right| \leq 4\nu^{-\frac{3}{2}+\epsilon},
\end{equation*}
and we can also find some absolute constant $C,K_0>0$ such that if $\nu \geq C$, then 
\begin{equation*}
	\left|\bbE(\log^2 (W'+1/2)) - \log^2 \nu\right| \leq 4/\nu,
\end{equation*}
\begin{equation*}
	\left\|\sqrt{\nu}\left(\log(W'+\frac{1}{2}) - \bbE \log(W'+\frac{1}{2})\right)\right\|_{\psi_2} \leq K_0.
\end{equation*}
\end{Lemma}
\textbf{Proof of Lemma \ref{lm-log-error}.} see \cite[Lemma 1, Lemma 2, Lemma 3]{shi2018high}.

\begin{Lemma}\label{lm-xi-Poisson}
Let $\cY_{ijk} \sim \text{Poisson}(I \exp(\cX_{ijk}))$ independently, and each entry of $\cX$ is bounded with $|\cX_{ijk}| \leq B$. Suppose $I \geq C\left(Be^B \log \overline p \vee 1\right)$, then with probability at least $1-\frac{c}{\overline p}$, we have
\begin{equation*}
 \sup_{\substack{\W_l \in \bbR^{p_l\times r_l}, \|\W_l\|\leq 1, l=1,2,3 \\ \cS \in \bbR^{r_1\times r_2 \times r_3},\left\|\cS\right\|_\tF \leq 1}} \left\langle \cZ, \S \times_1 \W_1 \times_2 \W_2 \times_3 \W_3 \right\rangle\leq C\sqrt{\frac{df \cdot e^B}{I}}.
\end{equation*}
where $\cZ:=\cY/I - \exp(\cX^*)$, and $df:=r_1r_2r_3 + \sum_{k=1}^3 p_kr_k$.
\end{Lemma}
\textbf{Proof of Lemma \ref{lm-xi-Poisson}.} 
For each Poisson random variable $\cY_{jkl}$, define 
\begin{equation*}
\cY'_{jkl} = \cY_{jkl} 1_{\{\cY_{jkl} \leq 3 I\exp(\cX_{jkl})\}},
\end{equation*}
and let $\cZ' = \cY'/I - \bbE \cY'/ I$. Notice that when $\cY = \cY'$, we have $\cZ = \cZ' + \bbE(\cY' - \cY)/I$. Now for any $t>0$, we have
\begin{equation}\label{eq-lm-xi-poisson-1}
	\begin{split}
		& \bbP  \left(\sup_{\substack{\W_l \in \bbR^{p_l\times r_l}, \|\W_l\|\leq 1 \\ \cS \in \bbR^{r_1\times r_2 \times r_3},\left\|\cS\right\|_\tF \leq 1}} \left| \left\langle \cZ, \cS \times_1 \W_1 \times_2 \W_2 \times_3 \W_3 \right\rangle\right| > t\right) \\
		& \leq 	\bbP\left(\sup_{\substack{\W_l \in \bbR^{p_l\times r_l}, \|\W_l\|\leq 1 \\ \cS \in \bbR^{r_1\times r_2 \times r_3},\left\|\cS\right\|_\tF \leq 1}} \left| \left\langle \cZ' + \bbE(\cY' - \cY)/I, \cS \times_1 \W_1 \times_2 \W_2 \times_3 \W_3 \right\rangle\right| > t\right) \\
		& + \bbP \left(\exists~(j,k,l), \cY_{jkl}' \neq \cY_{jkl} \right).
	\end{split}
\end{equation}
 We first provide an upper bound for the second probability:
\begin{equation}\label{eq-lm-xi-poisson-2}
	\begin{split}
		 \bbP (\exists~(j,k,l), \cY_{jkl}' & \neq \cY_{jkl})  = \bbP \left(\cup_{j,k,l}\left\{\cY_{jkl}' \neq \cY_{jkl}\right\} \right) \\
		& \leq \sum_{j,k,l} \bbP\left(\cY_{jkl}' \neq \cY_{jkl}\right) = \sum_{j,k,l} \bbP\left(\cY_{jkl} > 3I \exp(\cX_{jkl})\right) \\
		& \overset{\text{Lemma \ref{lm-Poisson-tail-bound}}}{\leq} \sum_{j,k,l}  \exp\left(-I\exp(\cX_{jkl})\right) \\
		& \leq p_1p_2p_3 \exp\left(-I e^{-B}\right) \leq  c/\overline{p}
	\end{split}
\end{equation}
Here in the last step we use the assumption that $I > Ce^B\log \overline p$. \\
Now we study the first probability in \eqref{eq-lm-xi-poisson-1}, first of all, we have:
\begin{equation}\label{eq-lm-xi-poisson-3}
	\begin{split}
		& \left| \bbE(\cY_{jkl}' - \cY_{jkl})/I \right| =  \bbE(\cY_{jkl}1_{\{\cY_{jkl} > 3I\exp(\cX_{jkl})\}})/I \\
		& \leq \frac{1}{I}\left(\bbE \cY_{jkl}^2\right)^{1/2} \left(\bbP\left(\cY_{jkl}>3I\exp(\cX_{jkl})\right)\right)^{1/2} \\
		& \overset{\text{Lemma \ref{lm-Poisson-tail-bound}}}{\leq} \frac{1}{I} \left(\sqrt{2} I\exp(\cX_{jkl})\right) \exp\left(-\frac{1}{2}I\exp(\cX_{jkl})\right)\\
		& \leq \sqrt{2} \exp\left(-\frac{1}{2}Ie^{-B} + B\right) \leq \sqrt{2}\exp\left(\frac{1}{2}B - \frac{1}{2}\log I\right) = \sqrt{\frac{2e^B}{I}}.
	\end{split}
\end{equation}
Next, we claim that $\{\cZ_{jkl}'\}$ are independent sub-Gaussian random variables. The following inequality bulids the tail bound of $\left|\cZ_{jkl}'\right|$:
\begin{equation*}
	\begin{split}
		&\bbP\left(|\cZ_{jkl}'| > t\right) = \bbP\left(\cZ_{jkl}' > t\right) + \bbP\left(\cZ_{jkl}' < -t\right)\\
		&  = \bbP\left(\cY'_{jkl} - \bbE \cY_{jkl}' > I t\right) + \bbP\left(\cY'_{jkl} - \bbE \cY_{jkl}' < - I t \right) \\
		& \overset{\eqref{eq-lm-xi-poisson-3}}{\leq} \bbP\left(\cY'_{jkl} > \bbE \cY_{jkl} - \sqrt{2Ie^B} + It\right) + \bbP\left(\cY'_{jkl} < \bbE \cY_{jkl} - I t\right)
	\end{split}
\end{equation*}
The upper tail probability can be bounded when $t$ takes different values:
\begin{itemize}
\item $2\sqrt{2I^{-1}e^B} < t < \sqrt{2I^{-1}e^B} + 2\exp(\cX_{jkl})$:
\begin{equation*}
	\begin{split}
		&\bbP\left(\cY'_{jkl} > \bbE \cY_{jkl} - \sqrt{2Ie^B} + It\right) \\
		& \overset{\text{Lemma \ref{lm-Poisson-tail-bound}}}{\leq} \exp\left(-\frac{1}{4} \frac{I^2(t-\sqrt{2I^{-1}e^B})^2}{I\exp(\cX_{jkl})}\right) \\
		& \leq \exp\left(-\frac{1}{4} \frac{I(t/2)^2}{\exp(\cX_{jkl})}\right)  \leq \exp\left(- \frac{I t^2}{16 e^B}\right).
	\end{split}
\end{equation*}
\item $t > \sqrt{2I^{-1}e^B} + 2\exp(\cX_{jkl})$:
\begin{equation*}
	\begin{split}
		\bbP\left(\cY'_{jkl} > \bbE \cY_{jkl} - \sqrt{2Ie^B} + I t\right) = 0.
	\end{split}
\end{equation*}
\end{itemize}
In conclusion, we have
\begin{equation*}
	\bbP\left(\cY'_{jkl} > \bbE \cY_{jkl} - \sqrt{2Ie^B} + It\right) \leq \exp\left(- \frac{It^2}{16 e^B}\right),\quad \forall t \geq 2\sqrt{2I^{-1}e^B} .
\end{equation*}
In the meantime, for any $t>0$ we have
\begin{equation*}
	\bbP\left(\cY'_{jkl} < \bbE \cY_{jkl} - It\right) \overset{\text{Lemma \ref{lm-Poisson-tail-bound}}}{\leq} \exp\left(-\frac{I^2 t^2}{2I\exp(\cX_{jkl})}\right) \leq \exp\left(-\frac{It^2}{2e^B}\right).
\end{equation*}
Thus
\begin{equation*}
	\begin{split}
		\bbP(|\cZ_{jkl}'| > t) \leq 2 \exp\left(-\frac{I t^2}{16e^B}\right), \quad \forall  t\geq 2\sqrt{2I^{-1}e^B}.
	\end{split}
\end{equation*}
Now for any $q \geq 1$,
\begin{equation*}
	\begin{split}
		& \bbE \left|\cZ_{jkl}'\right|^q = \int_0^\infty \bbP\left(\left|\cZ_{jkl}'\right| > x\right)qx^{q-1}dx  \\
		& \leq  \int_0^{2\sqrt{2I^{-1}e^B}} \bbP\left(\left|\cZ_{jkl}'\right| > x\right)qx^{q-1}dx + \int_{2\sqrt{2I^{-1}e^B}}^\infty 2\exp\left(-\frac{I x^2}{16e^B}\right)qx^{q-1}dx \\
		& \leq  \int_0^{2\sqrt{2I^{-1}e^B}} qx^{q-1}dx + \int_{0}^\infty 2\exp\left(-\frac{Ix^2}{16e^B}\right)qx^{q-1}dx \\
		& = \left(2\sqrt{2I^{-1}e^B}\right)^q + 4q\int_0^\infty e^{-y}\left(\frac{4e^{B/2}\sqrt{y}}{\sqrt{I}}\right)^{q-1} \frac{e^{B/2}}{\sqrt{Iy}}dy \\
		& = \left(2\sqrt{2I^{-1}e^B}\right)^q + 4^qq e^{qB/2} I^{-q/2} \Gamma(q/2) \\
		& \leq \left(2\sqrt{2I^{-1}e^B}\right)^q + 4^qq e^{qB/2} I^{-q/2} \left(\frac{q}{2}\right)^{q/2} \\
		& \leq \left(2\sqrt{2I^{-1}e^B} + 2\sqrt{2} q^{1/q}\sqrt{\frac{qe^B}{I}}\right)^q.
	\end{split}
\end{equation*}
Here the last inequality comes from the fact that $x^q + y^q \leq (x+y)^q$ for $q\geq 1$.
Then we can bound the $\psi_2$ norm of $\cZ_{jkl}'$:
\begin{equation*}
	\begin{split}
		\left\|\cZ_{jkl}'\right\|_{\psi_2} &= \sup_{q \geq 1} p^{-1/2} \left(\bbE |\cZ_{jkl}'|^q\right)^{1/p} \leq 2\sqrt{2}\sup_{q \geq 1} q^{-1/2} \left(\sqrt{I^{-1}e^B} + q^{1/q}\sqrt{\frac{qe^B}{I}}\right) \\
		& \leq 2\sqrt{2}\left(\sqrt{\frac{e^B}{I}} + e^{\frac{1}{e}}\sqrt{\frac{e^B}{I}}\right) \leq C\sqrt{\frac{e^B}{I}}.
	\end{split}
\end{equation*}
Now by Lemma \ref{lm-concentration-Gaussian-xi}, we have some universal constant $C$, such that
with probability $1-e^{c \cdot df}$, 
\begin{equation}\label{eq-lm-xi-poisson-4}
    \sup_{\substack{\W_l \in \bbR^{p_l\times r_l}, \|\W_l\|\leq 1 \\ \cS \in \bbR^{r_1\times r_2 \times r_3},\left\|\cS\right\|_\tF \leq 1}} \left\langle \cZ', \cS \times_1 \W_1 \times_2 \W_2 \times_3 \W_3 \right\rangle \leq C\sqrt{\frac{df \cdot e^B}{I}}.
\end{equation}

Combining \eqref{eq-lm-xi-poisson-3} and \eqref{eq-lm-xi-poisson-4}, we have
\begin{equation}\label{eq-lm-xi-poisson-5}
    \bbP\left(\sup_{\substack{\W_l \in \bbR^{p_l\times r_l}, \|\W_l\|\leq 1 \\ \cS \in \bbR^{r_1\times r_2 \times r_3},\left\|\cS\right\|_\tF \leq 1}} \left| \left\langle \cZ' + \bbE(\cY' - \cY)/I, \llbracket \cS;\W_1,\W_2,\W_3 \rrbracket \right\rangle\right| > (C+1) \sqrt{\frac{df \cdot e^B}{I}}\right) \leq e^{c\cdot df},
\end{equation}
and the conclusion follows by combining \eqref{eq-lm-xi-poisson-1}, \eqref{eq-lm-xi-poisson-2} and \eqref{eq-lm-xi-poisson-5}.\quad $\square$

\begin{Lemma}\label{lm-binomial-bias-variance}
    Suppose $Y \sim \text{Binomial}(N, s(x))$, where $s(x) = 1/(1+e^{-x})$. Let $\hat p' = \hat p1_{\left\{|\hat p - p| \leq \frac{p}{2}\right\}} + p1_{\left\{|\hat p - p|>\frac{p}{2}\right\}}$ and $X = \log\left(\frac{\hat p'+\frac{1}{2N}}{1-\hat p'+\frac{1}{2N}}\right)$ where $p = s(x)$, $\hat p = Y/N$. Then if $N \geq C|x|e^{3|x|}$, for any $\epsilon>0$, we have universal constants $K_0$, $C$, such that
    \begin{equation}\label{eq-binomial-bias-variance}
        \begin{split}
            \left|\bbE X - x\right| \leq Ce^{3|x|}N^{-\frac{3}{2}+\varepsilon}, \\
            \left\|e^{-|x|}\sqrt{N}\left(X - \bbE X \right)\right\|_{\phi_2} \leq K_0.
        \end{split}
    \end{equation}
\end{Lemma}
\textbf{Proof of Lemma \ref{lm-binomial-bias-variance}}: To simplify the proof, we only consider the case where $x \leq 0$. Let $g(t) = \log\left(\frac{t+\frac{1}{2N}}{1-t+\frac{1}{2N}}\right), t \in [0,1]$, then we have $X = g(\hat p')$. We first provide the derivatives of $g(t)$ up to order $3$:
\begin{equation}\label{eq-Binomial-derivatives}
    \begin{split}
        & g'(t) = \left(t+\frac{1}{2N}\right)^{-1} + \left(1-t+\frac{1}{2N}\right)^{-1},\\
        & g''(t) = -\left(t+\frac{1}{2N}\right)^{-2} + \left(1-t+\frac{1}{2N}\right)^{-2},\\
        & g'''(t) = 2\left(t+\frac{1}{2N}\right)^{-3} +  2\left(1-t+\frac{1}{2N}\right)^{-3}.
    \end{split}
\end{equation}
Let $p:=s(x)\leq \frac{1}{2}$, by Taylor's expansion of $g(t)$ at $t_0 = \frac{N+1}{N}p - \frac{1}{2N}$, we have
\begin{equation}\label{eq-Binomial-gt-taylor}
    g(t) = g(t_0) + g'(t_0)(t-t_0) + \frac{1}{2}g''(t_0)(t-t_0)^2 + \frac{1}{6}g'''(\xi)(\xi-t_0)^3,
\end{equation}
where $\xi$ is some number between $t_0$ and $t$. Now we let $f(t) = g(t) - g(t_0) - g'(t_0)(t-t_0) - \frac{1}{2}g''(t_0)(t-t_0)^2$, then one can see that as long as $-\frac{1}{2N} \leq t_0 - N^{-\frac{1}{2}+\epsilon} \leq t \leq t_0 + N^{-\frac{1}{2}+\epsilon} \leq 1+\frac{1}{2N}$ (the first and third inequality holds since $N > C|x|e^{3|x|} > Cp^{-3}$), we have some universal constant $C$, such that
\begin{equation}
    \begin{split}
        & |f(t)| \leq \frac{1}{3}\sup_{t: |t-t_0|\leq N^{-\frac{1}{2}+\epsilon}}\left(\left|t+\frac{1}{2N}\right|^{-3} + \left|1-t+\frac{1}{2N}\right|^{-3}\right)\left|t-t_0\right|^3 \\
        & \qquad \leq \frac{2}{3}\left|\frac{N+1}{N}p-N^{-\frac{1}{2}+\epsilon}\right|^{-3} N^{-3/2+3\epsilon} \leq Cp^{-3}N^{-\frac{3}{2}+\epsilon_1}
    \end{split}
\end{equation}
with $\epsilon_1 = 3\epsilon$. Here the last inequality comes from the following fact:
\begin{equation*}
	\begin{split}
		\frac{N+1}{N}p - N^{-\frac{1}{2}+\epsilon} & \geq  \frac{N+1}{N}p - p^{\frac{3}{2}-3\epsilon} \geq  cp.
	\end{split}
\end{equation*}

 Thus it follows that
\begin{equation}
    \begin{split}
        & \bbE \left|f(\hat p)1_{\{t_0-N^{-\frac{1}{2}+\epsilon}<\hat p<t_0-N^{-\frac{1}{2}+\epsilon}\}}\right| \leq Cp^{-3}N^{-\frac{3}{2}+\epsilon_1} = C(e^{|x|}+1)^3N^{-\frac{3}{2}+\epsilon_1}.
    \end{split}
\end{equation}
Then we have
\begin{equation}\label{eq-Binomial-bias-1}
    \begin{split}
        & \bbE \left|f(\hat p) \right| \leq \bbE \left|f(\hat p)1_{\{t_0-N^{-\frac{1}{2}+\epsilon}<\hat p<t_0-N^{-\frac{1}{2}+\epsilon}\}}\right| + \bbE \left|f(\hat p)1_{\{\hat p \in [0,t_0-N^{-\frac{1}{2}+\epsilon}] \cup [t_0+N^{-\frac{1}{2}+\epsilon}, 1]\}}\right| \\
        & \leq C(e^{|x|}+1)^3N^{-\frac{3}{2}+\epsilon_1} + \sup_{t\in [0,1]}|g'''(t)| \cdot \bbP\left(|\hat p-t_0|>N^{-\frac{1}{2}+\epsilon}\right) \\
        & \leq C(e^{|x|}+1)^3N^{-\frac{3}{2}+\epsilon_1} + \frac{32}{3}N^3\bbP\left(|\hat p-p|>N^{-\frac{1}{2}+\epsilon}-N^{-1}\right)\\
        & \leq C(e^{|x|}+1)^3N^{-\frac{3}{2}+\epsilon_1} + \frac{32}{3}N^3\bbP\left(|\hat p-p|>\frac{1}{2}N^{-\frac{1}{2}+\epsilon}\right)\\
        & \leq C(e^{|x|}+1)^3N^{-\frac{3}{2}+\epsilon_1} + CN^3e^{-\frac{1}{2}N^{2\varepsilon}} \leq Ce^{3|x|}N^{-\frac{3}{2}+\epsilon_1},
    \end{split}
\end{equation}
where the last but one inequality comes from applying the following concentration inequality for Binomial random variable:
\begin{equation}\label{eq-binomial-concentration}
    \bbP\left(\left|Y - p\right| \geq t\right) \leq 2\exp\left(-2Nt^2\right),\quad \forall t\geq 0.
\end{equation}
In the mean time, one can calculate that
\begin{equation}\label{eq-Binomial-bias-2}
    \begin{split}
        & \bbE g(\hat p) - \bbE f(\hat p) = g(t_0) +  \frac{2p-1}{(N+1)^2p(1-p)} \left(\frac{1}{8p(1-p)}-\frac{1}{2}\right)\\
        & = \log(\frac{p}{1-p}) + \frac{2p-1}{(N+1)^2p(1-p)} \left(\frac{1}{8p(1-p)}-\frac{1}{2}\right).
    \end{split}
\end{equation}
Combining \eqref{eq-Binomial-bias-1} and \eqref{eq-Binomial-bias-2}, we have
\begin{equation*}
    \begin{split}
        & \left|\bbE g(\hat p) - x\right| = \left| \bbE g(\hat p) - \log\left(\frac{p}{1-p}\right) - \bbE f(\hat p) + \bbE f(\hat p)\right| \\
        & \leq \left|\bbE g(\hat p) - \log\left(\frac{p}{1-p}\right) - \bbE f(\hat p)\right| + \left|\bbE f(\hat p)\right| \\
        & \leq \frac{2p-1}{(N+1)^2p(1-p)} \left(\frac{1}{8p(1-p)}-\frac{1}{2}\right) + Ce^{3|x|}N^{-\frac{3}{2}+\epsilon_1} \\
        & \leq Ce^{2|x|}N^{-2} + Ce^{3|x|}N^{-\frac{3}{2}+\epsilon_1} \\
        & \leq Ce^{3|x|}N^{-\frac{3}{2}+\epsilon_1}.
    \end{split}
\end{equation*}
Now we calculate the bias of $X = g(\hat p')$. First,
\begin{equation*}
	\begin{split}
		& \left|\bbE g(\hat p') - \bbE g(\hat p)\right| = \left|g(p)\bbP\left(|\hat p - p|\geq \frac{p}{2}\right) + \bbE g(\hat p)1_{\{|\hat p - p| \leq \frac{p}{2}\}} - \bbE g(\hat p)\right| \\
		\leq & \left|g(p)\bbP(\hat p \leq \frac{p}{2}) - \bbE g(\hat p)1_{\{\hat p \leq \frac{p}{2}\}}\right| + \left|g(p)\bbP(\hat p \geq \frac{3p}{2}) - \bbE g(\hat p)1_{\{\hat p \geq \frac{3p}{2}\}}\right| \\
		\leq & \max\left\{g(p)\bbP(\hat p \leq \frac{p}{2}), \bbE g(\hat p)1_{\{\hat p \leq \frac{p}{2}\}}\right\} + \max\left\{g(p)\bbP(\hat p \geq \frac{3p}{2}), \bbE g(\hat p)1_{\{\hat p \geq \frac{3p}{2}\}}\right\} \\
		\overset{(a)}{\leq} & \max\left\{g(p)e^{-\frac{1}{2}Np^2}, \left(\bbE g^2(\hat p)\right)^{1/2}\left(\bbP(\hat p \leq \frac{p}{2})\right)^{1/2}\right\} \\
		& + \max\left\{g(p)e^{-\frac{1}{2}Np^2}, \left(\bbE g^2(\hat p)\right)^{1/2}\left(\bbP(\hat p \geq \frac{3p}{2})\right)^{1/2}\right\} \\
		\leq & 2 \max\left\{g(p)e^{-\frac{1}{2}Np^2}, \left(\bbE g^2(\hat p)\right)^{1/2}e^{-\frac{1}{2}Np^2}\right\} \\
		& \overset{(b)}{\leq} 2\log(2N+1) e^{-\frac{1}{2}Np^2} \leq 2\log(2N+1) e^{-\frac{1}{2}N^{1/3}} \leq Ce^{3|x|}N^{-\frac{3}{2}+\epsilon_1}.
	\end{split}
\end{equation*}
Here, we use Cauchy-Schwarz inequality for (a), and apply the uniform bound $g(p) \leq \log(2N+1)$ for (b). Combining two inequalities above, we obtain
\begin{equation*}
	\left|\bbE g(\hat p') - x \right| \leq \left|\bbE g(p') - \bbE g(\hat p)\right| + \left|\bbE g(\hat p) -x\right| \leq Ce^{3|x|}N^{-\frac{3}{2}+\epsilon_1},
\end{equation*}
which gives the first inequality of \eqref{eq-binomial-bias-variance}.\\
Now we prove the second inequality of \eqref{eq-binomial-bias-variance}. We denote $D = g(\hat p') - \bbE g(\hat p')$ for convenience. By what we have proved, we have $\left|\bbE g(\hat p') - x\right| \leq Ce^{3|x|}N^{-\frac{3}{2}+\epsilon_1} \leq  Ce^{2|x|}N^{-1}$, thus we have
\begin{equation*}
    \begin{split}
        & \bbP\left(|e^{x}\sqrt ND| > t\right) = \bbP\left(e^{x}\sqrt ND>t\right) + \bbP\left(e^{x}\sqrt ND<-t\right) \\
        & = \bbP\left(g(\hat p') < \bbE g(\hat p') - \frac{te^{-x}}{\sqrt N} \right) + \bbP\left(g(\hat p') > \bbE g(\hat p') + \frac{te^{-x}}{\sqrt N} \right).
    \end{split}
\end{equation*}
Then
\begin{equation*}
	\begin{split}
		& \bbP\left(g(\hat p') < \bbE g(\hat p') - \frac{te^{-x}}{\sqrt{N}} \right) \leq \bbP\left(g(\hat p') < x - te^{-x}N^{-1/2} + Ce^{-2x}N^{-1}\right),\\
		& \bbP\left(g(\hat p') > \bbE g(\hat p') + \frac{te^{-x}}{\sqrt N} \right) \leq \bbP\left(g(\hat p') > x + te^{-x}N^{-1/2} - Ce^{-2x}N^{-1}\right).
	\end{split}
\end{equation*}
We first investigate the lower tail bound:
\begin{equation*}
	\begin{split}
		&\bbP\left(g(\hat p') < x - te^{|x|}N^{-1/2} + Ce^{2|x|}N^{-1}\right) \\
		& = \bbP\left(\hat p' < \frac{1+1/N}{\exp\left(-x+te^{|x|}N^{-1/2}-Ce^{2|x|}N^{-1}\right)+1} - \frac{1}{2N}\right)
	\end{split}
\end{equation*}
When $t>1$, since $te^{|x|}N^{-1/2} - Ce^{2|x|}N^{-1} > \frac{1}{2}t e^{|x|}N^{-1/2}$ by the assumption on $N$, we have
\begin{equation*}
	\begin{split}
		& \bbP\left(\hat p < \frac{1+1/N}{\exp\left(-x+te^{|x|}N^{-1/2}-Ce^{2|x|}N^{-1}\right)+1} - \frac{1}{2N}\right) \\
		& \leq \bbP\left(\hat p < \frac{1+1/N}{\exp\left(-x+\frac{1}{2}te^{|x|}N^{-1/2}\right)+1} - \frac{1}{2N}\right).
	\end{split}
\end{equation*}
When $t > 6e^xN^{1/2}$, we have
\begin{equation*}
	\begin{split}
		& \bbP\left(p' < \frac{1+1/N}{\exp\left(-x+\frac{1}{2}te^{|x|}N^{-1/2}\right)+1} - \frac{1}{2N}\right)  \leq \bbP\left(p' < \frac{1.1}{\exp\left(-x+\frac{1}{2}te^{|x|}N^{-1/2}\right)+1}\right) \\
		& \leq \bbP\left( p' \leq \frac{1.1}{\exp(-x + 3)+1}\right)  = \bbP\left( p' \leq \frac{1.1}{e^3\exp(-x )+1}\right) \\
		& \leq \bbP\left( p' \leq \frac{1}{2(\exp(-x)+1)}\right) =  \bbP\left( p' \leq \frac{p}{2}\right) = 0,
	\end{split}
\end{equation*}
where the last identity comes from the definition of $p'$. When $1 < t < 6e^xN^{1/2}$, we have
\begin{equation*}
	\begin{split}
		& \bbP\left(\hat p < \frac{1+1/N}{\exp\left(-x+\frac{1}{2}te^{|x|}N^{-1/2}\right)+1} - \frac{1}{2N}\right) \\
		& \leq \bbP\left(\hat p < \frac{1+1/N}{\exp\left(-x+\frac{1}{2}te^{|x|}N^{-1/2}\right)+1}\right) \\
		& = \bbP\left( \hat  p - p < \frac{1+1/N}{\exp\left(-x+\frac{1}{2}te^{|x|}N^{-1/2}\right)+1} - \frac{1}{1+e^{-x}}\right) \\
		& = \bbP\left( \hat  p - p < \frac{e^{-x}\left(1 - e^{\frac{1}{2}te^{|x|}N^{-1/2}}\right) + N^{-1}(1+e^{-x})}{\left(\exp\left(-x+\frac{1}{2}te^{|x|}N^{-1/2}\right)+1\right)(1+e^{-x})}\right) \\
		& \leq \bbP\left( \hat  p - p < \frac{-\frac{1}{2}\left(e^{\frac{1}{2}te^{|x|}N^{-1/2}} - 1\right) + N^{-1}}{\exp\left(-x+\frac{1}{2}te^{|x|}N^{-1/2}\right)+1}\right) \\
		& \leq \bbP\left( \hat  p - p < -\frac{1}{4}\frac{te^{|x|}N^{-1/2}}{\exp\left(-x+3\right)+1} + \frac{1}{N(e^{-x}+1)}\right) \\
		& \leq \bbP\left( \hat  p - p < -\frac{1}{4e^3}tN^{-1/2} + \frac{1}{N^{-1/2}(e^{-x}+1)}N^{-1/2}\right) \\
		& \leq \bbP\left( \hat  p - p < -\frac{1}{4e^3}tN^{-1/2} + ct N^{-1/2}\right) \leq \exp\left(-ct^2\right).
	\end{split}
\end{equation*}
Thus we have proved that $\forall t > 1$, there exists constant $c$, such that the lower tail bound
\begin{equation*}
 	\bbP\left(e^x\sqrt{N}D < -t\right) \leq e^{-ct^2}.
\end{equation*}
We can prove the similar result for upper tail bound. Thus, $\forall p \geq 1$, we have
\begin{equation*}
	\begin{split}
		&\bbE |e^{x}\sqrt{N}D|^q = \int_{0}^\infty \bbP\left(|e^{x}\sqrt{N}D|>t\right)qt^{q-1}dt \leq \int_{0}^1 qt^{q-1}dt + \int_1^\infty e^{-ct^2}qt^{q-1}dt \\ 
		\leq & \int_{0}^1 qt^{q-1}dt + \int_0^\infty e^{-ct^2}qt^{q-1}dt = 1 + c^{-\frac{q}{2}}q\Gamma\left(\frac{q}{2}\right)\\
		\leq & 1 + qc^{-\frac{q}{2}}\left(\frac{q}{2}\right)^{\frac{q}{2}} \leq \left(1+q^\frac{1}{q}c^{-1/2}\left(\frac{q}{2}\right)^{1/2}\right)^q, 
	\end{split}
\end{equation*}
where we use $x^q + y^q \leq (x+y)^q$ for all $x,y \geq 0, q \geq 1$. Then it follows that
\begin{equation*}
	\sup_{q \geq 1} q^{-\frac{1}{2}} \left(\bbE \left|e^x\sqrt{N}D\right|^q\right)^{1/q}\leq q^{-1/2} \left(1+q^\frac{1}{q}c^{-1/2}\left(\frac{q}{2}\right)^{1/2}\right) \leq (1+\frac{e^{1/e}c^{-1/2}}{\sqrt{2}}) =: K_0,
\end{equation*}
and we have proved the second inequality of \eqref{eq-binomial-bias-variance}. \quad\quad $\square$

\end{document}